\documentclass[mnsc]{arxiv_version} % current default for manuscript submission
\DoubleSpacedXI

\usepackage{multirow}
\usepackage{natbib}
 \bibpunct[, ]{(}{)}{,}{a}{}{,}%
 \def\bibfont{\small}%
\TheoremsNumberedThrough     % Preferred (Theorem 1, Lemma 1, Theorem 2)
\EquationsNumberedThrough    % Default: (1), (2), ...
\usepackage{pgfplots}
\pgfplotsset{compat=newest,
    width=7.5cm,
    height=5.5cm,
    scale only axis=true,
    max space between ticks=25pt,
    try min ticks=5,
    every axis/.style={
        axis y line=left,
        axis x line=bottom,
        axis line style={thick,->,>=latex, shorten >=-.4cm}
    },
    every axis plot/.append style={thick},
    tick style={black, thick}
}
\tikzset{
    semithick/.style={line width=0.8pt},
}
\usepgfplotslibrary{groupplots}
\usepgfplotslibrary{dateplot}

\usepackage{epstopdf}
\usepackage[ruled]{algorithm2e} % For algorithms

\usepackage{setspace}
\usepackage[noend]{algpseudocode}

\usepackage{xcolor}

\definecolor{blue}{HTML}{1F77B4}
\definecolor{orange}{HTML}{FC8D62}
\definecolor{teal}{HTML}{036c5f}
\definecolor{purple}{HTML}{9467BD}
\definecolor{yellow}{HTML}{E6AB02}
\definecolor{darkgreen}{HTML}{2CA02C}
\definecolor{brown}{HTML}{7B3F00}
\definecolor{line_blue}{HTML}{377eb8}
\definecolor{line_orange}{HTML}{ff7f00}
\definecolor{line_purple}{HTML}{9467BD}
\definecolor{line_yellow}{HTML}{E6AB02}
\definecolor{line_green}{HTML}{2CA02C}
\definecolor{line_brown}{HTML}{7B3F00}

\usepackage{tikz}
\usetikzlibrary{shapes.multipart,arrows,automata}
\usetikzlibrary{positioning,arrows.meta,bending}
\usetikzlibrary{patterns,calc,intersections,through,backgrounds,decorations.pathreplacing}

\usepackage{enumerate}
\usepackage{ifthen}
\usepackage{accents}
\usepackage{setspace}
\usepackage{hyperref}
\usepackage{cleveref}
\usepackage{bbm}

\newcommand\numberthis{\addtocounter{equation}{1}\tag{\theequation}}
\newcommand{\ubar}[1]{\underaccent{\bar}{#1}}
\newcommand{\pack}{p^\star} %battery pack size
\newcommand{\soc}{s} %state of charge
\newcommand{\bbZ}{\mathbb{Z}}
\newcommand{\bbR}{\mathbb{R}}
\newcommand{\1}[1]{\mathbf{1}\{#1\}} %indicator
\newcommand{\cp}{c_P} %cars picking up

\newcommand{\crt}{c_R} %cars on route
\newcommand{\cc}{c_C} %cars charging
\newcommand{\cd}{c_{DC}} %cars driving to charger
\newcommand{\ci}{c_I} %cars idling
\newcommand{\ecp}{\bar{c}_P} %equilibrium  cars picking up
\newcommand{\ecrt}{\bar{c}_R} % equilibrium cars on route
\newcommand{\ecc}{\bar{c}_C} %equilibrium  cars charging
\newcommand{\ecd}{\bar{c}_{DC}} %equilibrium  cars driving to charger
\newcommand{\eci}{\bar{c}_I} %equilibrium cars idling
\newcommand{\eq}{\bar{q}} %equilibrium customers

\newcommand{\pol}{\pi} 
\newcommand{\polSet}{\Pi} 
\newcommand{\lambdaeff}{\lambda_{\operatorname{eff}}} 
 
\newcommand{\alphasup}{\bar{\alpha}^\pi_{\operatorname{eff}}}
\newcommand{\alphainf}{\ubar{\alpha}^\pi_{\operatorname{eff}}}

\newcommand{\XC}{X^C} 
 
\newcommand{\dXC}{\dot{X}^C} 
\newcommand{\XB}{X^B}
\newcommand{\dXB}{\dot{X}^B} 
\newcommand{\bXC}{\bar{X}^C} 
\newcommand{\bXB}{\bar{X}^B} 
\newcommand{\tXC}{\tilde{X}^C}
\newcommand{\beq}{\begin{equation}} 
\newcommand{\eeq}{\end{equation}} 
\newcommand{\bit}{\begin{itemize}} 
\newcommand{\eit}{\end{itemize}} 
\newcommand{\Nt}{N_T} 
\newcommand{\amplitude}{a} 
\newcommand{\Tpeak}{T_p}
\newcommand{\Tvalley}{T_v}
\newcommand{\lbavg}{\lambda_{\operatorname{avg}}}

\tikzset{fontscale/.style = {font=\relsize{#1}}
    }
\newcommand{\revcolor}[1]{\textcolor{black}{#1}}

\begin{document}
%%%%%%%%%%%%%%%%

%\setcounter{page}{1}
%\DOI{10.1287/xxxx.0000.0000}%

\RUNAUTHOR{Varma, Castro, Maguluri}

\RUNTITLE{EV Capacity Planning}

% \TITLE{Bundling Information Goods of Decreasing Value}
\TITLE{Electric Vehicle Fleet and Charging Infrastructure Planning}

\ARTICLEAUTHORS{%
\AUTHOR{Sushil Mahavir Varma}
\AFF{Industrial and Operations Engineering, University of Michigan, Ann Arbor, \EMAIL{sushilv@umich.edu}}
\AUTHOR{Francisco Castro}
\AFF{Anderson School of Management, University of California, Los Angeles, \EMAIL{francisco.castro@anderson.ucla.edu}}
\AUTHOR{Siva Theja Maguluri}
\AFF{H.Milton School of Industrial and Systems Engineering, Georgia Institute of Technology, \EMAIL{siva.theja@gatech.edu}}
} % end of the block

\ABSTRACT{%
We study electric vehicle (EV) fleet and charging infrastructure planning in a spatial setting. For a centrally managed fleet that serves customer requests arriving continuously at a rate $\lambda$ throughout the day, we determine the minimum number of vehicles and chargers for a target service level, along with matching and charging policies. While non-EV systems require extra $\Theta(\lambda^{2/3})$ vehicles due to pickup times, EV systems differ. Charging increases nominal capacity, enabling pickup time reductions and allowing for an extra fleet requirement of only $\Theta(\lambda^{\nu})$ for $\nu \in (1/2, 2/3]$, depending on charging infrastructure and battery pack sizes.
We propose the Power-of-$d$ dispatching policy, which achieves this performance by selecting the closest vehicle with the highest battery level from $d$ options. We extend our results to accommodate time-varying demand patterns and discuss conditions for transitioning between EV and non-EV capacity planning. Simulations verify our scaling results, insights, and policy effectiveness. While long-range, fast-charging fleets resemble non-EV systems, short-range, low-cost fleets can still perform competitively---underscoring the need for EV-aware management policies.
}

%\KEYWORDS{deterministic inventory theory; infinite linear programming duality; 
%  existence of optimal policies; semi-Markov decision process; cyclic schedule}

%\KEYWORDS{}
% \HISTORY{}

\maketitle
\vspace{-1cm}
\section{Introduction}
In 2022, the United States federal government enacted the Inflation Reduction Act, which provides tax incentives for acquiring electric vehicles, encouraging individuals and businesses to transition to this technology. %Additionally, 
In the state of California, executive order N-79-20 mandated that, by 2035, all new vehicle sales within the state must be zero-emission vehicles, further solidifying the shift towards EVs in transportation \citep{link1}. Shifting consumer preferences---the electric vehicle global market share quadrupled from 2019 to 2021, reaching nearly 10\% in 2021 \citep{bibra2022global}---have prompted major automakers like Ford and General Motors to commit substantial resources to electric vehicle manufacturing \citep{link2,link3}.

In the on-demand transportation industry, several ride-hailing firms have commenced offering electric vehicle services and have put forward plans to transition to an entirely EV-based fleet \citep{link4,link5}. In 2021, Revel launched an employee-based all-electric ride-hailing service in New York City with the vision to continue growing its fleet and charging infrastructure \citep{link7}; BluSmart, an Indian ride-hailing service, started operating a fully EV fleet in 2019 \citep{link8}. Furthermore, recent advances in AI technology have enabled autonomous vehicle companies such as Waymo One to start providing a 24/7 service with its all-electric fleet in several cities in the US \citep{link9} through its own app but also within the Uber ecosystem \citep{link10}.

As companies in this industry expand their electric fleets, the charging infrastructure will need to scale accordingly, bringing to the forefront both economic and operational challenges. Operating with too few chargers can reduce costs\footnote{
Level 2 chargers can charge 100km in 1-2 hours with an installation cost between \$2,500-\$4,900. DC fast chargers can charge 100km in as little as 10 minutes with installation cost estimates ranging from \$20,000-\$100,000.} but may also reduce service quality if customer demand exceeds charging capacity. A smaller fleet can reduce capital expenses—a critical consideration for autonomous vehicle companies—but may also result in unfulfilled trip requests.
In addition, managing the complex charging and spatial dynamics innate to an all-EV on-demand transportation platform in the face of reduced charging availability  \citep{anderson2022electric} and longer refueling times \citep{lamonaca2022state} involves further operational tradeoffs. Postponing charging until late in the day can alleviate some of this complexity, but it risks vehicles being unable to serve some trips, can require a large charging infrastructure, and impose a sudden power demand increase \citep{garcia2014plug}. 
When dispatching vehicles, sending the nearest vehicle minimizes the distance traveled, but it may not be compatible with sustaining a fleet with a consistently high State of Charge (SoC) able to fulfill most customer requests.

To address these operational challenges while maintaining an appropriate service level and cost efficiency, a fully electric on-demand transportation platform should carefully design a system management policy, including charging and dispatching, and strike a balance between charging infrastructure and fleet size. We investigate infrastructure planning and management for an on-demand transportation platform that utilizes a fleet of electric vehicles operating continuously throughout the day, with full control over the fleet \revcolor{(e.g., autonomous vehicles)}. \label{response:setting_clarification_intro} Our objectives are (i) to establish guidelines for determining the minimum number of vehicles and charging stations needed to achieve a desired level of service as demand for rides grows, and (ii) to shed light on dispatch and charging policies that enable a cost-efficient system that appropriately balances workload and battery levels and performs well regardless of system conditions.

\subsection{Main Contributions}\label{intro:contributions}

{\it Infrastructure Planning.}  
We provide a policy-independent lower bound and a near-matching upper bound on the minimum number of vehicles and chargers needed to sustain a service level of $\alpha \in (0,1)$. These bounds depend on the arrival rate of requests $\lambda$, the average trip times $T_R$, and the ratio of charge and discharge rate $r$. In the baseline setting---constant arrival rate---the system operator can achieve any service level $\alpha \in (0,1)$ by setting the number of vehicles and chargers equal to no more and no less than  $(1+r) T_R \alpha \lambda + \Theta(\lambda^{1-\gamma})$ and 
$r T_R \alpha \lambda + \Theta(\lambda^{2\gamma})$, respectively, for any $\gamma\in [1/3, 1/2)$. That is, the system needs at least as many vehicles as the ones that are serving customers ($T_R\alpha \lambda$) and recharging ($rT_R \alpha \lambda $), and as many chargers to accommodate all EVs that need to charge due to the energy consumed during trips. Additionally, the system requires an extra $\Theta(\lambda^{1-\gamma})$ vehicles and $\Theta(\lambda^{2\gamma})$ chargers, showing a fundamental trade-off: a higher (lower) density of chargers implies a shorter (longer) drive-to-charger time, which reduces (increases) the fleet size requirement. 
However, there is a limit to the benefits of trading off EVs and chargers. The system cannot operate with less than $\Theta(\lambda^{1/2})$ extra vehicles or less than $\Theta(\lambda^{2/3})$ extra chargers. 

{\it EV vs. non-EV capacity planning.} 
In a non-EV system (i.e., the charging time is negligible), the results in \citet{besbes2022spatial} imply that the number of vehicles required to satisfy a given service level $\alpha \in (0, 1)$ scales approximately as $T_R\alpha \lambda + \Theta(\lambda^{2/3})$. The second-order term, $\lambda^{2/3}$, reflects the extra capacity needed (compared to a service system without spatial frictions) due to pickup times. We establish a related but different result that points to an important difference between the two systems.
In an EV system, the first-order term is higher because $\Theta(\lambda)$ vehicles are charging on average.  While this {\it extra} capacity is unavoidable, it offers an opportunity\revcolor{: EVs charging can be used to fulfill trip requests, reducing pickup times and} resulting in a second-order term that can be as low as $\Theta(\lambda^{1/2})$. The first-order term is higher, but the second-order term can be lower depending on the number of chargers and the matching policy.

We propose the \textit{Power-of-$d$ dispatch} (Po$d$) policy: for every trip request, choose the $d$ closest vehicles (including those charging) and dispatch the one with the highest SoC.
This policy effectively utilizes the extra vehicles. By selecting nearby, high-SoC vehicles, it attains low pickup times while balancing charge levels across the fleet, securing the availability of enough dispatchable EVs. This strategy enables the system operator to maintain a smaller fleet size, provided there are sufficient chargers. The resulting key insight is that optimal fleet capacity planning in an EV system is susceptible to the charging infrastructure, and it entails a dispatching policy that must balance fleet SoC.

{\it From EV to non-EV infrastructure planning insights.} Using the baseline results as a building block, we consider a time-varying arrival rate setting with demand peaks and valleys. \revcolor{This setting allows us to study when the fleet must be managed with EV‑aware policies and when it resembles a conventional (non‑EV) system.} Three different cases arise depending on the peak demand amplitude. For low and moderate demand amplitudes, there is not enough demand variability to justify shifting charging to the valleys. Vehicles charge during both peaks and valleys, maintaining an extra buffer of partially charged EVs at charging stations, thereby preserving the insights from the baseline model. At high amplitudes, with large battery packs and enough chargers, the fleet behavior approaches that of non-EVs, making the non-EV scaling potentially applicable. \revcolor{ In this case, while charging can be shifted to the valleys, an EV‑aware policy such as Po$d$---that replenishes SoC while still meeting all valley demand---may still be required.
}

{\it Simulations.} 
We corroborate our baseline setting theoretical findings with a simulator that explicitly accounts for the system's spatial dimension and the inherent stochasticity of trip requests and SoC evolution. Our simulation results confirm the scalings, the effectiveness of Po$d$, and the emerging trade-off between fleet size and chargers.

We then use Chicago trip data, which features a time-varying demand pattern, to simulate various systems with short-, medium-, or long-range EVs. We observe the same fleet size-charger trade-off and show how to use a dynamical system describing the system's evolution under Po$d$ to efficiently estimate the minimum fleet requirement. We demonstrate that Po$d$ performs better than or comparably to other matching policies that select the closest vehicle, the closest vehicle with sufficient SoC to complete a trip, or the highest-SoC vehicle within a radius. 
Consistent with our time-varying theoretical results, for a long-range fleet, a policy that performs well in a non-EV system (such as Closest Dispatch) continues to do so in our simulations, but it is still outperformed by EV-aware policies. Although longer ranges and higher charger power make the system behave more like a non-EV setting, a comparison of Po$d$ with an alternative night-only charging policy shows that Po$d$ can enable a short-range fleet with moderate (Level 2) charging to achieve a similar capacity requirement as a long-range fleet using DC fast charging under the alternative policy.

These results highlight the need for adaptable policies that smooth charging throughout the day and balance fleet workload---making short-range and medium-charge-power systems operationally competitive. Policies like Po$d$ and others that prioritize the highest-SoC vehicle within a radius can ensure adequate service quality with minimal infrastructure requirements.

\section{Literature Review} \label{sec:lit}
Our work lies in the intersection of studies about electric vehicles, capacity and infrastructure planning in service systems, and load balancing.

\subsection{Electric Vehicles}
Studies have explored challenges in electric vehicle charging, including grid strain and longer refueling times, proposing solutions like decentralized algorithms for night-time load flattening to avoid new demand peaks from EV commuters \citep{gan2012optimal, ma2011decentralized}, cost- and emission-based delayed charging \citep{wu2022smart}, and management of large-scale charging stations \citep{wang2016smart, mukherjee2014review}.
Other related research explores the optimal deployment of charging stations \citep{paganini2022optimization} and the impact and implementation of battery swapping technology \citep{qi2023scaling,avci2015electric}. Notably, \cite{mak2013infrastructure} develop a robust optimization model for cost-effective infrastructure deployment of EVs with battery swapping, and \cite{he2021charging} propose a methodology for optimizing the location, size, and repositioning of charging stations in an electric vehicle-sharing system like car2go. \cite{abouee2021adoption} investigate the viability and environmental impact of EVs in the vehicle-sharing market, highlighting the importance of charging speed and the number of charging stations. While these studies primarily focus on charging and repositioning aspects of EV-based systems, our work centers around infrastructure planning, specifically the analysis of charging station infrastructure and fleet size in a ride-hailing-like system. Other related applications include optimal repositioning and recharging policies for electric bike/scooter sharing \citep{akturk2022managing, osorio2021optimal, greening2021effective}, vehicle-to-grid electricity selling \citep{zhang2021values}, and the integration of solar power generation with battery storage \citep{solar_operations, kaps2022privately}.

\subsection{Capacity Planning}
The seminal work of \cite{halfin1981heavy}
on heavy traffic limits for $G/M/n$ queues 
implies that a system operator should staff $\lambda + \sqrt{\lambda}$ servers,  where $\lambda$ is the offered load (see \cite{reed2009g} for the $G/G/n$ case and \cite{gans2003telephone, aksin2007modern} for surveys). This result does not translate directly to a spatial service system in which the server needs to pick up the customer before service can start, as shown in \cite{besbes2022spatial}. We consider a related but considerably more complex system, as in addition to the spatial dimension of the problem, each vehicle has an evolving SoC, resulting in a larger state space. Moreover, EVs further introduce the challenge of balancing the SoC across the fleet. Our analysis carefully optimizes the tradeoff between SoC and pickup times as opposed to optimizing just for pickup times as in \cite{besbes2022spatial}. Our results and scalings are also different due to charging and the option of utilizing ``partially" charged vehicles to serve trips. For additional works in the space of spatial matching, we refer the interested reader to \cite{spatial_yash, spatial_stanford, benjaafar2022dimensioning}, and to \cite{daganzo2004bounds} and  \cite{lim2017agility} for spatial analysis in the context of supply chain and transportation.

\cite{deng2022fleet} propose a closed queueing network model to analyze infrastructure planning decisions for electric vehicle-sharing systems. They show that increasing the number of chargers improves vehicle availability and suggest that two slow chargers may be more effective than one fast charger when there is significant variability in charging times.
Our main result precisely characterizes the trade-off between the number of chargers and fleet size and sheds light on the forces behind that trade-off. Additionally, we explicitly model the pickup and drive-to-charger time and jointly minimize the fleet and charging infrastructure to achieve a given service level.

While we use an analytical methodology to guide infrastructure planning, simulation-based approaches have also been explored (see, e.g., \cite{levin2017general, bauer2018cost, loeb2019fleet, vosooghi2019shared, vosooghi2019robo, yang2023fleet}). These studies examine charging infrastructure planning \citep{bauer2018cost}; shared autonomous versus personal vehicle effectiveness \citep{levin2017general}; gasoline versus electric vehicle profitability \citep{loeb2019fleet}; and fleet sizing and charging infrastructure under congestion via simulations \citep{yang2023fleet}.

\subsection{Load Balancing}

Our analysis of the EV system suggests the Power-of-$d$ vehicles dispatch policy, which turns out to closely relate to the Power-of-$d$ choices policy in the load balancing literature \citep{vvedenskaya_power_of_2, mitzenmacher_thesis, mitzenmacher_power_of_d}. A standard observation in this literature is that the performance of a queueing system improves with $d$ as the load is balanced more effectively across its queues \citep{liu2020steady, liu2022steady, liu2022large, varma_sub_halfin_whitt}. In our system, even though a higher value of $d$ balances the SoC more effectively, the pickup times can also be higher. Thus, in contrast to the queueing literature, we must pick $d$ carefully to optimize this trade-off.

\section{Model Description} \label{sec:model1}
We study a deterministic electric vehicle (EV) ride-hailing system that continuously operates throughout the day. Customers arrive at the system at an exogenous rate $\lambda $ requesting trips with a uniformly distributed origin and destination over a bounded 2-dimensional region.\footnote{\Cref{sec:var_arrival_rate} considers time-varying arrivals,  and \Cref{sec: simulations} considers non-uniform origins and destinations.\label{footnote:r2_0_response}} We consider a loss system---an incoming customer leaves if not matched immediately. The system operator fully controls $n$ electric vehicles, each with a battery pack size of $\pack$ [kWh]. We let $r_d$ [kW] denote the discharge rate of the vehicles while driving. Hence, an EV can drive for $\pack/r_d$ hours before running out of charge. There is a total of $m$ chargers in the system uniformly distributed in space, which charge EVs at a rate of $r_c$ [kW]. We assume that $r_c>r_d$, and let $r=r_d/r_c$ denote the fraction of time a vehicle must spend at a charging station for every driven unit of time.\footnote{The charge rate is 10-20 [kW] for Level 2 chargers and more than 50 [kW] for the fast chargers. We consider constant charge and discharge rates that correspond to average values.} EVs have a dynamic State of Charge (SoC) denoted by $\soc_i(t) \in [0, 1]$ with $i\in[n]\triangleq\{1,\dots,n\}$. It represents the current battery level of a given vehicle at time $t$.

A vehicle can be available to serve a rider or not. Available vehicles are in one of three states: idling ($\ci(t)$), driving-to-charger ($\cd(t)$), or charging ($\cc(t)$). An EV that is picking up ($\cp(t)$) or driving with a customer ($\crt(t)$) is not available to serve new requests. The total number of EVs in the system satisfies 
\begin{equation}\label{eq:total_num_vehicles}
n = \cp(t) + \crt(t) +\cd(t) + \ci(t) + \cc(t), \quad \forall t.
\end{equation}
The system operator uses a system management policy to control the available EVs for pickup, the admission of customers, the routing of vehicles to chargers, the admission of EVs at chargers, and the length of charging sessions. We use $\pol$ to denote a policy and $\polSet$ to denote the set of admissible policies. We sometimes use the superscript $\pol$ to stress the dependence of a quantity on a policy.

The system operator's policy influences the evolution of the system and, consequently, determines the effective customer arrival rate, denoted by $\lambdaeff^{\pol}(t)$. For every match, the assigned vehicle picks up the customer and then drives to the customer's destination. Let $q(t)$ denote the total number of customers in the system at time $t$. We then have
\begin{equation}\label{eq:q_identity}
q(t) = \cp(t) + \crt(t),\quad \forall t.
\end{equation}
That is, the number of customers in the system equals the total number of vehicles picking up or driving with a customer. 
We use $T_R$ to denote the average trip time. For a given policy $\pol$, let the average pickup time be $T_P^{\pol}(t)$ at time $t$---the pickup time an average customer would experience at time $t$. We make the following assumption: there exists a constant $\tau_1>0$ such that
\begin{equation}\tag{A1} \label{eq:assumption1}
T_P^{\pol}(t) \geq \frac{\tau_1}{\sqrt{n-q(t)}},\quad \forall t,\quad \forall \: \pol \in \polSet.
\end{equation} 
The right-hand side above approximates the average closest-dispatch pickup time when the number of available vehicles is $n-q(t)$.\footnote{For similar modeling approaches, see \cite{besbes2022spatial, wang2022demand}, and for a quantization theory justification, see \cite{bucklew1982multidimensional}. A formal statement is provided in \Cref{lem:d-close}, Appendix~\ref{app:aux-tech-results}.\label{footnote:assumptions_low_bound}
} 
The assumption captures that for any policy $\pol$ and a given system state, the system operator may decide to (or be limited to) use a fraction of EVs to serve an arriving customer, thereby leading to a larger than the shortest possible average pickup time.

After serving a customer, a vehicle may be directed to a charging station. We let $T_{DC}^{\pol}(t)$ denote the system average time that it would take a vehicle to reach a charger at time $t$.  We assume that there exists a constant $\tau_2>0$ such that
\begin{equation}\tag{A2} \label{eq:assumption2}
T_{DC}^{\pol}(t) \geq \frac{\tau_2}{\sqrt{m-\cc(t)}},\quad \forall t,\quad \forall \: \pol \in \polSet.
\end{equation} 
Note that $m-\cc(t)$ corresponds to the total number of free chargers at time $t$. Hence, the right-hand side above approximates the fastest possible average closest-dispatch travel time to a charger. An arbitrary policy would have a higher average travel time to a charger, as it may route vehicles to chargers that are not necessarily the closest. Once an EV finishes charging, it becomes idle.

\subsection{Objective}\label{sec:model:objective}
The system operator aims to design a control policy 
$\pol \in \Pi$ that achieves a desired long-run service level—defined as the fraction of requested miles served or the fraction of trips completed—while minimizing the number of vehicles and chargers. We define a lower bound for the service level by
\begin{equation*}
\alphainf \overset{\Delta}{=} \frac{1}{\lambda} \liminf_{T\rightarrow \infty} \frac{1}{T}\int_0^T \lambdaeff^{\pol}(t)dt,
\end{equation*}
and we use $\alphasup$ to denote the $\limsup$ version of the quantity above.
Given a target service level $\alpha \in (0, 1)$, we seek a policy $\pi \in \Pi$ under which $\alphainf \geq \alpha$ and the system operates with as few vehicles and chargers as possible. 

In the subsequent analysis, we consider a sequence of systems in steady-state parameterized by $\lambda$ (e.g., \cite{halfin1981heavy}). We derive bounds that determine how the number of vehicles $n_\lambda$ and chargers $m_\lambda$ should scale to attain a certain service level as $\lambda$ increases. We omit the dependence of  $n_\lambda$ and $m_\lambda$ on $\lambda$ when clear from context.

\section{Minimum Requirement of EVs and Chargers}
We derive universal, policy-independent, lower bounds on the number of EVs and chargers in terms of a desired service level $\alpha \in [0, 1]$ and the exogenous arrival rate $\lambda$ for policies $\pol$ such that $\alphasup\geq \alpha$. 

The following proposition establishes first-order lower bounds on the number of vehicles and chargers required to sustain a service level of $\alphasup$, ignoring spatial frictions. 

\begin{proposition}[First Order Lower Bounds]\label{prop:lower-bound-n}
For any policy $\pi$, we have that $n\geq (1+r) T_R\alphasup \lambda$ and $
m\geq rT_R\alphasup \lambda$.
\end{proposition}
Ignoring the pickup time, each customer spends a total of $T_R$ in the system. Thus, at least $T_R\alphasup \lambda$ number of EVs are driving. To restore the SoC used while driving, EVs need $r T_R$ time at a charger (ignoring travel time to the charger). Thus, at least $r T_R \alphasup \lambda$ EVs are charging.
This implies a lower bound of $rT_R\alphasup \lambda$ for the number of chargers and of $T_R\alphasup \lambda + rT_R\alphasup \lambda$ for the number of vehicles. This highlights a fundamental difference between EV and non-EV systems: in an EV system, $\Theta(\lambda)$ extra vehicles are required due to charging needs.

The presence of spatial frictions (pickup times and drive-to-charger times) further increases the fleet size and the number of chargers required and gives rise to a natural trade-off between them. 
To study this trade-off, we analyze how the extra fleet size and number of chargers should scale as the demand increases by considering
 $\frac{n-(1+r)T_R\alpha\lambda}{\lambda^{\beta_1}}$ and $\frac{m- rT_R\alpha\lambda}{\lambda^{\beta_2}}$ (see bounds in \Cref{prop:lower-bound-n}) where $(\beta_1,\beta_2)\in(0,1)^2$ are design parameters. The next theorem characterizes the optimal choice of these parameters.\label{page:scaling-definition}

\begin{theorem}[Universal Lower Bounds]\label{thm:lower-bound}
Fix an $\alpha \in [\delta, 1]$ for some $\delta > 0$ and let $\pol$ be such that $1 \geq \alphasup \geq \alpha$. Then, for all $\lambda > 0$, there exists $\gamma \in [1/3, 1/2]$ such that
\begin{equation} \label{eq: takeaway_n_vs_m_1}
    n \geq (1+r)T_R\alpha\lambda + \Omega\left( \lambda^{1-\gamma}\right) \quad \text{and} \quad
    m \geq rT_R \alpha\lambda +\Omega\left(\lambda^{2\gamma}\right).
\end{equation}
\end{theorem}
\revcolor{To prove the theorem, we identify a set of differential equations that must be satisfied by any policy and then analyze their fixed points. We provide a proof overview in Appendix~\ref{sec:proof_over_thm_lower_bound} and a formal proof in Appendix~\ref{app: lower_bound_thm}.}
\Cref{thm:lower-bound} establishes that for any policy that achieves a service level of $\alpha \in (0, 1)$, the minimum number of vehicles and chargers should follow \eqref{eq: takeaway_n_vs_m_1} (i.e., $\beta_1=1-\gamma$ and $\beta_2=2\gamma$) for some design parameter $\gamma \in [1/3, 1/2]$.

The theorem also provides various insights related to fleet and infrastructure planning. First, the additional capacity requirements, compared to \Cref{prop:lower-bound-n}, emerge due to spatial frictions. The second-order term $\Omega(\lambda^{1-\gamma})$ for EVs corresponds to vehicles driving to pick up customers and reach chargers. The second-order term $\Omega(\lambda^{2\gamma})$ for chargers corresponds to the infrastructure needed to have enough spatial density so that vehicles do not spend too much time reaching a charger. Second, the parameter $\gamma$ captures a natural trade-off between $n$ and $m$. As $\gamma$ increases from 1/3 to 1/2, the fleet size requirement reduces at the expense of installing more chargers: more chargers reduce the drive-to-charger time, implying a lower fleet size requirement. 

Third, \Cref{thm:lower-bound} suggests a capacity planning prescription different from that in the context of non-EVs  \citep{besbes2022spatial}. For non-EVs, the number of vehicles should be at least $T_R \alpha\lambda + \Omega(\lambda^{2/3})$. Our result implies a minimum requirement of 
$ (1+r)T_R\alpha\lambda + \Omega\left( \lambda^{1-\gamma}\right)$ for EVs where $\gamma$ ranges from $1/3$ to $1/2$ depending on the number of chargers. The first-order term is higher due to charging needs; however, the second-order term can be of lower order with $1-\gamma \in [1/2, 2/3]$. The increase of the first-order term and potential reduction of the second-order term are an inherent byproduct of the physics of an EV system. If the fleet’s average SoC is high enough, the $\Theta(\lambda)$ charging vehicles can also serve arrivals, increasing available vehicles for pickups from $\Theta(\lambda^{2/3})$ in a non-EV system to $\Theta(\lambda)$ in an EV system. This boost in capacity increases EV density, reducing pickup times and thus the number of vehicles doing pickups. Importantly, these benefits depend on charger density: when chargers are sparse (low $\gamma$), drive-to-charger time dominates pickup time. However, the benefits of adding chargers eventually plateau at $\gamma=1/2$ at which point pickup times begin to dominate.

\Cref{thm:lower-bound} provides a best-case result, but it does not guarantee the existence of a policy that can attain a desired service level while minimizing the fleet and chargers, as in \eqref{eq: takeaway_n_vs_m_1}. In \Cref{sec:uppder_bound}, we prove that these bounds are near achievable.

\section{Achieving Optimal Performance}\label{sec:uppder_bound}

We analyze the system under the {\it Power-of-$d$ vehicles dispatch} policy and show that it nearly achieves the lower bounds in \Cref{thm:lower-bound}. We first introduce the policy, present the main theorem of this section, and discuss its implications. We then explain how to model and analyze the system's performance under this policy. 

The Power-of-$d$ dispatch policy selects $d$ closest vehicles and assigns the request to the one with the highest SoC (\Cref{alg:powerd}). Small $d$ ensures short pickup times but may lead to fleet SoC imbalance, while large $d$ yields longer pickup times but more balanced fleet SoC distribution. The parameter $d$ allows us to navigate this trade-off.

\begin{algorithm}[t]
\caption{Power-of-$d$ Vehicles Dispatch} \label{alg:powerd}
{\fontsize{10}{13}\selectfont
\SetAlgoLined
\SetKwInOut{Input}{Input}
\SetKwInOut{Initialize}{Initialization}
\Input{Fix $d \geq 1$}
\For{Every Arrival}{
\If{All vehicles are busy}{
The arrival is dropped\;
}
Choose the $d$ closest vehicles that are not driving\;
\If{None of the $d$ vehicles have enough SoC for the trip}{
The arrival is dropped\;
}
Match the arrival to the vehicle with the highest SoC among the $d$ closest\;
Once the trip ends, send the vehicle to the closest available charger\;
The vehicle becomes available for matching after it has reached a charger\;
}
}
\end{algorithm}

Analyzing the system under this policy requires a large state space and tracking of continuous SoC changes of each EV. Instead, we can define a modified yet tractable system whose analysis, under the Power-of-$d$ policy, implies the existence of a policy in the original system that nearly achieves the lower bound of \Cref{thm:lower-bound}.
To state our second main result, we define $\Nt$ as the number of trips an EV can complete with a 100\% SoC in the modified system (see \Cref{sec:thm2_state_space}).\label{page:r2_res_N_T}

\begin{theorem}[Near-Matching Upper Bound]\label{thm:informal_upper-bound}
Assume that the dynamical system in \eqref{eq:odes-power-d} for the modified system exhibits global stability. Fix an $\alpha \in [\delta, 1)$ for some $\delta > 0$ and $\gamma \in \left[\frac{1}{3}, \frac{1}{2+1/\Nt}\right)$ with $\Nt \geq 2$. Then, there exists a policy $\pi$ in the original system and $\lambda_0(\alpha)$ and $m, n$ with
\iffalse
Assume that the associated dynamical system exhibits global stability.\footnote{We present the dynamical system that describes the evolution of the system under power-of-$d$ vehicles dispatch in \Cref{sec:power_of_d_modeling}.}
Fix an $\alpha \in [\delta, 1)$ for some $\delta > 0$ and $\gamma \in \left[\frac{1}{3}, \frac{1}{2+1/\Nt}\right)$ with $\Nt \geq 2$. Then, there exists a policy $\pi$ and $\lambda_0(\alpha)$ and $m, n$ with
\fi
\begin{align*}
n = (1+r)T_R \alpha\lambda + \Theta\left(\lambda^{1-\gamma}\right), \quad \text{and,} \quad  m = rT_R \alpha\lambda + \Theta\left( \lambda^{2\gamma}\right),
\end{align*}
such that $\alphasup \geq \alphainf \geq \alpha,\: \forall \lambda \geq \lambda_0(\alpha)$. 
\end{theorem}
We defer the proof of the theorem and a discussion about the stability assumption to Appendix~\ref{sec:pod-analysis} (the assumption is satisfied in all our simulations in \Cref{sim:synthetic_data}). The result in Theorem \ref{thm:informal_upper-bound} nearly matches the lower bound in Theorem \ref{thm:lower-bound}, and the gap closes as the battery pack size increases.

The theorem also shows that the Power-of-$d$ policy effectively leverages the extra capacity introduced by charging. Serving a request with an EV that is both nearby and has a high SoC accomplishes two objectives: (i) it keeps pickup times short, which in turn supports vehicle availability; and (ii) it maintains a more uniform SoC across the fleet, preserving access to the extra capacity. Additionally, the analysis of the 
modified system under the Power-of-$d$ policy (see \Cref{lemma: coupling}) reveals insights into how the policy should be adjusted with the charger density and the battery pack size. The optimal choice of $d$ is $\tilde{\Theta}\left(\lambda^{\gamma / \Nt}\right)$.\footnote{A function $f(n)$ is $\tilde{\Theta}(g(n))$ if and only if there exists constants $C_1,C_2>0$ and $k\in(0,1]$  and $n_0>0$ such that for all $n\geq n_0$, we have $C_1\cdot g(n)\leq f(n) \leq C_2 \cdot g(n) \cdot \log(n)^k$.\label{rev:footnote_polylog}
} As the charger density increases (larger $\gamma$), the fleet size decreases. This smaller fleet serving the same arrival rate leads to a lower average SoC, requiring larger $d$ to find vehicles with sufficient charge. A higher $\Nt$ reduces the SoC constraint, allowing smaller $d$.

\begin{remark}\label{remark-constants}
We note that the theorem does not provide closed forms for the constants in front of $\lambda^{1-\gamma}$ and $\lambda^{2\gamma}$. \Cref{eq: in_terms_of_n} in the proof offers explicit values, but finding the optimal constants requires a more detailed analysis. Appendix~\ref{sec:sec-order-const-pod} outlines a procedure to approximate these values.
\end{remark}

\subsection{Power-of-\texorpdfstring{$d$}{d} Vehicles Dispatch: Modeling}\label{sec:power_of_d_modeling}

We model the system under Power-of-$d$ vehicles dispatch by first approximating the time a vehicle spends serving each customer request (\Cref{sec:busy_times_pod}). We then bound this time and use this bound to propose a modified system describing the fleet's SoC evolution (\Cref{sec:thm2_state_space,sec:trans_prob}). In \Cref{sec: dynamical_system_pod}, we present the associated dynamical system whose analysis establishes \Cref{thm:informal_upper-bound}.

\subsubsection{Vehicles' Busy Times}\label{sec:busy_times_pod}
In the Power-of-$d$ vehicles dispatch algorithm, an EV is dispatched to a charger right after finishing a trip. Hence, for each customer served, the corresponding EV is busy for the length 
of the drive-to-charger time, the pickup time, and the travel time.

Since there are $m$ chargers and a total of $\cc(t)$ vehicles charging, $m-\cc(t)$ chargers are available. We approximate $T_{DC}(t)$ by
\begin{equation}\label{eq: drive_to_charger_power_of_d}
T_{DC}(t) = \frac{\tau_2}{\sqrt{m-\cc(t)}}.
\end{equation}
The expression above is proportional to the average distance to the closest charger, assuming that there are many empty chargers uniformly distributed in space. 

Note that for $k$ uniformly distributed in space points, the  $d^{\text{th}}$ closest point to a given point scales as  $\Theta\left(\sqrt{\frac{d}{k}}\right)$ (cf. \Cref{lem:d-close} in Appendix~\ref{app:aux-tech-results}). We then approximation $T_P(t)$ by
\begin{equation}\label{eq: pickup_power_of_d}
T_P(t) = \tau_1 \sqrt{\frac{d}{n-(q(t) + \cd(t))}}.
\end{equation}
Recall that $n$ is the total number of vehicles in the system and $q(t) + \cd(t)$ vehicles are busy at time $t$, i.e., driving with a customer or driving to a charger. Hence, $n-(q(t) + \cd(t))$ vehicles are available to serve incoming customers. Note that $T_P(t)$ is consistent with \eqref{eq:assumption1} as $d\geq 1$. The rate $\mu(t)$ at which vehicles complete serving a customer and driving to a charger is then given by 
\begin{equation}\label{eq:power-d-service-rate}
\frac{1}{\mu(t)} = \tau_1 \sqrt{\frac{d}{n-(q(t) + \cd(t))}} +T_R+\frac{\tau_2}{\sqrt{m-\cc(t)}}.
\end{equation}
\subsubsection{State Space}\label{sec:thm2_state_space}
Our policy must identify the $d$ closest vehicles for each trip request and select the one with the highest SoC. Therefore, tracking the system state requires a detailed record of the fleet's SoC. Since monitoring individual vehicles' continuous SoC changes is intractable, we instead track the number of vehicles able to fulfill a specific number of trip requests at any time.
We accomplish this in three steps. We first introduce admission control to obtain an upper bound on the time a vehicle spends driving per customer, i.e., we upper bound $1/\mu(t)$. In turn, we obtain an upper bound on SoC spent per customer served, i.e., $r_d/\mu(t)$. Then, we introduce a modified system in which the SoC spent per trip is given by this upper bound, enabling us to discretize the SoC of each EV. Finally, our state space representation relies on this discretization to count the number of vehicles idling, charging, or busy that can fulfill a given number of trip requests. 

\begin{definition}[Admission Control] \label{def: admission_control}
    A trip request is not admitted if the number of busy cars, $q(t) + \cd(t)$, exceeds $\alpha T_R\lambda + \kappa\lambda^{1-\gamma}$ for some $\kappa>0$. A maximum of $A \overset{\Delta}{=} m - \lambda^{2\gamma}$ vehicles are allowed to charge at a time. Vehicles with lower SoC are given precedence to charge.
\end{definition}

\begin{remark}\label{remark:new_admission_control_discussion}
The choice of $\kappa$ follows from the analysis of \Cref{thm:informal_upper-bound} in Appendix~\ref{sec:pod-analysis}. While \Cref{def: admission_control} makes the system amenable to analysis, in our simulations in \Cref{sec: simulations}, we either omit admission control or impose a maximum allowable pickup time.
\label{res:new_admission_control_discussion}
\end{remark}

\Cref{def: admission_control} implies the following bound for the time an EV spends serving a request:
$$
\frac{1}{\mu(t)} \leq \tau_1 \sqrt{\frac{d}{n - \alpha T_R\lambda-\kappa\lambda^{1-\gamma}}} +T_R+\frac{\tau_2}{\sqrt{m-A}} \overset{\Delta}{=}T_B.
$$ 

For the remainder of the \Cref{sec:uppder_bound}, we consider a modified system in which each vehicle spends $\Delta\triangleq T_B r_d$ [kWh] of battery for each trip request (pickup and trip, followed by driving-to-charger). 
Since the modified system spends more SoC serving a trip request, any admissible decisions in the modified system are also admissible in the original system and, therefore, the bounds we derive in the former also hold in the latter. We make this observation precise in the following lemma.
\begin{lemma} \label{lemma: coupling}
    For a given $\lambda$, there exists a policy $\pi \in \Pi$ for the original system that serves the same sequence of arrivals as in the modified system. Thus, their effective arrival rates are identical.
\end{lemma}
We are ready to introduce the state space description. In the modified system, the state space can be reduced to a discrete set as the SoC of all EVs is an integer multiple of $\Delta$. Indeed, note that since the pack size of a vehicle is $\pack$, a fully charged vehicle can serve a total of $\Nt\triangleq \lfloor \pack/\Delta \rfloor$ trip requests. For $j\in\{0,\dots,\Nt\}$, we define: 
\begin{flalign*}
X^C_j(t) &\triangleq \text{\# of vehicles idling or charging that can serve no more than $j$ trip requests at time $t$,}\\
X^B_j(t) &\triangleq \text{\# of vehicles busy that can serve no more than $j$ trip requests at time $t$.}
\end{flalign*}
These state variables aggregate vehicles either charging/idling or busy that can serve a certain number of trips or fewer, but no more. Note that $X^C_{\Nt}(t)$ includes all vehicles that are idling and charging, and $
\XC_{\Nt}(t) + \XB_{\Nt}(t) = n
$ for all $t \in \bbR_+$.

\subsubsection{Transition Probabilities}\label{sec:trans_prob}
The evolution of $\XC_{j}(t)$ and  $\XB_{j}(t)$ depends on how the policy assigns vehicles with different SoC levels to trip requests. This assignment, in turn, depends on the SoC distribution of the fleet. For example, if vehicles that can serve only a small number of trip requests form a larger portion of the fleet, the policy will likely choose one of them to serve a trip request. To determine the assignment, we consider an approximation (validated by our numerical simulations in \Cref{sec: simulations}) that ignores the temporal correlations between consecutive trip requests: we treat the $d$ closest vehicles as if they were randomly sampled. Hence, their SoCs are as if they were randomly sampled without replacement from the pool of vehicles. Define the event 
$$
A_j(t) \triangleq \{\text{selecting at least one vehicle that can serve $j+1$ or more trips}\},\quad \forall j\in\{0,\dots,\Nt-1\}.
$$
Then, we have that
\begin{equation}\label{eq:prob-without-replacement}
p_j(t)\triangleq \mathbb{P}(A_j(t)) = 1 -\frac{\XC_{j}(t)}{\XC_{\Nt}(t)}\cdot \frac{\XC_{j}(t)-1}{\XC_{\Nt}(t)-1}\cdots\frac{\XC_{j}(t)-(d-1)}{\XC_{\Nt}(t)-(d-1)},\quad \forall j\in\{0,\dots,\Nt-1\},
\end{equation}
is the complement probability of choosing $d$ vehicles, uniformly at random without replacement, that can serve $j$ (but no more than $j$) or fewer trips. Intuitively, $p_j(t)$ determines the rate at which vehicles with enough SoC for $j+1$ or more trips start serving a customer. \revcolor{The expression above assumes that $\XC_j(t)\geq d-1$. For the case of $\XC_j(t) < d-1$,} we simply have $p_j(t) = 1$.

The expressions in  \eqref{eq:power-d-service-rate} and \eqref{eq:prob-without-replacement} elucidate the trade-off that the Power-of-$d$ vehicles dispatch policy balances. For a given state, as $d$ increases, the service time in the system increases because pickup times become larger. However, the probability of choosing at least one vehicle that can serve a certain number of trips increases with $d$. By varying $d$, we can induce small (large) pickup times but a low (high) likelihood of matching with a vehicle with high SoC.

\subsection{Power-of-\texorpdfstring{$d$}{d} Vehicles Dispatch: Dynamical System}\label{sec: dynamical_system_pod}
We are ready to introduce the differential equations describing the system's evolution under the Power-of-$d$ vehicles dispatch policy. Define the admission control variables $I(t) \triangleq \1{\XB_{\Nt}(t)\leq \alpha T_R \lambda + \kappa \lambda^{1-\gamma}}$ and $\tXC_j \triangleq \min\{\XC_j, A\}$. We have: 
%\vspace{-1cm}
\begin{subequations} \label{eq:odes-power-d}
\begin{flalign}\label{eq:powerd-1}
\dXB_j(t)&=\lambda \left(p_0(t)-p_j(t)\right) I(t) - \XB_{j}(t) \mu(t), \: 1\le j< \Nt,\\\label{eq:powerd-2}
\dXC_j(t)&=-\lambda \left(p_0(t)-p_j(t)\right) I(t) - \left(\tXC_j(t)-\tXC_{j-1}(t)\right) \frac{1}{rT_B} +\XB_{j+1}(t) \mu(t),\: 1\le j<\Nt,\\\label{eq:powerd-3}
\dXC_{\Nt}(t)&=-\lambda p_0(t) I(t) + \XB_{\Nt}(t) \mu(t), \\\label{eq:powerd-4}
\dXC_{0}(t)&=- \tXC_0(t)\frac{1}{rT_B} +\XB_{1}(t) \mu(t),
\end{flalign}
\end{subequations} 

The system state changes through three events: customer pickup, service completion, and charging progress. We next describe these transitions and their relation to \eqref{eq:odes-power-d}.

\emph{Pickup:} For $j \in \{1, \hdots, \Nt\}$, $\XC_j$ reduces by one, and $\XB_j$ increases by one if an EV with SoC at most $j$ picks up an incoming customer. %This transition corresponds to the down arrows in \Cref{fig: ODE}. 
The probability of this event is $p_0(t) - p_j(t)$. As the arrival rate of admitted customers is $\lambda I(t)$, $\XC_j$ reduces and, $\XB_j$ increases with rate $\lambda I(t)(p_0(t) - p_j(t)).$

\emph{Finish Service:} For $j \in \{1, \hdots, \Nt\}$, $\XB_j$ reduces by one, and $\XC_{j-1}$ increases by one if an EV with SoC at most $j$ finishes a trip request (pickup, travel to a destination, and drive-to-charger). As a total of $\XB_j$ EVs are finishing service with rate $\mu(t)$, these transitions happen with rate $\XB_j(t) \mu(t).$

\emph{Charge:} For $j \in \{0, \hdots, \Nt-1\}$, $\XC_j$ reduces by one if an EV with SoC equal to $j$ gains one unit of charge. Note that $r_d T_B$ corresponds to the total SoC in [kWh] that a vehicle expends during a trip. Hence, while charging, it takes $(r_d T_B)/r_c$ hours to gain the SoC required for a trip. Since there are $\tXC_j(t)- \tXC_{j-1}(t)$ charging vehicles that have SoC exactly $j$, $\XC_j(t)$ decreases at rate $(\tXC_j(t)-\tXC_{j-1}(t))/(rT_B).$

We obtain \eqref{eq:powerd-1} by combining the pickup and service rates. Similarly, we get \eqref{eq:powerd-2} by combining the pickup, service, and charge rate. Equations \eqref{eq:powerd-3} and \eqref{eq:powerd-4} are similar to Equations \eqref{eq:powerd-1} and \eqref{eq:powerd-2}, with the difference that fully charged vehicles do not gain more SoC, and fully discharged ones cannot serve incoming customers.

\section{Varying Arrival Rate: from EV to non-EV Infrastructure Planning}\label{sec:var_arrival_rate}
We consider the case of a time-varying arrival rate. We first present universal lower bounds that echo \Cref{prop:lower-bound-n} in our base model. We then discuss their implications and how they help us understand when EV and non-EV infrastructure planning apply.

Consider an arrival rate $\lambda(t)$ that alternates between $\lambda$ for valleys lasting $\Tvalley$ and $\amplitude\lambda$ for peaks lasting $\Tpeak$, with $\amplitude > 1$. The average arrival rate is $\lambda_{\operatorname{avg}} = \frac{\Tvalley + \amplitude\Tpeak}{\Tvalley+\Tpeak} \lambda$ and the goal is to meet $\alpha$ fraction of the demand.
Let $\alpha_v$ and $\alpha_p$ be the fraction of the demand served during the valley and peak:
\begin{align*}
\alpha_v = \frac{\lambda \Tvalley - (\lambda- \alpha\lambda_{\operatorname{avg}})^+\Tvalley}{\lambda \Tvalley} \quad  \text{and} \quad
\alpha_p =  \frac{\alpha\lambda_{\operatorname{avg}}\Tpeak + (\alpha\lambda_{\operatorname{avg}} -\alpha_v\lambda)^+\Tvalley}{\amplitude\lambda \Tpeak}, \numberthis \label{eq: peak_and_valley_service_level}
\end{align*}
where $(\lambda- \alpha\lambda_{\operatorname{avg}})^+\Tvalley$ is the demand excess over the average target demand during the valley, and $(\alpha\lambda_{\operatorname{avg}} -\alpha_v\lambda)^+\Tvalley$ is the excess of demand not served during the valley. 

We define the number of EVs driving during the peak ($n_{\operatorname{dp}}$), charging during the peak ($n_{\operatorname{cp},\eta}$) and charging during the valley ($n_{\operatorname{cv},\eta}$) by ($\eta$ is a parameter to be explained below)
\begin{align*}
n_{\operatorname{dp}} =  \alpha_{p}\amplitude\lambda T_R,\quad 
n_{\operatorname{cp}, \eta} = r\alpha \lambda_{\operatorname{avg}} T_R - \eta\frac{\Tvalley}{\Tpeak} \lambda_{\operatorname{avg}} T_R , \quad n_{\operatorname{cv}, \eta} = r\alpha \lambda_{\operatorname{avg}} T_R + \eta \lambda_{\operatorname{avg}} T_R,\quad \eta\geq 0.
\end{align*}
Since we meet a fraction $\alpha_p$ of the demand during the peak,  $\alpha_{p}\amplitude\lambda T_R$ vehicles drive during this period. Overall, we aim
to meet $\alpha$ fraction of the demand so $\alpha \lambda_{\operatorname{avg}} T_R$ vehicles drive and, thus, $r \alpha \lambda_{\operatorname{avg}} T_R$ charge on average. %to sustain the system operations. 
Charging can happen at different times of the day, e.g., it can be uniformly distributed during the peaks and valleys ($\eta= 0$), or it can be more concentrated during the valleys ($\eta>0$). Thus $\eta$ controls the allocation of charging throughout a day, and it must satisfy: $\frac{\Tvalley n_{\operatorname{cv}, \eta} + \Tpeak n_{\operatorname{cp}, \eta}}{\Tvalley + \Tpeak} = r \alpha \lambda_{\operatorname{avg}} T_R$.

Under any matching and charging policy, the system should operate with at least the peak vehicle and valley charging requirements: $n_{\operatorname{dp}} + n_{\operatorname{cp},\eta}$ and $n_{\operatorname{cv},\eta}$, respectively.
\begin{theorem}[First Order Lower Bounds] \label{thm: informal_varying_arrival_rate}
    Fix an $\alpha \in (0, 1)$,  let the arrival rate be given by $\lambda(t)$ with $r\Tpeak/\Tvalley < 1$.\footnote{
    Under the less realistic case of $r\Tpeak/\Tvalley \geq 1$, we can show that only Cases I and II exist. 
    } Let $\pi$ be such that $\alphasup \geq \alpha$ 
    then there exists $\eta\geq 0 $ and $w\geq0$ such that
    \begin{equation*}
            n \geq n_{\operatorname{dp}} + n_{\operatorname{cp}, \eta} - w\cdot \frac{T_R^2 \amplitude \lambda}{\Tpeak} \quad \text{and}\quad  m \geq n_{\operatorname{cv}, \eta}.
    \end{equation*}
    \begin{enumerate}
        \item [I.) ] If $\alpha\lbavg \leq \lambda$, then $\eta=0$ and $w=0$.
    \item [II.) ] If $\alpha\lbavg \in \left(\lambda, \frac{\lambda}{\left(1-r\Tpeak/\Tvalley\right)} \right]$, then  $\eta \in \left[0, \alpha - \frac{\lambda}{\lambda_{\operatorname{avg}}}\left(1+\frac{\amplitude T_R}{\Tvalley}\right) \right]$ and $w=1$.
    \item[III.)] If $\alpha\lbavg > \frac{\lambda}{\left(1-r\Tpeak/\Tvalley\right)}$, then  $\eta \in \left[0, \frac{\alpha r\Tpeak}{\Tvalley}-\frac{\amplitude \lambda T_R}{\Tvalley \lambda_{\operatorname{avg}}}\right]$ and $w=1$.
    \end{enumerate}
\end{theorem}
The theorem characterizes the minimum infrastructure requirement under $\lambda(t)$: the system needs vehicles equal to the sum of those charging and driving during the peak and chargers equal to the number of vehicles charging during the valley. The parameter $\eta$, determined by the policy in place, modulates the trade-off between the minimum number of EVs and chargers. As $\eta$ increases, fewer vehicles must charge during the peak but more must charge during the valley. The extent to which charging can be allocated to the valley while maintaining minimum infrastructure requirements depends on the peak's amplitude (the three cases in the theorem).

Next, we fix the service level $\alpha \in (0, 1)$ and peak-to-valley duration ratio $\Tpeak/\Tvalley$, and examine the three cases in \Cref{thm: informal_varying_arrival_rate}, discussing their implications for infrastructure planning as the peak amplitude $\amplitude$ varies. For low to medium amplitudes (I and II), the EV scaling and control (Power-of-$d$) are more appropriate. For large amplitudes (III), with high battery pack size and charger numbers, the non-EV scaling may be suitable, though this still requires EV-aware control to balance the fleet SoC. We explore this further in the subsequent sections and provide a summary in \Cref{tab:summary_of_time_varying}.

\begin{remark}
    We note that the terms $\frac{cT_R^2\lambda}{\Tpeak}$ and $\frac{cT_R}{\Tvalley}$ in the lower bound for $n$ and cases II and III are an edge effect induced by transitions between peaks and valleys.\footnote{ The number of EVs driving at the start of a valley will gradually decrease. A policy aware of this can potentially benefit, e.g., by admitting more customers at the end of the peak, resulting in a lower fleet size requirement. } 
These effects are negligible---and can be ignored---if peak and valley durations are much larger than the trip time.
\end{remark}

\begin{table}[]
    \centering
    \TABLE{Varying arrival rate scalings and system control.
    \label{tab:summary_of_time_varying}}{
    \begin{tabular}{|c||c|c|} \hline
      Peak Amplitude ($\amplitude$)  & \multicolumn{2}{c|}{System Behavior}  \\ \hline \hline
        Low & \multicolumn{2}{c|}{EV-scaling: Constant arrival rate results holds} \\ \hline
        Medium & \multicolumn{2}{c|}{EV-scaling: $\Theta(\lambda)$ EVs charging $\forall t$}\\ \hline
        \hline
     \multirow{2}{*}{High} &  \underline{Low $p^\star$ \textit{and/or} chargers} & \underline{High $p^\star$ \textit{and} chargers} \\ 
         & EV-scaling: $\Theta(\lambda)$ EVs charging $\forall t$ & Non-EV-scaling + EV-aware control \\ \hline
    \end{tabular}}{}
\end{table}

\subsection{EV Scaling for Small and Medium Demand Variability: Case I and II}

In Case I, the low peak amplitude implies a high valley arrival rate $\lambda$ compared to the average arrival rate $\lambda_{\operatorname{avg}}$. Here, fleet and charging requirements can be minimized by maintaining a constant effective arrival rate, disregarding demand peaks while serving at least an $\alpha$ fraction of the demand.

In Case II, the peak amplitude is large enough to require the system operator to meet a non-zero fraction of additional peak demand. However, the amplitude is not so large that the difference in vehicles driving during the peak and the valley is significant. Consequently, shifting too much charging to the valley would result in the number of EVs in the valley (driving and charging) dominating the number of EVs in peak and, in turn, in a total higher requirement of both vehicles and chargers. 
With $\eta\leq \alpha - \frac{\lambda}{\lbavg}(1+cT_R/\Tvalley)$, not all charging can be shifted to the valley, i.e., $n_{cp,\eta}>0$ for this range of $\eta$. 

In both cases, a non-zero fraction of EVs are charging, so $\Omega(\lambda)$ partially charged EVs are available to serve customers at all times. As in \Cref{thm:lower-bound}, these partially charged EVs can reduce the pickup times, and a similar EV-based scaling would be applicable for small to medium peak amplitudes.

\subsection{Large Demand Variability: Case III}\label{sec:case3}
In Case III, the peak amplitude is so large that one could shift all the charging to the valleys, and the entire fleet is employed to meet the peak demand. This is only possible if we have a sufficiently high number of chargers and EVs with a large battery pack size. The number of chargers has to be larger than $n_{\operatorname{cv}, \eta^*}$ where $\eta^*$ is such that no vehicles charge during the peak, i.e., $n_{\operatorname{cp}, \eta^*}=0$. The battery pack size has to be larger than the [kWh] spent during the peak, i.e., $\pack\geq r_d \Tpeak$. If either of these conditions is not satisfied, there will be $\Omega(\lambda)$
partially charged EVs available for dispatch and a similar EV-based scaling as in \Cref{thm:lower-bound} applies (see Appendix~\ref{sec:lim_char_lim_pack_size} for more details).

We next consider a setting with a sufficiently large battery pack size, $\pack \geq r_d \Tpeak$, and, for simplicity, an unlimited number of chargers (i.e., $m = \infty$). Since each EV has enough battery capacity to serve the peak demand, one might expect that an EV-agnostic policy that performs well in the non-EV setting would also perform well here. 
We demonstrate that this is not necessarily true. To make this point concrete, we first discuss the implications of a very simple EV-agnostic policy---Closest Dispatch ($d=1$). We then explain that an EV-aware policy that attempts to balance the load in the system can deliver performance meeting the minimal requirement in \Cref{thm: informal_varying_arrival_rate}, case III.

\subsubsection{The Limits of an EV-Agnostic Control}\label{sec:limits-cd-varying}
We start by bounding the service level and the number of vehicles with zero SoC under the Closest Dispatch policy in a constant arrival rate scenario. We then discuss the implications for varying arrival rates.

\begin{proposition} \label{prop: cd_varying_arrival_rate}
    \revcolor{Let the fleet size satisfy $n\in [(1+r)\lambda T_R,M\lambda T_R]$ for $M > 1+r$, and assume that the system under Closest Dispatch is globally stable. Then, there is a non-vanishing number of EVs with zero SoC ($\bXC_0 =\Omega (\lambda)$) and for large $\lambda$, the service level satisfies $\alphasup < 1$.}
\end{proposition}

The proof and additional details about the ODEs governing the system's evolution under Closest Dispatch are presented in Appendix~\ref{app: cd_varying_arrival_rate}.
The proposition uncovers the limits of this policy. Let us use
\Cref{prop: cd_varying_arrival_rate} to approximate the performance during the valley, assuming that we reached 
steady-state. We have $\Omega(\lambda)$ EVs with zero SoC at the end of the valley and, therefore, cannot serve any trips at the start of the peak. This implies that under Closest Dispatch, the system would need to operate with more EVs because it fails to meet 100\% of the demand during the valley $(\alphasup < 1)$, but also because it will have to charge some EVs during the peak due to poor load balancing. However, it is possible to verify (see \eqref{eq: alpha_eff_fixed_point_proof} in Appendix~\ref{app: cd_varying_arrival_rate}) that (i) as $\pack$ increases, the valley service level approaches 1, implying that larger battery packs reduce fleet size requirements, and (ii) higher demand amplitude accelerates this convergence, further decreasing the fleet size requirement. Thus, while Closest Dispatch requires more vehicles, its limitations lessen for fleets with larger battery packs or when demand variability is high (\Cref{sec:sim_chicago} provides empirical support).

\subsubsection{EV-Aware Control}\label{sec:pod-works-adaptive}
We first show that in a constant arrival rate setting with Power-of-$d$, in steady-state, almost all vehicles have a full charge if there are an extra $\Theta(\lambda)$ vehicles. We then use this result to argue that an EV-aware control policy is needed.

\begin{proposition} \label{prop: pod_var}
Suppose the fleet size is given by $n = (1 + r)\lambda T_R + r b \lambda T_R$ for some $b > 0$. Then, there exists a value $\lambda_v > 0$ and an admission control such that for all $\lambda \geq \lambda_v$, the Power-of-$d$ policy with
$d = (\log \lambda)^2$ has the following property: all fixed points of \eqref{eq:odes-power-d} fully meet demand ($p_0 I = 1$), and nearly all EVs remain close to full charge ($\XC_{\Nt-2} + \XB_{\Nt-2} \leq 1$).
\end{proposition}

The proof of the proposition is deferred to Appendix~\ref{app: pod_var}. The proposition suggests that the minimum infrastructure requirement from \Cref{thm: informal_varying_arrival_rate} case III can be met by employing an {\it adaptive} Power-of-$d$ policy. To see this, note that meeting $\alpha_p$ fraction of the demand during the peak without charging implies that only $\alpha_p a \lambda T_R (1-rT_p/T_v)$ EVs can be used to meet all of the demand during the valley.\footnote{$r_d T_p \alpha_p a \lambda T_R$ [kWh] are spent during the peak, so an average of $\frac{r_d T_p \alpha_p a \lambda T_R}{r_cT_v}$ EVs need to be charging during the entire valley to recover the lost SoC. %Thus, only $\alpha_p a \lambda T_R (1-rT_p/T_v)$ EVs can be used to meet all of the demand during the valley. 
} 
One can verify that (see Appendix~\ref{app: pod_var} for details)
\begin{align*}
    \alpha_p a \lambda T_R \left(1-\frac{rT_p}{T_v}\right) = (1+r)\lambda T_R + r b\lambda T_R,\:\: b>0. \numberthis \label{eq: fleet_pod_varying}
\end{align*}
Hence, if we approximate the performance in the valley by the steady-state then by \Cref{prop: pod_var}, we conclude that Power-of-$d$ with $d = (\log \lambda)^2$ is able to meet 100\% of the valley demand as $p_0 I = 1$. \Cref{prop: pod_var} also establishes that $n-1$ EVs would have sufficient SoC to serve $\Nt-1$ trips, i.e., almost all EVs have almost 100\% SoC by the end of the valley and are ready to serve trips during the next peak exclusively. Thus, Closest Dispatch can be employed during peak periods. Since the peak fleet size determines the total fleet requirement and Closest Dispatch is feasible, the system can operate like a non-EV system during peaks, making the non-EV scaling of \cite{besbes2022spatial} relevant. Importantly, policy adaptability is needed to balance workload across the fleet, which helps increase the SoC during valleys and avoids charging during peaks.

\section{Simulations} \label{sec: simulations}
In \Cref{sim:synthetic_data}, we simulate our baseline model and demonstrate the robustness of our theoretical scalings (\Cref{thm:informal_upper-bound}), their capacity planning implications, and we also test the assumptions on pickup and drive-to-charger times. In \Cref{sec:sim_chicago}, we use Chicago trip request data to demonstrate how the insights from the baseline model can translate into more complex settings. 

\subsection{Synthetic Data Study}\label{sim:synthetic_data}
We generate trip requests according to a Poisson process with rate $\lambda$, where the starting and ending locations are uniformly distributed in a 10-mile by 10-mile region, and $\lambda \in \left\{5, 10, 20, 40, 80, 160, 320\right\}$ arrivals per min. Chargers are set uniformly at random in the squared region with $m_p$ charging ports, and there are $m_l = \lfloor m / m_p \rfloor$ different locations of chargers where $m_p= 8$  \citep{no_of_posts}. Our simulations run for 1000 minutes, and we calculate outcomes like service level by considering only the second half of the simulation to ensure the system is in steady-state. We report the mean of outcomes over five randomly generated simulation runs.

For the system's evolution, we consider charge and discharge rates of 20[kW] and 5[kW], respectively (level 2 charging). Vehicles' average velocity is set to 20 [mi/hr] while their battery pack size is set to 40[kWh] (160-mile range). For each trip request, if the selected EV has sufficient battery to serve the trip and at least 20\% of SoC, then the customer is served. Otherwise, the customer leaves the system immediately. Once the trip ends, the EV becomes idle, and if its SoC is below 90\%,  the EV is dispatched to the closest available charger. While an EV is driving to a charger, another EV can reach the same charger and start charging earlier; in this case, the EV waits at the charger until a post becomes available. We also allow interrupting charging for dispatching, but if the charging is not interrupted and the SoC reaches 100\%, the EV becomes idle.

We provide an empirical evaluation of the theoretical results from \Cref{thm:informal_upper-bound} in \Cref{tab: asymptotic_sim}. 
Under the Power-of-$d$ policy, for a given $\beta$, we set $m = r T_R \alpha \lambda + \Theta(\lambda^\beta)$ and infer the number of EVs required to achieve 90\% service level via simulation (see Appendix~\ref{sec:app-sim-details}). The table
reports the fitted and theoretical scalings of the fleet, pickup time, and drive-to-charger time. In particular,
$$
\gamma = \min\left\{\frac{1}{2+1/\Nt}, \frac{\beta}{2}\right\} \implies 1-\gamma = \max\left\{0.508, 1-\frac{\beta}{2}\right\},\footnote{$\Nt = p^\star / (r_d \times T_B) \approx p^\star / (r_d \times T_R)$  = 32.}
$$
and the pickup time and the drive-to-charger time are approximately
$$
T_P \approx \frac{\tau_1\sqrt{d(1+r)}}{\sqrt{rn}} \sim n^{\gamma/(2\Nt) - 1/2} = n^{\gamma/64 - 1/2}, \quad \text{and}\quad
T_C \approx \frac{\tau_2}{\sqrt{m-A}} \sim n^{-\beta/2}.
$$

\begin{table}[bth!]
    \centering
    \TABLE{\revcolor{Theoretical versus empirical $1-\gamma$, average pickup time fit, and average drive-to-charger time.}% fit as a function of arrival rate.
    \label{tab: asymptotic_sim}}{
    \begin{tabular}{|c||c||c|c|c||c|c|c||c|c|c||} \hline
        Series & \multicolumn{1}{c||}{$\beta$} & \multicolumn{3}{c||}{$1-\gamma$} & \multicolumn{3}{c||}{Pickup time Scaling} &\multicolumn{3}{c||}{Drive-To-Charger time Scaling} \\ \hline \hline
         &  & Fit & Theo & Error \% & Fit & Theo & Error \% & Fit & Theo  & Error \% \\ \hline
         A & 1 & 0.525 & 0.508 & 3.3\% & -0.474 & -0.492 & 3.7\% & -0.496 & -0.5  & 0.8\% \\ \hline
         B & 0.906 & 0.546 & 0.547 & 0.2\% & -0.461 & -0.493 & 6.5\% & -0.458 & -0.453 & 1.1\% \\ \hline
         C & 0.803 & 0.577 & 0.599 & 3.7\% & -0.442 & -0.494 & 10.5\% & -0.412 & -0.402 & 2.5\% \\ \hline
         D & 0.714 & 0.605 & 0.643 & 5.9\% & -0.414 & -0.494 & 16.2\% & -0.386 & -0.357 & 8.1\% \\ \hline
    \end{tabular}}{}
\end{table}

Observe that the error percentage for $\gamma$ is consistently below \revcolor{6\%}, providing empirical verification of \Cref{thm:informal_upper-bound}. 
In all cases, the fleet scaling is below 2/3, confirming the fundamental difference between EV and non-EV systems and showcasing how Power-of-$d$ leverages partially charged vehicles to reduce fleet size requirements.
In Appendix~\ref{app: matching_algo_comparison}, we perform an extensive comparison of Power-of-$d$ with other policies in this setting, and in the \Cref{sec:sim_chicago}, we show how it compares to other policies in the context of time-varying and non-uniform spatial demand patterns.

\begin{figure}
\centering
\scalebox{0.8}{
\FIGURE{\begin{tikzpicture}    
\node at (0, -6.25) {\scalebox{0.5}{\includegraphics{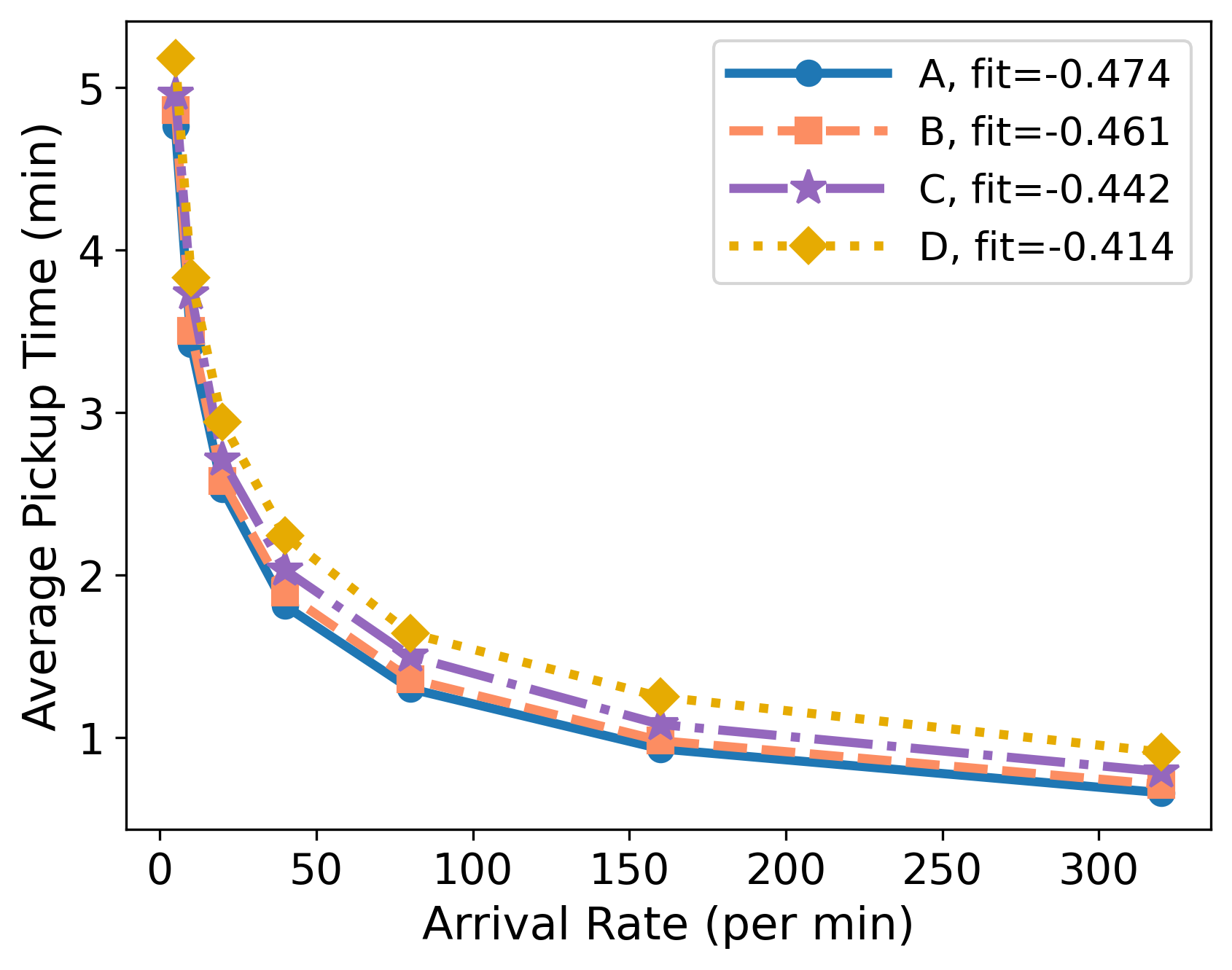}}};
\node at (8, -6.25) {\scalebox{0.5}{\includegraphics{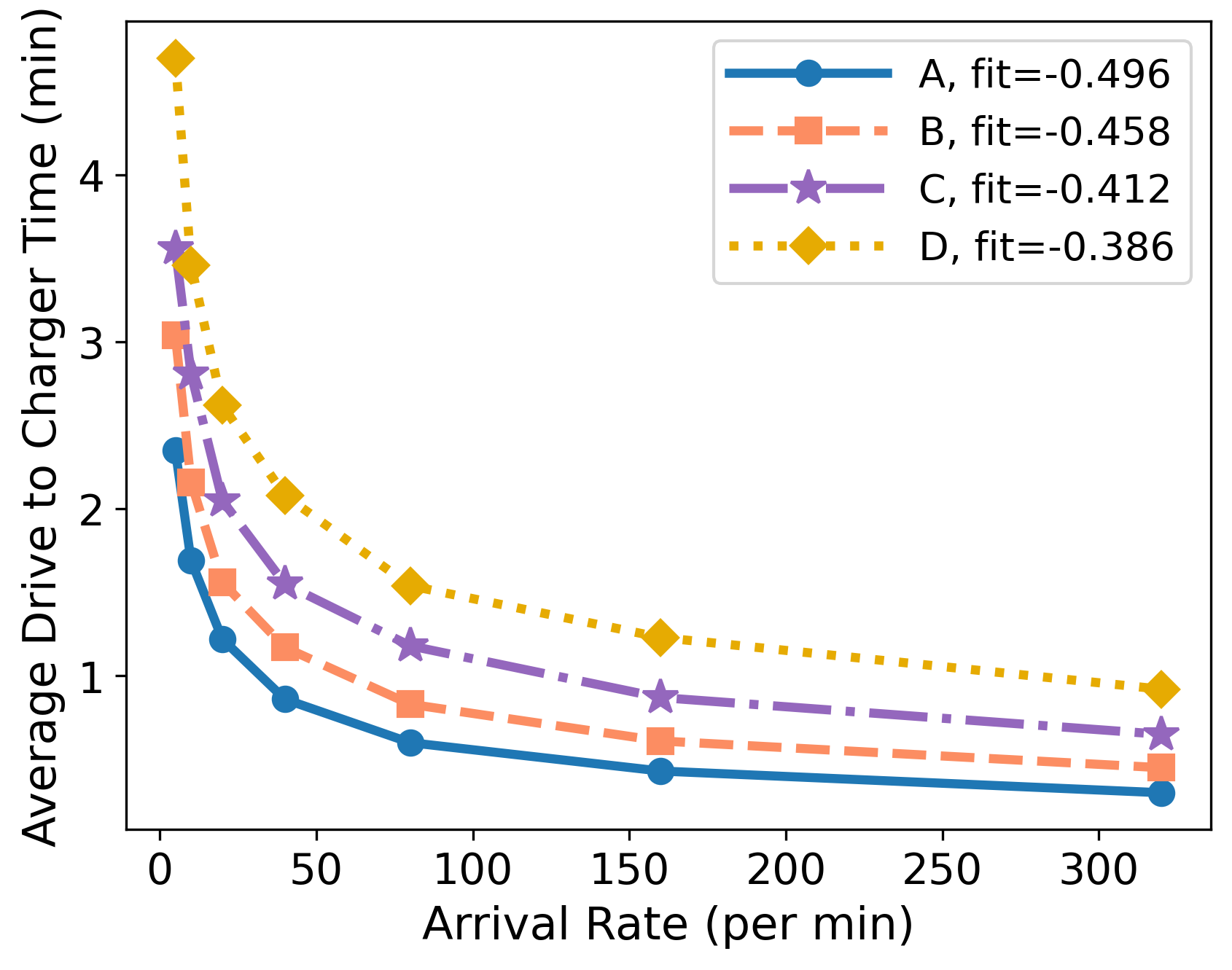}}};
\end{tikzpicture}}{
    \revcolor{Average pickup time (bottom left) and average drive to charger time (bottom right) for Po2.} 
    \label{fig: asymptotic_sim_pick_driveto}}{}
    }
\vspace{-0.25cm}
\end{figure}

The error percentage for the exponent of drive-to-charger is also small (below 8\%), which verifies the functional form in \eqref{eq: drive_to_charger_power_of_d}. The error percentage for pickup time is larger for small $\beta$. The error is due to the departure from the assumption that EVs considered for a match are uniformly distributed in space: the considered EVs are mostly at a charger, and a smaller $\beta$ (a smaller number of chargers) exacerbates this. However, the theory suggests that the pickup time should be insensitive to the values of $\beta$ and $\gamma$ for large $p^\star$, which we indeed observe in \Cref{fig: asymptotic_sim_pick_driveto}.  

\begin{figure}[htb!]
    \FIGURE{
    \begin{minipage}[c]{0.48\textwidth}
    \includegraphics[scale=0.35]{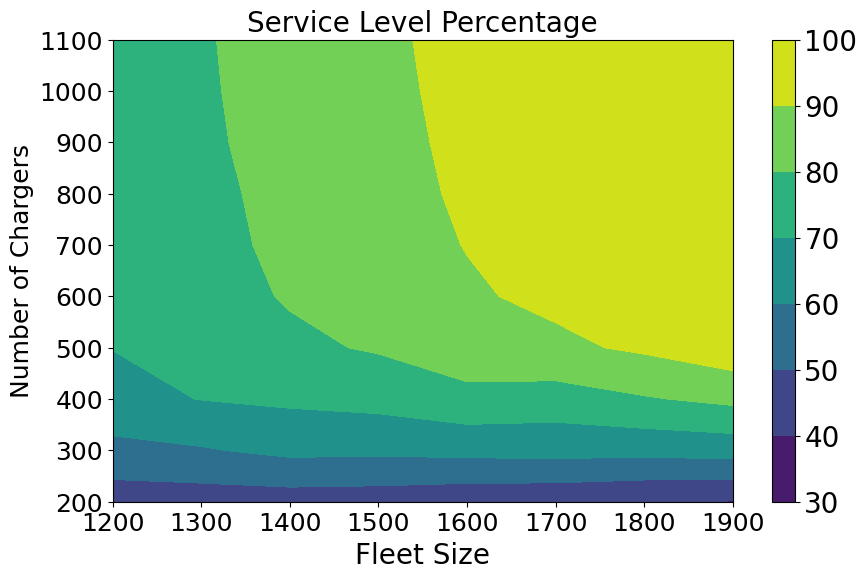}
    \end{minipage}
    \begin{minipage}[c]{0.48\textwidth}
    \centering
\includegraphics[scale=0.35]{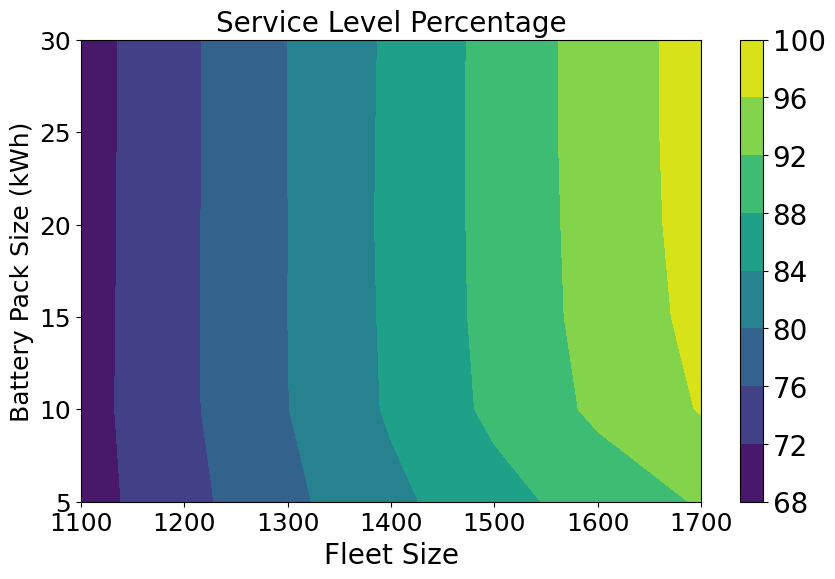}
    \end{minipage}}{\revcolor{Contours of service level (averaged over four runs) as a function of the fleet size and the number of chargers with $p^\star=40$ [kWh] (left), and fleet size and battery pack size (right) with $m_p=1$ and $m_l=1500$. We set $\lambda = 80$ per min.} \label{fig: contour_service_level}}{}
    \vspace{-0.25cm}
\end{figure}

Finally, our results support the idea of a trade-off between chargers and EVs. The left panel of \Cref{fig: contour_service_level} shows how more vehicles are needed to reach a certain service level as the number of chargers decreases, and vice versa. Additionally, in agreement with \Cref{thm:informal_upper-bound}, the right panel of the figure shows that for larger values of $p^\star$, the fleet size requirement remains largely constant.

\subsection{Chicago Trip Data Study}\label{sec:sim_chicago}
\begin{figure}
    \FIGURE{
    \includegraphics[scale=0.4]{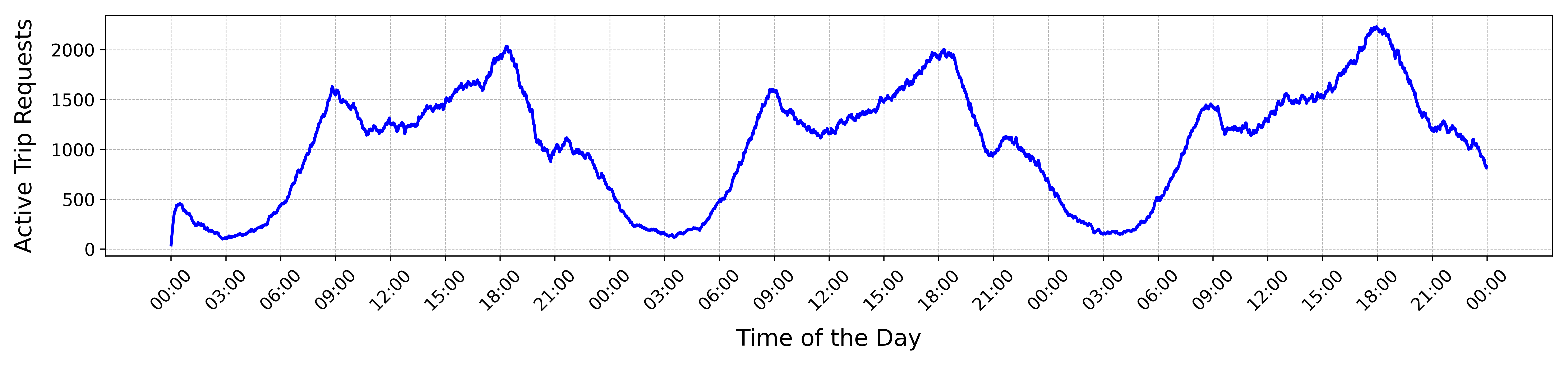}}{Chicago data:
    Active trip requests at different times of the day.
    \label{fig:raw_demand_data}}{}
\end{figure}

We start with a brief description of the system's parameters and data; a more extensive description can be found in Appendix~\ref{sec:sim_chicago_additional_setup}. We use trip-level data from ride-sharing companies operating in Chicago, IL \citep{chicago_dataset}. For demonstration purposes, we fixed the date range to June 14, 2022 (Tuesday) to June 17, 2022 (Thursday) for all our simulations (our results remain qualitatively the same for other dates). \Cref{fig:raw_demand_data} shows the active trip requests during this period. We consider four representative EVs: 
Nissan Leaf ($\sim$\$30,000), Tesla Model 3 Standard Range ($\sim$\$40,000),
Mustang Mach-E SR RWD ($\sim$\$50,000), and Hyundai IONIQ 5 ($\sim$\$50,000).\footnote{\url{https://ev-database.org/} accessed on 07/23/2024.} After accounting for a nominal 10\% battery degradation typically experienced in 4-5 years of usage and considering an average consumption, we use a range of 130, 200, 260, \revcolor{and 290 miles,} respectively. For most of our simulations, we consider a fixed charge rate of 20 [kW] (level 2 chargers), while in \Cref{sec:chicago-perf-pod}, we include 100 [kW] (DC fast charging). We set the charging station locations to be uniformly at random within the city of Chicago, each endowed with 4 charging ports (i.e., $m_p = 4$). 
 
We simulate five main policies: (i) Power-of-$d$ (Po$d$),\footnote{We allow fractional values of $d$ by sampling (independently) $\lfloor d \rfloor$ EVs with probability $\lceil d \rceil - d$ and $\lceil d \rceil$ EVs otherwise.} (ii) Closest Dispatch (CD), (iii) Closest Available Dispatch (CAD), (iv) Power-of-$T_{P,\text{max}}$ (Po$T_{P,\text{max}}$), and (v) Charge-at-Night (CaN). For policies (i) to (iv), after serving a customer, an EV is dispatched to the nearest available charging station, where it stays until fully charged or matched with a customer. All EVs charging, idling, waiting for a charger, or driving to a charger are potential candidates to serve customers. A vehicle is dispatched only if it has enough SoC\footnote{To reach the customer's pickup location, destination, and the nearest charger with at least 5\% SoC remaining.}
 and if the pickup time does not exceed $T_{P,\text{max}} \in \{30, 45, 60, \infty\}$; otherwise, the customer is dropped. The policies differ in their selection criteria: (ii) considers only the closest EV to the customer request; (iii) considers the closest EV with sufficient SoC to serve the customer request; and (iv) selects the EV with the highest SoC among those with a pickup time of at most $T_{P,\text{max}}$. Policy (v) can be combined with any of the other four matching policies. During the day (5 am to 12 am), EVs charge only if their SoC falls below 10\% or 20\% (optimally chosen), replenishing just enough to operate until 12 am. During the night (12 am to 5 am), all EVs charge to 100\%. Under this policy, charging is not interrupted, so only idle EVs are dispatched and, after serving a customer, an EV remains at the drop-off location (see Appendix~\ref{sec: detailed_policy_description} for additional details).

We start by considering infrastructure planning, then show the importance of adaptivity to non-uniform demand and load-balancing, and compare the performance of the Po$d$ policy with the aforementioned policies.

\subsubsection{Infrastructure Planning} We demonstrate how to estimate the minimum fleet and charging infrastructure using our baseline model. Then, we discuss the infrastructure tradeoff in the current context.

A first, data-driven approach to infrastructure planning is to use simulation to find the minimum number of EVs and chargers that achieve a target served workload.
While effective, this approach can be time-consuming. However, it is possible to efficiently approximate the results from the simulator using the Power-of-$d$ ODEs given by \eqref{eq:odes-power-d}. In particular, the arrival rate $\lambda(t)$ and the mean trip time $T_R(t)$ are allowed to vary with $t$ and are fitted from the data.
The pickup time $T_P(t)$ 
and the drive-to-charger time $T_{DC}(t)$ are estimated from simulations. We set a nominal admission control for incoming customers and charging EVs to ensure that $T_P(t)$ and $T_{DC}(t)$ are bounded (by $10$ min and $15$ min). The maximum service time $T_B$ is approximated by the mean cumulative trip, pickup, and drive-to-charger time, which turns out to be $30.5$ min (see 
Appendix~\ref{app:infra_plan_odes} for a detailed description of the estimated parameters). Because $T_P(t), T_{DC}(t)$, and $T_B$ are estimated values, we also test the performance of the ODEs by considering a parameter misspecification of $\pm 20\%$.

We simulate the ODEs to determine the fleet size that serves 90\% of requested miles and compare the results from the detailed simulations (see \Cref{tab: ode_fleet_sizing}). When the parameters $(T_P(t), T_{DC}(t), T_B)$ are not misspecified, the ODE outputs fleet sizes with less than 2\% error. Moreover, even with 20\% parameter misspecification, the error remains within \revcolor{5-20\%}. While the ODEs abstract away the system's stochasticity and its spatial nature, they provide prescriptions aligned with more detailed simulations, serving as a simple first approximation for fleet and charging infrastructure planning.

\begin{table}[tbh!]
    \centering
    \TABLE{\revcolor{Comparing the fleet sizing corresponding to 90\% of service level (served workload) for $m \in \{125, 175, 225\} \times 4$ under the simulator and the ODEs given by \eqref{eq:odes-power-d} for Po$d$ with Tesla.} 
    \label{tab: ode_fleet_sizing}}{
    \begin{tabular}{|c|c|c|c|c|c|c|} \hline
       \multirow{2}{*}{Chargers}  & \multirow{2}{*}{Optimal $d$} & \multicolumn{4}{c|}{90\% Fleet Size (\% Error)}  \\ \cline{3-6}
       & & Sim & ODE & ODE $+20\%$ Param Err & ODE $-20\%$ Param Err \\ \hline \hline
       500 & 2.6 & 2405 & 2503 (4\%) & 2947 (22\%) & 2208 (8\%) \\ \hline
       700 & 1.6 & 2290 & 2349 (3\%) & 2573 (12\%) & 2150 (6\%)  \\ \hline
       900 & 1.4 & 2254 & 2312 (3\%) & 2494 (11\%) & 2133 (5\%)  \\ \hline
\end{tabular}\label{tab:ode_planning_chicago}}{}
\end{table}

\begin{figure}[htb!]
    \centering
    \FIGURE{
    \begin{minipage}[c]{0.48\textwidth}
    \includegraphics[scale=0.36]{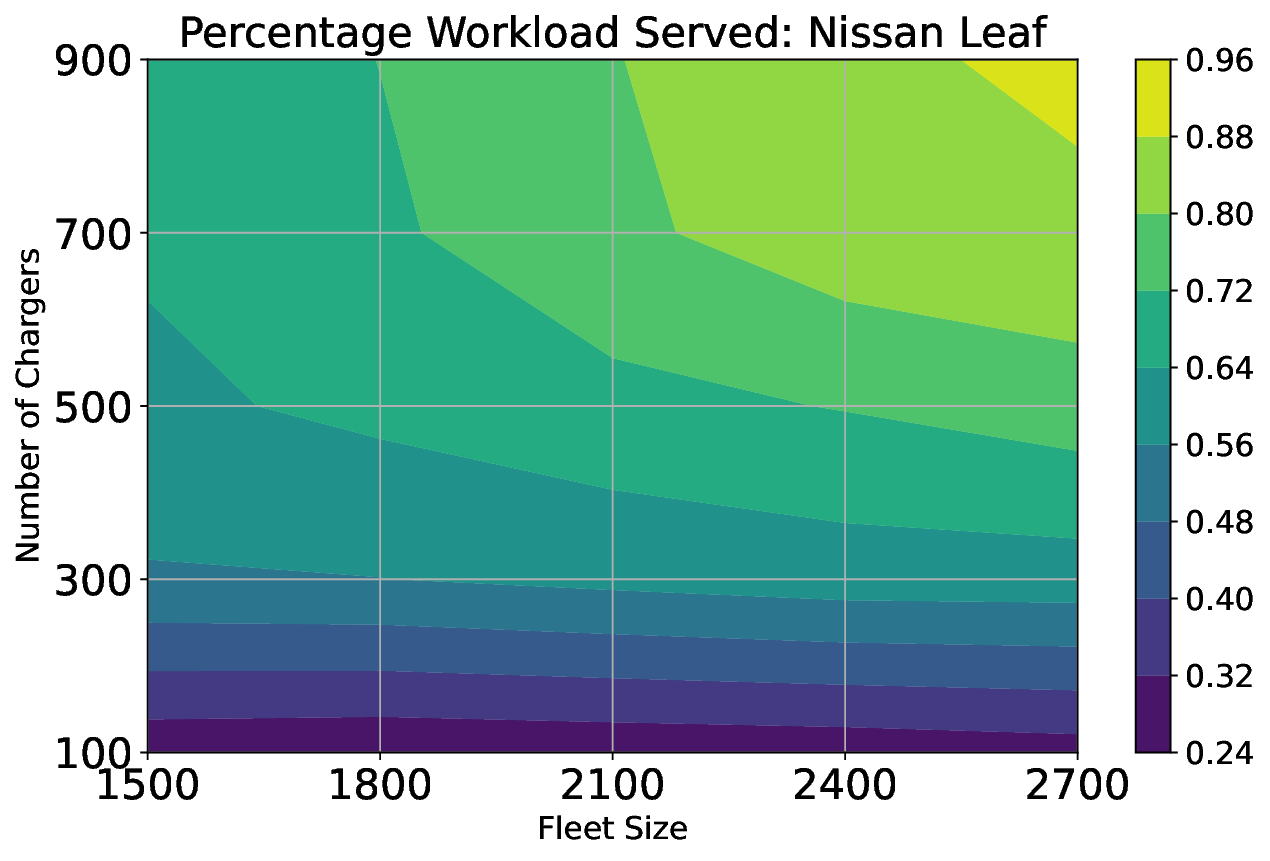}
    \end{minipage}
    \begin{minipage}[c]{0.48\textwidth}
    \includegraphics[scale=0.36]{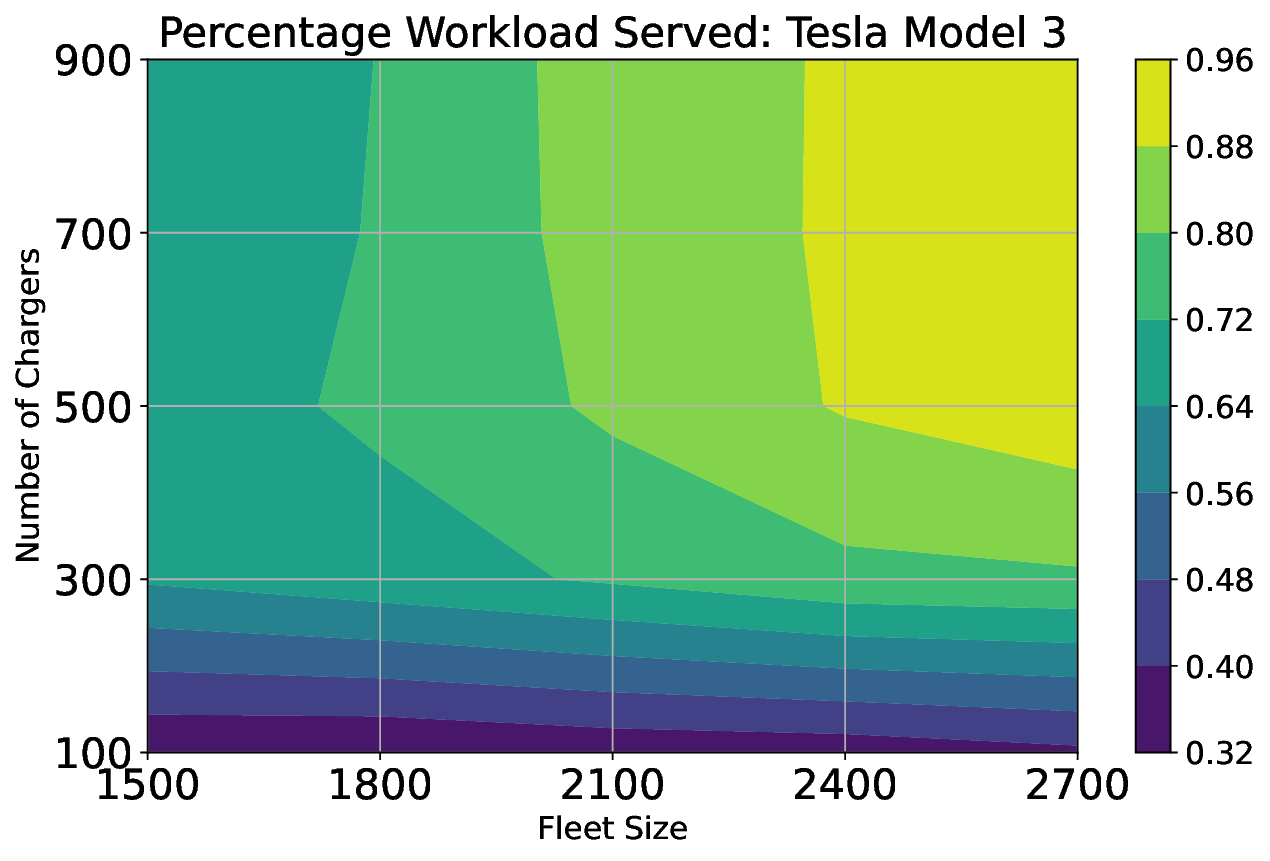}
    \end{minipage}
    }{Contours of service level as a function of the fleet size and the number of chargers for Nissan (left) and Tesla (right) under Po$2$ in the detailed simulation. 
    \label{fig:infra_planning_chicago}
    }{}
\end{figure}

Finally, consistent with our observations in \Cref{sim:synthetic_data} and our findings in \Cref{thm:lower-bound}, \Cref{fig:infra_planning_chicago} illustrates the trade-off between fleet size and the number of chargers. The figure demonstrates that under the Po$2$ policy, this fundamental trade-off is still evident in the detailed simulation, reinforcing our earlier theoretical results and their applicability to more complex, realistic scenarios.

\subsubsection{Performance of Power-of-\texorpdfstring{$d$}{d} 
}\label{sec:chicago-perf-pod}
We next discuss mechanisms that drive the performance of Po$d$: adaptivity to non-uniform demand---adaptive charging and adaptive relocation---and load balancing. We will then compare the various policies.

\subsubsection*{Adaptive Charging.} %\label{sec: adaptive_charging}
\begin{figure}[tbh!]
    \centering
    \FIGURE{\scalebox{0.75}{
    \begin{tikzpicture}
        \node at (0, 0) {\includegraphics[scale=0.47]{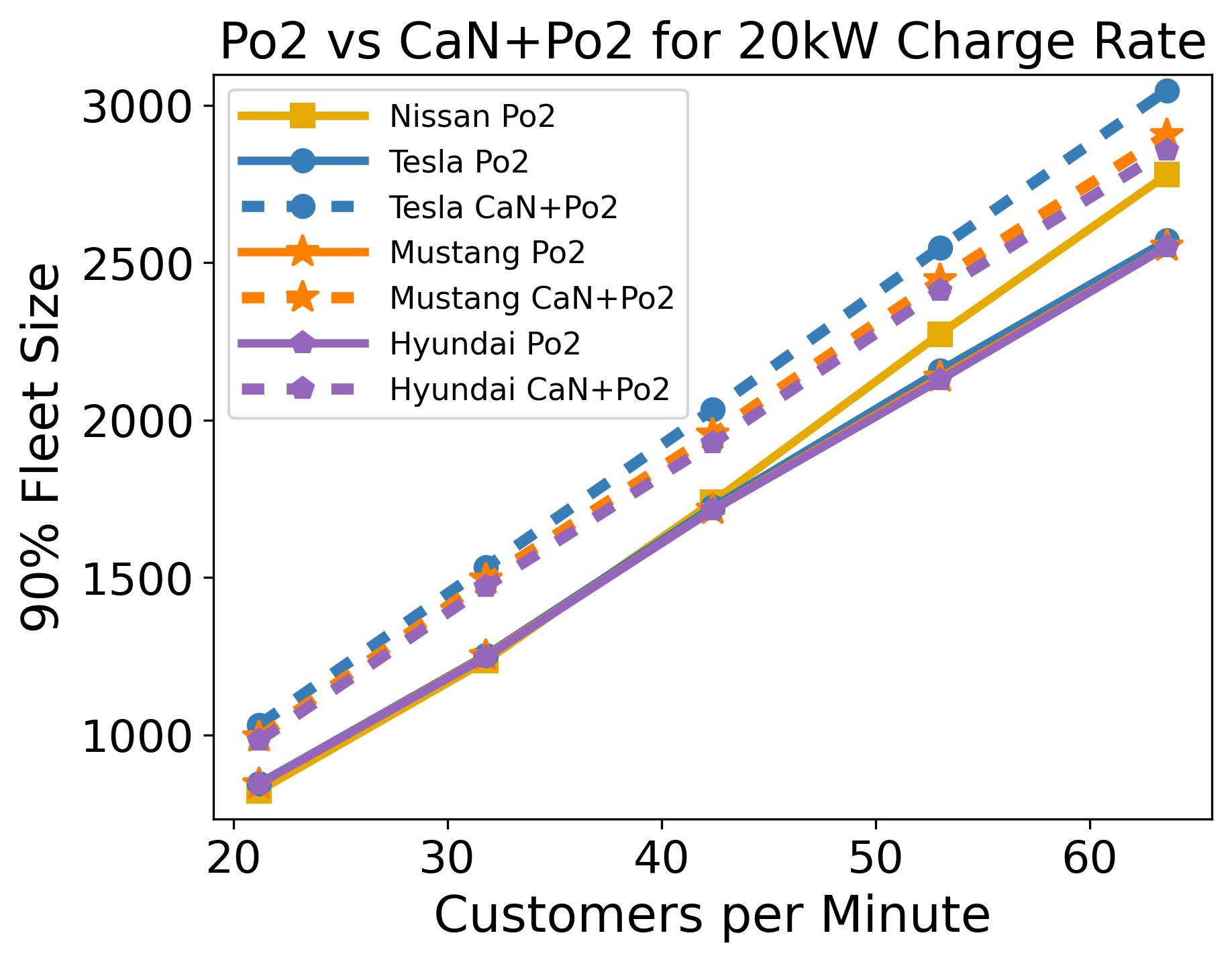}
        };
        \node at (7.7, 0) {\includegraphics[scale=0.47]{
        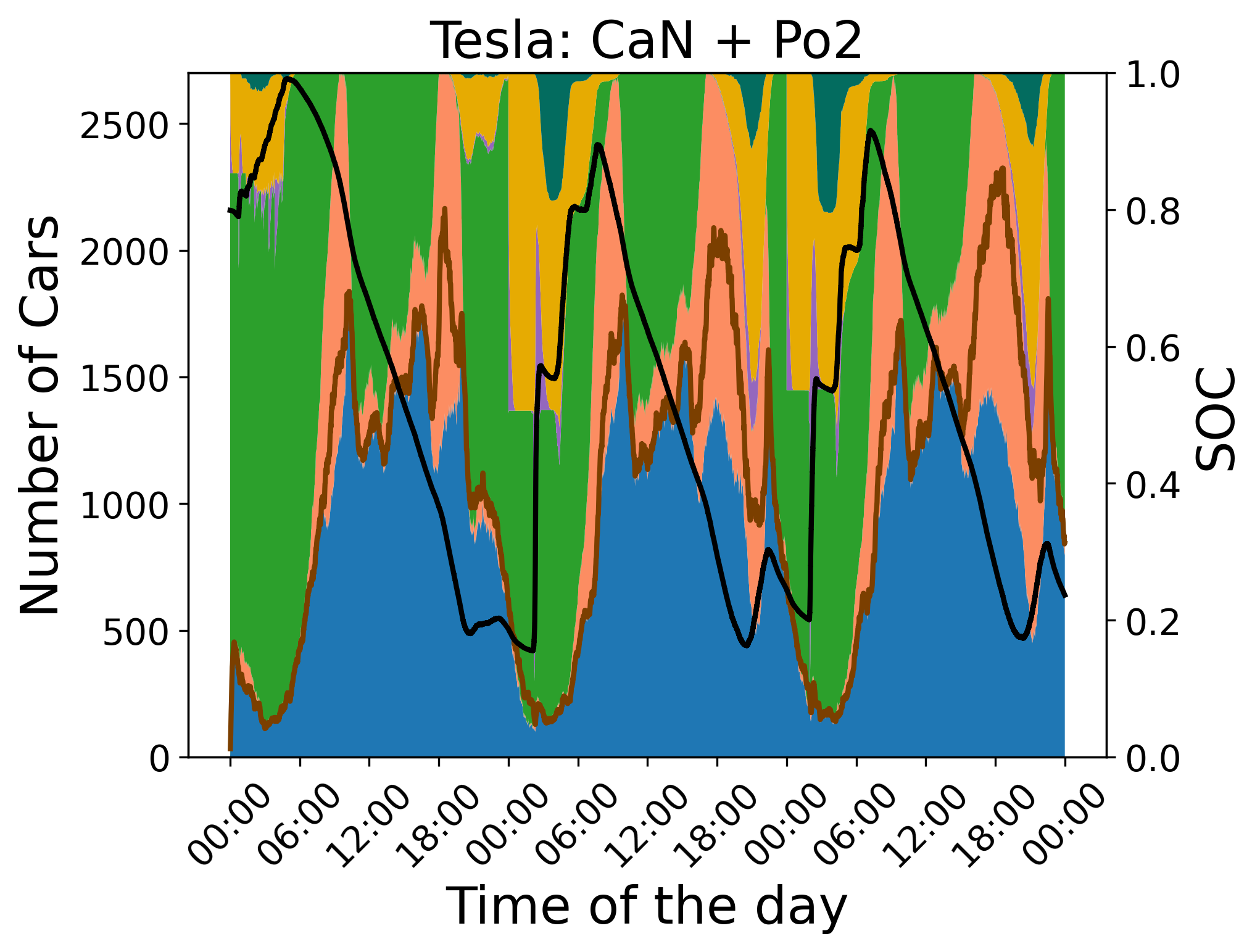}};
        \node at (15.8, 0) {\includegraphics[scale=0.47]{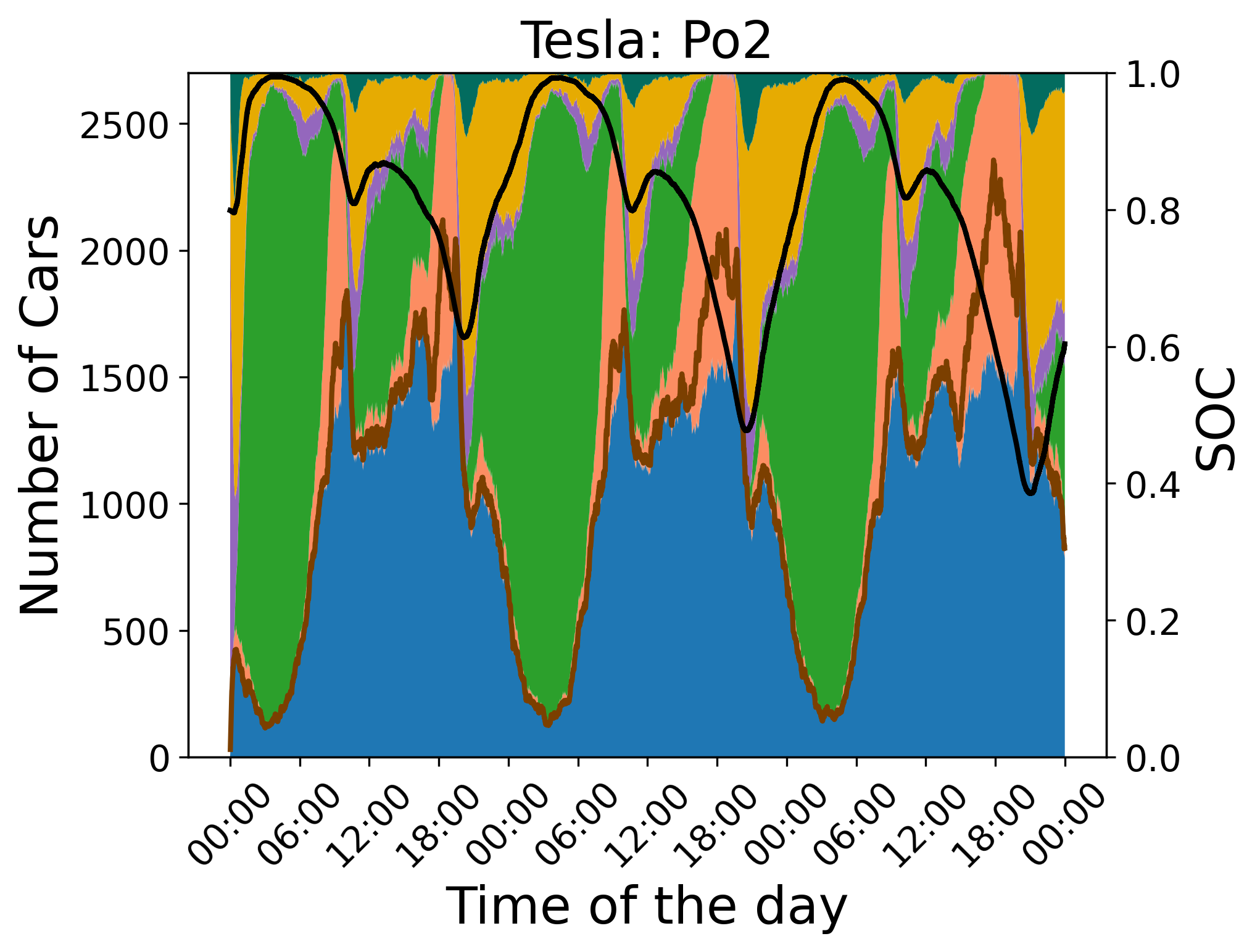}};
        \def \sw {0.5}
        \def \dbs {0.2}
        \begin{scope}[yshift=0.3cm]
        \draw[rectangle] (-1.2, -3.8) -- (18, -3.8) -- (18, -5.4) -- (-1.2, -5.4) -- (-1.2, -3.8);
        \filldraw[rectangle, blue] (-1, -4) -- (-0.5, -4) -- (-0.5, -4.5) -- (-1, -4.5) -- (-1, -4) node[black, right] at (-0.3, -4.25) {Driving with Customer};
        \filldraw[rectangle, orange] (4.2, -4) -- (4.7, -4) -- (4.7, -4.5) -- (4.2, -4.5) -- (4.2, -4)  node[black, right] at (4.9, -4.25) {Picking up Customer};
        \filldraw[rectangle, darkgreen] (10, -4) -- (10.5, -4) -- (10.5, -4.5) -- (10, -4.5) -- (10, -4)  node[black, right] at (10.7, -4.25) {Idle};
        \filldraw[rectangle, purple] (12.8, -4) -- (13.3, -4) -- (13.3, -4.5) -- (12.8, -4.5) -- (12.8, -4) node[black, right] at (13.5, -4.25) {Driving to Charger};
        \filldraw[rectangle, yellow] (-1, -4.7) -- (-0.5, -4.7) -- (-0.5, -5.2) -- (-1, -5.2) -- (-1, -4.7) node[black, right] at (-0.3, -4.95) {Charging};
        \filldraw[rectangle, teal] (4.2, -4.7) -- (4.7, -4.7) -- (4.7, -5.2) -- (4.2, -5.2) -- (4.2, -4.7) node[black, right] at (4.9, -4.95) {Waiting for Charger};
        \draw[black, ultra thick] (10, -4.95) -- (10.5, -4.95) node[black, right] at (10.7, -4.95) {Avg SoC};
        \draw[brown, ultra thick] (12.8, -4.95) -- (13.3, -4.95) node[black, right] at (13.5, -4.95){(Potential) Active Trips};;
        \end{scope}
    \end{tikzpicture}}}{\revcolor{Pod vs. CaN. Fleet size serving 90\% of the workload with $m \in \{150, 225, 300, 375, 450\} \times 4$ for five different arrival rates (left). System state for CaN + Po2 (center) and Po2 (right) for Tesla with $m = 450\times 4$.}
    \label{fig: pod_vs_can_stackplot}}{}
\end{figure}

The demand variations observed in \Cref{fig:raw_demand_data} suggest that, given enough range, a fleet operating without charging during the day can attain a desired service level without necessarily increasing the fleet size requirement. We compare the performance of Po2 with $T_{P,max}=$ 45 min with that of Charge-at-Night (with the best-performing choice between Po2, CAD, and $T_{P,max}$). We compare the performance of the policies by calculating the fleet size corresponding to a target 90\% workload served (when charging stations are not a limiting factor). The resultant fleet size requirements are summarized in the left panel of \Cref{fig: pod_vs_can_stackplot}.

We observe that Po2 outperforms CaN for all simulated EVs, achieving the target served workload with a smaller fleet size. Moreover, Po2 achieves this target with a smaller and more cost-efficient battery pack size: the fleet size requirement for Nissan under Po2 is smaller than that for other vehicles under CaN (Nissan under CaN is excluded due to a high fleet requirement). We note that the inferior performance of the longer-range EVs is partly due to CaN's lack of adaptive charging and---as discussed in the next subsection---partly due to the impact of relocation.
\vspace{-0.1cm}
To understand how adaptive charging affects the results, we plot the system state evolution as a stack plot in \Cref{fig: pod_vs_can_stackplot}. Po$d$ naturally adapts by charging more EVs during low-demand periods (and not charging during high-demand periods): note the increased yellow area in the right panel from 11 am to 3 pm and 11 pm to 5 am. Unlike CaN, Po$d$ leverages low-demand periods in the afternoon to maintain a consistently high SoC of the fleet (the black curve stays above 0.4).

Under CaN, Tesla's range is insufficient, forcing some EVs to charge around 7-8 pm as seen by the small peak in SoC at that time, reducing the number of available EVs, increasing pickup times (larger orange area around 7-8 pm during the second and third days), and worsening system performance.\footnote{A similar but less pronounced effect is observed with the Mustang fleet in \Cref{fig: pod_vs_can_stackplot_full} in Appendix~\ref{sec:sim_chicago_additional_setup}.} As we explain next, while short-range vehicles unavoidably miss some demand due to insufficient range and the lack of charging adaptivity, the poor performance of Mustang and Hyundai is largely driven by relocation.

\subsubsection*{Adaptive Relocation.} 

Because the Chicago dataset exhibits {\it spatially imbalanced demand}, EVs under CaN tend to concentrate in destination-heavy regions after drop-offs. As a result, vehicles face long pickup times, reducing their availability and driving up the capacity requirement. In contrast, under Po$d$, EVs are dispatched to charging stations after drop-off, which are distributed across the city, inducing a natural relocation of the fleet.

To shed light on the effect of relocation, we simulate a variation of CaN, called CaN-R, which mimics Pod’s relocation behavior: after drop-off, vehicles are relocated to a charging station and occupy it for the duration they would take to charge, but without actually increasing their SoC. Note that such relocation has the benefit of reducing pickup times, but it also runs the risk of losing demand due to a resulting reduced fleet SoC and the lack of adaptive charging.

\begin{figure}[tbh!]
    \centering
    \FIGURE{\scalebox{0.75}{
    \begin{tikzpicture}
        \node at (0, 0) {\includegraphics[scale=0.5]{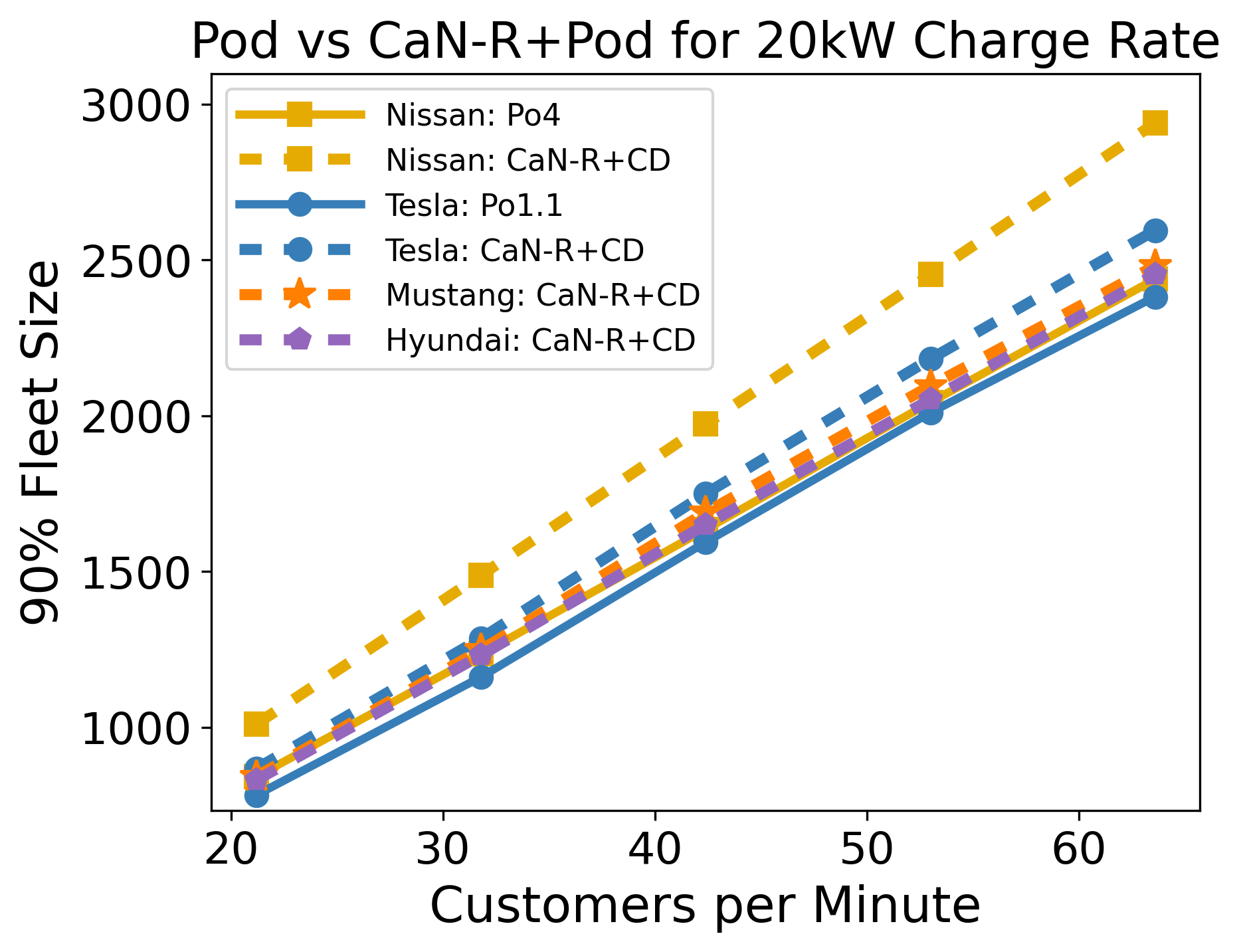}
        };
        \node at (7.8, 0) {\includegraphics[scale=0.5]{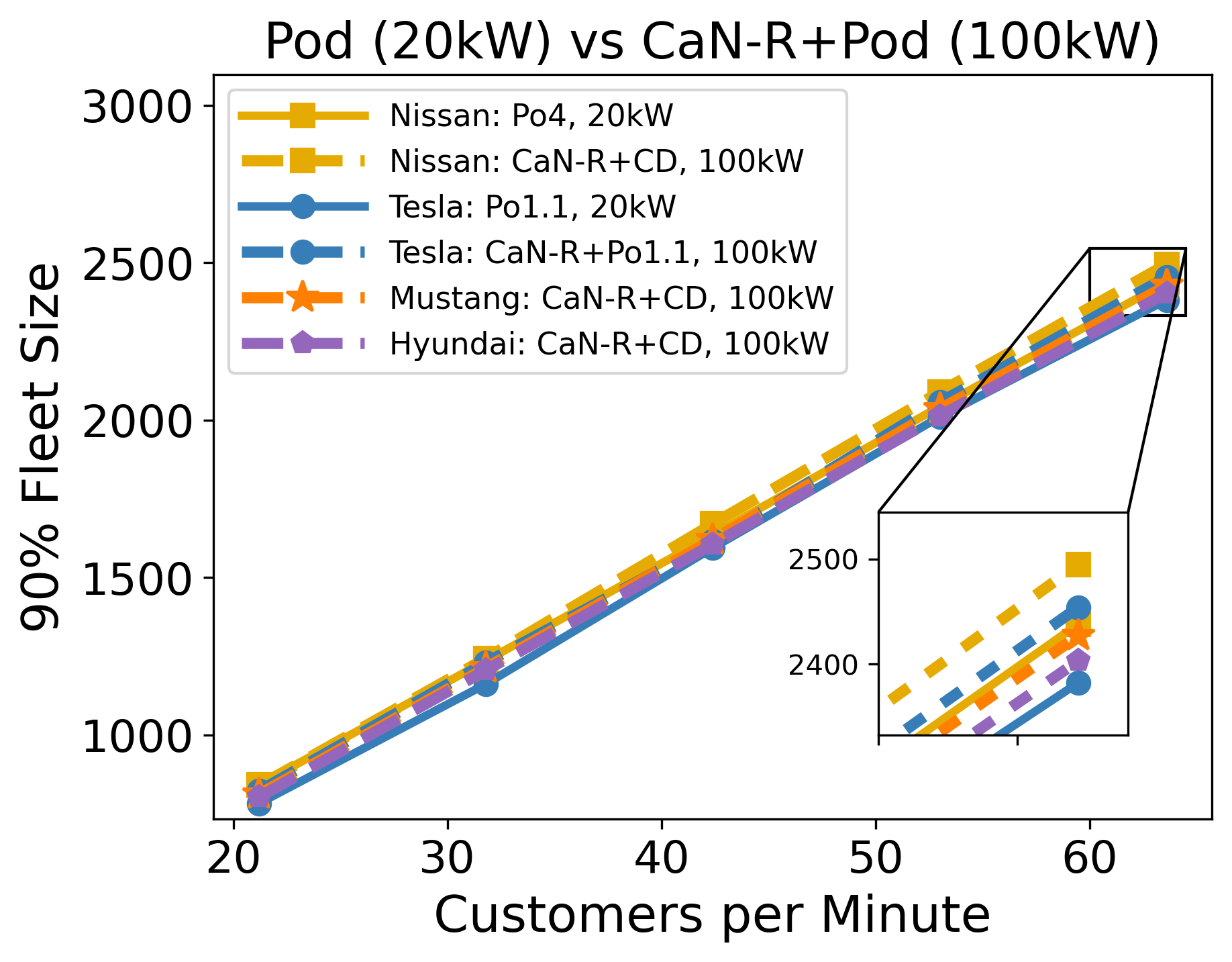}};
        \node at (15.85, 0) {\includegraphics[scale=0.5]{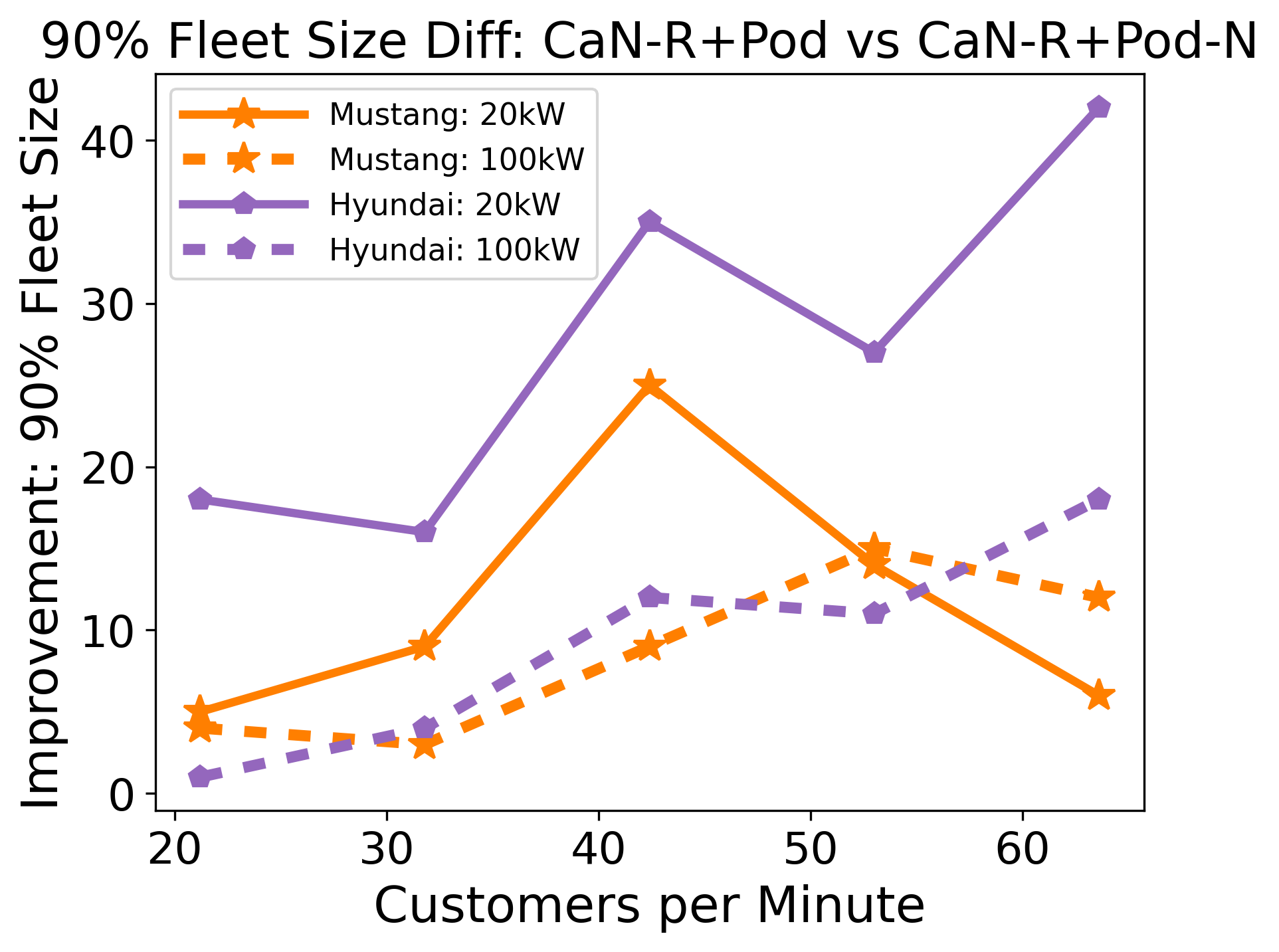}};
    \end{tikzpicture}}}
    {\revcolor{Pod vs. CaN. Fleet size serving 90\% of the workload with $r_c=20$ [kW] (left) and $r_c=100$ [kW] (center), fleet size improvement (decrease) for Hyundai and Mustang under CaN-R with Pod at night (right);
    $d$ is optimally chosen and $m \in \{150, 225, 300, 375, 450\} \times 4$ for five different arrival rates.} 
    \label{fig: pod_vs_can_20_100_opt_d}}{}
\end{figure}
\Cref{fig: pod_vs_can_20_100_opt_d} together with the left panel of \Cref{fig: pod_vs_can_stackplot} show that 
(i) adding relocation to the baseline CaN policy (CaN‑R) sharply improves performance: it reduces pickup times, boosts availability, and lowers capacity requirements because long-range EVs such as the Hyundai---and, to a lesser extent, the Mustang---can absorb the SoC loss incurred while relocating and still operate all day without recharging; 
(ii) relocation alone is not enough for shorter‑range vehicles with Level 2 chargers: under CaN-R, Nissan and Tesla fleets still deplete their batteries during the day and miss demand; 
(iii) augmenting CaN‑R with an EV‑aware night‑time Po$d$ policy that can interrupt charging for high‑SoC vehicles further increases availability and leaves the fleet with a higher, more uniform SoC for the next day (cf. \Cref{sec:pod-works-adaptive}), benefiting even long‑range vehicles (see right panel in \Cref{fig: pod_vs_can_20_100_opt_d} and Appendix~\ref{sec:additional_result_sim} for relevant stack plots); and (iv) although DC fast charging (middle panel) benefits every type of fleet, careful SoC management on Level 2 infrastructure can perform competitively: with optimally chosen $d$, Po$d$ enables shorter‑range models like the Nissan or Tesla to match the performance of longer‑range vehicles or systems equipped with fast chargers.

When demand is spatially imbalanced, both SoC management and vehicle relocation matter for system performance. The relative importance of each depends on vehicle and infrastructure characteristics---particularly range and charging rate.
For short‑range fleets, optimally selecting $d$ to keep SoC levels high is critical, whereas long‑range fleets or those with fast chargers gain more from relocation and can be managed like a non‑EV system during the day (with an EV‑aware policy at night). Po$d$ already manages both aspects effectively; however, further optimizing relocation---or any other mechanism for keeping workload low---in both Po$d$ and CaN‑R, particularly for long‑range fleets, could further enhance performance (Appendix~\ref{sec:additional_result_sim} supplements this discussion).

\subsubsection*{Load Balancing.} 
\begin{figure}[hbt!]
    \centering
    
    \FIGURE{\scalebox{0.75}{
    \begin{tikzpicture}
        \node at (0, 0) {\includegraphics[scale=0.48]{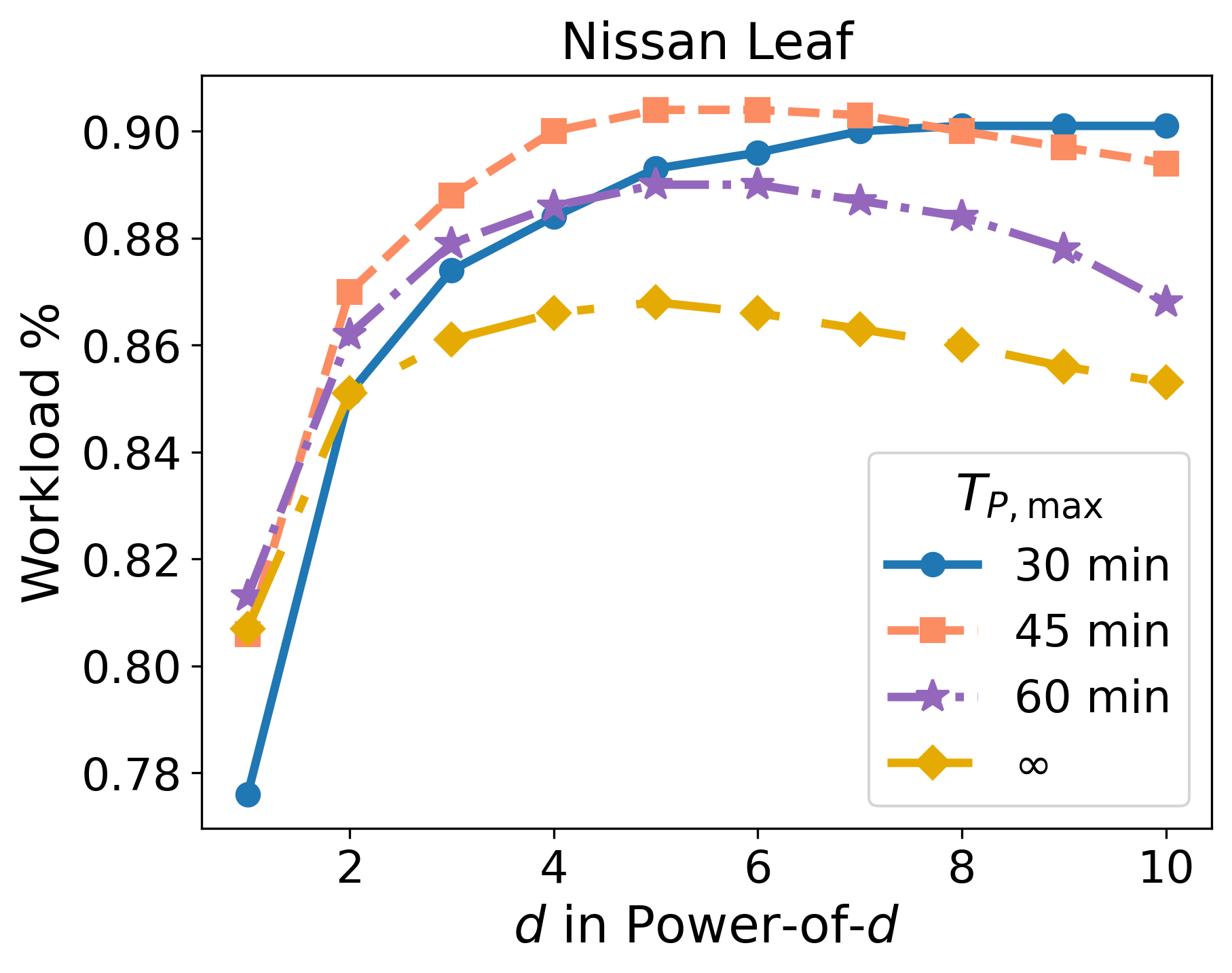}};
        \node at (7.5, 0) {\includegraphics[scale=0.48]{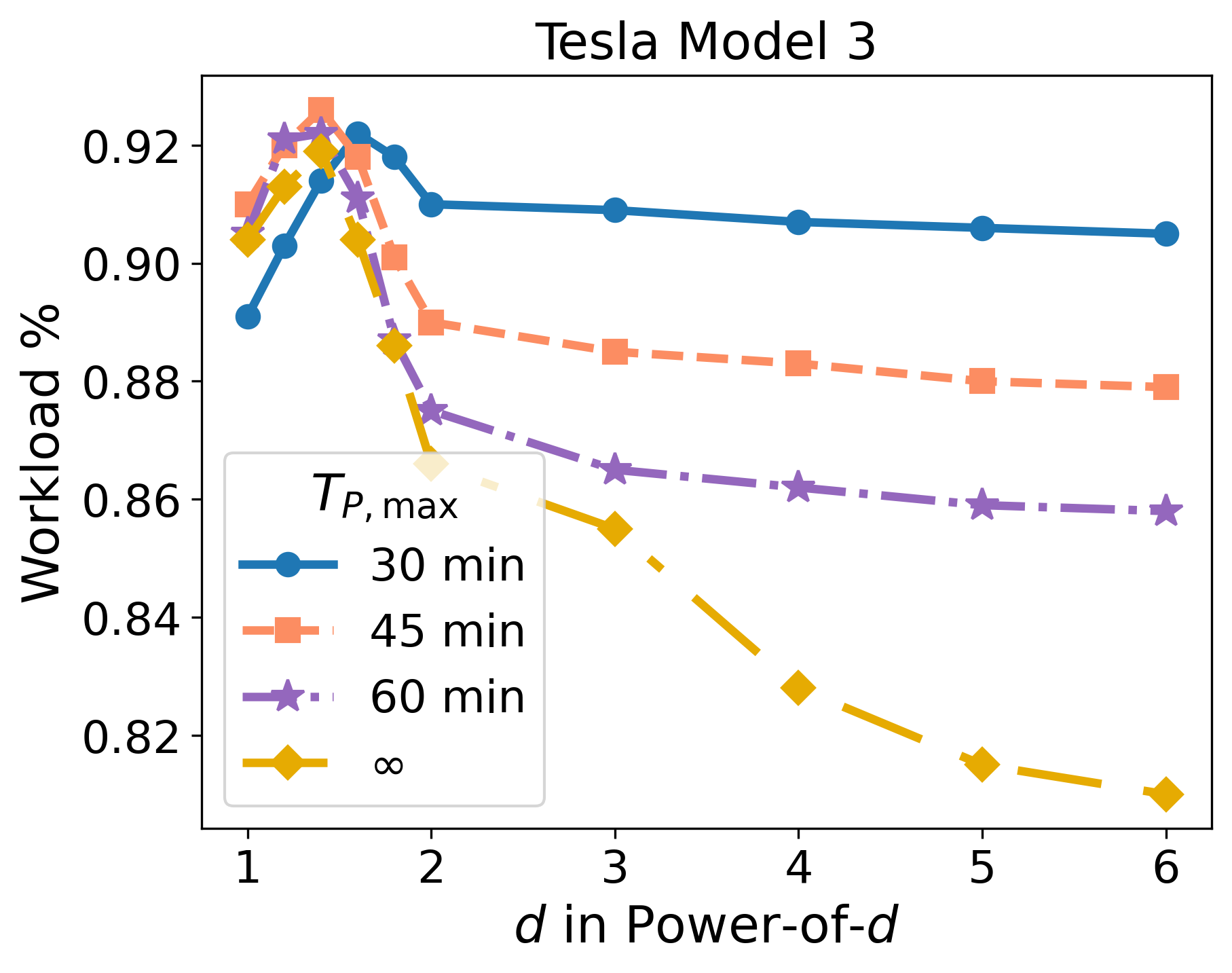}};
        \node at (15, 0) {\includegraphics[scale=0.25]{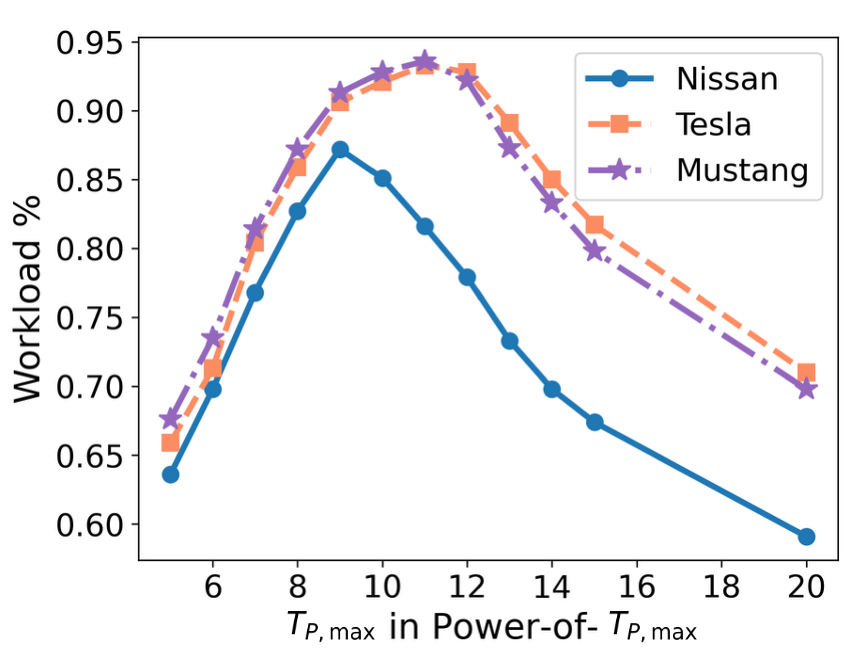}};
    \end{tikzpicture}}
    }
    {
     \revcolor{Served workload (\%) (averaged over four runs) for Po$d$ and Po$T_{p,max}$; $n = 2400$ and $m = 225 \times 4$.}
    \label{fig: real_sim_load_balancing}}{}
\end{figure}

We illustrate the need for load balancing in \Cref{fig: real_sim_load_balancing}, which shows the served workload for Po$d$ varying $d$ and $T_{P, \max}$ for Nissan and Tesla fleets. The optimal $d$ value greater than 1 demonstrates the necessity of load balancing. The difference between CD and Po$d$ becomes more pronounced with smaller battery sizes, owing to worsening fleet SoC balancing.\footnote{
\Cref{prop:closest-fails} in
Appendix \ref{sec: closest_dispatch} shows that CD requires a larger fleet than Po$d$ to achieve the same service level.}
Additionally, aligned with Theorem~\ref{thm:informal_upper-bound}, for $T_{P, \max} = 45$, the optimal $d$ decreases from \revcolor{$d = 6$} for Nissan to $d = 1.4$ for Tesla: with larger $p^\star$, closer vehicles are more likely to have a sufficiently high SoC (see also \Cref{tab: real_sim_cad_vs_pod}).

This shows the importance of balancing the need for short pickup times (low $d$) and a fleet with a balanced SoC (high $d$). Achieving the right trade-off avoids having EVs drive long distances to pick up customers, saving SoC and reducing the likelihood of dropping customers. \revcolor{Values of 
$d$ in the range of 2–4 for Nissan are particularly effective in managing this trade-off}\label{R1-round-2:comment-8-2}. Balancing the load in this way enables low SoC EVs to recover during the valley and be available for the next peak (cf. \Cref{sec:pod-works-adaptive}), as shown by the SoC increase at 12 pm in the right panel of \Cref{fig: pod_vs_can_stackplot}.

\subsubsection{Policy Comparison} \label{sec:sim_chicago_policy_comp} 
\Cref{tab: real_sim_cad_vs_pod} shows that Po$d$ consistently delivers a competitive performance.\footnote{\revcolor{Hyundai and Mustang obtain a similar performance.}}
It outperforms CD and CAD and is generally better or comparable to Po$T_{P, \max}$.\footnote{See \Cref{tab:chicago_served_requests} in Appendix~\ref{sec:sim_chicago_additional_setup} for similar results using the percentage of customers served.} For the optimal $T_{P, \max}$, Po$T_{P, \max}$ can result in a better served workload but higher pickup times compared to Po$d$. Po$T_{P, \max}$ balances SoC, but it does so in an adaptive manner---when many (few) EVs are idle, it effectively operates like Po$d$ with large (low) $d$---explaining its good performance (cf. \Cref{sec:pod-works-adaptive}). However, because this policy selects the EV with the highest SoC within $T_{P, \max}$, the pickup times are skewed towards $T_{P, \max}$ and are, therefore, longer than those under Po$d$. 

\begin{table}[hbt!]
    \centering
    \TABLE{\revcolor{Served workload (\%) (averaged over four runs) and pickup time (min) for CD, CAD, Po$T_{P, \max}$, and Po$d$ with $d$ optimal; $n=2400$ and $m = 225 \times 4$.} \label{tab: real_sim_cad_vs_pod}}{
    \begin{tabular}{|c||c|c|c|c|c|c|c|c|} \hline
       \multirow{2}{*}{$T_{P, \max} = \infty$}  &  \multicolumn{2}{c|}{CD} & \multicolumn{2}{c|}{CAD}  & \multicolumn{2}{c|}{Po$T_{P, \max}$} & \multicolumn{2}{c|}{Po$d$ ($d = \cdot$)} \\ \cline{2-9} 
       &  Workload & Pickup & Workload & Pickup & Workload & Pickup & Workload & Pickup 
         \\ \hline \hline
         Nissan  & 81.0\% & 2.97  & 83.8\% & 6.36 & $\leq$65\% & 12$\geq$ & 86.8\% (5) &  5.22 \\ \hline
         Tesla  & 90.5\% & 3.76  & 88.5\% & 6.68 & $\leq$75\% & 12$\geq$ & 91.4\% (1.2) & 4.00  \\ \hline
         Mustang  & 91.5\% & 4.43 & 88.5\% & 6.71 & $\leq$75\% & 12$\geq$& 91.5\% (1) &  4.43 \\ \hline \hline
         best $T_{P, \max}$ &  \multicolumn{2}{c|}{CD ($T_{P, \max} = \cdot$)} & \multicolumn{2}{c|}{CAD ($T_{P, \max} = \cdot$)}  & \multicolumn{2}{c|}{Po$T_{P, \max}$ ($T_{P, \max} = \cdot$)} & \multicolumn{2}{c|}{Po$d$ ($T_{P, \max} = \cdot, d = \cdot$)} \\ \hline \hline
         Nissan  & 81.1\% (60) &  2.88 & 87.1\% (30) & 4.85 & 87.1\% (9) & 5.58 & 90.3\% (45, 6) & 4.45 \\ \hline
         Tesla  & 91.0\% (45)& 3.34 & 92.4\% (30) & 4.40  & 93.3\% (11) & 6.76 & 92.4\% (45, 1.4) & 3.99 \\ \hline
         Mustang  & 92.5\% (45) & 4.26 & 92.5\% (30) & 4.42 & 93.6\% (11) & 6.77 & 92.7\% (45, 1.05) & 3.86 \\ \hline
    \end{tabular}}{}
\end{table}

\begin{figure}
    \centering
    \FIGURE{\scalebox{0.75}{
    \begin{tikzpicture}
        \node at (7.8, 0) {\includegraphics[scale=0.5]{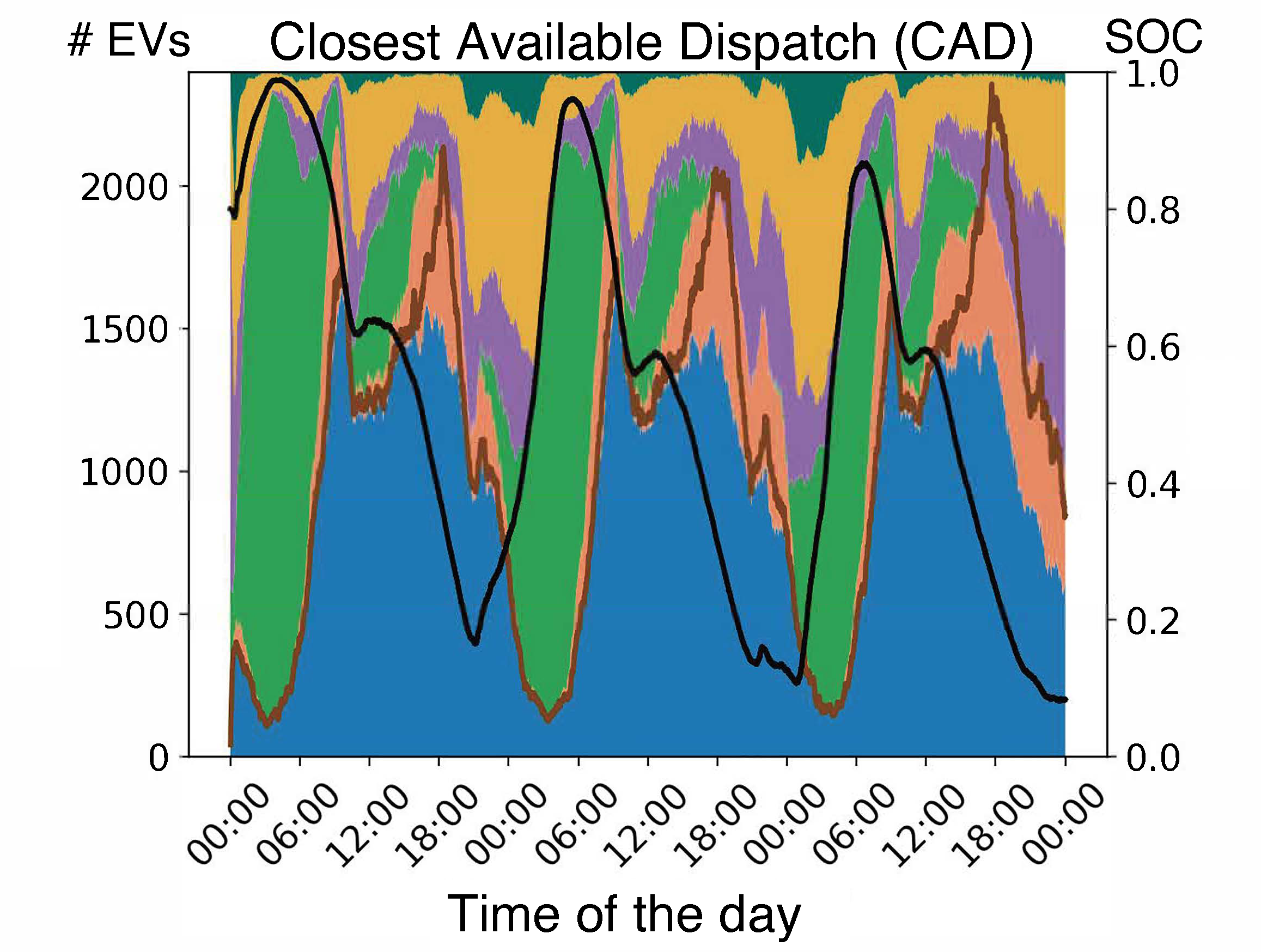}};
        \node at (16, 0) {\includegraphics[scale=0.5]{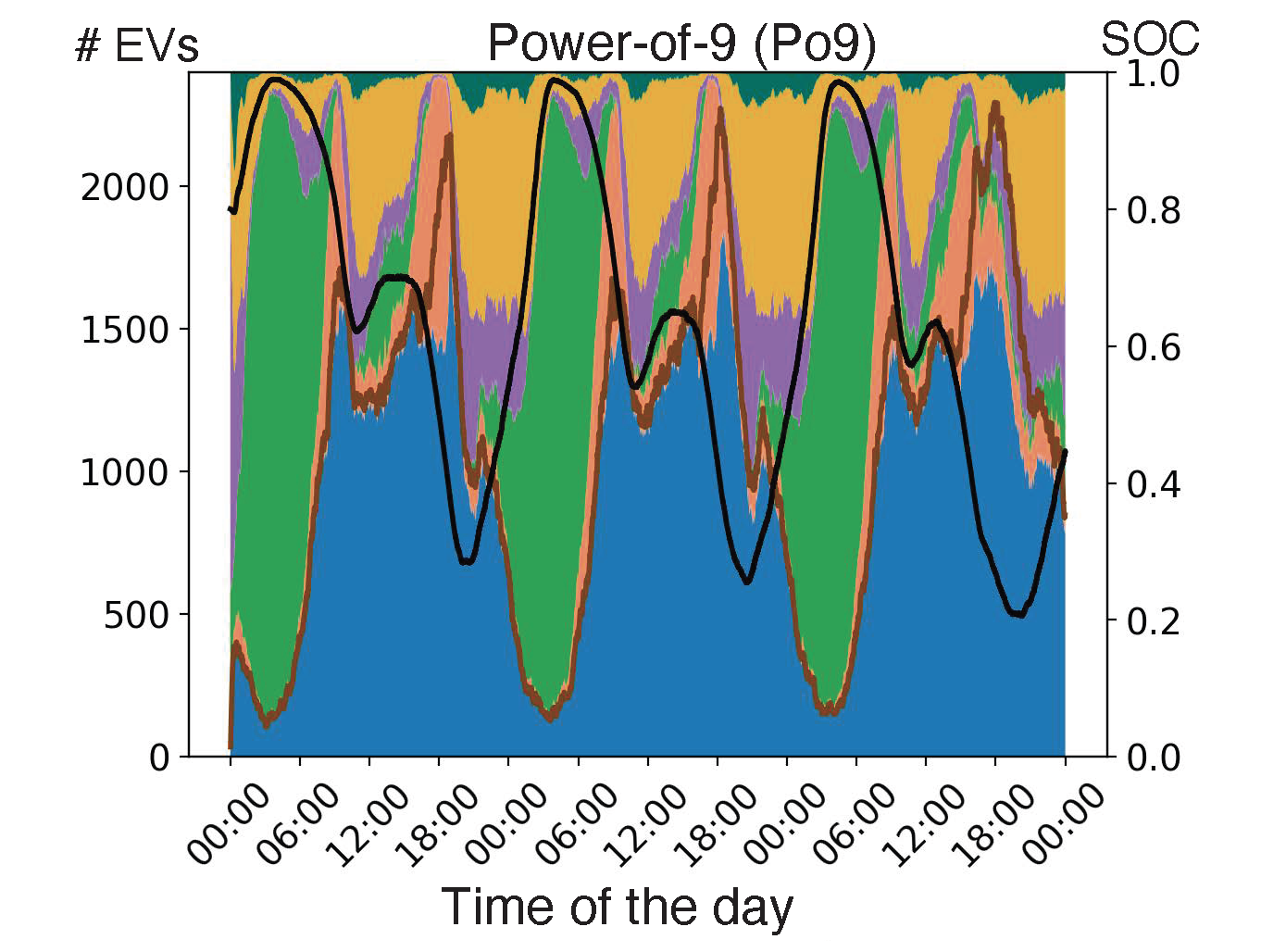}};
        \node at (0, 0) {\includegraphics[scale=0.5]{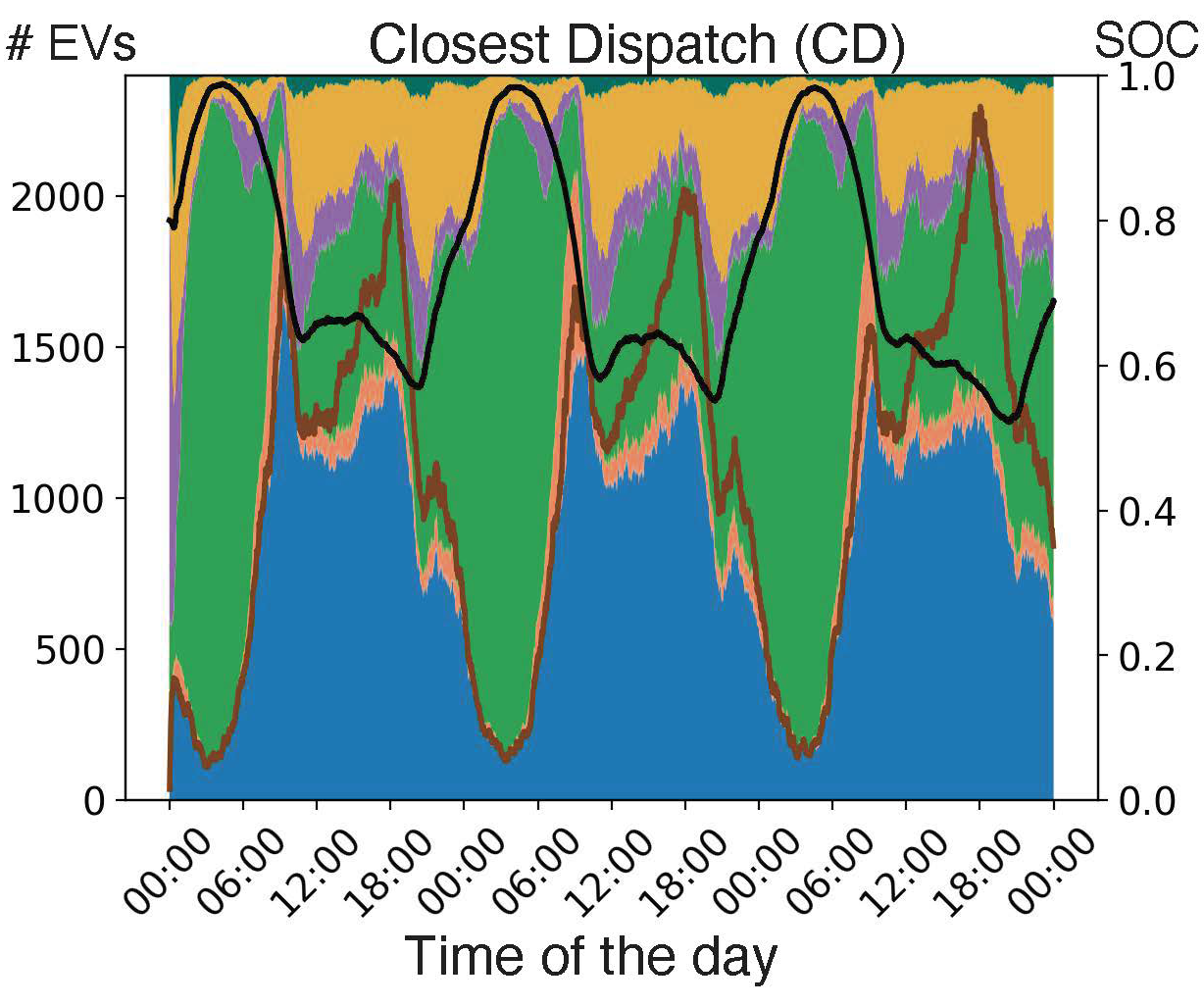}};
        \def \sw {0.5}
        \def \dbs {0.2}
        \begin{scope}[yshift=0.3cm]
        \draw[rectangle] (-1.2, -3.8) -- (18, -3.8) -- (18, -5.4) -- (-1.2, -5.4) -- (-1.2, -3.8);
        \filldraw[rectangle, blue] (-1, -4) -- (-0.5, -4) -- (-0.5, -4.5) -- (-1, -4.5) -- (-1, -4) node[black, right] at (-0.3, -4.25) {Driving with Customer};
        \filldraw[rectangle, orange] (4.2, -4) -- (4.7, -4) -- (4.7, -4.5) -- (4.2, -4.5) -- (4.2, -4)  node[black, right] at (4.9, -4.25) {Picking up Customer};
        \filldraw[rectangle, darkgreen] (10, -4) -- (10.5, -4) -- (10.5, -4.5) -- (10, -4.5) -- (10, -4)  node[black, right] at (10.7, -4.25) {Idle};
        \filldraw[rectangle, purple] (12.8, -4) -- (13.3, -4) -- (13.3, -4.5) -- (12.8, -4.5) -- (12.8, -4) node[black, right] at (13.5, -4.25) {Driving to Charger};
        \filldraw[rectangle, yellow] (-1, -4.7) -- (-0.5, -4.7) -- (-0.5, -5.2) -- (-1, -5.2) -- (-1, -4.7) node[black, right] at (-0.3, -4.95) {Charging};
        \filldraw[rectangle, teal] (4.2, -4.7) -- (4.7, -4.7) -- (4.7, -5.2) -- (4.2, -5.2) -- (4.2, -4.7) node[black, right] at (4.9, -4.95) {Waiting for Charger};
        \draw[black, ultra thick] (10, -4.95) -- (10.5, -4.95) node[black, right] at (10.7, -4.95) {Avg SoC};
        \draw[brown, ultra thick] (12.8, -4.95) -- (13.3, -4.95) node[black, right] at (13.5, -4.95){(Potential) Active Trips};;
        \end{scope}
    \end{tikzpicture}}}{
    System state for CD (left), CAD (center), and Po$9$ (right) for Nissan; $n = 2400$, $m = 225 \times 4$, and $r_c=20$ [kW].
    \label{fig: real_sim_stackplot_cad}}{}
\end{figure}

To further understand the good performance of Po$d$, \Cref{fig: real_sim_stackplot_cad} shows the time series of CD, CAD, and Po$9$ with $T_{P, \max} = 30$ min for Nissan. For CD, a request is lost if the closest vehicle does not have enough battery to serve the request. This allows some vehicles to fully charge at the expense of losing customers. In the figure, CD shows more idling EVs (green) and fewer EVs driving with customers (blue) compared to Po$9$.
The CAD policy performs worse than Po$9$ because it is a greedy policy and does not balance the fleet's SoC. This behavior results in fewer EVs charging (less yellow) between 11 am and 3 pm, as some of the fleet is already fully charged due to the larger imbalance (more green). Consequently, the SoC increases less under CAD during this time, leading to very low SoC later (7-8 pm) and forcing the algorithm to send EVs to charge during peak hours.

Finally, note that most policies tend to perform well for vehicles with larger battery packs\revcolor{---load balancing is more important for shorter-range vehicles.}\label{page:rev-r2-load-important-for-shor}
In particular, CD offers competitive performance (cf., \Cref{sec:limits-cd-varying}).
Regardless, Po$d$ offers consistent and robust performance due to its adaptability and load-balancing features. 
In practice, our results suggest testing Po$1$ \revcolor{to Po$4$} with bounded $T_{P, \max}$ (e.g., 30 \revcolor{or 45 mins}), or Po$T_{P,\max}$ (for optimally chosen $T_{P,\max}$) given its adaptive Po$d$-like nature. %(cf. \Cref{sec:pod-works-adaptive}).

\section{Conclusions}\label{sec:conlusion}

We study fleet and infrastructure planning for electric vehicle ride-hailing systems that are fully controlled by a platform (e.g., autonomous vehicles), providing universal lower and matching upper bounds on the minimum number of vehicles and chargers required to achieve target service levels. Our findings reveal a critical trade-off between fleet size and the number of chargers: having insufficient charging stations can increase the workload in the system (drive-to-charger time) and lead to a higher fleet size requirement. Importantly, the system operator can take advantage of partially charged vehicles to significantly reduce the workload in the system due to pickups. However, to realize the benefits of this ``extra" capacity, a load-balancing type policy is necessary. We verify this both theoretically and via extensive simulations. We also discuss the scope of these insights by examining how demand amplitude, vehicle battery pack size, charging rate, and relocation affect EV capacity planning and when non-EV thinking is more appropriate. Our simulations suggest that a fleet with long-range EVs or with DC fast charging technology can operate as a non-EV system. Nevertheless, our findings reveal the viability of employing a cost-effective fleet with short-range vehicles and Level 2 charging technology when managed appropriately.

Our methodology and results open several avenues for future research.  The characterization of the trade-off between fleet size and the number of chargers provides a framework for optimizing resource allocation based on factors such as fleet maintenance and repair costs, as well as charging infrastructure investments and availability. 
More granular insights could be gleaned from formal analysis of more complex settings—including non-uniform spatial demand patterns (optimizing relocation), general time-varying demand and electricity price variations, and optimal pricing strategies. As these systems scale, their interaction with the electricity grid becomes increasingly critical. For example, charging fleets during low-demand periods could impose prohibitive stress on the grid. For larger fleets, balancing electricity costs and demand fluctuations can be a priority for a company's bottom line. 
Our insights and the Power-of-$d$ policy provide a starting point, but further research is needed to devise optimal management strategies that are aware of the broader electricity supply chain in which these systems operate.

\bibliographystyle{informs2014}
\bibliography{references}

\begin{thebibliography}{59}
\providecommand{\natexlab}[1]{#1}
\providecommand{\url}[1]{\texttt{#1}}
\providecommand{\urlprefix}{URL }

\bibitem[{Abouee-Mehrizi et~al.(2021)Abouee-Mehrizi, Baron, Berman,
  \protect\BIBand{} Chen}]{abouee2021adoption}
Abouee-Mehrizi H, Baron O, Berman O, Chen D (2021) Adoption of electric
  vehicles in car sharing market. \emph{Production and Operations Management}
  30(1):190--209.

\bibitem[{Akbarpour et~al.(2022)Akbarpour, Alimohammadi, Li, \protect\BIBand{}
  Saberi}]{spatial_stanford}
Akbarpour M, Alimohammadi Y, Li S, Saberi A (2022) The value of excess supply
  in spatial matching markets. \emph{Proceedings of the 23rd ACM Conference on
  Economics and Computation}, 62, EC '22 (New York, NY, USA: Association for
  Computing Machinery), ISBN 9781450391504,
  \urlprefix\url{http://dx.doi.org/10.1145/3490486.3538375}.

\bibitem[{Aksin et~al.(2007)Aksin, Armony, \protect\BIBand{}
  Mehrotra}]{aksin2007modern}
Aksin Z, Armony M, Mehrotra V (2007) The modern call center: A
  multi-disciplinary perspective on operations management research.
  \emph{Production and operations management} 16(6):665--688.

\bibitem[{Akturk et~al.(2025)Akturk, Candogan, \protect\BIBand{}
  Gupta}]{akturk2022managing}
Akturk D, Candogan O, Gupta V (2025) Managing resources for shared
  micromobility: Approximate optimality in large-scale systems.
  \emph{Management Science} 71(7):5676--5695.

\bibitem[{{Alto}(2024)}]{link4}
{Alto} (2024) \url{https://ridealto.com/ev}, [Online; accessed 23-May-2024].

\bibitem[{Anderson et~al.(2022)Anderson, Bhargava, Boehm, \protect\BIBand{}
  Parker}]{anderson2022electric}
Anderson EG, Bhargava HK, Boehm J, Parker G (2022) Electric vehicles are a
  platform business: What firms need to know. \emph{California Management
  Review} 64(4):135--154.

\bibitem[{Avci et~al.(2015)Avci, Girotra, \protect\BIBand{}
  Netessine}]{avci2015electric}
Avci B, Girotra K, Netessine S (2015) Electric vehicles with a battery
  switching station: Adoption and environmental impact. \emph{Management
  Science} 61(4):772--794.

\bibitem[{Bauer et~al.(2018)Bauer, Greenblatt, \protect\BIBand{}
  Gerke}]{bauer2018cost}
Bauer GS, Greenblatt JB, Gerke BF (2018) Cost, energy, and environmental impact
  of automated electric taxi fleets in manhattan. \emph{Environmental science
  \& technology} 52(8):4920--4928.

\bibitem[{Bellan(2021)}]{link7}
Bellan R (2021) Revel launches an all-electric rideshare service with a fleet
  of 50 {T}eslas. \url{https://bit.ly/RevelTesla}, [Online; accessed
  05-Feb-2023].

\bibitem[{Benjaafar et~al.(2022)Benjaafar, Wu, Liu, \protect\BIBand{}
  Gunnarsson}]{benjaafar2022dimensioning}
Benjaafar S, Wu S, Liu H, Gunnarsson EB (2022) Dimensioning on-demand vehicle
  sharing systems. \emph{Management Science} 68(2):1218--1232.

\bibitem[{Besbes et~al.(2022)Besbes, Castro, \protect\BIBand{}
  Lobel}]{besbes2022spatial}
Besbes O, Castro F, Lobel I (2022) Spatial capacity planning. \emph{Operations
  Research} 70(2):1271--1291.

\bibitem[{Bibra et~al.(2022)Bibra, Connelly, Dhir, Drtil, Henriot, Hwang,
  Le~Marois, McBain, Paoli, \protect\BIBand{} Teter}]{bibra2022global}
Bibra EM, Connelly E, Dhir S, Drtil M, Henriot P, Hwang I, Le~Marois JB, McBain
  S, Paoli L, Teter J (2022) Global ev outlook 2022: Securing supplies for an
  electric future.

\bibitem[{{BluSmart}(2024)}]{link8}
{BluSmart} (2024) \url{https://en.wikipedia.org/wiki/BluSmart}, [Online;
  accessed 23-May-2024].

\bibitem[{Boudette \protect\BIBand{} Davenport(2021)}]{link3}
Boudette N, Davenport C (2021) {G.M.} will sell only zero-emission vehicles by
  2035.
  \url{https://www.nytimes.com/2021/01/28/business/gm-zero-emission-vehicles.html},
  [Online; accessed 05-Feb-2023].

\bibitem[{Bucklew \protect\BIBand{} Wise(1982)}]{bucklew1982multidimensional}
Bucklew J, Wise G (1982) Multidimensional asymptotic quantization theory with r
  th power distortion measures. \emph{IEEE Transactions on Information Theory}
  28(2):239--247.

\bibitem[{{City of Chicago}(2023)}]{chicago_dataset}
{City of Chicago} (2023) Transportation network providers - trips (2018 -
  2022).
  \url{https://data.cityofchicago.org/Transportation/Transportation-Network-Providers-Trips-2018-2022-/m6dm-c72p/about_data}.

\bibitem[{Daganzo \protect\BIBand{} Smilowitz(2004)}]{daganzo2004bounds}
Daganzo CF, Smilowitz KR (2004) Bounds and approximations for the
  transportation problem of linear programming and other scalable network
  problems. \emph{Transportation science} 38(3):343--356.

\bibitem[{Deng et~al.(2022)Deng, Gupta, \protect\BIBand{}
  Shroff}]{deng2022fleet}
Deng Y, Gupta A, Shroff NB (2022) Fleet sizing and charger allocation in
  electric vehicle sharing systems. \emph{IFAC Journal of Systems and Control}
  22:100210.

\bibitem[{{Ford Media Center}(2021)}]{link2}
{Ford Media Center} (2021) \url{https://bit.ly/FordEVShift}, [Online; accessed
  05-Feb-2023].

\bibitem[{Gan et~al.(2012)Gan, Topcu, \protect\BIBand{} Low}]{gan2012optimal}
Gan L, Topcu U, Low SH (2012) Optimal decentralized protocol for electric
  vehicle charging. \emph{IEEE Transactions on Power Systems} 28(2):940--951.

\bibitem[{Gans et~al.(2003)Gans, Koole, \protect\BIBand{}
  Mandelbaum}]{gans2003telephone}
Gans N, Koole G, Mandelbaum A (2003) Telephone call centers: Tutorial, review,
  and research prospects. \emph{Manufacturing \& Service Operations Management}
  5(2):79--141.

\bibitem[{Garc{\'\i}a-Villalobos et~al.(2014)Garc{\'\i}a-Villalobos, Zamora,
  San~Mart{\'\i}n, Asensio, \protect\BIBand{} Aperribay}]{garcia2014plug}
Garc{\'\i}a-Villalobos J, Zamora I, San~Mart{\'\i}n JI, Asensio FJ, Aperribay V
  (2014) Plug-in electric vehicles in electric distribution networks: A review
  of smart charging approaches. \emph{Renewable and Sustainable Energy Reviews}
  38:717--731.

\bibitem[{Greening \protect\BIBand{} Erera(2021)}]{greening2021effective}
Greening LM, Erera AL (2021) Effective heuristics for distributing vehicles in
  free-floating micromobility systems. Technical report, Working Paper.

\bibitem[{Halfin \protect\BIBand{} Whitt(1981)}]{halfin1981heavy}
Halfin S, Whitt W (1981) Heavy-traffic limits for queues with many exponential
  servers. \emph{Operations research} 29(3):567--588.

\bibitem[{He et~al.(2021)He, Ma, Qi, \protect\BIBand{} Wang}]{he2021charging}
He L, Ma G, Qi W, Wang X (2021) Charging an electric vehicle-sharing fleet.
  \emph{Manufacturing \& Service Operations Management} 23(2):471--487.

\bibitem[{Inside{EV}s(2022)}]{no_of_posts}
Inside{EV}s (2022) Tesla supercharging network: Almost 400 stations added in q4
  2022.
  \url{https://insideevs.com/news/633377/tesla-supercharging-network-2022q4/},
  [Online; accessed 15-Jun-2023].

\bibitem[{Kanoria(2022)}]{spatial_yash}
Kanoria Y (2022) Dynamic spatial matching. \emph{Proceedings of the 23rd ACM
  Conference on Economics and Computation}, 63--64.

\bibitem[{Kaps et~al.(2022)Kaps, Marinesi, \protect\BIBand{}
  Netessine}]{solar_operations}
Kaps C, Marinesi S, Netessine S (2022) When should the off-grid sun shine at
  night? optimum renewable generation and energy storage investments.
  \emph{Optimum Renewable Generation and Energy Storage Investments.(January 5,
  2022)} .

\bibitem[{Kaps \protect\BIBand{} Netessine(2022)}]{kaps2022privately}
Kaps C, Netessine S (2022) Privately-owned battery storage-reshaping the way we
  do electricity. \emph{Available at SSRN} .

\bibitem[{LaMonaca \protect\BIBand{} Ryan(2022)}]{lamonaca2022state}
LaMonaca S, Ryan L (2022) The state of play in electric vehicle charging
  services--a review of infrastructure provision, players, and policies.
  \emph{Renewable and sustainable energy reviews} 154:111733.

\bibitem[{Levin et~al.(2017)Levin, Kockelman, Boyles, \protect\BIBand{}
  Li}]{levin2017general}
Levin MW, Kockelman KM, Boyles SD, Li T (2017) A general framework for modeling
  shared autonomous vehicles with dynamic network-loading and dynamic
  ride-sharing application. \emph{Computers, Environment and Urban Systems}
  64:373--383.

\bibitem[{Lim et~al.(2017)Lim, Mak, \protect\BIBand{} Shen}]{lim2017agility}
Lim MK, Mak HY, Shen ZJM (2017) Agility and proximity considerations in supply
  chain design. \emph{Management Science} 63(4):1026--1041.

\bibitem[{Liu et~al.(2022{\natexlab{a}})Liu, Gong, \protect\BIBand{}
  Ying}]{liu2022large}
Liu X, Gong K, Ying L (2022{\natexlab{a}}) Large-system insensitivity of
  zero-waiting load balancing algorithms. \emph{ACM SIGMETRICS Performance
  Evaluation Review} 50(1):101--102.

\bibitem[{Liu et~al.(2022{\natexlab{b}})Liu, Gong, \protect\BIBand{}
  Ying}]{liu2022steady}
Liu X, Gong K, Ying L (2022{\natexlab{b}}) Steady-state analysis of load
  balancing with coxian-2 distributed service times. \emph{Naval Research
  Logistics (NRL)} 69(1):57--75.

\bibitem[{Liu \protect\BIBand{} Ying(2020)}]{liu2020steady}
Liu X, Ying L (2020) Steady-state analysis of load-balancing algorithms in the
  sub-halfin--whitt regime. \emph{Journal of Applied Probability}
  57(2):578--596.

\bibitem[{Loeb \protect\BIBand{} Kockelman(2019)}]{loeb2019fleet}
Loeb B, Kockelman KM (2019) Fleet performance and cost evaluation of a shared
  autonomous electric vehicle (saev) fleet: A case study for austin, texas.
  \emph{Transportation Research Part A: Policy and Practice} 121:374--385.

\bibitem[{Ma et~al.(2011)Ma, Callaway, \protect\BIBand{}
  Hiskens}]{ma2011decentralized}
Ma Z, Callaway DS, Hiskens IA (2011) Decentralized charging control of large
  populations of plug-in electric vehicles. \emph{IEEE Transactions on control
  systems technology} 21(1):67--78.

\bibitem[{Mak et~al.(2013)Mak, Rong, \protect\BIBand{}
  Shen}]{mak2013infrastructure}
Mak HY, Rong Y, Shen ZJM (2013) Infrastructure planning for electric vehicles
  with battery swapping. \emph{Management science} 59(7):1557--1575.

\bibitem[{Mitzenmacher(1996)}]{mitzenmacher_thesis}
Mitzenmacher M (1996) Load balancing and density dependent jump markov
  processes. \emph{Proceedings of 37th Conference on Foundations of Computer
  Science}, 213--222 (New York, NY, USA: IEEE),
  \urlprefix\url{http://dx.doi.org/10.1109/SFCS.1996.548480}.

\bibitem[{Mitzenmacher(2001)}]{mitzenmacher_power_of_d}
Mitzenmacher M (2001) The power of two choices in randomized load balancing.
  \emph{IEEE Transactions on Parallel and Distributed Systems}
  12(10):1094--1104.

\bibitem[{Mukherjee \protect\BIBand{} Gupta(2014)}]{mukherjee2014review}
Mukherjee JC, Gupta A (2014) A review of charge scheduling of electric vehicles
  in smart grid. \emph{IEEE Systems Journal} 9(4):1541--1553.

\bibitem[{Osorio et~al.(2021)Osorio, Lei, \protect\BIBand{}
  Ouyang}]{osorio2021optimal}
Osorio J, Lei C, Ouyang Y (2021) Optimal rebalancing and on-board charging of
  shared electric scooters. \emph{Transportation Research Part B:
  Methodological} 147:197--219.

\bibitem[{Padilla(2020)}]{link1}
Padilla A (2020) {EXECUTIVE ORDER N-79-20}. \url{https://bit.ly/ExecOrderCali},
  [Online; accessed 05-Feb-2023].

\bibitem[{Paganini et~al.(2022)Paganini, Esp{\'\i}ndola, Marvid,
  \protect\BIBand{} Ferragut}]{paganini2022optimization}
Paganini F, Esp{\'\i}ndola E, Marvid D, Ferragut A (2022) Optimization of
  spatial infrastructure for ev charging. \emph{2022 IEEE 61st Conference on
  Decision and Control (CDC)}, 5035--5041 (Cancun, Mexico: IEEE).

\bibitem[{Qi et~al.(2023)Qi, Zhang, \protect\BIBand{} Zhang}]{qi2023scaling}
Qi W, Zhang Y, Zhang N (2023) Scaling up electric-vehicle battery swapping
  services in cities: A joint location and repairable-inventory model.
  \emph{Management Science} 69(11):6855--6875.

\bibitem[{Reed(2009)}]{reed2009g}
Reed J (2009) The g/gi/n queue in the halfin--whitt regime. \emph{The Annals of
  Applied Probability} 19(6):2211 -- 2269.

\bibitem[{Srinivasa \protect\BIBand{} Haenggi(2009)}]{srinivasa2009distance}
Srinivasa S, Haenggi M (2009) Distance distributions in finite uniformly random
  networks: Theory and applications. \emph{IEEE Transactions on Vehicular
  Technology} 59(2):940--949.

\bibitem[{{Uber}(2020)}]{link5}
{Uber} (2020) Millions of trips a day, zero emissions and a shift to
  sustainable packaging. \url{https://bit.ly/UberSustainability}, [Online;
  accessed 05-Feb-2023].

\bibitem[{{Uber}(2023)}]{link10}
{Uber} (2023)
  \url{https://investor.uber.com/news-events/news/press-release-details/2023/Waymo-and-Uber-Partner-to-Bring-Waymos-Autonomous-Driving-Technology-to-the-Uber-Platform/default.aspx},
  [Online; accessed 02-November-2024].

\bibitem[{Varma et~al.(2022)Varma, Castro, \protect\BIBand{}
  Maguluri}]{varma_sub_halfin_whitt}
Varma SM, Castro F, Maguluri ST (2022) Power-of-$d$ choices load balancing in
  the sub-halfin whitt regime.

\bibitem[{Vosooghi et~al.(2019{\natexlab{a}})Vosooghi, Kamel, Puchinger,
  Leblond, \protect\BIBand{} Jankovic}]{vosooghi2019robo}
Vosooghi R, Kamel J, Puchinger J, Leblond V, Jankovic M (2019{\natexlab{a}})
  Robo-taxi service fleet sizing: assessing the impact of user trust and
  willingness-to-use. \emph{Transportation} 46:1997--2015.

\bibitem[{Vosooghi et~al.(2019{\natexlab{b}})Vosooghi, Puchinger, Jankovic,
  \protect\BIBand{} Vouillon}]{vosooghi2019shared}
Vosooghi R, Puchinger J, Jankovic M, Vouillon A (2019{\natexlab{b}}) Shared
  autonomous vehicle simulation and service design. \emph{Transportation
  Research Part C: Emerging Technologies} 107:15--33.

\bibitem[{Vvedenskaya et~al.(1996)Vvedenskaya, Dobrushin, \protect\BIBand{}
  Karpelevich}]{vvedenskaya_power_of_2}
Vvedenskaya ND, Dobrushin RL, Karpelevich FI (1996) Queueing system with
  selection of the shortest of two queues: An asymptotic approach.
  \emph{Problemy Peredachi Informatsii} 32(1):20--34.

\bibitem[{Wang et~al.(2024)Wang, Zhang, \protect\BIBand{}
  Zhang}]{wang2022demand}
Wang G, Zhang H, Zhang J (2024) On-demand ride-matching in a spatial model with
  abandonment and cancellation. \emph{Operations Research} 72(3):1278--1297.

\bibitem[{Wang et~al.(2016)Wang, Liu, Du, \protect\BIBand{}
  Kong}]{wang2016smart}
Wang Q, Liu X, Du J, Kong F (2016) Smart charging for electric vehicles: A
  survey from the algorithmic perspective. \emph{IEEE Communications Surveys \&
  Tutorials} 18(2):1500--1517.

\bibitem[{{Waymo One}(2024)}]{link9}
{Waymo One} (2024)
  \url{https://support.google.com/waymo/answer/9059119?hl=en&ref_topic=9210304&sjid=14856197520396367784-NC},
  [Online; accessed 23-May-2024].

\bibitem[{Wu et~al.(2022)Wu, Y{\"u}cel, \protect\BIBand{} Zhou}]{wu2022smart}
Wu OQ, Y{\"u}cel {\c{S}}, Zhou Y (2022) Smart charging of electric vehicles: An
  innovative business model for utility firms. \emph{Manufacturing \& Service
  Operations Management} 24(5):2481--2499.

\bibitem[{Yang et~al.(2023)Yang, Levin, Hu, Li, \protect\BIBand{}
  Jiang}]{yang2023fleet}
Yang J, Levin MW, Hu L, Li H, Jiang Y (2023) Fleet sizing and charging
  infrastructure design for electric autonomous mobility-on-demand systems with
  endogenous congestion and limited link space. \emph{Transportation Research
  Part C: Emerging Technologies} 152:104172.

\bibitem[{Zhang et~al.(2021)Zhang, Lu, \protect\BIBand{}
  Shen}]{zhang2021values}
Zhang Y, Lu M, Shen S (2021) On the values of vehicle-to-grid electricity
  selling in electric vehicle sharing. \emph{Manufacturing \& Service
  Operations Management} 23(2):488--507.

\end{thebibliography}

\newpage

\renewcommand{\theHsection}{A\arabic{section}}
\begin{APPENDICES}
\setcounter{page}{1}   % restart numbering
\section*{\centering{\large Online Appendix: Electric Vehicle Fleet and Charging Infrastructure Planning}}
\addcontentsline{toc}{section}{Online Appendix: Electric Vehicle Fleet and Charging Infrastructure Planning}

\begin{center}
\textbf{Sushil Mahavir Varma}$^{1}$ \quad \textbf{Francisco Castro}$^{2}$\quad \textbf{Siva Theja Maguluri}$^{3}$\par
\normalsize
$^1$Industrial and Operations Engineering, University of Michigan Ann Arbor \\
$^{2}$Anderson School of Management, University of California, Los Angeles \\
$^{3}$H.Milton School of Industrial and Systems Engineering, Georgia Institute of Technology 
\end{center}
\section{Universal Lower Bounds} 
\label{appx:proof_sec:min}

We prove \Cref{thm:lower-bound} by analyzing a subset of state variables and describing their evolution through a set of differential equations that are valid under any policy. In Section~\ref{sec:proof_over_thm_lower_bound}, we present a proof overview. In particular, we present the set of differential equations and explain how its analysis leads to \Cref{thm:lower-bound}. In Appendix~\ref{app:lower_bound_ode}, we show how to analyze the set of ODEs, and in Appendix~\ref{app: lower_bound_thm}, we provide a formal proof of \Cref{thm:lower-bound}.

\subsection{Proof Overview for Theorem 1}\label{sec:proof_over_thm_lower_bound}
%We start by introducing some notations. 
Consider an arbitrary policy $\pi \in \Pi$ and let $p^\pol_{dc}(t) \in [0,1]$ be the portion of EVs that are driving to a charger at time $t$, and that eventually reach a charger. Let $T^\pol_C(t)$ be the average time an EV charges at time $t$ (note this is bounded by $\pack/r_c$). The evolution of the number of customers ($q(t)$), the number of EVs charging ($\cc(t)$), and the aggregate SoC ($s(t)$) is:
\begin{subequations} \label{eq:partial_ode}
\begin{align}\label{eq:ode-q}
\dot{q}(t) &=\lambdaeff^{\pol}(t) - \frac{q(t)}{T_R+T_P^{\pol}(t)},\quad  q(0)=q_0, \\ \label{eq:ode-cc}
\dot{c}_C(t) &= p^\pol_{dc}(t)\cdot\frac{\cd(t)}{T^\pol_{DC}(t)} - \frac{\cc(t)}{T^\pol_C(t)} ,\quad  \cc(0)=c_{C0}, \\ \label{eq:ode-soc}
\dot{\soc}(t) &= \cc(t)\cdot r_c - (n-\cc(t)-\ci(t))\cdot r_d, \quad  \soc(0)= \soc_0,
\end{align}
\end{subequations}
for an initial condition $(q_0, c_{C0}, s_0)$. 
Customers exit the system at a rate of $1/(T_R+T_P^{\pol}(t))$. Then $q(t)$ decreases at rate $q(t)/(T_R+T_P^{\pol}(t))$ and increases at rate $\lambdaeff^{\pol}(t)$.
EVs finish charging at rate $1/T^\pol_C(t)$, decreasing $c_C(t)$ at rate $c_C(t)/T^\pol_C(t)$. Additionally, EVs driving to a charger arrive at rate $p_{dc}^\pol(t)/T^\pol_{DC}(t)$, increasing $c_C(t)$ by $p_{dc}^\pol(t) c_{DC}(t)/T^\pol_{DC}(t)$. For $s(t)$, $\cc(t)$ vehicles charging at rate $r_c$ increase the SoC, while $(n-\cc(t)-\ci(t))$ vehicles driving decrease it at rate $r_d$.

The set of equations \eqref{eq:ode-q} to \eqref{eq:ode-soc} partially characterizes the system evolution for an arbitrary policy $\pi$. It turns out that these equations already provide enough information to derive the bounds in \Cref{thm:lower-bound}. We next analyze the equilibrium points of these differential equations. 
\begin{proposition}\label{prop:equilibrium-points-min}
For any policy $\pi$, there exists $(\eq,\ecc,\eci,\ecd)$ such that
\begin{subequations} \label{eq:stability-gen}
\begin{flalign}\label{eq:stability-1}
 \alphasup \lambda \left(\frac{\tau_1}{\sqrt{n-\eq}} + T_R \right)\leq \eq,\\\label{eq:stability-2}
\frac{\ecc}{\sqrt{m-\ecc}} \cdot \frac{r_c}{\pack}\leq \frac{\ecd}{\tau_2},\\\label{eq:stability-3}
\ecc r_c= \left(n-\ecc-\eci\right) r_d.
\end{flalign}
\end{subequations} 
\end{proposition}
We prove this result by analyzing the long-run average of Equations \eqref{eq:ode-q}--\eqref{eq:ode-soc} and using the convexity of the lower bounds for $T_P^{\pol}(t)$ and $T_{DC}^{\pol}(t)$ in \eqref{eq:assumption1} and \eqref{eq:assumption2}. \Cref{eq:stability-1} is analogous to Little's law: the number of customers in the system is equal to the product of the effective arrival rate and the amount of time each customer spends in the system (pickup time plus trip time).
\Cref{eq:stability-2} says that the minimum rate of EVs departing from a charger has to be less than the maximum rate at which they arrive at a charger. 
\Cref{eq:stability-3} balances the aggregate charge and discharge rates. 
We next explain, using approximations, how to derive \Cref{thm:lower-bound} from \Cref{prop:equilibrium-points-min}.
Because of the extra $\Theta(\lambda)$ vehicles charging, the pickup time is $\Omega(\lambda^{-1/2})$ (cf. \Cref{eq:assumption1}). Thus, 
\begin{equation}
    n \geq  \ecrt + \ecc + \ecp +  \ecd \geq T_R\alphasup \lambda + rT_R\alphasup \lambda + \Omega(\lambda^{1/2})  + \ecd \label{eq: n_cars_lb}, 
\end{equation}
    where we have used that the number of EVs driving with a customer is equal to the number of customers on-trip ($T_R\alphasup \lambda$); the number of EVs charging is at least $r$ times $\ecrt$; and the number of EVs picking up customers is $\lambda^{-1/2} \times \lambda = \lambda^{1/2}$. We have
\begin{equation}
    m \stackrel{\eqref{eq:stability-2}}{\geq} \ecc + \left(\frac{r_c\tau_2}{\pack}\right)^2\cdot \left(\frac{\ecc}{\ecd}\right)^2 
    \geq  r\eq +r\ecd +\left(\frac{r_d\tau_2}{\pack}\right)^2\cdot \left(\frac{\eq}{\ecd}\right)^2 
    \stackrel{\eqref{eq:stability-1}}{\geq}
    rT_R\alphasup \lambda + r\ecd + \left(\frac{\lambda}{\ecd}\right)^2, \label{eq: n_chargers_lb}
\end{equation}
where in the second-to-last inequality we used $\ecc = (\eq +\ecd)r$ (by \eqref{eq:stability-3} and \eqref{eq:q_identity}).
Parameterizing $\ecd$ by $\lambda^{1-\gamma}$ (or, equivalently, setting $\beta_2= 2\gamma$), and noting that $\alphasup \geq \alpha$, \eqref{eq: n_cars_lb} and \eqref{eq: n_chargers_lb} yield
$$
n \geq (1+r)T_R \alpha \lambda + \Omega\left(\lambda^{\max\left\{\frac{1}{2}, 1-\gamma\right\}}\right)\quad \text{and}\quad
m \geq rT_R \alpha \lambda + \Omega\left(\lambda^{\max\left\{2\gamma, 1-\gamma\right\}}\right).
$$
If we consider $\gamma > 1/2$ or $\gamma < 1/3$, we increase the infrastructure requirement (increase the lower bounds above) without reducing the minimum fleet size nor the minimum number of chargers.

\subsection{Analyzing the ODE}\label{app:lower_bound_ode}
\proof{Proof of \Cref{prop:equilibrium-points-min}}
For this proof, we restrict ourselves to the sequence that achieves the limit supremum of $\frac{1}{T} \int_0^T \lambdaeff(t) dt$. Now, note that $(q(t),\cc(t),\ci(t),\cd(t))$ are all bounded in $[0,n]$. Hence, if we denote any of these quantities by $y(t)$, we have that $\frac{1}{T}\int_0^T y(t)dt$ is a bounded sequence in a compact set. Therefore, it has a subsequence that converges. We can therefore create a subsequence for which all of the integrals $\frac{1}{T}\int_0^T y(t)dt$  with $y(t)\in \{q(t),\cc(t),\ci(t),\cd(t)\}$ converge. Note that here we are abusing notation and using $T$ instead of explicitly writing the subsequence. We note that this does not affect what we want to prove. Let us use $(\eq,\ecc,\eci,\ecd)$ to denote the limit of each of these integrals (of the subsequences).

We first show \cref{eq:stability-1}. Recall that from \cref{eq:ode-q} we have that
$$
\dot{q}(t) =\lambdaeff^{\pol}(t) - \frac{q(t)}{T_R+T_P^{\pol}(t)}.
$$
Integrating the equation above from $0$ to $T$ and dividing by $T$ yields
\begin{flalign*}
\frac{q(T)-q_0}{T} &=\frac{1}{T}\int_0^T \lambdaeff^{\pol}(t)dt - \frac{1}{T}\int_0^T\frac{q(t)}{T_R+T_P^{\pol}(t)}dt\\
&\stackrel{\eqref{eq:assumption1}}{\geq}\frac{1}{T}\int_0^T \lambdaeff^{\pol}(t)dt - \frac{1}{T}\int_0^T\frac{q(t)}{T_R+\frac{\tau_1}{\sqrt{n-q(t)}}}dt\\
&\geq\frac{1}{T}\int_0^T \lambdaeff^{\pol}(t)dt - \frac{\frac{1}{T}\int_0^T q(t)dt}{T_R+\frac{\tau_1}{\sqrt{n-\frac{1}{T}\int_0^Tq(t)dt}}}
\end{flalign*}
where the last inequality follows from the concavity of $h(x) = x/(T_R +\tau_1/\sqrt{n-x})$ (we provide a proof of this at the end of the current proof) and Jensen's inequality. Taking the 
limit when $T \uparrow \infty$ (and noting that $h(x)$ is a continuous function) we have 
$$
\eq \geq \alphasup\lambda\cdot \left(T_R + \frac{\tau_1}{\sqrt{n-\eq}}\right),
$$
where we have used the definition of $\alphasup$. Now we prove \cref{eq:stability-2}. From \cref{eq:ode-cc} we have 
\begin{flalign*}
\dot{c}_C(t) &= p^\pol_{dc}(t)\cdot\frac{\cd(t)}{T^\pol_{DC}(t)} - \frac{\cc(t)}{T^\pol_C(t)} \\
&\stackrel{(a)}{\leq}  \frac{\cd(t)}{T^\pol_{DC}(t)} - \cc(t)\cdot \frac{r_c}{\pack} \\
&\stackrel{\eqref{eq:assumption2}}{\leq}  \frac{\cd(t)\cdot \sqrt{m-\cc(t)}}{\tau_2} - \cc(t)\cdot \frac{r_c}{\pack}
\end{flalign*}
where in $(a)$ we have used that $p^\pol_{dc}(t)\leq 1$ and that $T^\pol_C(t)\leq \pack/r_c$. Dividing both sides above by $\sqrt{m-\cc(t)}$ and integrating yields 
\begin{flalign}\nonumber
\frac{1}{T}(-2\sqrt{m-\cc(T)} + 2\sqrt{m-\cc(0)})&\leq \frac{1}{T}\int_0^T\frac{\cd(t)}{\tau_2}dt - \frac{1}{T}\int_0^T\frac{\cc(t)}{\sqrt{m-\cc(t)}}\cdot \frac{r_c}{\pack}dt\\\label{eq:ode-cc-proof-1}
&\leq \frac{1}{T}\int_0^T\frac{\cd(t)}{\tau_2}dt - \frac{\frac{1}{T}\int_0^T\cc(t)dt}{\sqrt{m-\frac{1}{T}\int_0^T\cc(t)dt}}\cdot \frac{r_c}{\pack}
\end{flalign}
where we have used that 
$$
\int_{0}^T \frac{\dot{c}_c(t)}{\sqrt{m-\cc(t)}}dt = -2\sqrt{m-\cc(T)} + 2\sqrt{m-\cc(0)},
$$
and the convexity of $x/\sqrt{m-x}$ together with Jensen's inequality. Observe that \cref{eq:ode-cc-proof-1} implies $\frac{1}{T} \int_0^T \cc(t) dt$ is uniformly bounded away from $m$. Noting that $x/\sqrt{m-x}$ is continuous on this restricted domain and taking the limit as $T\uparrow\infty$ in \cref{eq:ode-cc-proof-1} yields
$$
\frac{\ecc}{\sqrt{m-\ecc}}\cdot \frac{r_c}{\pack}\leq \frac{\ecd}{\tau_2}.
$$

Finally, from \eqref{eq:ode-soc} we have
$$
\dot{\soc}(t) = \cc(t)\cdot r_c - (n-\cc(t)-\ci(t))\cdot r_d.
$$
As before we integrate from 0 to $T$, divide by $T$ and take the limit. By doing so, because $\soc(t)$ is bounded, we obtain that
$$
\ecc \cdot r_c = (n-\ecc-\eci)\cdot r_d.
$$
To conclude the proof we show that $h(x) = x/(T_R +\tau_1/\sqrt{n-x})$ is concave. We have
$$
\frac{d^2}{dx^2}h(x) = \frac{\tau_1 \left(x \left(T_R \sqrt{n-x}+3 \tau_1\right)-4 n \left(T_R \sqrt{n-x}+\tau_1\right)\right)}{4 (n-x)^{3/2} \left(T_R \sqrt{n-x}+\tau_1\right)^3},
$$
the denominator is positive, while the numerator is non-positive because
$$
x \left(T_R \sqrt{n-x}+3 \tau_1\right) \leq n \left(T_R \sqrt{n-x}+3 \tau_1\right)\leq 4n \left(T_R \sqrt{n-x}+ \tau_1\right).
$$
\endproof
\proof{Proof of \Cref{prop:lower-bound-n}}
By \cref{eq:stability-3} we have 
$$
\ecc = (n - \ecc -\eci) r \Longleftrightarrow \ecc = (n-\eci) \frac{r}{1+r}.
$$
Therefore, 
$$
n = \eq+\ecc+ \eci+\ecd \geq \eq+\ecc+ \eci= \eq+ n\frac{r}{1+r} + \frac{\eci}{1+r}\geq \eq+ n \frac{r}{1+r}.
$$
Now, from \cref{eq:stability-1} we have that $\eq\geq T_R\alphasup \lambda$. Hence, we have 
$$\frac{n}{1+r} = n -  n \frac{r}{1+r}\geq \eq \geq T_R\alphasup \lambda \Longrightarrow n\geq (1+r)T_R\alphasup \lambda.$$
Now, by \eqref{eq:stability-2}, we get
$$
m \geq \ecc + \left(\frac{r_c\tau_2\ecc}{\pack \ecd }\right)^2 \geq \ecc \overset{\eqref{eq:stability-3}}{=} (n-\ecc-\eci)r \overset{\eqref{eq:total_num_vehicles}}{=} (\eq+\ecd)r \geq r\eq \overset{\eqref{eq:stability-1}}{\geq} rT_R \alphasup \lambda.
$$
This completes the proof.
\endproof

\subsection{Proof of Theorem \ref{thm:lower-bound}} \label{app: lower_bound_thm}
\proof{Proof of Theorem \ref{thm:lower-bound}}
We prove the contrapositive. In particular, we show that if for all $\kappa_1, \kappa_2 > 0$ there exists $\lambda > 0$ such that for all $\gamma \in [1/3, 1/2]$, we have
\begin{align} \label{eq: lower_bound_contrapositive}
n \leq (1+r)T_R \alpha \lambda + \kappa_1 \lambda^{1-\gamma} \quad \text{OR} \quad m \leq rT_R \alpha \lambda + \kappa_2 \lambda^{2\gamma},
\end{align}
then, $\alphasup < \alpha$. Now, we fix the following expressions for $\kappa_1, \kappa_2$:
\begin{align*}
    \kappa_1 = \frac{\tau_1}{2(\tau_1+T_R)}, \quad \kappa_2 = \min\left\{\frac{3}{2}r_d^{4/3}\left(\frac{\delta T_R\tau_2}{r_c p^\star}\right)^{2/3}, r(r+1)\left(\frac{r_d T_R^2 \tau_2 \delta}{p^\star \tau_1}\right)^2\right\}.
\end{align*}
For the above defined $\kappa_1, \kappa_2 > 0$, there exists $\lambda > 0$ such that \eqref{eq: lower_bound_contrapositive} is satisfied for all $\gamma \in [1/3, 1/2]$. Fix such a $\lambda > 0$ for the rest of the proof. Now, let $\Gamma_1 = \left\{\gamma \in \left[\frac{1}{3}, \frac{1}{2}\right] : n \leq (1+r)T_R \alpha \lambda + \kappa_1\lambda^{1-\gamma}\right\}$ and $\Gamma_2 = \left\{\gamma \in \left[\frac{1}{3}, \frac{1}{2}\right] : m \leq rT_R \alpha \lambda + \kappa_2\lambda^{2\gamma}\right\}$. Note that $\Gamma_1 \cup \Gamma_2 = [1/3, 1/2]$. Now, we consider three cases.

\textbf{Case I ($\Gamma_2 = \emptyset$):} In this case, we have $\Gamma_1 = [1/3, 1/2]$ and so $n \leq (1+r)T_R \alpha \lambda + \kappa_1 \lambda^{1-\gamma}$ for all $\gamma \in [1/3, 1/2]$. By using \eqref{eq:stability-1}, we get
\begin{align*}
    \lambda\alphasup &\leq \frac{\eq}{\frac{\tau_1}{\sqrt{n-\eq}}+T_R} \overset{(a)}{\leq}   \frac{n/(1+r)}{\frac{\tau_1}{\sqrt{n-\eq}}+T_R} - \frac{\ecd}{\frac{\tau_1}{\sqrt{n-\eq}} + T_R} \overset{(b)}{\leq} \frac{n/(1+r)}{\frac{\tau_1}{\sqrt{n}}+T_R} - \frac{\ecd}{\frac{\tau_1\sqrt{1+r}}{\sqrt{r}} + T_R} \\
    &= \frac{n}{(1+r)T_R}\left(1+\frac{\tau_1}{T_R\sqrt{n}}\right)^{-1}- \frac{\ecd}{\frac{\tau_1\sqrt{1+r}}{\sqrt{r}} + T_R}\\
    &\overset{(c)}{\leq} \frac{n}{(1+r)T_R}\left(1-\frac{\tau_1}{\tau_1+T_R}\frac{1}{\sqrt{n}}\right)- \frac{\ecd}{\frac{\tau_1\sqrt{1+r}}{\sqrt{r}} + T_R} \\
    &=\frac{n}{(1+r)T_R} - \frac{\tau_1}{T_R(1+r)(\tau_1+T_R)}\sqrt{n}- \frac{\ecd}{\frac{\tau_1\sqrt{1+r}}{\sqrt{r}} + T_R} \\
    &\leq \alpha \lambda + \frac{\kappa_1}{(1+r)T_R} \lambda^{1-\gamma}- \frac{\tau_1}{T_R(1+r)(\tau_1+T_R)}\sqrt{n}- \frac{\ecd}{\frac{\tau_1\sqrt{1+r}}{\sqrt{r}} + T_R} \numberthis \label{eq: cdc_upper_bound} \\
    &\overset{(d)}{<} \alpha\lambda
\end{align*}
where $(a)$ follows by evaluating $\eq$ as follows:
\begin{align}
    \eq &= n - \ecc - \eci - \ecd \leq n - \ecc - \eci \frac{r}{1+r} - \ecd\overset{\eqref{eq:stability-3}}{=}  \frac{n}{1+r}- \ecd. \label{eq:upper_bound_q}
\end{align}
Next, $(b)$ follows by lower bounding $\eq \geq 0$ in the first term and upper bounding $\eq \leq n/(1+r)$ along with $n\geq 1$ in the second term. Further, $(c)$ follows as 
$$
\frac{1}{1+\frac{C}{\sqrt{n}}} \leq 1 - \frac{C}{C+1}\frac{1}{\sqrt{n}} \quad \forall n \geq 1.
$$
Lastly, $(d)$ follows by picking $\gamma = 1/2$ and noting that $\kappa_1 = \frac{\tau_1}{2(\tau_1+T_R)}$. This completes Case I as the last inequality implies $\alphasup < \alpha$.

\textbf{Case II ($\Gamma_1 = \emptyset$):} In this case, we have $\Gamma_2 = [1/3, 1/2]$ and so $m \leq rT_R\alpha \lambda + \kappa_2\lambda^{2\gamma}$ holds for all $\gamma \in [1/3, 1/2]$. By using \eqref{eq:stability-2}, we get
\begin{align*}
m & \geq \ecc + \left(\frac{r_c\tau_2}{\pack}\right)^2\left(\frac{\ecc}{\ecd}\right)^2 \nonumber \\
    &\overset{(a)}{\geq} r \eq + r \ecd  + \left(\frac{r_c\tau_2}{\pack}\right)^2 \left(\frac{r \eq}{\ecd}\right)^2 \nonumber \\
    &\overset{(b)}{\geq} rT_R \alphasup \lambda + r\ecd  + \left(\frac{r_d T_R\tau_2}{\pack}\right)^2\left(\frac{\alphasup \lambda}{\ecd}\right)^2 \numberthis \label{eq: cdc_lower_bound} \\
    &\overset{(c)}{\geq} rT_R \alphasup \lambda + \frac{3}{2^{2/3}} r_d^{4/3} \left(\frac{T_R\tau_2}{\pack r_c}\right)^{2/3} \left(\alphasup \lambda\right)^{2/3},\nonumber \\
    &\overset{(d)}{\geq} rT_R \alphasup \lambda +  \frac{3}{2^{2/3}} r_d^{4/3} \left(\frac{T_R\tau_2}{\pack r_c}\right)^{2/3} \delta^{2/3} \lambda^{2/3}
\end{align*}
where $(a)$ follows by \eqref{eq:stability-3} and \eqref{eq:total_num_vehicles}. Specifically, $\ecc =  r(n-\ecc-\eci)= r(\eq + \ecd)$. Next, $(b)$ follows by \eqref{eq:stability-1} as $\eq \geq  \alphasup \lambda T_R + \alphasup \lambda \frac{\tau_1}{\sqrt{n-\eq}} \geq  \alphasup \lambda T_R$. Further, $(c)$ follows by noting that $c_1x + c_2/x^2$ is minimized at $x = (2c_2/c_1)^{1/3}$, which implies
$$
c_1x +\frac{c_2}{x^2} \geq \frac{3}{2^{2/3}} c_1^{2/3}c_2^{1/3}.
$$
Lastly, $(d)$ follows as $\alpha \geq \delta$. Now, by substituting $m \leq rT_R \alpha \lambda + \kappa_2 \lambda^{2\gamma}$, picking $\gamma = 1/3$, and noting that $\kappa_2 < \frac{3}{2^{2/3}} r_d^{4/3} (T_R\tau_2 / (r_c\pack))^{2/3}\delta^{2/3}$, we obtain that $\alphasup < \alpha$.

\textbf{Case III ($\Gamma_1, \Gamma_2 \neq \emptyset$):} Let $\gamma_1 = \sup\{\gamma \in \Gamma_1\}$ and $\gamma_2 = \inf\{\gamma \in \Gamma_2\}$. As $\Gamma_1, \Gamma_2 \neq \emptyset$, we have $\gamma_1, \gamma_2 \in [1/3, 1/2]$. First note that $(\gamma_1 + \gamma_2)/2 \in [1/3, 1/2]$. Now, if $\gamma_1 < \gamma_2$, then by the definition of $\gamma_1, \gamma_2$, we have $(\gamma_1+\gamma_2)/2 \notin \Gamma_1 \cup \Gamma_2$ as $\gamma_1 < (\gamma_1+\gamma_2)/2 < \gamma_2$. This is a contradiction as $\Gamma_1 \cup \Gamma_2 = [1/3, 1/2]$ and $(\gamma_1 + \gamma_2)/2 \in [1/3, 1/2]$. Thus, we must have $\gamma_1 \geq \gamma_2$.

The rest of the proof is by contradiction. Let $\alphasup \geq \alpha$, then by \eqref{eq: cdc_upper_bound} and \eqref{eq: cdc_lower_bound}, we get
\begin{align*}
\ecd &\leq \frac{\kappa_1}{(1+r)T_R}\left(\frac{\tau_1\sqrt{1+r}}{\sqrt{r}} + T_R\right)\lambda^{1-\gamma}  \leq \frac{\tau_1}{2T_R \sqrt{r(r+1)}}\lambda^{1-\gamma} \quad \forall \gamma \in \Gamma_1, \\
\ecd &\geq \frac{ r_dT_R \tau_2 \delta}{\pack\sqrt{\kappa_2}} \lambda^{1-\gamma} \quad \forall \gamma \in \Gamma_2.
\end{align*}
Now, by taking the sup over $\Gamma_1$ and inf over $\Gamma_2$ in the above equation, we get
\begin{align*}
   \frac{ r_dT_R \tau_2 \delta}{\pack\sqrt{\kappa_2}} \lambda^{1-\gamma_2} \leq  \ecd \leq \frac{\tau_1}{2T_R \sqrt{r(r+1)}}\lambda^{1-\gamma_1} \leq \frac{\tau_1}{2T_R \sqrt{r(r+1)}}\lambda^{1-\gamma_2},
\end{align*}
where the last inequality follows as $\gamma_1 \geq \gamma_2$. Now, as $\kappa_2 \leq r(r+1)\left(\frac{r_d T_R^2 \tau_2 \delta}{p^\star \tau_1}\right)^2$ we get a contradiction. This completes the proof of the theorem. \hfill $\square$
\endproof

\section{Matching Upper Bounds} \label{app: matching_upper_bound}
\subsection{Proof of Lemma \ref{lemma: coupling}}
\proof{Proof of Lemma \ref{lemma: coupling}}
     Consider a coupled modified and original model, i.e., they both have an identical arrival sequence and an identical initial state. Now, let $\pi_a$ be such that the dispatching decisions in the original model are exactly the same as in the modified model operating under $\pi_s$, where $\pi_s$ is an admissible policy for the modified model. In particular, let $\pi_a$ be such that for any $i \in [n]$, EV $i$ serves the same customer and follows it up by driving to the same charger as in the modified model under $\pi_s$. Under such a coupling, the location of all the EVs is identical under the two models. In addition, the following is true.
     \begin{itemize}
        \item EV $i$ is picking up customer $k$ in the original model if and only if it is picking up customer $k$ in the modified model.
        \item EV $i$ is driving with customer $k$ in the original model if and only if it is driving with customer $k$ in the modified model.
        \item EV $i$ is driving to charger $j \in [m]$ in the original model if and only if it is driving to charger $j$ in the modified model.
        \item If charger $j$ is available in the modified model, then it is available in the original model.
        \item In the modified model, EV $i$ expends at least as much SoC as in the original model by picking up and dropping the customer, followed by driving to the charger.
         \item If EV $i$ is charging in the modified model, then it is either charging in the original model or it is idling at the charger location with 100\% SoC in the original model. On the other hand, if EV $i$ is charging in the original model, then it is charging in the modified model.  
     \end{itemize}
Using the above information, we conclude that $\operatorname{SoC}_i^a(t) \geq \operatorname{SoC}_i^s(t)$ for all $t \in \bbR_+$, $i \in [n]$. That is, the SoC of each EV is path-wise larger in the original model. Lastly, if $\pi_s$ is admissible for the modified model, then we have $\operatorname{SoC}_i^s(t) \geq 0$ for all $t \in \bbR_+$, $i \in [n]$. Thus, $\pi_a$ is admissible for the original model as $\operatorname{SoC}_i^a(t) \geq \operatorname{SoC}_i^s(t) \geq 0$ for all $t \in \bbR_+$, $i \in [n]$. This completes the proof. \hfill $\square$
\endproof

\subsection{Power-of-\texorpdfstring{$d$}{d} Vehicles Dispatch: Analysis} \label{sec:pod-analysis}

We are ready to prove \Cref{thm:informal_upper-bound}. The proof has three parts. First, we analyze the system of differential equations in \eqref{eq:odes-power-d}. We characterize its equilibrium points and show existence and uniqueness. Then, under a stability assumption and a given selection of the admission control and infrastructure parameters, we analyze the unique equilibrium of the system and obtain
a lower bound on the number of busy vehicles in equilibrium. Finally, we translate this bound into a lower bound for $\alphainf$ which, in turn, delivers the achievability result in \Cref{thm:informal_upper-bound}.

We set the fleet size, number of chargers, and the admission control parameters $\kappa$ and $A$ to be
\begin{equation} \label{eq: in_terms_of_n}
n = \alpha (1+r) T_R \lambda + \kappa_1 \lambda^{1-\gamma}, \quad m = \alpha r T_R \lambda + \kappa_2\lambda^{2\gamma}, \quad \kappa = \frac{2+r}{2(1+r)}\kappa_1, \quad A = m -  \lambda^{2\gamma},
\end{equation}
with $\kappa_1 = 32\tau_2$ and $\kappa_2 = r\kappa_1 + 4\tau_2 + 1$. Note that these choices of $m$ and $n$ are indeed consistent with \Cref{thm:informal_upper-bound}. Now, we present two propositions that are key in the proof of \Cref{thm:informal_upper-bound}. 
\begin{proposition}\label{prop:scalings-infty}  Let $\kappa$, $n$, $m$, and $A$ be given by \eqref{eq: in_terms_of_n} with $\Nt \geq 2$. There exist $\lambda_{p1} > 0$ such that if 
$$
    d =  \lambda^{\gamma/\Nt} (\log \lambda)^{1/\Nt}  \quad\text{with}\quad
    \gamma \in \left[\frac{1}{3}, \frac{1}{2+1/\Nt}\right),
$$
then for all $\lambda \geq \lambda_{p1}$, the dynamical system \eqref{eq:odes-power-d} admits a unique equilibrium $(\mathbf{\XC}, \mathbf{\XB})$  such that:
$$
\XC_{\Nt} \leq \alpha r T_R \lambda + \frac{(1+2r)\kappa_1}{2(1+r)}\lambda^{1-\gamma},
 \quad \text{and} \quad 
  \XB_{\Nt} \geq \alpha T_R \lambda + \frac{\kappa_1}{2(1+r)}\lambda^{1-\gamma}.
$$
\end{proposition}
The above proposition proves that \eqref{eq:odes-power-d} exhibits a unique equilibrium, and provides bounds on
$\XC_{\Nt}$  and $\XB_{\Nt}$. To prove the proposition, we set $(\mathbf{\dXC}, \mathbf{\dXB})$ to zero in \eqref{eq:odes-power-d} and analyze the corresponding system of non-linear equations. Using this system of equations, we obtain a non-linear equation in one variable, which we show to have a unique solution, thereby implying the existence of a unique fixed point. The analysis of this solution delivers the bounds
%Then, by further manipulations and by carefully choosing $\kappa_i : i \in [3]$, we obtain the required bounds 
on $\XC_{\Nt}$ and $\XB_{\Nt}$. In particular, we show that the number of busy vehicles is at least $\alpha T_R\lambda + \Theta(\lambda^{1-\gamma})$. The choice of $d$ corresponds to the minimum value of $d$ such that the lower bound on $\XB_{\Nt}$ holds.
%Recall that $\XB_{\Nt}$ corresponds to all the busy vehicles in the system. 
Now, in equilibrium, \eqref{eq:powerd-3} results in $\tilde{\lambda}p_0I = \XB_{\Nt} \mu$, where $p_0, I, \mu$ are the equilibrium quantities corresponding to $p_0(t), I(t), \mu(t)$, respectively. Note that $\tilde{\lambda}p_0 I$ is exactly the effective arrival rate of customers, i.e., the arriving customers that are served (cf. \Cref{def: admission_control}). Indeed, out of the admitted customers, a $p_0$ fraction is not served by the Power-of-$d$ vehicles dispatch policy due to the low SoC of the selected EVs. Thus, the effective rate of customers is proportional to the number of busy vehicles, and so, the lower bound on $\XB_{\Nt}$ translates into a lower bound on $\alphainf$. We obtain said lower bound on $\alphainf$ in the following proposition. 
\begin{proposition} \label{prop: lower_bound_lambda} Under the same setting as in Proposition \ref{prop:scalings-infty}, if the dynamical system \eqref{eq:odes-power-d} is globally stable, then there exists $\lambda_{p_2} > 0$ such that for all $\lambda \geq \lambda_{p_2}$ and for any initial condition $(\XC_0, \XB_0)$, we have $\alphainf \geq \alpha$. In fact, there exists $T_0(\lambda)>0$ such that for all $T\geq T_0$, we have
$$
\frac{1}{T}\int_0^T \lambdaeff^{\pol(d)}(t)dt \geq \alpha \lambda.
$$
\end{proposition}
As a proof sketch of the above proposition, consider the following back-of-the-envelope calculation: $$\lambda \alphasup = \tilde{\lambda}p_0 I = \XB_{\Nt} \mu \approx \frac{\XB_{\Nt}}{T_R} \geq \alpha \lambda.$$
Observe that \Cref{thm:informal_upper-bound} follows immediately by \Cref{prop: lower_bound_lambda} by setting $\lambda_0 = \lambda_{p2}$. Furthermore, \Cref{prop: lower_bound_lambda} goes beyond the limiting characterization of the service level $(\alphainf)$ and provides a finite time lower bound on the service level.

Finally, note that while we have established the existence and uniqueness of an equilibrium point in the dynamical system \eqref{eq:odes-power-d}, we have assumed global stability. In \Cref{sec: simulations}, we conduct hundreds of simulations with different initial conditions and observe that the system converges to a fixed point in all instances. Furthermore, we note that the ODE for Power-of-$d$ vehicles dispatch policy is similar to the ODE for Power-of-$d$ choices in load balancing (cf. \Cref{app: matching_algo_comparison})
which is known to be globally stable \citep{mitzenmacher_power_of_d, mitzenmacher_thesis}. However, a central difference between our setting and load balancing is that we have two types of vehicles---busy and idle/charging. This distinction, in addition to highly non-linear ODEs, makes the stability analysis challenging. Establishing global stability is then an interesting open problem that may have broader applications for load balancing in queueing systems. 

\subsection{Fixed point for Power-of-d Vehicles Dispatch}
\label{pg:rev2-prop5}
\proof{Proof of \Cref{prop:scalings-infty}}
Consider the system of differential equations in \eqref{eq:odes-power-d}. Any $(\XC,\XB)$ of this dynamical system is a fixed point if and only if the following are satisfied:
\begin{subequations} \label{eq:odes-power-d-eq}
\begin{flalign}\label{eq:powerd-1-eq}
\lambda \left(p_0-p_j\right) I = \XB_{j}\mu, \: 1\le j< \Nt,\\\label{eq:powerd-2-eq}
\lambda \left(p_0-p_j\right) I + \left(\tXC_j-\tXC_{j-1}\right)/(rT_B) =\XB_{j+1} \mu,\: 1\le j<\Nt,\\\label{eq:powerd-3-eq}
\lambda p_0 I = \XB_{\Nt} \mu, \\\label{eq:powerd-4-eq}
\tXC_0 /(rT_B) =\XB_{1} \mu, \\
0 \leq \XC_0 \leq \XC_1 \leq \hdots \XC_{\Nt}, \quad 0\leq \XB_1 \leq \XB_2 \leq \hdots \XB_{\Nt}, \label{eq:powerd-5-eq} \\
\XC_{\Nt} +\XB_{\Nt} = n \label{eq:powerd-6-eq}
\end{flalign}
\end{subequations} 
where $I = \1{\XB_{\Nt}\leq \alpha T_R \lambda + \frac{2+r}{2(1+r)}\kappa_1\lambda^{1-\gamma}}$. To establish the proposition, we start with a proof preamble in which we state and prove facts that will be used at various points during the proof. After the preamble, we split the proof into four main steps: 
\begin{enumerate}
\item [\textbf{(a)}] (\textbf{System simplification}) We first write \eqref{eq:odes-power-d-eq} exclusively in terms of $\XC$ by showing $\tXC_j=r T_B \mu \XB_{j+1}$ for $j\in\{0,\dots,\Nt-1\}$ and $\tXC_j = \XC_j$ for all $j \in \{0, 1, \hdots, \Nt-1\}$.
\item [\textbf{(b)}] (\textbf{Reduction to a single variable}) From part $(a)$, we obtain $\Nt+1$ equations and $\Nt+1$ unknowns. We use $\Nt$ of these equations to show that there exists continuous and strictly decreasing functions $h_{j}$ such that $\XC_{j} = h_{j}(\XC_{\Nt})$ for all $j \in \{0, 1, \hdots, \Nt-1\}$. The last unused equation can then be written only in terms of $\XC_{\Nt}$: the system reduces to finding the root of a function in one dimension.
\item [\textbf{(c)}] (\textbf{Existence and Uniqueness}) We employ the intermediate value theorem to conclude the existence of a solution of the last unused equation. We also exploit the monotonicity properties of $h_j$ to establish uniqueness. In particular, we conclude the existence of a unique $X^\star \leq \alpha r T_R \lambda + \frac{(1+2r)\kappa_1}{2(1+r)}\lambda^{1-\gamma}$ such that $\XC_j = h_j(X^\star)$ for $j \in \{0, 1, \hdots, \Nt-1\}$ and $\XC_{\Nt}=X^\star$ solves \eqref{eq:odes-power-d-eq} up to \eqref{eq:powerd-5-eq} (along with the $(\XB, \tilde{\XC})$ obtained in part (a)).
\item [\textbf{(d)}] (\textbf{Final verifications}) We end by showing that $X^\star$ verifies the inequalities in \eqref{eq:powerd-5-eq}.
\end{enumerate}
\underline{\textbf{Proof preamble:}} Note that $p_j$ is given by \eqref{eq:prob-without-replacement} when $\XC_j \geq d-1$ and $p_j=1$ otherwise. So, we have
\begin{equation} \label{eq: bounds_on_p_j}
 p_j = 1- \frac{\Pi_{i=0}^{d-1} \left[\XC_j-i\right]^+}{\Pi_{i=0}^{d-1} (\XC_{\Nt}-i)} \geq 1-\left(\frac{X_j^C}{X_{\Nt}^C}\right)^d
\end{equation}
where $\frac{0}{0}$ should be understood as $0$. Note that the lower bound corresponds to the complement probability of choosing $d$ vehicles, uniformly at random \emph{with replacement}, that can serve $j$ (but no more than $j$) or fewer trips.
And
$$
\frac{1}{\mu} = \tau_1 \sqrt{\frac{d}{n-X_{\Nt}^B}} + T_R + \frac{\tau_2}{\sqrt{m-\min\left\{X_{\Nt-1}^C, A\right\}}}.
$$
Next, we present upper and lower bounds on $T_B$ that will be useful later in the proof. First note that $I=1$ as if $I=0$ then $\XB_{\Nt} > \alpha T_R \lambda + \frac{2+r}{2(1+r)}\kappa_1\lambda^{1-\gamma}$ and from \eqref{eq:powerd-3-eq}, we would get that $\XB_{\Nt}=0$, a contradiction. Thus, $\XB_{\Nt} \leq \alpha T_R \lambda + \frac{2+r}{2(1+r)}\kappa_1\lambda^{1-\gamma}$ which implies
\begin{align}
    T_R \leq T_B &\leq \tau_1 \sqrt{\frac{d}{n-\alpha T_R\lambda-\frac{2+r}{2(1+r)}\kappa_1\lambda^{1-\gamma}}} + T_R + \frac{\tau_2}{\sqrt{m-A}} \nonumber \\
    &\overset{*}{\leq}  \frac{\tau_1}{\sqrt{\alpha rT_R}}\lambda^{\frac{\gamma/\Nt-1}{2}}\log\lambda + T_R + \tau_2\lambda^{-\gamma} \overset{**}{\leq} T_R + 2\tau_2\lambda^{-\gamma}
    \label{eq: upper_bound_t_b}
\end{align}
where $(*)$ follows as $n = \alpha(1+r)T_R \lambda+\kappa_1\lambda^{1-\gamma}$ and $A = m - \lambda^{2\gamma}$. Next, $(**)$ follows for all $\lambda \geq \lambda_{00}$ for some $\lambda_{00} > 0$ as $-\gamma > \frac{\gamma/\Nt-1}{2}$.

\underline{\textbf{Proof of (a):}} \textbf{System simplification}

We first show (a). Combining \eqref{eq:powerd-1-eq} and \eqref{eq:powerd-2-eq}, we deduce that 
$$
\tXC_j-\tXC_{j-1} = r T_B \mu\left(\XB_{j+1}-\XB_j\right) \quad \forall j\in\{1,\dots,\Nt-1\}.
$$
By carrying out telescopic sum from $j$ to $1$, we get
$$
\tXC_j-\tXC_{0} = rT_B \mu \left(\XB_{j+1}-\XB_1\right) \quad \forall j\in\{1,\dots,\Nt-1\}.
$$
Using \eqref{eq:powerd-4-eq} we obtain that $\tXC_{j}=r T_B \mu \XB_{j+1}$ for $j\in\{0, 1,\dots,\Nt-1\}$. Now, we establish $\XC_{\Nt-1} \leq A$. Assume not, then as $\XB_{\Nt} = \tXC_{\Nt-1}/(rT_B\mu)$, we have 
\begin{align*}
    \XB_{\Nt} = \frac{A}{rT_B \mu } &\geq \frac{T_R}{rT_B}\left(\alpha rT_R\lambda + (\kappa_2-1)\lambda^{2\gamma}\right) \overset{\eqref{eq: upper_bound_t_b}}{\geq} \frac{\alpha T_R\lambda + \frac{\kappa_2-1}{r}\lambda^{2\gamma}}{1+\frac{2\tau_2}{T_R}\lambda^{-\gamma}} \overset{*}{\geq} \left(1-\frac{2\tau_2}{T_R}\lambda^{-\gamma}\right)\left(\alpha T_R\lambda + \frac{\kappa_2-1}{r}\lambda^{2\gamma}\right) \\
    &\overset{**}{\geq} \alpha T_R\lambda+\frac{\kappa_2-1}{r}\lambda^{2\gamma}-4\alpha\tau_2\lambda^{1-\gamma} > \alpha T_R \lambda + \frac{2+r}{2(1+r)}\kappa_1 \lambda^{1-\gamma},
\end{align*}
where the first inequality follows as $1/\mu \geq T_R$. Next, $(*)$ follows as $1/(1+x) \geq 1-x$ for $x \geq 0$. Further, $(**)$ follows for $\lambda \geq \lambda_{01}$ for some $\lambda_{01} > 0$ as $2\gamma < 1$, implying $\frac{2\tau_2}{T_R}\lambda^{-\gamma}\left(\alpha T_R\lambda + \frac{\kappa_2-1}{r}\lambda^{2\gamma}\right) \leq 4\alpha \tau_2 \lambda^{1-\gamma}$. Finally, the last inequality follows by noting that $\alpha<1$, $\gamma \geq 1/3$, and $\kappa_2 = 1 + 4\tau_2 + r\kappa_1 > 1+4r\tau_2+r\frac{2+r}{2(1+r)}\kappa_1$ as $r<1$. We now arrive at a contradiction as $I=1$ implying $\XB_{\Nt} \leq \alpha T_R \lambda+\frac{2+r}{2(1+r)}\kappa_1\lambda^{1-\gamma}$. Thus, we have $\XC_{j} \leq \XC_{\Nt-1} \leq A$ for all $j \in \{0, 1, \hdots, \Nt-1\}$, which implies $\tXC_j = \XC_j$ for all $j \in \{0, 1, \hdots, \Nt-1\}$. So, we conclude that \eqref{eq:odes-power-d-eq} is equivalent to solving the following set of equations:
\begin{subequations} \label{eq:eq-XC-a}
\begin{flalign}\label{eq:eq-XC-a-1}
\XC_{j-1}=\XC_{\Nt-1}-rT_B \lambda p_j , \: 1\le j< \Nt,\\\label{eq:eq-XC-a-2}
\XC_{\Nt-1}=rT_B\lambda p_0, \\\label{eq:eq-XC-a-3}
\XC_{\Nt} + \frac{\XC_{\Nt-1}}{rT_B \mu}=n, \\
0 \leq \XC_0 \leq \XC_1 \leq \hdots \leq \XC_{\Nt-1} \leq \min\{\XC_{\Nt}, A\}, \label{eq:eq-XC-a-4} \\
0 \leq X^B_{\Nt}=(n-X^C_{\Nt}) \leq \alpha T_R\lambda + \frac{2+r}{2(1+r)}\kappa_1 \lambda^{1-\gamma} \label{eq:eq-XC-a-5}
\end{flalign}
\end{subequations} 
In addition to replacing $(\XB, \tXC)$ in terms of $\XC$, note that we also replaced $p_0$ with $\XC_{\Nt-1}$ in the first equation owing to \eqref{eq:eq-XC-a-2}. Such a substitution simplifies the analysis in the rest of the proof. We can also write $\mu$ in terms of $\XC$ as follows:
\begin{align*}
    \frac{1}{\mu} &= \tau_1 \sqrt{\frac{d}{n-X_{\Nt}^B}} + T_R + \frac{\tau_2}{\sqrt{m-\min\left\{X_{\Nt-1}^C, A\right\}}} = \tau_1 \sqrt{\frac{d}{X_{\Nt}^C}} + T_R + \frac{\tau_2}{\sqrt{m-X_{\Nt-1}^C}}
\end{align*}
The expression of $\mu$ along with \eqref{eq:eq-XC-a} concludes part $(a)$ of the proof.

\underline{\textbf{Proof of (b)}:} \textbf{
Reduction to a single variable
}

The proof strategy for this part of the proof is as follows: We start with \eqref{eq:eq-XC-a-3} to show that $\XC_{\Nt-1}$ is a continuous, decreasing function of $\XC_{\Nt}$ and then work our way backwards using \eqref{eq:eq-XC-a-1} to show $\XC_j$ is a continuous, decreasing functions of $\XC_{\Nt}$ for all $j \in \{\Nt-2, \hdots, 0\}$. %Lastly, we use \eqref{eq:eq-XC-a-2} to conclude the existence and uniqueness using the intermediate value theorem. 
First, by substituting the expression of $\mu$ in \eqref{eq:eq-XC-a-3}, we get
\begin{align*}
    f(\XC_{\Nt}, \XC_{\Nt-1}) = \XC_{\Nt} + \frac{\XC_{\Nt-1}}{rT_B} \left(\tau_1 \sqrt{\frac{d}{X_{\Nt}^C}} + T_R + \frac{\tau_2}{\sqrt{m-X_{\Nt-1}^C}}\right) - n = 0
\end{align*}
In the following claim, we show that $\XC_{\Nt-1}$ is a continuous, decreasing function of $\XC_{\Nt}$.
\begin{claim} \label{claim1}
    There exists $\lambda_c>0$ such that the following holds for $\lambda \geq \lambda_c$: there exists a unique solution of $f(x, x) = 0, x \in [0, A]$; denote it by $L$. Then, $L \in \left[\alpha rT_R\lambda+\frac{r\kappa_1}{2(1+r)}\lambda^{1-\gamma}, \alpha rT_R\lambda+\frac{2.25r\kappa_1}{2+r}\lambda^{1-\gamma}\right]$. Also, there exists a unique, continuous, strictly decreasing, and bijective function $h_{\Nt-1}: [L, n] \rightarrow [0, L]$ such that all $(\XC_{\Nt}, \XC_{\Nt-1})$ satisfying
    \begin{align*}
        f(\XC_{\Nt}, \XC_{\Nt-1}) = 0, \quad  \XC_{\Nt} \geq \XC_{\Nt-1}, \quad \XC_{\Nt} \in (0, n], \quad \XC_{\Nt-1} \in [0, A], \numberthis \label{eq: claim_hypothesis}
    \end{align*}
are given by $\{(x, h_{\Nt-1}(x)) : x \in [L, n]\}$.
\end{claim}
We provide a proof of Claim~\ref{claim1} right after the present proof. Now, we recursively define $h_{j}: [L, n] \rightarrow \bbR$ as follows: 
\begin{align*}
    h_{j-1}(x) = h_{\Nt-1}(x) - r T_B \lambda\left(1 - \frac{\Pi_{i=0}^{d-1} \left[h_{j}(x)-i\right]^+}{\Pi_{i=0}^{d-1} (x-i)}\right) \quad \forall j \in \{0, 1, \hdots, \Nt-1\} \numberthis \label{eq: h_functions}
\end{align*}
As $h_{\Nt-1}(\cdot)$ is a continuous and strictly decreasing function (by Claim~\ref{claim1}), it is easy to observe that $h_j(\cdot)$ is a continuous and strictly decreasing function for all $j \in \{-1, 0, \hdots, \Nt-1\}$. Also, by \eqref{eq:eq-XC-a-1}, we conclude that any $\XC$ satisfying \eqref{eq:eq-XC-a} must be such that $\XC_j = h_j(\XC_{\Nt})$ for all $j \in \{0, \hdots, \Nt-1\}$ and $\XC_{\Nt} \in [L, n]$. Thus, we conclude part $(b)$ of the proof. Note that we have now reduced \eqref{eq:eq-XC-a} (except \eqref{eq:eq-XC-a-4}) to solving a single variable equation: \eqref{eq:eq-XC-a-2} with $\XC_j = h_j(\XC_{\Nt})$ for $j \in \{0, \hdots, \Nt-1\}$, which is identical to $h_{-1}(\XC_{\Nt})=0$.

\underline{\textbf{Proof of (c)}:} \textbf{Existence and Uniqueness}

To establish part $(c)$, we consider \eqref{eq:eq-XC-a-2} as a function of $\XC_{\Nt}$ which is the same as $h_{-1}(\XC_{\Nt})=0$. Note that $h_{-1}(\cdot)$ is a strictly decreasing function, which establishes uniqueness. Now, we invoke the intermediate value theorem for existence. First, note that $h_{\Nt-1}(L) = L$ by Claim~\ref{claim1} which then inductively implies $h_{-1}(L) = L > 0$ by \eqref{eq: h_functions}. Next, we show $h_{-1}\left(U\right) < 0$ where $U = \alpha r T_R \lambda + \frac{(1+2r)\kappa_1}{2(1+r)}\lambda^{1-\gamma}$. First note that $n> U > L$ by \eqref{eq: in_terms_of_n} and noting that $\frac{1+2r}{2(1+r)} \geq \frac{2.25r}{2+r}$ for $r \in (0, 1)$. Now, by Claim~\ref{claim1}, we have $f(U, h_{\Nt-1}(U)) = 0$, which implies
\begin{align*}
    h_{\Nt-1}(U) &= r T_B (n-U) \frac{1}{\tau_1 \sqrt{\frac{d}{U}} + T_R + \frac{\tau_2}{\sqrt{m-h_{\Nt-1}(U)}}} \leq \frac{r T_B}{T_R} (n-U) \\
     &\overset{\eqref{eq: upper_bound_t_b}}{\leq}  r\left(1 + \frac{2\tau_2}{T_R}\lambda^{-\gamma}\right)\left(\alpha T_R\lambda + \frac{\kappa_1}{2(1+r)}\lambda^{1-\gamma}\right)\overset{*}{\leq} \alpha rT_R\lambda  + \frac{r\kappa_1}{2(1+r)}\lambda^{1-\gamma} + 4\alpha r \tau_2 \lambda^{1-\gamma} \\
    &\overset{**}{\leq} \alpha rT_R\lambda  + \frac{3r\kappa_1}{4(1+r)}\lambda^{1-\gamma}. \numberthis \label{eq: upper_bound_nt_minus_one}
\end{align*}
where $(*)$ follows for all $\lambda \geq \lambda_{02}$ for some $\lambda_{02} > 0$ such that $\frac{\kappa_1}{2(1+r)}\lambda^{1-\gamma} \leq \alpha T_R \lambda$. Further, $(**)$ holds because $\kappa_1 = 32\tau_2 \geq 16(1+r)\tau_2$ as $r<1$. Now, we inductively show that $h_j(U) \leq \alpha r T_R \lambda - \lambda^{1-(j+1)\gamma/\Nt}$ for $j \in \{\Nt-2, \hdots, 0\}$ and $h_{-1}(U) < 0$. The base case for $j = \Nt-2$ is verified as follows:
\begin{align*}
    h_{\Nt-2}(U) &\overset{\eqref{eq: bounds_on_p_j}}{\leq} h_{\Nt-1}(U) - r T_R \lambda\left(1 - \left(\frac{h_{\Nt-1}(U)}{U}\right)^d\right) \quad \text{(as $T_B \geq T_R$)} \\
    &\overset{\eqref{eq: upper_bound_nt_minus_one}}{\leq} \alpha rT_R\lambda  + \frac{3r\kappa_1}{4(1+r)}\lambda^{1-\gamma} - rT_R\lambda \left(1 - \left(\frac{1+\frac{3\kappa_1}{4\alpha T_R (1+r)}\lambda^{-\gamma}}{1+\frac{(1+2r)\kappa_1}{2\alpha r T_R(1+r)}\lambda^{-\gamma}}\right)^d\right) \\
    &\overset{*}{\leq} \alpha rT_R\lambda  + \frac{3r\kappa_1}{4(1+r)}\lambda^{1-\gamma} - rT_R\lambda \left(1 - \left(1 - \frac{(2+r)\kappa_1}{8\alpha r T_R(1+r)}\lambda^{-\gamma}\right)^d\right) \\
    &\overset{**}{\leq} \alpha rT_R\lambda  + \frac{3r\kappa_1}{4(1+r)}\lambda^{1-\gamma} - \frac{(2+r)\kappa_1}{16\alpha  (1+r)}\lambda^{1-\gamma}d \\
    &\overset{***}{\leq} \alpha rT_R\lambda - \lambda^{1-(\Nt-1)\gamma/\Nt}.
\end{align*}
The last three inequalities  holds for all $\lambda \geq \lambda_{03}$ for some $\lambda_{03} > 0$. In particular, $(*)$ holds as $\frac{1+\frac{3\kappa_1}{4\alpha T_R(1+r)}\lambda^{-\gamma}}{1+\frac{(1+2r)\kappa_1}{2\alpha r T_R(1+r)}\lambda^{-\gamma}} = 1-\frac{\frac{(2+r)\kappa_1}{4\alpha rT_R (1+r)}\lambda^{-\gamma}}{1+\frac{(1+2r)\kappa_1}{2\alpha r T_R(1+r)}\lambda^{-\gamma}} \leq 1-\frac{(2+r)\kappa_1}{8\alpha rT_R (1+r)}\lambda^{-\gamma}$ for $\lambda $ large enough. Next, $(**)$ holds as $\lim_{\lambda \rightarrow \infty} \frac{1-(1-c\lambda^{-\gamma})^d}{c\lambda^{-\gamma}d} = 1$ for any constant $c>0$ as $\lambda^{-\gamma}d = \lambda^{-\gamma+\gamma/\Nt} (\log \lambda)^{1/\Nt} \rightarrow 0$ since $\Nt \geq 2$. Thus, we have 
$(1-c\lambda^{-\gamma})^d \leq 1-c\lambda^{-\gamma}d/2$ for $\lambda$ large enough and any constant $c>0$. %as $\lambda^{\gamma/\Nt} \leq d \leq \lambda^{\gamma/\Nt} \log \lambda$ and $\Nt \geq 2$. 
Lastly, the inequality $(***)$ holds as $\lambda^{1-\gamma} d = \lambda^{1-(\Nt-1)\gamma/\Nt} (\log \lambda)^{1/\Nt}$ which is bigger order than $\lambda^{1-(\Nt-1)\gamma/\Nt}$. Now, we proceed to show the induction step: we assume the bound $h_{j+1}(U) \leq \alpha r T_R \lambda - \lambda^{1-(j+2)\gamma/\Nt}$ and show the required bound on $h_j(U)$ below.
% \textbf{Induction Step:}
\begin{align*}
    h_{j}(U) &\overset{\eqref{eq: bounds_on_p_j}}{\leq} h_{\Nt-1}(U) - r T_R \lambda \left(1-\left(\frac{h_{j+1}(U)}{U}\right)^d\right) \quad \text{(as $T_B \geq T_R$)} \\
    &\overset{\eqref{eq: upper_bound_nt_minus_one}}{\leq} \alpha rT_R\lambda  + \frac{3r\kappa_1}{4(1+r)}\lambda^{1-\gamma} - rT_R\lambda \left(1 - \left(\frac{1-\frac{1}{\alpha r T_R}\lambda^{-(j+2)\gamma/\Nt}}{1+\frac{(1+2r)\kappa_1}{2\alpha r T_R(1+r)}\lambda^{-\gamma}}\right)^d\right) \quad \text{(Using IH)} \\
    &\overset{*}{\leq} \alpha rT_R\lambda  + \frac{3r\kappa_1}{4(1+r)}\lambda^{1-\gamma} - rT_R\lambda \left(1 - \left(1 - \frac{1}{2 \alpha r T_R}\lambda^{-(j+2 )\gamma/\Nt}\right)^d\right) \\
    &\overset{**}{\leq} \alpha rT_R\lambda  + \frac{3r\kappa_1}{4(1+r)}\lambda^{1-\gamma} - \min\left\{\frac{1}{4\alpha}\lambda^{1-(j+2)\gamma/\Nt}d, \frac{(1+\alpha)r T_R \lambda}{2} \right\} \\
    &\overset{***}{\leq}
    \begin{cases}
        \alpha rT_R\lambda - \lambda^{1-(j+1)\gamma/\Nt} &\textit{if } j \in \{0, 1, 2, \hdots, \Nt-3\} \\
        \frac{(\alpha-1)r T_R \lambda}{4} & \textit{if } j = -1
    \end{cases}
\end{align*}
Similar to the base case, the last three inequalities above hold for $\lambda \geq \lambda_{04}$ for some $\lambda_{04} > 0$. In particular, $(*)$ holds as $\frac{1 - \frac{1}{\alpha r T_R}\lambda^{-(j+2)\gamma/\Nt}}{1+\frac{(1+2r)\kappa_1}{2\alpha r T_R(1+r)}\lambda^{-\gamma}} \leq 1 - \frac{ \frac{1}{\alpha r T_R}\lambda^{-(j+2)\gamma/\Nt}}{1+\frac{(1+2r)\kappa_1}{2\alpha r T_R(1+r)}\lambda^{-\gamma}}\leq 1 - \frac{1}{2\alpha r T_R}\lambda^{-(j+2)\gamma/\Nt}$ for $\lambda$ large enough. Next, for $j \in \{0, 1, \hdots, \Nt-3\}$, %$(**)$ holds as $(1-c\lambda^{-(j+2)\gamma/\Nt})^d \leq 1-c\lambda^{-(j+2)\gamma/\Nt} d + c^2 \lambda^{-2(j+2)\gamma/\Nt} d^2/2 \leq 1-c\lambda^{-(j+2)\gamma}d/2$ for $\lambda$ large enough and any constant $c>0$ as $\lambda^{\gamma/\Nt} \leq d \leq \lambda^{\gamma/\Nt} \log \lambda$. 
$\lim_{\lambda \rightarrow \infty} \frac{1-(1-c\lambda^{-(j+2)\gamma/\Nt})^d}{c\lambda^{-(j+2)\gamma/\Nt}d} = 1$ for any constant $c>0$ as $\lambda^{-(j+2)\gamma}d = \lambda^{-(j+1)\gamma/\Nt} (\log \lambda)^{1/\Nt} \rightarrow 0$ since $\Nt \geq 2$. Thus, we have 
$(1-c\lambda^{-(j+2)\gamma})^d \leq 1-c\lambda^{-(j+2)\gamma}d/2$ for $\lambda$ large enough and any constant $c>0$. In addition, for $j=-1$, we have $(1-c\lambda^{-\gamma/\Nt})^d \rightarrow 0$ as $\lambda \rightarrow \infty$ and so $(1-c\lambda^{-\gamma/\Nt})^d \leq (1-\alpha)/2$ for $\lambda$ large enough as $\alpha<1$. Note that we also use that $\min\left\{\frac{1}{4\alpha}\lambda^{1-\gamma/\Nt}d, \frac{(1+\alpha)r T_R \lambda}{2} \right\} = \frac{(1+\alpha)r T_R \lambda}{2}$ for $\lambda$ large enough as $\lambda^{1-\gamma/\Nt}d = \lambda (\log \lambda)^{1 / \Nt}$. Lastly, for $j \in \{0, 1, \hdots, \Nt-3\}$, the inequality $(***)$ holds as $\lambda^{1-(j+2)\gamma/\Nt} d = \lambda^{1-(j+1)\gamma/\Nt} (\log \lambda)^{1/\Nt}$ which is bigger order than $\lambda^{1-\gamma}$ and $\lambda^{1-(j+1)\gamma/\Nt}$. Moreover, for $j=-1$, we simply use that $\frac{(1-\alpha)rT_R \lambda}{4} \geq \frac{3r\kappa_1}{4(1+r)} \lambda^{1-\gamma}$ for $\lambda$ large enough as $\gamma>0$.

The above completes the induction step and shows that $h_{-1}(U) < 0$ as $\alpha<1$. Thus, by the intermediate value theorem, there exists $X^\star \in [L, U]$ that satisfies $h_{-1}(X^\star) = 0$. Thus, $\XC_{j} = h_{j}(X^\star)$ for $j \in \{0, \hdots, \Nt-1\}$ along with $\XC_{\Nt}=X^\star$ is a unique solution of \eqref{eq:eq-XC-a-1}, \eqref{eq:eq-XC-a-2}, \eqref{eq:eq-XC-a-3}, \eqref{eq:eq-XC-a-5}.

\underline{\textbf{Proof of (d):}} \textbf{Final verifications}

Now, it remains to verify \eqref{eq:eq-XC-a-4} which we prove by induction with hypothesis $h_k(X^\star) \leq h_{k+1}(X^\star)$ for all $k \in \{j, \hdots, \Nt-1\}$. The base case of $h_{\Nt-1}(X^\star) \leq X^\star$ is satisfied by Claim~\ref{claim1} as $h_{\Nt-1}(x) \leq x$ for any $x \in [L, U]$. Now, for the induction step, we have
\begin{align*}
    h_{j-1}(X^\star) &= h_{\Nt-1}(X^\star) - r T_B \lambda\left(1 - \frac{\Pi_{i=0}^{d-1} \left[h_{j}(X^\star)-i\right]^+}{\Pi_{i=0}^{d-1} (X^\star-i)}\right) \\
    &\leq h_{\Nt-1}(X^\star) - r T_B \lambda\left(1 - \frac{\Pi_{i=0}^{d-1} \left[h_{j+1}(X^\star)-i\right]^+}{\Pi_{i=0}^{d-1} (X^\star-i)}\right) \\
    &= h_{j}(X^\star).
\end{align*}
Lastly, we conclude the proof by showing $h_0(X^\star) \geq 0$. If not, then, $h_{-1}(X^\star) \leq U - rT_B \lambda < 0$ for $\lambda \geq \lambda_{05}$ for some $\lambda_{05}>0$ as $\alpha<1$,
which is a contradiction. This completes the proof by setting $\lambda_{p1} = \max_{k \in \{0, \hdots, 5\}}\{\lambda_{0k}, \lambda_{c}\}$.
\hfill $\square$

\endproof

\proof{Proof of Claim \ref{claim1}}
First we show that $f(x, x) = 0, x \in [0, A]$ has a unique solution $L \in \left[\alpha rT_R\lambda+\frac{r\kappa_1}{2(1+r)}\lambda^{1-\gamma}, \alpha rT_R\lambda+\frac{2.25r\kappa_1}{2+r}\lambda^{1-\gamma}\right]$. Note that $g : [0, A] \rightarrow \bbR$ defined as
\begin{align*}
    g(x) = f(x, x) = x\left(1+\frac{T_R}{rT_B}+\frac{\tau_2}{rT_B\sqrt{m-x}}\right)+\frac{\tau_1\sqrt{dx}}{rT_B}-n
\end{align*}
is a strictly increasing and continuous function, so there can exist at most one root of the function. Now we show that $g\left(\alpha rT_R\lambda+\frac{r\kappa_1}{2(1+r)}\lambda^{1-\gamma}\right) \leq 0$ and $g\left(\alpha rT_R\lambda+\frac{2.25r\kappa_1}{2+r}\lambda^{1-\gamma}\right) \geq 0$. First, we have
\begin{align*} &g\left(\alpha rT_R\lambda+\frac{r\kappa_1}{2(1+r)}\lambda^{1-\gamma}\right) \\
={}& \left(\alpha rT_R\lambda+\frac{r\kappa_1}{2(1+r)}\lambda^{1-\gamma}\right)\left(1+\frac{T_R}{rT_B}+\frac{\tau_2}{rT_B\sqrt{m-\alpha rT_R\lambda-\frac{r\kappa_1}{2(1+r)}\lambda^{1-\gamma}}}\right)\\
&+\frac{\tau_1}{rT_B}\sqrt{d\left(\alpha rT_R\lambda+\frac{r\kappa_1}{2(1+r)}\lambda^{1-\gamma}\right)}-n \\
\overset{(a)}{\leq}{}& \left(\alpha rT_R\lambda+\frac{r\kappa_1}{2(1+r)}\lambda^{1-\gamma}\right)\left(1+\frac{1}{r}+\frac{\tau_2\lambda^{-\gamma}}{rT_R}\right)+\frac{\tau_1}{rT_R}\sqrt{d\left(\alpha rT_R\lambda+\frac{r\kappa_1}{2(1+r)}\lambda^{1-\gamma}\right)}-n  \\
\overset{(b)}{\leq}{}& \left(\alpha rT_R\lambda+\frac{r\kappa_1}{2(1+r)}\lambda^{1-\gamma}\right)\left(1+\frac{1}{r}+\frac{\tau_2\lambda^{-\gamma}}{rT_R}\right)+\frac{2\tau_1}{\sqrt{rT_R}}\lambda^{\frac{1}{2}+\frac{\gamma}{2\Nt}}\log \lambda-\alpha(1+r)T_R\lambda-\kappa_1\lambda^{1-\gamma} \\
\overset{(c)}{\leq}{}& -\left(\frac{\kappa_1}{2}-2\tau_2\right)\lambda^{1-\gamma}+\frac{2\tau_1}{\sqrt{rT_R}}\lambda^{\frac{1}{2}+\frac{\gamma}{2\Nt}}\log \lambda \overset{(d)}{\leq} 0, \numberthis \label{eq: claim1_left_point}
\end{align*}
where $(a)$ follows as $T_B \geq T_R$ and $\frac{r\kappa_1}{2(1+r)}\lambda^{1-\gamma} \leq (\kappa_2-1)\lambda^{1-\gamma} \leq (\kappa_2-1)\lambda^{2\gamma}$ as $\gamma \geq 1/3$, $r>0$, and $\kappa_2 = r\kappa_1+4\tau_2+1$. Next, $(b)$ follows as $\alpha<1$, $d \leq \lambda^{\gamma/\Nt}(\log \lambda)^2$, and $\alpha rT_R\lambda + \frac{r\kappa_1}{2(1+r)}\lambda^{1-\gamma} \leq 4rT_R\lambda$ for $\lambda \geq \lambda_{c0}$ for some $\lambda_{c0} > 0$. Further, $(c)$ follows as $\left(\alpha rT_R\lambda + \frac{r\kappa_1}{2(1+r)}\lambda^{1-\gamma}\right) \tau_2\lambda^{-\gamma}/rT_R \leq 2\tau_2 \lambda^{1-\gamma}$ for $\lambda \geq \lambda_{c1}$ for some $\lambda_{c1} > 0$ since $\alpha<1$. Lastly, $(d)$ follows as $\kappa_1 = 32\tau_2$ and $14\tau_2\lambda^{1-\gamma} \geq \frac{2\tau_1}{\sqrt{rT_R}} \lambda^{\frac{1}{2}+\frac{\gamma}{2\Nt}}\log\lambda$ for $\lambda \geq \lambda_{c2}$ for some $\lambda_{c2} > 0$ as $1-\gamma > \frac{1}{2}+\frac{\gamma}{2\Nt}$. This shows that $g\left(\alpha rT_R\lambda+\frac{r\kappa_1}{2(1+r)}\lambda^{1-\gamma}\right) \leq 0$. Now, we analyze $g\left(\alpha rT_R\lambda+\frac{2.25r\kappa_1}{2+r}\lambda^{1-\gamma}\right)$ below. We have
\begin{align*}
    g\left(\alpha rT_R\lambda+\frac{2.25r\kappa_1}{2+r}\lambda^{1-\gamma}\right) &\geq \left(\alpha rT_R\lambda+\frac{2.25r\kappa_1}{2+r}\lambda^{1-\gamma}\right)\left(1+\frac{T_R}{rT_B}\right)-n \\
    &\overset{\eqref{eq: upper_bound_t_b}}{\geq} \left(\alpha rT_R\lambda+\frac{2.25r\kappa_1}{2+r}\lambda^{1-\gamma}\right)\left(1+\frac{1}{r\left(1 + \frac{2\tau_2}{T_R}\lambda^{-\gamma}\right)}\right)-n \\
     &\overset{(a)}{\geq} \left(\alpha rT_R\lambda+\frac{2.25r\kappa_1}{2+r}\lambda^{1-\gamma}\right)\left(1+\frac{1}{r} - \frac{2\tau_2}{rT_R}\lambda^{-\gamma}\right)-n \\
     &\overset{(b)}{\geq} \frac{2.25(1+r)\kappa_1}{2+r}\lambda^{1-\gamma} - 4\tau_2\lambda^{1-\gamma}-\kappa_1\lambda^{1-\gamma} \\
     &\geq \frac{0.25(1+5r)\kappa_1}{2+r}\lambda^{1-\gamma} - 4\tau_2\lambda^{1-\gamma}  \overset{(c)}{\geq} 0, \numberthis \label{eq: claim1_right_point}
\end{align*}
where $(a)$ follows as $1/(1+x) \geq 1-x$ for $x \geq 0$. Next, $(b)$ follows as $\alpha<1$ as well as $\frac{2.25r\kappa_1}{2+r}\lambda^{1-\gamma} \leq rT_R\lambda$ for $\lambda \geq \lambda_{c3}$ for some $\lambda_{c3} > 0$. Lastly, $(c)$ follows as $\kappa_1 = 32\tau_2 \geq 16(2+r)\tau_2/(1+5r)$ as $r>0$. Thus, by using the intermediate value theorem along with \eqref{eq: claim1_left_point} and \eqref{eq: claim1_right_point}, and noting that $g$ is strictly monotonic, we conclude that there exists a unique $L \in \left[\alpha rT_R\lambda+\frac{r\kappa_1}{2(1+r)}\lambda^{1-\gamma}, \alpha rT_R\lambda+\frac{2.25r\kappa_1}{2+r}\lambda^{1-\gamma}\right]$ that solves $g(x)=0$, completing the first part of the claim.

\revcolor{
Next, we show that there exists a bijective mapping $h_{\Nt-1} : [L, n] \rightarrow [0, L]$ such that all solutions satisfying \eqref{eq: claim_hypothesis} are given by $\{(x, h_{\Nt-1}(x)): x \in [L, n]\}$. The idea of the proof is illustrated in Fig.~\ref{fig:illustration_f}. Multiplying $f(x, y)$ by $\sqrt{x}$ and replacing $\sqrt{x}$ with $z$, we get the following cubic polynomial of $z$:
\begin{align*}
    \tilde{f}_y(z) = z^3 + z\left(\frac{T_Ry}{rT_B} + \frac{\tau_2 y}{rT_B\sqrt{m-y}}-n\right) + \frac{\tau_1 y\sqrt{d}}{rT_B}.
\end{align*}
Thus, for any $y \in [0, A]$, there exist at most three values of $x$ such that $f(x, y) = 0$. It is easy to verify that one of the solutions of the above cubic polynomial is negative for any $y \in [0, A]$. In particular, if $y \in (0, A]$, then $\tilde{f}_y(0) > 0$ and $\lim_{z \rightarrow -\infty}\tilde{f}_y(z) = -\infty$. Otherwise, if $y=0$, then $z = -\sqrt{n}$ is one of the roots. As one of the roots of $\tilde{f}_y(z)$ is negative, there exist at most two values of $x$ for which $f(x, y) = 0$ for any $y \in [0, A]$. Now we show that only one of these solutions satisfies $x \geq y$, showing the existence of a unique $x$ such that $f(x, y) = 0$ and $x \geq y$ for any $y \in [0, A]$. We consider the following four cases depending on the value of $y \in [0, A]$.
}

\begin{figure}
    \centering
\begin{tikzpicture}[scale=1.0]

  %%% Define some example coordinates (tweak as needed):
  \def\A{3}  % A on the y-axis
  \def\L{2}  % L on both axes (smaller than A)
  \def\n{5}  % n on the x-axis (bigger than L)

  % Compute coefficients for the cubic: f(x)= a*x^3 + b*x^2 + c*x.
  % These satisfy f(0)=0, f(2)=2, f(5)=0 and f'(2)=0.
  \pgfmathsetmacro{\a}{1/18}    % a = 1/18
  \pgfmathsetmacro{\b}{-13/18}  % b = -13/18
  \pgfmathsetmacro{\c}{20/9}    % c = 20/9

  % 1) Axes
  \draw[->, thick] (0,0) -- (6,0) node[right] {$x$};
  \draw[->, thick] (0,0) -- (0,4) node[above] {$y$};

  % 2) Axis labels
  \node[left]  at (0,0) {0};
  \node[left]  at (0,\L) {$L$};
  \node[left]  at (0,\A) {$A$};
  \node[below] at (\L,0) {$L$};
  \node[below] at (\n,0) {$n$};

  % 3) Diagonal line x = y (draw it up to y=A)
  \draw[black, thick] (0,0) -- (\A,\A) node[near end, above right, xshift=3mm] {$x=y$};

  % 4) Shaded region: 0 <= x <= L, from y=x up to y=A
  %    This forms a trapezoid: (0,0)->(0,A)->(L,A)->(L,L)
  \fill[red!30, opacity=0.6]
    (0,0) -- (0,\A) -- (\A,\A) -- cycle;

  % 5) Horizontal lines at y=L and y=A
  \draw[blue] (0,\L) -- (\n,\L);
  \draw[blue] (0,\A) -- (\n,\A);

  % 6) Vertical lines at x=L and x=n
  \draw[blue] (\L,0) -- (\L,\A);
  \draw[blue] (\n,0) -- (\n,\A);

  % 7) The green curve f(x,y) = 0, from (0,L) to (n,0)
  \draw[green!70!black, thick, domain=0:\n, samples=200] 
    plot (\x, {\a*(\x)^3 + \b*(\x)^2 + \c*(\x)})
    node[right, xshift=-18mm, yshift=16mm] at (\n,0) {$f(x,y)=0$};

  % 8) Mark key points
  \fill (0,\L) circle (1.5pt);
  \fill (\n,0) circle (1.5pt);

\end{tikzpicture}
\caption{Illustration of all $x \in (0, n]$ and $y \in [0, A]$ satisfying $f(x, y) = 0$. The red shaded area depicts the region that does not satisfy $x \geq y$.}
    \label{fig:illustration_f}
\end{figure}

\textbf{Case I:} If $y=0$, then $x=n$ is the only root of $f(x, y) = 0$.

\textbf{Case II:} We show that for any $y \in (0, L)$, there exists a unique $x \in (L, n)$ satisfying $f(x, y)= 0, x \geq y$ and vice versa. Fix a $y \in (0, L)$. Observe the following:
\begin{itemize}
    \item $\lim_{x \rightarrow 0}f(x, y) = \infty$ and $f(n, y) > 0$.
    \item $f(y, y) = g(y)<0$: The first equality holds by the definition of $g$, and the last one holds as $L$ is the unique root of $g$, $g$ is a strictly increasing function, and $y < L$.
    \item $f(L, y) < f(L, L) = g(L) = 0$: the first equality holds as $f(L, \cdot)$ is a strictly increasing function.
\end{itemize}
Noting that $f(\cdot, y)$ is a continuous function and invoking the intermediate value theorem, we conclude that the two roots of $f(x, y)$ lie in the intervals $(0, y)$ and $(L, n)$ respectively. Thus, for any $y \in (0, L)$ there exists a unique $x \in (L, n)$ such that $f(x,y)=0, x \geq y$ and no $x \in [y, L]$ satisfies $f(x,y)=0, x \geq y$. Now, we show that for any $x \in (L, n)$, there exists a unique $y \in (0, L)$ satisfying $f(x,y)=0, x \geq y$ and no $y \in [0, A] \backslash (0, L)$ satisfies $f(x,y)=0, x \geq y$. Uniqueness is readily established by noting that $f(x, \cdot)$ is a strictly increasing function for any $x \in (L, n)$. To show existence, fix an $x \in (L, n)$ and note that $f(x, 0) = x - n < 0$. Next, we show that $f(x, L) > 0$. We analyze the derivative of $f$ as follows:
\begin{align*}
    \frac{df(x, L)}{dx} &= 1 - \frac{\tau_1 \sqrt{d} L}{2r T_B x^{3/2}} \geq 1 - \frac{\tau_1 \sqrt{d}}{2r T_B \sqrt{L}} \geq 1 - \frac{\tau_1}{2r T_B \sqrt{\alpha r T_R}} \lambda^{\frac{\gamma}{2\Nt}-\frac{1}{2}} \log \lambda > 0 \quad \forall x \in (L, n), \numberthis \label{eq: derivative_claim}
\end{align*}
where the last inequality follows for $\lambda \geq \lambda_{c4}$ for some $\lambda_{c4} > 0$ as $\gamma/\Nt < 1$. In addition, we have $f(L, L) = g(L) = 0$ which combined with \eqref{eq: derivative_claim} shows that $f(x, L) > 0$ for all $x \in (L, n)$.  Thus, by noting that $f(x, \cdot)$ is a continuous function, by the intermediate value theorem, for any $x \in (L, n)$, there exists a unique $y \in (0, L)$ satisfying $f(x, y) = 0$.

\textbf{Case III:} For $y = L$, we have $f(L, L) =g(L)= 0$. In addition by \eqref{eq: derivative_claim}, we have $f(x, L) > 0$ for all $x \in (L, n]$. Thus,  $L$ is the unique solution of $f(x, L) = 0, x \in [L, n]$.

\textbf{Case IV:} For $y \in (L, A]$, we have $f(x, y) > f(x, L)$ for all $x \in [L, n]$, as $f(x, \cdot)$ is a strictly increasing function for any fixed $x \in [L, n]$. In addition, $f(x, L) \geq 0$ as previously argued. Thus, there exists no $x \in [L, n]$ satisfying $f(x, y) = 0, x \geq y, y \in (L, A]$.

By Case I-III, we conclude that there exists a unique bijective mapping $h_{\Nt-1}:[L,n] \rightarrow [0, L]$ such that $\{(x, h_{\Nt-1}(x)): x \in [L, n]\}$ solves \eqref{eq: claim_hypothesis}. In addition, by Case IV, we conclude that these are the only possible solutions.

Now, we show that $h$ is a strictly decreasing function. For any $y_1 < y_2 \in [0, L]$, let $x_1, x_2 \in [L, n]$ be the corresponding unique $x$ satisfying the hypothesis of the claim. We now show that $x_1 > x_2$. First note that $f(x, \cdot)$ is a strictly increasing function for any $x \in [L, n]$. So, as $f(x_1, y_1) = 0$ and $y_2 > y_1$, we conclude that $f(x_1, y_2) > 0$.
% \begin{align*}
%      f(x_1, y_2) = f(x_1, y_2) - f(x_1, y_1) = \frac{1}{rT_B} \left(\tau_1\sqrt{\frac{d}{x_1}} + T_R\right)(y_2-y_1) + \frac{\tau_2}{rT_B}\left(\frac{y_2}{\sqrt{m-y_2}} - \frac{y_1}{\sqrt{m-y_1}}\right) > 0,
% \end{align*}
% where the last inequality follows as $ y/\sqrt{m-y}$ is a strictly increasing function of $y \in \left[0, L\right]$. 
Similarly, as $f(L, L) = 0$ and $y_2 \leq L$, we have $f(L, y_2) \leq 0$.
% get $f(L, y_2) = f(L, y_2) - f(L, L) \leq 0$ as $y_2 \leq L$. 
As $f(L, y_2) \leq 0$ and $f(x_1, y_2) > 0$, we must have $x_2 \in [L, x_1)$ which shows that $x_2 < x_1$. Thus, $h$ is a strictly decreasing function. Now, we show that $h$ is a continuous function. As $h$ is strictly decreasing, the only possible discontinuity is a jump discontinuity. Assume there exists an $x_0 \in (L, n)$ and $c > 0$ such that $h_{\Nt-1}(x_0^+) < h_{\Nt-1}(x_0^-)$ where $h_{\Nt-1}(x_0^s) = \lim_{x \rightarrow x_0^{s}} h_{\Nt-1}(x)$ for $s \in \{+, -\}$. Then, for any $y \in (h_{\Nt-1}(x_0^+), h_{\Nt-1}(x_0^-))$, there exists no $x \in (L, n)$ such that $y = h_{\Nt-1}(x)$ which is a contradiction as $h_{\Nt-1}$ is a bijective function. One can similarly handle the boundary $x \in \{L, n\}$. This shows that $h_{\Nt-1}$ is a continuous function, which completes the proof by setting $\lambda_{c} = \max_{k \in [4]} \lambda_{ck}$. \hfill $\square$
\endproof

\subsection{Service Level for Power-of-d Vehicles Dispatch}
\proof{Proof of \Cref{prop: lower_bound_lambda}}
First, consider $\lambda \geq \lambda_{p_1}$. Then, by global stability and Proposition \eqref{prop:scalings-infty}, there exists $t_0(\lambda) > 0$ such that for all $t \geq t_0(\lambda)$, we have
\begin{align*}
\XB_{\Nt}(t) \geq \alpha T_R\lambda +\frac{\kappa_1}{4(1+r)}\lambda^{1-\gamma} - 1,
\end{align*}
Using this and \eqref{eq:powerd-3}, we get
\begin{align*}
\lambda p_0(t) I(t) &= \XB_{\Nt}(t)\mu(t) - \dXC_{\Nt}(t) \\
&\geq \left(\alpha T_R\lambda +\frac{\kappa_1}{4(1+r)}\lambda^{1-\gamma}-1\right)\mu(t) - \dXC_{\Nt}(t) \\
&\overset{(a)}{\geq} \frac{1}{T_B}\left(\alpha T_R\lambda +\frac{\kappa_1}{4(1+r)}\lambda^{1-\gamma}\right) - \frac{1}{T_R} - \dXC_{\Nt}(t) \\
&\overset{\eqref{eq: upper_bound_t_b}}{\geq} \frac{1}{ 1 + \frac{2\tau_2}{T_R}\lambda^{-\gamma}}\left(\alpha \lambda +\frac{\kappa_1}{4(1+r)}\lambda^{1-\gamma}\right) - \frac{1}{T_R} - \dXC_{\Nt}(t) \\
&\overset{(b)}{\geq}  \left(1- \frac{2\tau_2}{T_R}\lambda^{-\gamma}\right)\left(\alpha \lambda +\frac{\kappa_1}{4(1+r)T_R}\lambda^{1-\gamma}\right)-\frac{1}{T_R} - \dXC_{\Nt}(t) \\
&\overset{(c)}{\geq} \alpha \lambda - \frac{2\tau_2}{T_R}\lambda^{1-\gamma} + \frac{\kappa_1}{4(1+r)T_R}\lambda^{1-\gamma} - \dXC_{\Nt}(t) \\
&\overset{(d)}{\geq} \alpha \lambda+ \frac{\kappa_1}{8(1+r)T_R}\lambda^{1-\gamma} - \dXC_{\Nt}(t) \numberthis \label{eq: admitted_arrival_rate}
\end{align*}
where $(a)$ follows as $1/T_B \leq \mu(t) \leq 1/T_R$. Next, $(b)$ follows as $1/(1+x) \geq 1-x$ for $x \geq 0$. \revcolor{Next, $(c)$ follows for all $\lambda \geq \lambda_{11}$ for some $\lambda_{11} > 0$ as $\frac{2\tau_2 \lambda^{-\gamma}}{T_R}\left(\alpha \lambda +\frac{\kappa_1}{4(1+r)T_R}\lambda^{1-\gamma}\right)+\frac{1}{T_R} \leq \frac{2\tau_2 \lambda^{1-\gamma}}{T_R}$ for $\lambda$ large enough as $\alpha<1$.} Lastly, $(d)$ follows by noting that $\kappa_1 = 32\tau_2 \geq 16(1+r)\tau_2$ as $r<1$. Now, by integrating both sides of \eqref{eq: admitted_arrival_rate}, we get
\begin{align*}
\frac{1}{T}\int_0^T \lambdaeff^{\pol(d)}(t)dt &= \frac{1}{T} \int_0^T \lambda p_0(t) I(t)dt \\
&\geq \frac{1}{T} \int_{t_0}^T \lambda p_0(t) I(t)dt\\
&\geq \frac{T-t_0(\lambda)}{T} \left(\alpha \lambda+ \frac{\kappa_1}{8(1+r)T_R}\lambda^{1-\gamma}\right) - \frac{\XC_{\Nt}(T) - \XC_{\Nt}(t_0)}{T} \\
&\geq \frac{T-t_0(\lambda)}{T} \left(\alpha \lambda+ \frac{\kappa_1}{8(1+r)T_R}\lambda^{1-\gamma}\right) - \frac{n}{T} \\
&= \alpha \lambda+ \frac{\kappa_1}{8(1+r)T_R}\lambda^{1-\gamma} - \frac{1}{T}\left(t_0(\lambda)\left(\alpha \lambda+ \frac{\kappa_1}{8(1+r)T_R}\lambda^{1-\gamma}\right)+n\right) \\
&\geq \alpha \lambda \quad \forall T \geq T_0(\lambda) = \frac{t_0(\lambda)\left(\alpha \lambda+ \frac{\kappa_1}{8(1+r)T_R}\lambda^{1-\gamma}\right)+n}{\frac{\kappa_1}{8(1+r)T_R}\lambda^{1-\gamma}}.
\end{align*}
So, by taking the limit infimum as $T \rightarrow \infty$ on both sides, we get
\begin{align*}
 \alphainf \lambda = \liminf_{T\rightarrow \infty}\frac{1}{T}\int_0^T \lambdaeff^{\pol(d)}(t)dt \geq \alpha \lambda \implies \alphainf \geq \alpha.
\end{align*}
This completes the proof with $\lambda_{p_2} = \max\{\lambda_{p_1}, \lambda_{11}\}$.
\hfill $\square$
\endproof

\subsection{Second Order Terms Constants for Power-of-d Vehicles Dispatch }\label{sec:sec-order-const-pod}
In this section, we discuss how to approximate the optimal values for the constants in front of the second-order terms in \Cref{thm:informal_upper-bound}. 
To obtain approximate values for the constants, we can follow the following procedure:
\begin{enumerate}
    \item Fix a desired service level $\alpha$.
    \item Set $m(\lambda)= r T_R \alpha \lambda + c_1\lambda^{2\gamma}$, for some constant $c_1>0$. The constant $\kappa>0$ in the admission control in \Cref{def: admission_control}, can be set as in \Cref{eq: in_terms_of_n}.
    \item For different values of $n$:
    \begin{enumerate}
        \item Simulate the system \eqref{eq:odes-power-d} for $T>0$ large.
        \item Compute $\alpha_{\text{eff}}(n,\lambda) = \frac{1}{T}\int_0^T \lambda p_0(t)I(t) dt$.
        \item If $n$ is such that $\alpha_{\text{eff}}(n,\lambda) =\alpha$ (or close, up to some tolerance), stop. 
    \end{enumerate}
    \item The output of step 3 gives $n(\lambda)$. 
    \item Doing these steps for $\lambda$ large provides an approximated value for the constant in front of $\lambda^{1-\gamma}$ for the number of vehicles, which can be estimated by: $\frac{n(\lambda)-(1+r)T_R \alpha \lambda}{\lambda^{1-\gamma}}$.
\end{enumerate}
We note that the constant $c_1$ for the number of chargers and admission control are choice parameters. As such, they will influence the value of the constant for the fleet size.

\section{Extensions and Additional Insights}\label{sec:extension_add_insights}
\subsection{Failure of Closest Dispatch} \label{sec: closest_dispatch}
We now present the performance of \Cref{alg:powerd} when $d=1$. For the admission control policy as defined in \Cref{def: admission_control}, fix an arbitrary $\kappa > 0$ and $\gamma \in [1/3, 1/2)$ and limit the number of busy EVs by $\min\{\alpha T_R\lambda+\kappa\lambda^{1-\gamma}, n-1\}$. Note that the minimum with $n-1$ is to ensure regularity of the pickup times. In this section, we only focus on the case of $m=\infty$, so we have $A = \infty$, i.e., we don't implement admission control for the chargers. Now, similar to Power-of-$d$, we write down a set of ODEs that will provide an upper bound on the performance of the Closest Dispatch. We make a similar simplification by considering a modified system in which each vehicle spends $\Delta_{\operatorname{CD}}\triangleq T_R r_d$ [kWh] of battery for each trip request (pickup and trip). Note that the drive to the charger time is zero as $m=\infty$. Since the modified system spends less SoC serving a trip request (as $\mu(t) \leq 1/T_R$), we expect that the original system cannot achieve a better service level than the modified system. However, formally establishing a result equivalent to \Cref{lemma: coupling} is challenging. In particular, as the modified system expels smaller kWh serving a trip, an incoming trip that is served in the modified system could be infeasible to serve in the original system. After such an event occurs, the SoC in the modified system could become smaller compared to the original system. Such non-monotonicity is a challenge in implementing a coupling argument similar to \Cref{lemma: coupling}. Nonetheless, as the modified system assumes that less SoC is spent serving trips, we analyze Closest Dispatch under this modified system. We define $\Nt^{\operatorname{CD}} \overset{\Delta}{=} p^\star / \Delta_{\operatorname{CD}}$ and consider the following set of ODEs:
\begin{subequations} \label{eq:odes-cd}
\begin{flalign}\label{eq:cd-1}
\dXB_j(t)&=\lambda I(t) \frac{\XC_j(t) - \XC_0(t)}{\XC_{\Nt^{\operatorname{CD}}}(t)} - \XB_{j}(t) \mu(t), \: 1\le j< \Nt^{\operatorname{CD}},\\\label{eq:cd-2}
\dXC_j(t)&=-\lambda I(t) \frac{\XC_j(t) - \XC_0(t)}{\XC_{\Nt^{\operatorname{CD}}}(t)} - \left(\XC_j(t)-\XC_{j-1}(t)\right) \frac{1}{rT_R} +\XB_{j+1}(t) \mu(t),\: 1\le j<\Nt^{\operatorname{CD}},\\\label{eq:cd-3}
\dXC_{\Nt^{\operatorname{CD}}}(t)&=-\lambda I(t) \left(1-\frac{\XC_0(t)}{\XC_{\Nt^{\operatorname{CD}}}(t)}\right) + \XB_{\Nt^{\operatorname{CD}}}(t) \mu(t), \\\label{eq:cd-4}
\dXC_{0}(t)&=- \XC_0(t)\frac{1}{rT_R} +\XB_{1}(t) \mu(t),
\end{flalign}
\end{subequations} 
Note that the above set of ODEs is the same as \eqref{eq:odes-power-d} with $d=1$, $T_B$ replaced with $T_R$, $\Nt$ replaced with $\Nt^{\operatorname{CD}}$, and $\tXC$ replaced with $\XC$. We now analyze this set of ODEs to obtain a lower bound on the fleet size requirement in the following proposition.
\begin{proposition}[Failure of Closest Dispatch] \label{prop:closest-fails}
Fix an $\alpha \in [\delta, 1]$ for some $\delta>0$. Also, let $m = \infty$ and $1 \geq \alphasup \geq \alpha$, where $\pol$ is the Closest Dispatched policy governed by \eqref{eq:odes-cd}. Then, there exists $\lambda_0, \epsilon >0$, such that 
$$
    n \geq (1+r)T_R \alpha \lambda + \epsilon \lambda \quad \forall  \lambda \geq \lambda_0.
$$
\end{proposition}
Proposition~\ref{prop:closest-fails} establishes that the Closest Dispatch algorithm requires more vehicles to satisfy the demand, even with infinitely many chargers. Under this policy, we greedily minimize the pickup distance. However, in terms of battery level, we are less likely to choose an EV with a high SoC because the closest distance essentially picks a vehicle at random regarding its battery levels. Thus, the variance of SoC across the fleet is higher, resulting in a fraction of EVs being unable to serve customers due to low SoC. To compensate for unavailable EVs, a larger fleet size is warranted, as otherwise, an incoming arrival is more likely to be rejected due to the closest vehicle having a low SoC.

\proof{Proof of \Cref{prop:closest-fails}}
For ease of notation, we denote $\Nt^{\operatorname{CD}}$ simply by $\Nt$ for the rest of the proof. Note that $(\XB_j(t), \XC_j(t))$ are all bounded in $[0,n]$. Hence, if we denote any of these quantities by $y(t)$, we have that $\frac{1}{T}\int_0^T y(t)dt$ is a bounded sequence in a compact set. Therefore, it has a subsequence that converges. We can therefore create a subsequence for which all of the integrals $\frac{1}{T}\int_0^T y(t)dt$  with $y(t)\in \{\XB_j(t), \XC_j(t)\}$ converge. Without loss of generality, let this sequence correspond to the limit supremum of $\frac{1}{T}\int_0^T \XB_{\Nt}(t) dt$. Note that here we are abusing notation and using $T$ instead of explicitly writing the subsequence. We note that this does not affect what we want to prove. Let us use $(\bXB_j(t), \bXC_j(t))$ to denote the limit of each of these integrals (of the subsequences). The proof is divided into the following steps:
\begin{itemize}
    \item[(a)] We first prove the following coarse upper and lower bounds on $\mu(t)$ for all $t \in \bbR_+$:
    \begin{align*}
        \frac{1}{T_R+\tau_1} \leq \mu(t) \leq \frac{1}{T_R}.
    \end{align*}
    \item[(b)] Next, we use these coarse bounds along with \eqref{eq:cd-3} to obtain a coarse lower bound on $n$. Such a lower bound then allows to obtain the following stronger bound on $\mu(t)$:
    \begin{align*}
        \frac{1}{T_R+\frac{CT_R}{\sqrt{\lambda}}} \leq \mu(t) \leq \frac{1}{T_R} \quad \forall t \in \bbR_+,
    \end{align*}
    where $C = \frac{2\tau_1\sqrt{T_R+\tau_1}}{T_R^2\sqrt{\delta r}}$.
    \item[(c)] Analyzing long term behavior of \eqref{eq:cd-1}, \eqref{eq:cd-2}, and \eqref{eq:cd-4} allows us to obtain
    \begin{align*}
        \frac{r}{1+\frac{C}{\sqrt{\lambda}}}\bXB_{j+1} \leq \bXC_j  \leq r\bXB_{j+1} \quad \forall j \in \{0, 1, \hdots, \Nt-1\}.
    \end{align*}
    \item[(d)] We next show that $f(n)^j\left(1+\frac{C}{\sqrt{\lambda}}\right)^{\min\{1, j\}}\left(\bXC_{j} - \bXC_{j-1}\right)$ with $f(n) =\frac{rT_R\lambda}{n - \alpha T_R\lambda-\kappa\lambda^{1-\gamma}}$ is a non-decreasing function of $j \in \{0, 1, \hdots, \Nt\}$ where $\bXC_{-1} = 0$, by analyzing the long term behavior of \eqref{eq:cd-1}, \eqref{eq:cd-2}, and \eqref{eq:cd-4}.
    \item[(e)] We use $n = \bXB_{\Nt} + \bXC_{\Nt}$ to obtain
    \begin{align*}
    \bXC_{\Nt} - \bXC_{\Nt-1} \geq \frac{n r}{r+1+\frac{C}{\sqrt{\lambda}}}\left(1-\frac{C}{\sqrt{\lambda}}\right) \left(\frac{1}{\Nt+1}\mathbbm{1}\left\{f(n) = 1\right\} + \frac{1-f(n)}{1-f(n)^{\Nt+1}}\mathbbm{1}\left\{f(n) \neq 1\right\}\right)
\end{align*}
    \item[(f)] Using $n = \bXB_{\Nt} + \bXC_{\Nt}$ again allows us to translate the above lower bound to an upper bound on $\bXB_{\Nt}$.
    \item[(g)] Lastly, the upper bound on $\bXB_{\Nt}$ combined with \eqref{eq:cd-3} provides a lower bound on $n$ which completes the proof.
\end{itemize}
We start by proving part $(a)$. As $m=\infty$, we have
$$
\mu(t) = \frac{1}{T_R+\tau_1/\sqrt{n-q(t)}} \leq \frac{1}{T_R} \quad \forall t \in \bbR_+.
$$
By the admission control policy, we have $q(t) \leq \min\left\{\alpha \lambda T_R + \kappa \lambda^{1-\gamma}, n-1\right\}$ for all $t \in \bbR_+$. Thus, we have
$$
\mu(t) = \frac{1}{T_R+\tau_1/\sqrt{n-q(t)}} \geq \frac{1}{T_R + \tau_1} \quad \forall t \in \bbR_+.
$$
This completes part $(a)$ of the proof. We now show part $(b)$ by bootstrapping this coarse lower bound on $\mu(t)$ to obtain a tighter lower bound. By adding \eqref{eq:cd-1} and \eqref{eq:cd-2}, we get
\begin{align*}
    \dXB_j(t) + \dXC_j(t) = \left(\XB_{j+1}(t) - \XB_j(t)\right) \mu(t) - \frac{1}{rT_R}\left(\XC_j(t) - \XC_{j-1}(t)\right).
\end{align*}
By integrating both sides wrt $t \in [0, T]$ and using the upper and lower bounds on $\mu(t)$, we get
\begin{align*}
   & \frac{1}{T_R + \tau_1}\left(\int_0^T \XB_{j+1}(t) dt - \int_0^T \XB_j(t) dt \right) - \frac{1}{rT_R}\left(\int_0^T \XC_j(t) dt - \int_0^T \XC_{j-1}(t) dt\right) \\
   \leq{}& \XB_j(T) + \XC_j(T)-\XB_j(0)-\XC_j(0) \\
   \leq{}& \frac{1}{T_R}\left(\int_0^T \XB_{j+1}(t) dt - \int_0^T \XB_j(t) dt \right) - \frac{1}{rT_R}\left(\int_0^T \XC_j(t) dt - \int_0^T \XC_{j-1}(t) dt\right).
\end{align*}
Next, by dividing by $T$ and taking the limit $T \rightarrow \infty$, we get
\begin{align*} \numberthis \label{eq: x_c_b_j}
    \frac{r}{1 + \frac{\tau_1}{T_R}}\left( \bXB_{j+1} - \bXB_j(t) \right) \leq \bXC_j(t) - \bXC_{j-1}(t) \leq r\left( \bXB_{j+1} - \bXB_j(t) \right) \quad \forall j \in \{1, \hdots, \Nt-1\}.
\end{align*}

Following the same steps for \eqref{eq:cd-4}, we get
\begin{equation} \label{eq: x_c_b_1}
\frac{r}{1 + \frac{\tau_1}{T_R}} \bXB_1 \leq \bXC_0 \leq r\bXB_1.
\end{equation}
By adding \eqref{eq: x_c_b_j} for $1 \leq j \leq k$ along with \eqref{eq: x_c_b_1}, we get
\begin{align*}
    \frac{r}{1 + \frac{\tau_1}{T_R}}\bXB_{k+1} \leq \bXC_k \leq r \bXB_{k+1} \quad \forall k \in \left\{0, 1, \hdots, \Nt-1\right\}. \numberthis \label{eq: coarse_bounds_on_xc}
\end{align*}
We use the above inequality to obtain a coarse lower bound on $n$. We have
\begin{align*}
    n = \bXB_{\Nt}+\bXC_{\Nt} \geq \bXB_{\Nt} + \bXC_{\Nt-1}  \geq \bXB_{\Nt} \left(1+\frac{r}{1+\frac{\tau_1}{T_R}}\right) \geq \alpha T_R\lambda \left(1+\frac{r}{1+\frac{\tau_1}{T_R}}\right), \numberthis \label{eq: coarse_bound_n}
\end{align*}
where the last inequality follows by integrating \eqref{eq:cd-3}. In particular, we have
\begin{align*}
     \alpha \lambda \leq \alphasup \lambda = \limsup_{T\rightarrow \infty}\frac{1}{T}\int_0^T \XB_{\Nt}(t) \mu(t) dt \leq \frac{\bXB_{\Nt}}{T_R} \implies \bXB_{\Nt} \geq \alpha T_R\lambda, \numberthis \label{eq: x_b_nt_cd}
\end{align*}
where the last inequality holds as $\mu(t) \leq 1/T_R$ and noting that we are working with a subsequence for which the limit supremum of $\frac{1}{T}\int_0^T \XB_{\Nt}(t) dt$ is achieved.
Now, using the coarse bound on $n$ given by \eqref{eq: coarse_bound_n}, we obtain a tighter lower bound on $\mu(t)$. We have
\begin{align*}
    \mu(t) \geq \frac{1}{T_R + \frac{\tau_1}{\sqrt{\frac{\alpha rT_R\lambda}{1+\frac{\tau_1}{T_R}}-\kappa\lambda^{1-\gamma}}}} \overset{*}{\geq} \frac{1}{T_R + \frac{2\tau_1\sqrt{1+\frac{\tau_1}{T_R}}}{\sqrt{\delta r T_R\lambda}}} \overset{\Delta}{=} \frac{1}{T_R + \frac{C T_R}{\sqrt{\lambda}}},
\end{align*}
where $(*)$ follows for all $\lambda \geq \lambda_{00}$ as $1-\gamma < 1$ and $\alpha \geq \delta$. This completes part $(b)$ of the proof. Now, using the above lower bound on $\mu(t)$, we use \eqref{eq:cd-1} and \eqref{eq:cd-2} to tighten the bound \eqref{eq: coarse_bounds_on_xc}. In particular, by repeating the same steps as in \eqref{eq: x_c_b_j}-\eqref{eq: coarse_bounds_on_xc}, we get
\begin{align*}
    \frac{r}{1 + \frac{C}{\sqrt{\lambda}}}\bXB_{k+1} \leq \bXC_k \leq r \bXB_{k+1} \quad \forall k \in \left\{0, 1, \hdots, \Nt-1\right\}. \numberthis \label{eq: correct_bound_xc}
\end{align*}
The above inequality completes part $(c)$ of the proof. Next to show part $(d)$, we add \eqref{eq:cd-1} for $j+1$ and \eqref{eq:cd-2} for $j$, to get
\begin{align*}
\dXB_{j+1}(t) + \dXC_j(t) 
&= \frac{\lambda}{\XC_{\Nt}(t)} \left(\XC_{j+1}(t)-\XC_j(t)\right)I(t) - \left(\XC_j(t) - \XC_{j-1}(t)\right)\frac{1}{rT_R} \\
&\leq \frac{\lambda}{n - \alpha T_R \lambda-\kappa\lambda^{1-\gamma}} \left(\XC_{j+1}(t)-\XC_j(t)\right) - \left(\XC_j(t) - \XC_{j-1}(t)\right)\frac{1}{rT_R},
\end{align*}
where the last inequality is followed by admission control. In particular, $\XB_{\Nt}(t) \leq \alpha \lambda T_R + \kappa \lambda^{1-\gamma}$ which implies $X_{\Nt}^C(t) \geq n- \alpha\lambda T_R - \kappa \lambda^{1-\gamma}$ for all $t \in \bbR_+$. Now, by integrating both sides from $0$ to $T$, dividing by $T$, and then taking the limit supremum, we get
$$
\bXC_j - \bXC_{j-1} \leq \frac{rT_R\lambda}{n - \alpha T_R\lambda-\kappa\lambda^{1-\gamma}}  \left(\bXC_{j+1} - \bXC_j\right) \quad \forall j \in \left\{1, \hdots, \Nt-2\right\}.
$$
Similarly, by subtracting \eqref{eq:cd-3} from \eqref{eq:cd-2} for $j = \Nt-1$, we get
\begin{align*}
    \dXC_{\Nt-1}(t) - \dXC_{\Nt}(t)  = \frac{\lambda}{\XC_{\Nt}(t)} \left(\XC_{\Nt}(t) - \XC_{\Nt-1}(t)\right)I(t) - \left(\XC_{\Nt-1}-\XC_{\Nt-2}\right)\frac{1}{rT_R}.
\end{align*}
Thus, we can repeat the steps above, to obtain
$$
\bXC_j - \bXC_{j-1} \leq \frac{rT_R\lambda}{n - \alpha T_R\lambda-\kappa\lambda^{1-\gamma}}  \left(\bXC_{j+1} - \bXC_j\right) \quad \forall j \in \left\{1, \hdots, \Nt-1\right\}.
$$
For the ease of notation, we define $f(n) =\frac{rT_R\lambda}{n - \alpha T_R\lambda-\kappa\lambda^{1-\gamma}}$. Using the above inequality multiple times, we get
\begin{align*}
    \bXC_j - \bXC_{j-1} \leq f(n) \left(\bXC_{j+1} - \bXC_j\right) \leq f(n)^{\Nt-j}\left(\bXC_{\Nt} - \bXC_{\Nt-1}\right) \quad \forall j \in \left\{1, \hdots, \Nt-1\right\}. \numberthis \label{eq: telescopic_sum_bound_j}
\end{align*}
Now, by taking $j = 1$ in \eqref{eq:cd-1}, we get
\begin{align*}
    \dXB_1(t) &= \frac{\lambda}{\XC_{\Nt}} \left(\XC_1(t) - \XC_0(t)\right)I(t) - \XB_1(t)\mu(t) \\
    &\leq \frac{\lambda}{n-\alpha T_R\lambda-\kappa\lambda^{1-\gamma}} \left(\XC_1(t) - \XC_0(t)\right) - \XB_1(t)\frac{1}{T_R+\frac{CT_R}{\sqrt{\lambda}}}.
\end{align*}
Now, by integrating both sides from $0$ to $T$, dividing by $T$, and then taking the limit supremum, we get
\begin{align*}
    \bXC_0 &\leq r\bXB_1 \leq f(n)\left(1+\frac{C}{\sqrt{\lambda}}\right)\left(\bXC_{1} - \bXC_0\right) \\
    &\leq f(n)^j\left(1+\frac{C}{\sqrt{\lambda}}\right)\left(\bXC_{j} - \bXC_{j-1}\right) \quad \forall j \in \{1, \hdots, \Nt\}. \numberthis \label{eq: telescopic_sum_bound_0}
\end{align*}
The above inequality combined with \eqref{eq: telescopic_sum_bound_j} completes part $(d)$ of the proof. Now, we use the upper bounds on $\bXC_j-\bXC_{j-1}$ to obtain a lower bound on $\bXC_{\Nt}-\bXC_{\Nt-1}$ as follows:
\begin{align*}
    n & = \bXC_{\Nt} + \bXB_{\Nt} \\
    & \overset{\eqref{eq: correct_bound_xc}}{\leq} \bXC_{\Nt} + \frac{\bXC_{\Nt-1}}{r}\left(1+\frac{C}{\sqrt{\lambda}}\right) \\
    &\overset{*}{\leq} \frac{r+1+\frac{C}{\sqrt{\lambda}}}{r} \bXC_{\Nt}  \\
    &= \frac{r+1+\frac{C}{\sqrt{\lambda}}}{r}\left(\sum_{j=1}^{\Nt} \left(\bXC_j - \bXC_{j-1}\right) + \bXC_0\right) \\
    &\overset{**}{\leq} \frac{r+1+\frac{C}{\sqrt{\lambda}}}{r}\left(1+\frac{C}{\sqrt{\lambda}}\right)\left(\bXC_{\Nt} - \bXC_{\Nt-1}\right) \sum_{j=0}^{\Nt} f(n)^{\Nt-j} \\
    &= \frac{r+1+\frac{C}{\sqrt{\lambda}}}{r}\left(1+\frac{C}{\sqrt{\lambda}}\right)\left(\bXC_{\Nt} - \bXC_{\Nt-1}\right) \left((\Nt+1)\mathbbm{1}\left\{f(n) = 1\right\} + \frac{1-f(n)^{\Nt+1}}{1-f(n)}\mathbbm{1}\left\{f(n) \neq 1\right\}\right),
\end{align*}
where $(*)$ follows as $\bXC_{\Nt-1} \leq \bXC_{\Nt}$. Next $(**)$ follows by \eqref{eq: telescopic_sum_bound_j} and \eqref{eq: telescopic_sum_bound_0}. Now, the above inequality implies the following:
\begin{align*}
    \bXC_{\Nt} - \bXC_{\Nt-1} \geq \frac{n r}{r+1+\frac{C}{\sqrt{\lambda}}}\left(1-\frac{C}{\sqrt{\lambda}}\right) \left(\frac{1}{\Nt+1}\mathbbm{1}\left\{f(n) = 1\right\} + \frac{1-f(n)}{1-f(n)^{\Nt+1}}\mathbbm{1}\left\{f(n) \neq 1\right\}\right). \numberthis \label{eq: result_of_telescopic_bounds}
\end{align*}
Note that we used the bound $1/(1+x) \geq 1-x$ for $x \geq 0$ to conclude the above inequality. This bound on $\bXC_{\Nt}-\bXC_{\Nt-1}$ completes part $(e)$ of the proof. Now, using the above inequality, we get
\begin{align*}
    n &= \bXB_{\Nt} + \bXC_{\Nt} \\
    &\overset{*}{\geq} \bXB_{\Nt} + r \bXB_{\Nt}\left(1-\frac{C}{\sqrt{\lambda}}\right) + \bXC_{\Nt} - \bXC_{\Nt-1} \\
    &\overset{**}{\geq} \left(r+1-\frac{C}{\sqrt{\lambda}}\right)\bXB_{\Nt} + \frac{n r}{r+1+\frac{C}{\sqrt{\lambda}}}\left(1-\frac{C}{\sqrt{\lambda}}\right) \left(\frac{1}{\Nt+1}\mathbbm{1}\left\{f(n) = 1\right\} + \frac{1-f(n)}{1-f(n)^{\Nt+1}}\mathbbm{1}\left\{f(n) \neq 1\right\}\right),
\end{align*}
where, $(*)$ follows by adding and subtracting $\bXC_{\Nt-1}$ and using \eqref{eq: correct_bound_xc}  and noting that $1/(1+x) \geq 1-x$ for $x \geq 0$. Next, $(**)$ follows by using $r<1$ to bound the first term and \eqref{eq: result_of_telescopic_bounds} to bound the second term. The above inequality then provides us with the following upper bound on $\bXB_{\Nt}$, completing part $(f)$ of the proof.
\begin{align*}
    \bXB_{\Nt} \leq \frac{n}{r+1-\frac{C}{\sqrt{\lambda}}} - \frac{n r}{(r+1)^2-\frac{C^2}{\lambda}}\left(1-\frac{C}{\sqrt{\lambda}}\right) \left(\frac{1}{\Nt+1}\mathbbm{1}\left\{f(n) = 1\right\} + \frac{1-f(n)}{1-f(n)^{\Nt+1}}\mathbbm{1}\left\{f(n) \neq 1\right\}\right).
\end{align*}
Now, we move on to showing part $(g)$, which is the last part of the proof. By \eqref{eq: x_b_nt_cd}, we get
\begin{align*}
    \lambda \alpha &\leq \lambda \alphasup \leq \frac{\bXB_{\Nt}}{T_R} \\
    &\leq \frac{n}{T_R\left(r+1-\frac{C}{\sqrt{\lambda}}\right)} - \frac{n r}{T_R\left((r+1)^2-\frac{C^2}{\lambda}\right)}\left(1-\frac{C}{\sqrt{\lambda}}\right) \left(\frac{1}{\Nt+1}\mathbbm{1}\left\{f(n) = 1\right\} + \frac{1-f(n)}{1-f(n)^{\Nt+1}}\mathbbm{1}\left\{f(n) \neq 1\right\}\right) \numberthis \label{eq: loose_bound_lambda} \\
    &\overset{*}{\leq} \frac{n}{T_R\left(r+1-\frac{C}{\sqrt{\lambda}}\right)} - \frac{n r}{T_R\left((r+1)^2-\frac{C^2}{\lambda}\right)} \left(1-\frac{C}{\sqrt{\lambda}}\right)  \min\left\{\frac{1}{\Nt+1}, \frac{\frac{2}{\alpha}-1}{\left(\frac{2}{\alpha}\right)^{\Nt+1}-1}\right\} \\
    &\overset{**}{\leq} \frac{n}{T_R\left(r+1\right)}\left(1+\frac{C}{\sqrt{\lambda}}\right) - \frac{n r}{2T_R(r+1)^2} \min\left\{\frac{1}{\Nt+1}, \frac{\frac{2}{\alpha}-1}{\left(\frac{2}{\alpha}\right)^{\Nt+1}-1}\right\} \\
    &\overset{***}{\leq} \frac{n}{T_R\left(r+1\right)} - \frac{n r}{4T_R(r+1)^2} \min\left\{\frac{1}{\Nt+1}, \frac{\frac{2}{\alpha}-1}{\left(\frac{2}{\alpha}\right)^{\Nt+1}-1}\right\} \\
    &\leq \frac{n}{T_R(1+r)}\left(1-\frac{\epsilon}{T_R(1+r)\delta}\right), \numberthis \label{eq: lambda_tight_bound}
\end{align*}
where $(*)$ follows as $n \geq (1+r)\alpha \lambda T_R - \alpha CT_R\sqrt{\lambda}$ by \eqref{eq: loose_bound_lambda}, which further implies 
$$
f(n) =\frac{rT_R\lambda}{n - \alpha T_R\lambda-\kappa\lambda^{1-\gamma}} \leq \frac{rT_R\lambda}{r\alpha T_R\lambda- \alpha  CT_R\sqrt{\lambda}-\kappa\lambda^{1-\gamma}} \leq \frac{1}{\alpha-\frac{2\kappa}{r T_R}\lambda^{-\gamma}} \leq \frac{2}{\alpha}
$$
where the last two inequalities follows for all $\lambda \geq \lambda_{01}$ for some $\lambda_{01} > 0$ as $\gamma \in [1/3, 1/2)$. Thus, $(*)$ follows as $\frac{x-1}{x^{\Nt+1}-1}$ is a decreasing function of $x \geq 0$. Next, $(**)$ follows for all $\lambda \geq \lambda_{02}$ for some $\lambda_{02} > 0$ as $1-C/\sqrt{\lambda} \geq 1/2$, $(r+1)^2-C^2/\lambda \leq (r+1)^2$, and $\frac{1}{1-C/(r+1)/\sqrt{\lambda}} \leq 1+C/\sqrt{\lambda}$ as $r>0$. Further, $(***)$ follows for all $\lambda \geq \lambda_{03}$ for some $\lambda_{03} > 0$ as the last term is $\Theta(\lambda)$ while $\frac{C n}{T_R(r+1)\sqrt{\lambda}}$ is $\Theta(\sqrt{\lambda})$. The last inequality follows by setting
$$
\epsilon \overset{\Delta}{=} \frac{r\delta T_R}{4}\min\left\{\frac{1}{\Nt+1}, \frac{\frac{2}{\alpha}-1}{\frac{2^{\Nt+1}}{\alpha^{\Nt+1}}-1}\right\},
$$
Now, using \eqref{eq: lambda_tight_bound}, we get the following lower bound on $n$:
\begin{align*}
n &\geq (1+r)T_R \alpha \lambda \left(1-\frac{\epsilon}{T_R(1+r)\delta}\right)^{-1} \overset{*}{\geq} (1+r)T_R \alpha \lambda \left(1+\frac{\epsilon}{T_R(1+r)\delta}\right) \overset{**}{\geq}  (1+r)T_R \alpha \lambda + \epsilon \lambda,
\end{align*}
where $(*)$ follows by noting that $1/(1-x) \geq 1+x$ for $x \in [0, 1]$ and noting that $\epsilon/(T_R(1+r)\delta) \leq 1$. Lastly, $(**)$ follows as $\alpha \geq \delta$. This completes the proof by setting $\lambda_0 = \max_{k \in [4]}\{\lambda_{0k}\}$. \hfill $\square$
\endproof

\subsection{Pack Size versus Fleet Size}
By \Cref{thm:informal_upper-bound}, to obtain a service level of $\alpha$ under the Power-of-$d$ vehicles dispatch algorithm, we have
\begin{equation} \label{eq: pack_vs_fleet}
n \asymp (1+r)T_R \lambda + \lambda^{\frac{1+1/\Nt}{2+1/\Nt}},
\end{equation}
for large enough $m$. Note that the right-hand side of the above equation is a decreasing function of $\Nt$. Thus, a smaller fleet size is sufficient, given a larger battery. As evident from the discussion on ICE cars, the exponent of the second-order term is less than 2/3 due to the double counting of partially charged cars. Now, if the pack size is large, the fraction of partially charged cars will be higher, resulting in increased benefits. \Cref{eq: pack_vs_fleet} quantifies this phenomenon. In short, a trade-off exists between the fleet and the pack size quantified by \eqref{eq: pack_vs_fleet}. It is also interesting to note that the benefits from a larger pack size exhibit the law of diminishing marginal utility. In particular, define $f(\Nt) = (1+r)T_R \lambda + \lambda^{\frac{1+1/\Nt}{2+1/\Nt}}$ and note that $f^{\prime\prime}(\Nt) \geq 0$.
Thus, $n$ is a convex function of $\Nt$, which implies that the benefit of increasing the pack size reduces for larger values of the pack size.

\subsection{Battery Charging versus Battery Swapping}
Battery swapping is an alternate approach to replenishing the SoC of an EV. Rather than charging the battery pack, a spare one with a higher SoC replaces it. In particular, the spare battery pack charges while the EV is driving, which later replaces the original battery pack. One can interpret a fleet of EVs equipped with battery swapping technology as a fleet of ICE vehicles. In particular, the time required to replenish the charge by swapping the battery would be of the same order as refueling an ICE vehicle. Thus, by \cite{besbes2022spatial}, we can conclude that the fleet size $(n_{bs})$ and battery $(b_{bs})$ requirements for EV with battery swapping is
\begin{equation}
    n_{bs} = T_R \lambda + \Theta\left(\lambda^{2/3}\right),  \quad b_{bs} = (1+r)T_R \lambda + \Theta\left(\lambda^{2/3}\right). \label{eq: battery_swapping} 
\end{equation}
In particular, the fleet size requirement follows from \cite{besbes2022spatial} and for each EV, at least $(1+r)$ batteries are required to balance the aggregate discharge rate with the aggregate charge rate. On the other hand, by Theorem~\ref{thm:lower-bound} and \ref{thm:informal_upper-bound}, the fleet size $(n_{bc})$ and battery requirements $(b_{bc})$ for EV with battery charging is
\begin{equation}
    n_{bc} = b_{bc} = (1+r)T_R \lambda + \Theta\left(\lambda^{1-\gamma}\right) \quad  \text{for} \  \gamma \in \left[\frac{1}{3}, \frac{1}{2+1/\Nt}\right). \label{eq: battery_charging}
\end{equation}
Thus, even though the fleet size requirement under battery charging is higher, the number of battery packs required is smaller compared to battery swapping, especially with a high density of chargers. In today's world, battery packs are the bottleneck in manufacturing electric vehicles. Thus, battery swapping is a viable option for ride-hailing systems only if battery manufacturing is accelerated in the future. 

\section{Varying Arrival Rate} \label{app: varying_arrival_rate}
\subsection{Proof of Theorem~\ref{thm: informal_varying_arrival_rate}}
Recall that $\alphasup$ is the long-run fraction of trip requests met, given by
\begin{align*}
    \alphasup = \frac{1}{\lambda_{\operatorname{avg}}}\limsup_{T \rightarrow \infty} \frac{1}{T} \int_0^T \lambdaeff^{\pol}(t) dt. \numberthis \label{eq:alpha_eff_var}
\end{align*}
In addition, we define $\bar{\alpha}^{\pi}_v$ and $\bar{\alpha}^{\pi}_p$ as the long-run fraction of trips met during the valleys and peaks, respectively. Mathematically, we have
\begin{subequations} \label{eq:def_alpha_var}
\begin{align}
   \lambda \bar{\alpha}^{\pi}_v &= \frac{T_v+T_p}{T_v}\limsup_{T \rightarrow \infty} \frac{1}{T} \sum_{i=0}^{\left\lfloor \frac{T-T_v}{T_v+T_p} \right\rfloor-1} \int_{i(T_v+T_p)}^{(i+1)T_v+iT_p} \lambdaeff^{\pol}(t) dt \label{eq:def_alpha_v} \\
    a\lambda \bar{\alpha}^{\pi}_p &= \frac{T_v+T_p}{T_p}\limsup_{T \rightarrow \infty} \frac{1}{T} \sum_{i=0}^{\left\lfloor \frac{T}{T_v+T_p} \right\rfloor-1} \int_{(i+1)T_v+iT_p}^{(i+1)(T_v+T_p)} \lambdaeff^{\pol}(t) dt. \label{eq:def_alpha_p}
\end{align} 
\end{subequations}
Now, note that the following equations are satisfied:
\begin{subequations} \label{eq:def_var_arrival_rate}
    \begin{align}
        q(t) &\geq \int_{\max\{0, t-T_R\}}^t \lambdaeff^\pi(s)ds \label{eq:def_var-1} \\
        \ecc r_c &\geq \bar{q} r_d \label{eq:def_var-2} \\
        \cc(t) &\leq \min\{m, n - q(t)\}, \label{eq:def_var-3}
    \end{align}
\end{subequations}
for all $t>0$. Additionally, we assume that $T_v$ is divisible by $T_R$. Equation \eqref{eq:def_var-1} holds as any customer spends at least $T_R$ amount of time in the system. Thus, any admitted customer within $[\max\{t-T_R, 0\}, t]$ time must be currently present in the system. Next, \eqref{eq:def_var-2} ensures that the long-run charge rate is equal to the long-run discharge rate. Lastly, \eqref{eq:def_var-3} bounds the number of EVs charging by the minimum of $m$ and the number of vehicles not serving a trip. Now, we show that these equations imply that \Cref{thm: informal_varying_arrival_rate} holds. We first present a lemma that relates the service level with the fleet size and number of chargers under an arbitrary policy $\pi$. Such a result is analogous to analyzing the ODE \eqref{eq:partial_ode} for the constant arrival rate setting.
\begin{lemma} \label{lemma: stability_equations_varying_arrival_rate}
    For any arbitrary policy $\pi$ satisfying \eqref{eq:def_var_arrival_rate}, the following relations hold.
\begin{subequations}
\begin{align}
\alphasup\lambda_{\operatorname{avg}} &= \frac{\bar{\alpha}^{\pi}_v \lambda T_v + \bar{\alpha}^{\pi}_p a\lambda T_p}{T_v+T_p} \label{eq: service_level_balance} \\
n &\geq \bar{\alpha}^{\pi}_p a \lambda T_R \label{eq: peak_demand_bound} \\
n &\geq (1+r)\alphasup \lambda_{\operatorname{avg}}T_R, \quad m \geq r\alphasup\lambda_{\operatorname{avg}}T_R \label{eq: stable_soc_bound_weak} \\
r_d \alphasup\lambda_{\operatorname{avg}}T_R &\leq \frac{T_vr_c}{T_v+T_p}\min\left\{m, n-T_R\lambda\bar{\alpha}_v^\pi+\frac{T_R^2\lambda}{T_v}\right\} + \frac{T_pr_c}{T_v+T_p}\min\left\{m, n-T_R a\lambda\bar{\alpha}_p^\pi+\frac{aT_R^2\lambda}{T_p}\right\} \label{eq: stable_soc_bound}
\end{align}
\end{subequations}
\end{lemma}
We defer the proof of \Cref{lemma: stability_equations_varying_arrival_rate} to the end of the section. First, we present the proof of \Cref{thm: informal_varying_arrival_rate} below.
\proof{Proof of \Cref{thm: informal_varying_arrival_rate}}

\textbf{Case I:} The first case of the theorem follows directly from \eqref{eq: stable_soc_bound_weak} in \Cref{lemma: stability_equations_varying_arrival_rate}.

\textbf{Case II:} In this case, we have
\begin{align*}
    \alpha_p = \frac{\alpha \lambda_{\operatorname{avg}}(T_v+T_p)}{a\lambda T_p} - \frac{T_v}{aT_p}.
\end{align*}
Now, using the above, we have
\begin{align*}
   n_{\operatorname{dp}} + n_{\operatorname{cp}, \eta} &= \alpha_p a \lambda T_R + \alpha r \lambda_{\operatorname{avg}}T_R - \eta \frac{T_v}{T_p}\lambda_{\operatorname{avg}}T_R \\
   &= \alpha \lambda_{\operatorname{avg}}T_R \left(1+\frac{T_v}{T_p}\right) - \lambda T_R \frac{T_v}{T_p}+\alpha r \lambda_{\operatorname{avg}}T_R - \eta \frac{T_v}{T_p}\lambda_{\operatorname{avg}}T_R \\
   &= (1+r)\alpha \lambda_{\operatorname{avg}}T_R + \frac{T_v}{T_p}\left(\alpha - \frac{T_v+T_p}{T_v+aT_p}-\eta\right)\lambda_{\operatorname{avg}}T_R.
\end{align*}
Similarly, by substituting $\lambda_{\operatorname{avg}} = \frac{T_v + aT_p}{T_v+T_p} \lambda$, we get that the range of $\eta$ is $\left[0, \alpha - \frac{T_v+T_p+aT_R(1+T_p/T_v)}{T_v+aT_p}\right]$. We prove the contrapositive. In particular, we show that if for all $\eta \in \left[0, \alpha - \frac{T_v+T_p+aT_R(1+T_p/T_v)}{T_v+aT_p}\right]$, we have
\begin{equation} \label{eq: n_or_m}
    n < (1+r)\alpha \lambda_{\operatorname{avg}}T_R + \frac{T_v}{T_p}\left(\alpha - \frac{T_v+T_p}{T_v+aT_p}-\eta\right)\lambda_{\operatorname{avg}}T_R-\frac{aT_R^2\lambda}{T_p}, \quad \text{OR} \quad m < r\alpha\lambda_{\operatorname{avg}}T_R+\eta\lambda_{\operatorname{avg}}T_R,
\end{equation}
then $\alphasup < \alpha$. First, consider the case when the former holds for all $\eta$ in \eqref{eq: n_or_m}. Then, we have $n < (1+r)\alpha\lambda_{\operatorname{avg}}T_R$. Thus, by \eqref{eq: stable_soc_bound_weak}, we get $\alphasup < \alpha$. Now, consider the case when the latter holds for all $\eta$ in \eqref{eq: n_or_m}. Then, we have $m < r\alpha\lambda_{\operatorname{avg}}T_R$. Thus, by \eqref{eq: stable_soc_bound_weak}, we get $\alphasup < \alpha$. Lastly, consider the case when there exists a non-zero subset of $\left[0, \alpha - \frac{T_v+T_p+aT_R(1+T_p/T_v)}{T_v+aT_p}\right]$ for which the former holds in \eqref{eq: n_or_m}. Let $\eta^{(1)}$ be the supremum of such constants. Similarly, a non-zero subset for which the latter holds and let $\eta^{(2)}$ be the infimum of such constants. To ensure that \eqref{eq: n_or_m} holds for all $\eta \in \left[0, \alpha - \frac{T_v+T_p+aT_R(1+T_p/T_v)}{T_v+aT_p}\right]$, we must have $\eta^{(1)} \geq \eta^{(2)}$. In addition, either the supremum or the infimum is attained. Thus, there exists a $\eta \in [\eta^{(2)}, \eta^{(1)}]$ such that
\begin{equation} \label{eq: n_or_m_simplified_1}
    n \leq (1+r)\alpha \lambda_{\operatorname{avg}}T_R + \frac{T_v}{T_p}\left(\alpha - \frac{T_v+T_p}{T_v+aT_p}-\eta\right)\lambda_{\operatorname{avg}}T_R-\frac{aT_R^2 \lambda}{T_p}, \quad m < r\alpha\lambda_{\operatorname{avg}}T_R+\eta\lambda_{\operatorname{avg}}T_R,
\end{equation}
or
\begin{equation} 
    n < (1+r)\alpha \lambda_{\operatorname{avg}}T_R + \frac{T_v}{T_p}\left(\alpha - \frac{T_v+T_p}{T_v+aT_p}-\eta\right)\lambda_{\operatorname{avg}}T_R-\frac{aT_R^2\lambda}{T_p}, \quad m \leq r\alpha\lambda_{\operatorname{avg}}T_R+\eta\lambda_{\operatorname{avg}}T_R. \label{eq: n_or_m_simplified_2}
\end{equation}
Now, if \eqref{eq: n_or_m_simplified_1} is satisfied, and $\alphasup \geq \alpha$, then we obtain a contradiction. By \eqref{eq: stable_soc_bound}, we have
\begin{align*}
    n &\geq r\alphasup \lambda_{\operatorname{avg}}T_R \left(1+\frac{T_v}{T_p}\right) - \frac{mT_v}{T_p} + T_R a\lambda \bar{\alpha}_p^{\pi} - \frac{aT_R^2\lambda}{T_p} \\
    &\overset{(a)}{>} r\alphasup \lambda_{\operatorname{avg}}T_R \left(1+\frac{T_v}{T_p}\right) - \frac{T_v}{T_p}\left(r\alpha \lambda_{\operatorname{avg}}T_R+\eta\lambda_{\operatorname{avg}}T_R\right) + T_R a\lambda \bar{\alpha}_p^{\pi} - \frac{aT_R^2\lambda}{T_p} \\
    &\overset{(b)}{\geq} r\alphasup \lambda_{\operatorname{avg}}T_R \left(1+\frac{T_v}{T_p}\right) - \frac{T_v}{T_p}\left(r\alpha \lambda_{\operatorname{avg}}T_R+\eta\lambda_{\operatorname{avg}}T_R\right) + T_R a\lambda\left(\alphasup \left(1+\frac{T_v}{aT_p}\right)-\frac{T_v}{aT_p}\right) - \frac{aT_R^2\lambda}{T_p} \\
    &\overset{(c)}{\geq} r\alpha \lambda_{\operatorname{avg}}T_R - \frac{T_v}{T_p}\eta\lambda_{\operatorname{avg}}T_R + T_R a\lambda_{\operatorname{avg}}\frac{T_v+T_p}{T_v+aT_p}\left(\alpha \left(1+\frac{T_v}{aT_p}\right)-\frac{T_v}{aT_p}\right) - \frac{aT_R^2\lambda}{T_p} \\
    &= (1+r)\alpha \lambda_{\operatorname{avg}}T_R + \frac{T_v}{T_p}\left(\alpha - \frac{T_v+T_p}{T_v+aT_p}-\eta\right)\lambda_{\operatorname{avg}}T_R-\frac{aT_R^2\lambda}{T_p}, \numberthis \label{eq: contradiction_n_m_varying_arrival_rate}
\end{align*}
where $(a)$ follows by using the upper bound on $m$ given by \eqref{eq: n_or_m_simplified_1}. Next, $(b)$ follows by substituting $\bar{\alpha}_v^{\pi}\leq1$ in \eqref{eq: service_level_balance} to obtain an upper bound on $\bar{\alpha}_p^{\pi}$. Lastly, $(c)$ follows as we assume $\alphasup \geq \alpha$. Note that the last equation contradicts \eqref{eq: n_or_m_simplified_1}. Similarly, we can obtain a contradiction if instead, \eqref{eq: n_or_m_simplified_2} were satisfied. This completes the proof of Case II.

\textbf{Case III:} Similarly to the previous case, we prove the contrapositive for this case as well. In particular, we show that if for all $\eta \in \left[0, \frac{\alpha rT_p}{T_v}-\frac{aT_R(1+T_p/T_v)}{T_v+aT_p}\right]$, we have
\begin{equation} \label{eq: n_or_m_case_3}
    n < (1+r)\alpha \lambda_{\operatorname{avg}}T_R + \frac{T_v}{T_p}\left(\alpha - \frac{T_v+T_p}{T_v+aT_p}-\eta\right)\lambda_{\operatorname{avg}}T_R-\frac{aT_R^2\lambda}{T_p}, \quad \text{OR} \quad m < r\alpha\lambda_{\operatorname{avg}}T_R+\eta\lambda_{\operatorname{avg}}T_R,
\end{equation}
then, $\alphasup < \alpha$. First, consider the case when the former in \eqref{eq: n_or_m_case_3} holds for all $\eta$. Then, we have
\begin{equation*}
    n < \alpha \lambda_{\operatorname{avg}}T_R + \frac{T_v}{T_p}\left(\alpha - \frac{T_v+T_p}{T_v+aT_p}\right)\lambda_{\operatorname{avg}}T_R.
\end{equation*}
In addition, by \eqref{eq: peak_demand_bound} and \eqref{eq: service_level_balance}, we get
\begin{align*}
    n \geq \bar{\alpha}_p^{\pi} a \lambda T_R \geq  a\lambda T_R\left(\alphasup\left(1+\frac{T_v}{aT_p}\right)-\frac{T_v}{aT_p}\right) = \alphasup \lambda_{\operatorname{avg}}T_R +\frac{T_v}{T_p}\left(\alphasup -  \frac{T_v+T_p}{T_v+aT_p}\right)\lambda_{\operatorname{avg}}T_R.
\end{align*}
Thus, the lower and upper bounds on $n$ imply that $\alphasup < \alpha$. 

Now, consider the case when the latter in \eqref{eq: n_or_m_case_3} holds for all $\eta$. Then, we have $m < r\alpha \lambda_{\operatorname{avg}}T_R$. Thus, \eqref{eq: stable_soc_bound_weak} implies $\alphasup < \alpha$. Lastly, consider the case when there exists some $\eta$ for which the former in \eqref{eq: n_or_m_case_3} holds and some $\eta$ for which the latter holds. Thus, similar to Case II, there exists $\eta \in \left[0, \frac{\alpha rT_p}{T_v}-\frac{aT_R(1+T_p/T_v)}{T_v+aT_p}\right]$, such that either
\begin{equation} \label{eq: n_or_m_case_3_simplified_1}
    n \leq (1+r)\alpha \lambda_{\operatorname{avg}}T_R + \frac{T_v}{T_p}\left(\alpha - \frac{T_v+T_p}{T_v+aT_p}-\eta\right)\lambda_{\operatorname{avg}}T_R-\frac{aT_R^2\lambda}{T_p}, \quad m < r\alpha\lambda_{\operatorname{avg}}T_R+\eta\lambda_{\operatorname{avg}}T_R,
\end{equation}
or
\begin{equation} \label{eq: n_or_m_case_3_simplified_2}
    n < (1+r)\alpha \lambda_{\operatorname{avg}}T_R + \frac{T_v}{T_p}\left(\alpha - \frac{T_v+T_p}{T_v+aT_p}-\eta\right)\lambda_{\operatorname{avg}}T_R-\frac{aT_R^2\lambda}{T_p}, \quad m \leq r\alpha\lambda_{\operatorname{avg}}T_R+\eta\lambda_{\operatorname{avg}}T_R.
\end{equation}
As $\alpha$ is large enough for this case, we can verify that $\left[0, \frac{\alpha rT_p}{T_v}\right] \subseteq \left[0, \alpha-\frac{T_v+T_p}{T_v+aT_p}\right]$. Thus, by \eqref{eq: contradiction_n_m_varying_arrival_rate}, we conclude that \eqref{eq: n_or_m_case_3_simplified_1} implies $\alphasup < \alpha$. Similarly, \eqref{eq: n_or_m_case_3_simplified_2} also implies $\alphasup < \alpha$. This completes the proof of Case III.
\hfill $\square$
\endproof
\proof{Proof of \Cref{lemma: stability_equations_varying_arrival_rate}}
We start by proving the first equation. By the definition of $\alphasup$, we get
\begin{align*}
    \lambda_{\operatorname{avg}}\alphasup &\overset{\eqref{eq:alpha_eff_var}}{=} \limsup_{T \rightarrow \infty} \frac{1}{T} \int_{0}^{T} \lambdaeff^{\pol}(t) dt \\
    &= \limsup_{T \rightarrow \infty} \frac{1}{T} \sum_{i=0}^{\left\lfloor \frac{T-T_v}{T_v+T_p} \right\rfloor} \int_{i(T_v+T_p)}^{(i+1)T_v+iT_p} \lambdaeff^{\pol}(t) dt + \limsup_{T \rightarrow \infty} \frac{1}{T} \sum_{i=0}^{\left\lfloor \frac{T}{T_v+T_p} \right\rfloor-1} \int_{(i+1)T_v+iT_p}^{(i+1)(T_v+T_p)} \lambdaeff^{\pol}(t) dt \\ 
    &\overset{\eqref{eq:def_alpha_var}}{=} \frac{\bar{\alpha}^{\pi}_v \lambda T_v + \bar{\alpha}^{\pi}_p a\lambda T_p}{T_v+T_p}. \numberthis \label{eq:alpha_relation}
\end{align*}
Note that, for the second equality, we divide the integral into two parts, corresponding to the peaks and valleys. Next, we use the definition of $\bar{\alpha}_v^\pi$ and $\bar{\alpha}_p^\pi$ given by \eqref{eq:def_alpha_var}. Now, we prove \eqref{eq: peak_demand_bound} below.
 \begin{align*}
     n \geq q(T+T_R) &\overset{\eqref{eq:def_var-1}}{\geq}  \int_{T}^{T+T_R} \lambdaeff^\pi(t) dt \quad\forall T \ : \ \frac{T}{T_v+T_p}-\left\lfloor\frac{T}{T_v+T_p}\right\rfloor \in \left[\frac{T_v}{T_v+T_p}, 1-\frac{T_R}{T_v+T_p}\right] \\
      &\overset{(a)}{\geq} \frac{T_R}{T_p} \int_{(i+1)T_v+iT_p}^{(i+1)(T_v+T_p)} \lambdaeff^\pi(t) dt \quad \forall i \in \bbZ_+ \\
     &\overset{(b)}{\geq} \frac{T_R}{kT_p} \sum_{i=0}^{k-1}\int_{(i+1)T_v+iT_p}^{(i+1)(T_v+T_p)} \lambdaeff^\pi(t) dt \quad \forall k \in \bbZ_+ \\
     &\geq \bar{\alpha}_p^{\pi}a\lambda T_R,
 \end{align*}
where $(a)$ follows by taking the average of the previous expression for $T \in \{(i+1)T_v+iT_p+kT_R : k \in \left[T_p/T_R\right]\}$ and using the fact that $T_p$ is divisible by $T_R$. Similarly, $(b)$ follows by taking the average over the first $k$ peaks in the arrival rate. Now, the last inequality follows by first taking the limit supremum as $k \rightarrow \infty$, and then, using the definition of $\bar{\alpha}_p^\pi$ given by \eqref{eq:def_alpha_p}. Next, to prove \eqref{eq: stable_soc_bound_weak}, we use \eqref{eq:def_var-2}, to get $\ecc r_c \geq \bar{q}r_d$. To further simplify, we start by lower-bounding the long-run average number of customers in the system $(\bar{q})$. We have
\begin{align*}
    \bar{q} &= \limsup_{T\rightarrow \infty} \frac{1}{T} \int_0^T q(t)dt \overset{\eqref{eq:def_var-1}}{\geq} \limsup_{T\rightarrow \infty} \frac{1}{T} \int_0^T \int_{\max\{t-T_R, 0\}}^{t} \lambdaeff^{\pi}(s) ds dt  \\
    &\overset{(a)}{=} \limsup_{T\rightarrow \infty} \frac{1}{T} \int_{0}^T \int_{s}^{\min\{s+T_R, T\}} \lambdaeff^{\pi}(s) dt ds \\
    &= \limsup_{T\rightarrow \infty} \frac{1}{T} \int_{0}^T\min\{T_R, T-s\}\lambdaeff^{\pi}(s) ds \\
    &\geq T_R\limsup_{T\rightarrow \infty} \frac{1}{T} \int_{0}^{T-T_R}\lambdaeff^{\pi}(s) ds \\
    &\overset{(b)}{=} T_R\limsup_{T\rightarrow \infty} \frac{1}{T} \int_{0}^{T}\lambdaeff^{\pi}(s) ds \\
    &\overset{\eqref{eq:def_alpha_var}}{=} \alphasup \lambda_{\operatorname{avg}} T_R \overset{\eqref{eq:alpha_relation}}{=} \frac{\lambda T_R}{T_v+T_p} \left(T_v\bar{\alpha}_v^{\pi} + aT_p\bar{\alpha}_p^{\pi}\right), \numberthis \label{eq:q_bar_var_bound}
\end{align*}
where $(a)$ follows by interchanging the order of the integrals. Next, $(b)$ follows by noting that $0\leq \int_{T-T_R}^T \lambdaeff^{\pi}(s)/T ds \leq T_Ra\lambda/T \rightarrow 0$ as $T \rightarrow \infty$. Now, we upper bound the long-run average number of EVs charging in the system.
\begin{align*}
   r\bar{q} \overset{\eqref{eq:def_var-2}}{\leq} \ecc \overset{\eqref{eq:def_var-3}}{\leq}{}& \liminf_{T \rightarrow \infty} \frac{1}{T}\int_0^T \min\left\{m, n - q(t)\right\} dt
    \leq{} \min\left\{m, n - \bar{q}\right\},
\end{align*}
where the inequality follows using Jensen's inequality as $\min$ is a concave function. By the above equation, we immediately get
\begin{align*}
    m \geq r\bar{q} \overset{\eqref{eq:q_bar_var_bound}}{\geq} r\alphasup\lambda_{\operatorname{avg}}T_R, \quad n \geq (1+r)\bar{q} \overset{\eqref{eq:q_bar_var_bound}}{\geq} (1+r)\alphasup\lambda_{\operatorname{avg}}T_R.
\end{align*}
Lastly, to prove \eqref{eq: stable_soc_bound}, we obtain a more fine tuned bound on $\ecc$ by dividing the integral over $t$ into peaks and valleys as follows:
\begin{align*}
    \ecc \overset{\eqref{eq:def_var-3}}{\leq}{}& \liminf_{T \rightarrow \infty} \frac{1}{T}\int_0^T \min\left\{m, n - q(t)\right\} dt \\
    \overset{\eqref{eq:def_var-1}}{\leq}{}& \liminf_{T \rightarrow \infty} \frac{1}{T}\int_0^T \min\left\{m, n - \int_{\max\{t-T_R, 0\}}^t \lambdaeff^{\pi}(s)ds\right\} dt \\
    ={}& \liminf_{T \rightarrow \infty} \frac{1}{T} \sum_{i=0}^{\left\lfloor \frac{T-T_v}{T_v+T_p} \right\rfloor} \int_{i(T_v+T_p)}^{(i+1)T_v+iT_p} \min\left\{m, n - \int_{\max\{t-T_R, 0\}}^t \lambdaeff^{\pi}(s) ds\right\} dt \\
    &+ \liminf_{T \rightarrow \infty} \frac{1}{T} \sum_{i=0}^{\left\lfloor \frac{T}{T_v+T_p} \right\rfloor-1} \int_{(i+1)T_v+iT_p}^{(i+1)(T_v+T_p)} \min\left\{m, n - \int_{\max\{t-T_R, 0\}}^t \lambdaeff^{\pi}(s) ds\right\} dt 
\end{align*}
Now, we simplify each of the above integrals separately. For the valleys, we have
\begin{align*}
    &\liminf_{T \rightarrow \infty} \frac{1}{T} \sum_{i=0}^{\left\lfloor \frac{T-T_v}{T_v+T_p} \right\rfloor} \int_{i(T_v+T_p)}^{(i+1)T_v+iT_p} \min\left\{m, n - \int_{\max\{t-T_R, 0\}}^t \lambdaeff^{\pi}(s) ds\right\} dt \\
    \overset{(a)}{\leq}{}& \min\left\{\frac{mT_v}{T_v+T_p}, \frac{T_vn}{T_v+T_p} - \limsup_{T \rightarrow \infty} \frac{1}{T} \sum_{i=0}^{\left\lfloor \frac{T-T_v}{T_v+T_p} \right\rfloor} \int_{i(T_v+T_p)}^{(i+1)T_v+iT_p} \int_{\max\{t-T_R, 0\}}^t \lambdaeff^{\pi}(s) ds dt \right\} \\
    \overset{(b)}{=}{}& \min\left\{\frac{mT_v}{T_v+T_p}, \frac{T_vn}{T_v+T_p} - \limsup_{T \rightarrow \infty} \frac{1}{T} \sum_{i=0}^{\left\lfloor \frac{T-T_v}{T_v+T_p} \right\rfloor} \int_{\max\{i(T_v+T_p)-T_R, 0\}}^{(i+1)T_v+iT_p} \int_{\max\left\{s, i(T_v+T_p)\right\}}^{\min\left\{s+T_R, (i+1)T_v+iT_p\right\}} \lambdaeff^{\pi}(s) dt ds \right\} \\
    \overset{(c)}{=}{}& \min\left\{\frac{mT_v}{T_v+T_p}, \frac{T_vn}{T_v+T_p} - T_R\lambda \bar{\alpha}_v^{\pi}\frac{T_v}{T_v+T_p} - 
    \limsup_{T \rightarrow \infty} \frac{1}{T} \sum_{i=0}^{\left\lfloor \frac{T-T_v}{T_v+T_p} \right\rfloor} \int_{\max\{i(T_v+T_p)-T_R, 0\}}^{i(T_v+T_p)} \left(s+T_R-i(T_v+T_p)\right) \lambdaeff^{\pi}(s) ds \right. \\
    &\left.+ \liminf_{T \rightarrow \infty} \frac{1}{T} \sum_{i=0}^{\left\lfloor \frac{T-T_v}{T_v+T_p} \right\rfloor} \int_{(i+1)T_v+iT_p-T_R}^{(i+1)T_v+iT_p} \left(s+T_R-(i+1)T_v-iT_p\right) \lambdaeff^{\pi}(s) ds \right\} \\
    \leq{}& \min\left\{\frac{mT_v}{T_v+T_p}, \frac{T_vn}{T_v+T_p} - T_R\lambda \bar{\alpha}_v^\pi \frac{T_v}{T_v+T_p} + \frac{T_R^2\lambda}{T_v+T_p} \right\} = \frac{T_v}{T_v+T_p}\min\left\{m, n-T_R\lambda\bar{\alpha}_v^\pi+\frac{T_R^2\lambda}{T_v}\right\}
\end{align*}
where, $(a)$ follows by Jensen's inequality as $\min$ is a concave function. Next, $(b)$ follows by the interchange of integrals. Now, $(c)$ follows by splitting the outer integral into three parts - $[i(T_v+T_p)-T_R, i(T_v+T_p)]$, $[i(T_v+T_p), (i+1)T_v+iT_p-T_R]$, and $[(i+1)T_v+iT_p-T_R, (i+1)T_v+iT_p]$. Similarly, we can upper-bound the integral corresponding to the peaks to get
\begin{equation*}
    \ecc \leq \frac{T_v}{T_v+T_p}\min\left\{m, n-T_R\lambda\bar{\alpha}_v^\pi+\frac{T_R^2\lambda}{T_v}\right\} + \frac{T_p}{T_v+T_p}\min\left\{m, n-T_R a\lambda\bar{\alpha}_p^\pi+\frac{aT_R^2\lambda}{T_p}\right\}.
\end{equation*}
Now, by substituting the bounds on $\bar{q}$ and $\ecc$ in \eqref{eq:def_var-2}, we get \eqref{eq: stable_soc_bound}. This completes the proof.
\hfill $\square$
\endproof
\subsection{Additional Details for Case III, Theorem~\ref{thm: informal_varying_arrival_rate}: Limited Chargers or Limited Pack Size}\label{sec:lim_char_lim_pack_size}
\textbf{EV Scaling with Limited Chargers:}
For the Case III, after ignoring the edge effects, $n_{\operatorname{cp}, \eta} = 0$ if and only if $\eta = \frac{\alpha r \Tpeak}{\Tvalley}$. Thus, if the number of chargers is at most $m< n_{\operatorname{cv}, \frac{\alpha r \Tpeak}{\Tvalley}} = \alpha r \lambda_{\operatorname{avg}} T_R \left(1 + \frac{\Tpeak}{\Tvalley}\right)$, then, we have $n_{\operatorname{cp}, \eta} = \Omega(\lambda)$ implying that a non-zero fraction of partially charged EVs are available to dispatch to pick up an incoming customer at all times. Similar to \Cref{thm:lower-bound}, one can benefit from these partially charged EVs to reduce the pickup times. Thus, we conclude that EV-based scaling is applicable for large peak amplitudes whenever the number of chargers is less than $\alpha r \lambda_{\operatorname{avg}} T_R \left(1 + \frac{\Tpeak}{\Tvalley}\right)$. Now, for the rest of the section, we assume that the number of chargers is greater than this threshold.

\textbf{EV Scaling with Small Battery Pack Size:}
We obtain a necessary condition on the battery pack size $p^\star$ to ensure sustained operations without needing to charge during the peaks. Let $\Delta \operatorname{soc}_{\operatorname{valley}}$ be the average SoC gained during a demand valley. Thus, we can only budget to spend a total of $\Delta \operatorname{soc}_{\operatorname{valley}}$ SoC during the peak. Thus, $p^\star$ should be such that $\Delta \operatorname{soc}_{\operatorname{valley}} \geq r_d T_2/p^\star$. Note that $r_dT_2$ is the kWh expelled during the peak, and so $r_d T_2/p^\star$ is the SoC lost during the peak. Now, we calculate the $\Delta \operatorname{soc}_{\operatorname{valley}}$ to obtain a lower bound on $p^\star$. Note that $\lambda T_R$ EVs are driving and $n-\lambda T_R$ EVs are charging during the valley. Thus, $\Delta \operatorname{soc}_{\operatorname{valley}}$ is given by 
\begin{align*}
\Delta \operatorname{soc}_{\operatorname{valley}}= \min\left\{\frac{(n-\lambda T_R)r_c T_1 - \lambda T_R r_d T_1}{n p^\star}, 1\right\}. \numberthis \label{eq: soc_valley}
\end{align*}
Now substituting $\Delta \operatorname{soc}_{\operatorname{valley}} \geq r_d T_2/p^\star$ in the above equation, we get
\begin{align*}
    p^\star \geq r_d T_2, \quad n(T_1 - rT_2) \geq (1+r)\lambda T_R T_1. \numberthis \label{eq: nec_condition_pack_size}
\end{align*}
By substituting $n = n_{\operatorname{dp}}$, which corresponds to Case III of \Cref{thm: informal_varying_arrival_rate} with the largest value of $c$ and noting that $\alpha \geq \frac{\lambda}{\lambda_{\operatorname{avg}}(1-rT_2/T_1)}$, one can verify the second inequality above is satisfied. Thus, we simply have $p^\star \geq r_d T_2$. In other words, if $p^\star < r_d T_2$, then it is not possible (under any policy) to sustain driving during the peak with no charging. Thus, in this case, we will have a constant fraction of EVs charging at all times. So, we conclude that EV-based scaling similar to \Cref{thm:lower-bound} is applicable whenever $p^\star < r_d T_2$. In the rest of the section, we now fix $p^\star \geq r_d T_2$.
\subsection{Details of Closest Dispatch for Varying Arrival Rate} \label{app: cd_varying_arrival_rate}
We start this section by substantiating the discussion after \Cref{prop: cd_varying_arrival_rate} in \Cref{sec:limits-cd-varying}, we then provide the proposition's formal proof. Denote by $\bar{\alpha}_v$ and $\bar{\alpha}_p$ the fraction of the demand met during the valley and peak under Closest Dispatch. From the proposition, $\bar{\alpha}_v<1$ and $\bar{\alpha}_p$ is determined as in \eqref{eq: peak_and_valley_service_level} by substituting $\bar{\alpha}_v$ in place of $\alpha_v$. Note that $\bar{\alpha}_p>\alpha_p$. Now, at the beginning of a peak, there are at least $\bXC_0$ EVs that need to charge during the peak for $T^\prime $ hours so that they can only drive for $(\Tpeak-T^\prime)$ hours, where $r_cT^\prime = r_d(T_p-T^\prime)$. Hence, to serve $\bar{\alpha}_p$ fraction of the demand during the peak, the average number of available vehicles must satisfy $(n-\bXC_0)T^\prime/\Tpeak + n(\Tpeak-T^\prime)/\Tpeak\geq \bar{\alpha}_p \amplitude \lambda T_R$, or equivalently,   
\begin{align*}
    n \geq \bar{\alpha}_p a \lambda T_R + \frac{\bXC_0}{(1+r)} \quad \textit{with } \bar{\alpha}_p = \frac{\alpha\lambda_{\operatorname{avg}}(T_v+T_p)}{\amplitude\lambda T_p} - \frac{\bar{\alpha}_v T_v}{\amplitude T_p}  \numberthis. \label{eq: lower_bound_fleet_size_cd_var} 
\end{align*}
The fleet size requirement as given by \eqref{eq: lower_bound_fleet_size_cd_var} is strictly greater than the lower bound of \Cref{thm: informal_varying_arrival_rate}: first, $\bar{\alpha}_p - \alpha_p$ is the additional service level during the peak to compensate for lost demand during the valley ($\bar{\alpha}_v<1$); second, $\bXC_0/(1+r)$ is the additional fleet size to compensate for EVs with zero SoC at the start of the peak due to poor load balancing by Closest Dispatch.

\proof{Proof of \Cref{prop: cd_varying_arrival_rate}}
For ease of notation, we denote $\Nt^{\operatorname{CD}}$ simply by $\Nt$ as we work with the Closest Dispatch policy throughout the proof. Similar to the proof of \cref{prop:scalings-infty}, one can show that the Closest Dispatch ODE given by \eqref{eq:odes-cd} has a unique fixed point $(\bXC, \bXB)$. We omit this proof for brevity. Now, assuming global stability of the Closest Dispatch ODE, for all $\epsilon > 0$, there exists $t_0 > 0$ such that for all $t \geq t_0$, we have
\begin{align*}
    |\bXC_0 - \XC_0(t)| \leq \epsilon, |\bXC_{\Nt} - \XC_{\Nt}(t)| \leq \epsilon, \XB_{\Nt}(t) > 0
\end{align*}
Thus, we have
\begin{align*}
    \frac{1}{T}\int_0^T I(t)\left(1-\frac{\XC_0(t)}{\XC_{\Nt}(t)}\right) dt \leq \frac{t_0}{T} + \left(1 - \frac{\bXC_0-\epsilon}{\bXC_{\Nt}+\epsilon}\right) \frac{T-t_0}{T}.
\end{align*}
Now, by taking the limit supremum as $T \rightarrow \infty$ on both sides, we get
\begin{align*}
     \alphasup \leq  1 - \frac{\bXC_0-\epsilon}{\bXC_{\Nt}+\epsilon}.
\end{align*}
As the above is true for an arbitrary $\epsilon > 0$, we have
\begin{align*}
    \alphasup \leq 1 - \frac{\bXC_0}{\bXC_{\Nt}}. \numberthis \label{eq: alpha_eff_var_proof}
\end{align*}
Now, we analyze the fixed point $(\bXC, \bXB)$ of the Closest Dispatch ODE. First note that $I = 1$ for such a fixed point as otherwise \eqref{eq:cd-3} with $\dXC_{\Nt}(t)=0$ implies $\bXB_{\Nt} = 0$ which is a contradiction. Next, as $I=1$, we have $\bXB_{\Nt} \leq \lambda T_R + \kappa \lambda^{1-\gamma}$ for some $\kappa > 0$, $\gamma \in [1/3, 1/2]$. As $n \geq (1+r)\lambda T_R$, we have 
\begin{align*}
    \frac{1}{T_R} \geq \mu(t) = \frac{1}{T_R + \frac{\tau_1}{\sqrt{n-\XB_{\Nt}(t)}}} \geq \frac{1}{T_R + \frac{\tau_1}{\sqrt{r \lambda T_R - \kappa \lambda^{1-\gamma}}}} \overset{*}{\geq} \frac{1}{T_R + \frac{\tau_1\sqrt{2}}{\sqrt{rT_R \lambda}}} \overset{\Delta}{=} \frac{1}{T_R\left(1+\frac{C}{\sqrt{\lambda}}\right)}, \numberthis \label{eq: mu_bound_var}
\end{align*}
where $(*)$ follows for all $\lambda \geq \lambda_{00}$ for some $\lambda_{00} > 0$ as $rT_R\lambda > 2\kappa \lambda^{1-\gamma}$. In addition, the last equality follows by setting $C \overset{\Delta}{=} \frac{\sqrt{2}\tau_1}{T_R\sqrt{rT_R}}$. Now, by setting $\dXC_{\Nt}=0$ and $I = 1$ in \eqref{eq:cd-3} and using that $\mu(t) \leq 1/T_R$, we get $\bXB_{\Nt} \geq \alphasup\lambda T_R$. So, we get $\bXC_{\Nt} \leq n - \alphasup\lambda T_R$. Next, we set $\dXB_{j+1}(t) = 0$ in \eqref{eq:cd-1} and $\dXC_j(t) = 0$ in \eqref{eq:cd-2} and add the resultant equations to get
\begin{align*}
    \bXC_{j+1} - \bXC_j = \frac{\bXC_{\Nt}}{rT_R\lambda} \left(\bXC_j-\bXC_{j-1}\right) \leq \frac{n-\alphasup\lambda T_R}{rT_R\lambda}\left(\bXC_j-\bXC_{j-1}\right) \quad \forall j \in \{1, \hdots, \Nt-2\} \numberthis \label{eq: gs_cd_1}
\end{align*}
Similarly, setting $\dXB_{1}(t) = 0$ in \eqref{eq:cd-1} and $\dXC_0(t) = 0$ in \eqref{eq:cd-4} and adding the resultant equations, we get
\begin{align*}
    \bXC_1 - \bXC_0 = \frac{\bXC_{\Nt}}{rT_R\lambda} \bXC_0 \leq \frac{n-\alphasup\lambda T_R}{rT_R\lambda}\bXC_0 \numberthis \label{eq: gs_cd_2}
\end{align*}
Now, by using \eqref{eq: gs_cd_1} multiple times along with \eqref{eq: gs_cd_2}, we get
\begin{align*}
    \bXC_{j+1} - \bXC_j \leq \left(\frac{n-\alphasup\lambda T_R}{rT_R\lambda}\right)^{j+1}\XC_0 \quad \forall j \in \{0, \hdots, \Nt-2\}. \numberthis \label{eq: upto_nt_minus_two}
\end{align*}
Lastly, setting $\dXC_{\Nt}(t) = 0$ in \eqref{eq:cd-3} and $\dXC_{\Nt-1}(t)=0$ in \eqref{eq:cd-2} and subtracting the resultant equations, we get
\begin{align*}
    \bXC_{\Nt} - \bXC_{\Nt-1} = \frac{\bXC_{\Nt}}{rT_R\lambda} \left(\bXC_{\Nt-1} - \bXC_{\Nt-2}\right) \leq \frac{n-\alphasup\lambda T_R}{rT_R\lambda}\left(\bXC_{\Nt-1} - \bXC_{\Nt-2}\right).
\end{align*}
Combining the above equation with \eqref{eq: upto_nt_minus_two}, we get
\begin{align*}
    \bXC_{j+1} - \bXC_j \leq \left(\frac{n-\alphasup\lambda T_R}{rT_R\lambda}\right)^{j+1}\XC_0 \quad \forall j \in \{0, \hdots, \Nt-1\}. \numberthis \label{eq: upto_nt_minus_one}
\end{align*}
Now, we have
\begin{align*}
    \bXC_{\Nt} - \bXC_0 =  \sum_{j=0}^{\Nt-1} \bXC_{j+1} - \bXC_j &\overset{\eqref{eq: upto_nt_minus_one}}{\leq} \bXC_0 \sum_{j=0}^{\Nt-1} \left(\frac{n-\alphasup\lambda T_R}{rT_R\lambda}\right)^{j+1} 
\end{align*}
Now, by simplifying the above equation, we get
\begin{align*}
    \bXC_{\Nt} \leq \bXC_0 \sum_{j=0}^{\Nt} \left(\frac{n-\alphasup\lambda T_R}{rT_R\lambda}\right)^{j} = \bXC_0 \frac{\left(\frac{n-\alphasup\lambda T_R}{rT_R\lambda}\right)^{\Nt+1} -1}{\left(\frac{n-\alphasup\lambda T_R}{rT_R\lambda}\right)-1} \numberthis \label{eq:lower_bound_x0}
\end{align*}
The above equation now implies the following bound on $\alphasup$. We have
\begin{align*}
    \alphasup \leq 1 - \frac{\bXC_0}{\bXC_{\Nt}} \leq 1 - \frac{\left(\frac{n-\alphasup\lambda T_R}{rT_R\lambda}\right) -1}{\left(\frac{n-\alphasup\lambda T_R}{rT_R\lambda}\right)^{\Nt+1}-1}. \numberthis \label{eq: alpha_eff_fixed_point_proof}
\end{align*}

Now, consider the following fixed-point equation:
\begin{align*}
    \bar{\alpha}  = 1 - \frac{\frac{n-\bar{\alpha}\lambda T_R}{rT_R\lambda}-1}{\left(\frac{n-\bar{\alpha}\lambda T_R}{rT_R\lambda}\right)^{\Nt+1}-1}. \numberthis \label{eq: fixed_point_alpha_proof}
\end{align*}
First, observe that the RHS is a continuous and strictly decreasing function of $\bar{\alpha}$. When $\bar{\alpha} = 0$, the RHS is positive as 
\begin{align*}
    \frac{\frac{n}{rT_R\lambda}-1}{\left(\frac{n}{rT_R\lambda}\right)^{\Nt+1}-1} < 1 \quad \forall n > 0.
\end{align*}
In addition, for all $\bar{\alpha} \in [0, 1]$, we have
\begin{align*}
    \sup_{\lambda > 0} \frac{\frac{n-\bar{\alpha}\lambda T_R}{rT_R\lambda}-1}{\left(\frac{n-\bar{\alpha}\lambda T_R}{rT_R\lambda}\right)^{\Nt+1}-1} \geq \sup_{\lambda > 0} \frac{\frac{n}{rT_R\lambda}-1}{\left(\frac{n}{rT_R\lambda}\right)^{\Nt+1}-1} \geq \frac{\frac{M}{r}-1}{\left(\frac{M}{r}\right)^{\Nt+1}-1} > 0 \numberthis \label{eq: lower_bound_geo_sum}
\end{align*}
as $n \leq M\lambda T_R$ for some fixed constant $M > 1$. Thus, by the intermediate value theorem, there exists a unique $\bar{\alpha} \in \left[0, 1-\frac{\frac{M}{r}-1}{\left(\frac{M}{r}\right)^{\Nt+1}-1}\right]$ such that \eqref{eq: fixed_point_alpha_proof} is satisfied. As $1-\frac{\frac{M-1}{r}-1}{\left(\frac{M-1}{r}\right)^{\Nt+1}-1} < 1$, we have $\sup_{\lambda > \lambda_{00}} \bar{\alpha} < 1$. Also, by \eqref{eq: alpha_eff_fixed_point_proof}, we conclude that $\alphasup \leq \bar{\alpha}$. Also, as discussed in Section~\ref{sec:limits-cd-varying}, note that $\bar{\alpha}$ approaches 1 as $\Nt$ increases and/or $M$ increases.

Now, to complete the proof, note that \eqref{eq:lower_bound_x0} implies 
\begin{align*}
    \bXC_0 \geq \bXC_{\Nt} \frac{\left(\frac{n-\alphasup\lambda T_R}{rT_R\lambda}\right)-1}{\left(\frac{n-\alphasup\lambda T_R}{rT_R\lambda}\right)^{\Nt+1} -1} \overset{\eqref{eq: lower_bound_geo_sum}}{\geq} \bXC_{\Nt} \frac{\frac{M}{r}-1}{\left(\frac{M}{r}\right)^{\Nt+1}-1}.
\end{align*}
Lastly, as $I=1$, we have $\bXB_{\Nt} \leq \lambda T_R + \kappa \lambda^{1-\gamma}$. Now, as $n \geq (1+r)\lambda T_R$ and $\gamma \in [1/3, 1/2]$, there exists $\lambda_{01} > 0$ such that for all $\lambda \geq \lambda_{01}$, we get $\bXC_{\Nt} = n - \bXB_{\Nt} \geq r \lambda T_R / 2$.
This completes the proof by setting $\lambda_0 \overset{\Delta}{=} \max\{\lambda_{00}, \lambda_{01}\}$.
\hfill $\square$
\endproof
\subsection{Power-of-\texorpdfstring{$d$}{d} for Varying Arrival Rate} \label{app: pod_var}
\textbf{Justification for \eqref{eq: fleet_pod_varying}:} We simplify $\bar{\alpha}_p a \lambda T_R(1-rT_p/T_v)$ below. We have
\begin{align*}
    &\bar{\alpha}_p a \lambda T_R \left(1-\frac{rT_p}{T_v}\right) \\
    ={}& \left(\frac{\alpha \lambda_{\operatorname{avg}}(T_v+T_p)}{a\lambda T_p}-\frac{\bar{\alpha}_vT_v}{aT_p}\right)a \lambda T_R \left(1-\frac{rT_p}{T_v}\right) \\
    ={}& \alpha \lambda_{\operatorname{avg}} T_R \left(1+\frac{T_v}{T_p}\right)\left(1-\frac{rT_p}{T_v}\right) - \frac{\bar{\alpha}_v\lambda T_R T_v}{T_p} + \bar{\alpha}_v r \lambda T_R \\
    ={}&\alpha \lambda_{\operatorname{avg}} T_R \left(1+\frac{T_v}{T_p}\right)\left(1-\frac{rT_p}{T_v}\right) - \bar{\alpha}_v\lambda T_R\left(1+\frac{T_v}{T_p}\right) + (1+r)\bar{\alpha}_v \lambda T_R \\
    ={}& (1+r)\bar{\alpha}_v\lambda T_R + (1-\bar{\alpha}_v)\left(1+\frac{T_v}{T_p}\right) \lambda T_R + \lambda_{\operatorname{avg}}T_R\left(1+\frac{T_v}{T_p}\right)\left(1-\frac{rT_p}{T_v}\right)\left(\alpha - \frac{\lambda}{\lambda_{\operatorname{avg}}\left(1-\frac{rT_p}{T_v}\right)}\right) \\
    ={}& (1+r)\bar{\alpha}_v\lambda T_R + (1-\bar{\alpha}_v)\left(1+\frac{T_v}{T_p}\right) \lambda T_R +  \lambda T_R\left(a+\frac{T_v}{T_p}\right)\left(1-\frac{rT_p}{T_v}\right)\left(\alpha - \frac{\lambda}{\lambda_{\operatorname{avg}}\left(1-\frac{rT_p}{T_v}\right)}\right) \\
    \overset{\Delta}{=}{}& (1+r)\bar{\alpha}_v\lambda T_R + (1-\bar{\alpha}_v)\left(1+\frac{T_v}{T_p}\right) \lambda T_R + b r\lambda T_R,
\end{align*}
where we define
\begin{align*}
    b \overset{\Delta}{=} \frac{1}{r}\left(a+\frac{T_v}{T_p}\right)\left(1-\frac{rT_p}{T_v}\right)\left(\alpha - \frac{\lambda}{\lambda_{\operatorname{avg}}\left(1-\frac{rT_p}{T_v}\right)}\right).
\end{align*}
Note that $b>0$ whenever we are in Case III of \Cref{thm: informal_varying_arrival_rate}. 

\proof{Proof of \Cref{prop: pod_var}}
\revcolor{We consider $m = \infty$ and admission control such that $A = \infty$ and the maximum number of busy vehicles is $(1+r)\lambda T_R + rb \lambda T_R/2$.}
Consider the system of differential equations in \eqref{eq:odes-power-d}. Any equilibrium $(\XC,\XB)$ of this dynamical system must satisfy \eqref{eq:odes-power-d-eq}. Now, we show that there exists a solution of \eqref{eq:odes-power-d-eq} and that all solutions of \eqref{eq:odes-power-d-eq} satisfy the hypothesis of the proposition. We start by quantifying some relations of $(\XC, \XB)$ that must be satisfied for any such fixed point of \eqref{eq:odes-power-d-eq}.

If $I = 0$ then $\XB_{\Nt} > T_R \lambda + rb\lambda T_R/2$ and from \eqref{eq:powerd-3-eq}, we would get that $\XB_{\Nt} = 0$, a contradiction. Thus, we must have $I = 1$. Next, as $m = \infty$ and $A=\infty$, we have $\tXC_j = \XC_j$ for all $j \in \{1, \hdots, \Nt\}$. Now, combining \eqref{eq:powerd-1-eq} and \eqref{eq:powerd-2-eq}, we deduce that 
$$
\XC_j-\XC_{j-1} = r T_B \mu\left(\XB_{j+1}-\XB_j\right) \quad \forall j\in\{1,\dots,\Nt-1\}.
$$
By carrying out telescopic sum from $j$ to $1$, we get
$$
\XC_j-\XC_{0} = rT_B \mu \left(\XB_{j+1}-\XB_1\right) \quad \forall j\in\{1,\dots,\Nt-1\}.
$$
Using \eqref{eq:powerd-4-eq} we obtain that $\XC_{j}=r T_B \mu \XB_{j+1}$ for $j\in\{0, 1,\dots,\Nt-1\}$. Substituting this back into \eqref{eq:odes-power-d-eq}, any $(\XC, \XB)$ is a fixed point of \eqref{eq:odes-power-d-eq} if and only if $\XC$ satisfies the following: 
\begin{subequations} \label{eq:eq-XC-a-new}
\begin{flalign}\label{eq:eq-XC-a-1-new}
\XC_{j-1}=rT_B\lambda \left(p_0-p_j\right) , \: 1\le j< \Nt,\\\label{eq:eq-XC-a-2-new}
\XC_{\Nt-1}=rT_B\lambda p_0, \\\label{eq:eq-XC-a-3-new}
\XC_{\Nt} + \frac{\XC_{\Nt-1}}{rT_B \mu}=n, \\
0\leq \XC_0 \leq \XC_1 \leq \hdots \leq \XC_{\Nt-1} \leq \XC_{\Nt}, \label{eq:eq-XC-a-4-new} \\
X^B_{\Nt}=(n-X^C_{\Nt}) \leq (1+r)\lambda T_R + \frac{rb \lambda T_R}{2}. \label{eq:eq-XC-a-5-new}
\end{flalign}
\end{subequations} 
Now, we verify that there exists $\lambda_{v} > 0$ such for all $\lambda \geq \lambda_{v}$ the following is a solution of \eqref{eq:eq-XC-a-new}: $\XC_j = 0$ for all $j \in \{1, \hdots, \Nt-3\}$, $\XC_{\Nt-2} = rT_B\lambda \left(\frac{rT_B\lambda}{n-\lambda/\mu}\right)^d$, $\XC_{\Nt-1} = rT_B\lambda$, and $\XC_{\Nt} = n - \lambda/\mu$. We start by simplifying the expression of $T_B$ below. We have
\begin{align*}
    T_B = T_R + \frac{\tau_1\sqrt{d}}{\sqrt{n-(1+r)\lambda T_R - rb \lambda T_R/2}} = T_R + \frac{\tau_1\sqrt{2} \log \lambda}{\sqrt{rbT_R \lambda}}.
\end{align*}
Then we have
\begin{align*}
    \XC_{\Nt-1} = rT_B \lambda = rT_R\lambda \left(1 + \frac{\tau_1\sqrt{2}\log \lambda}{T_R\sqrt{rbT_R \lambda}}\right) \leq r T_R \lambda + \frac{r^2 b T_R \lambda}{2}, \numberthis \label{eq: upper_bound_nt_minus_1}
\end{align*}
where the last inequality follows for all $\lambda \geq \lambda_{v0}$ for some $\lambda_{v0} > 0$ as $rbT_R > 0$. Now, we show that \eqref{eq:eq-XC-a-5-new} is satisfied by lower bounding the choice of $\XC_{\Nt}$. We have
\begin{align*}
    \XC_{\Nt} = n - \frac{\lambda}{\mu} \overset{*}{\geq} n - \lambda T_B \overset{\eqref{eq: upper_bound_nt_minus_1}}{\geq} r T_R \lambda + \frac{rbT_R}{2}\lambda, \numberthis \label{eq: lower_bound_nt}
\end{align*}
where $(*)$ follows as $\mu \geq 1/T_B$ by definition. Thus, \eqref{eq:eq-XC-a-5-new} is satisfied. Next, we upper bound our choice of $\XC_{\Nt-2}$ below. We have
\begin{align*}
    \XC_{\Nt-2} = r T_B \lambda \left(\frac{rT_B\lambda}{n-\lambda/\mu}\right)^d \overset{*}{\leq} 2r T_R \lambda \left(\frac{1+rb/2}{1+b/2}\right)^{(\log \lambda)^2} \overset{**}{\leq} \frac{r}{2}, \numberthis \label{eq: upper_bound_nt_minus_two}
\end{align*}
where $(*)$ follows by upper bounding $rT_B \lambda$ by \eqref{eq: upper_bound_nt_minus_1} and lower bounding $n-\lambda/\mu$ by \eqref{eq: lower_bound_nt}. In addition, there exists $\lambda_{v1} > 0$ such that for all $\lambda \geq \lambda_{v1}$, we have $T_B \leq 2T_R$. Next, $(**)$ follows by first noting that the base of $(\log \lambda)^2$ is strictly less than 1 as $r<1$ and $b>0$. Thus, $\left(\frac{1+rb/2}{1+b/2}\right)^{(\log \lambda)^2}$ is of the order $\lambda^{-\log \lambda}$ which implies that $\lambda\left(\frac{1+rb/2}{1+b/2}\right)^{(\log \lambda)^2} \rightarrow 0$. Thus, there exists $\lambda_{v2} > 0$ such that for all $\lambda \geq \lambda_{v2}$, $(**)$ is satisfied. Thus, by \eqref{eq: lower_bound_nt}, \eqref{eq: upper_bound_nt_minus_1}, and \eqref{eq: upper_bound_nt_minus_two}, we conclude that \eqref{eq:eq-XC-a-4-new} is satisfied. Next, observe that \eqref{eq:eq-XC-a-3-new} is satisfied for the choice of $(\XC_{\Nt}, \XC_{\Nt-1})$. Next, we have $p_0=1$ as $\XC_0 = 0$, and so, \eqref{eq:eq-XC-a-2-new} is satisfied. Now, we complete the proof by verifying \eqref{eq:eq-XC-a-1-new} below. First note that $p_j = 1$ for all $j \in \{1, \hdots, \Nt-3\}$ as $\XC_j=0$ for all $j \in \{1, \hdots, \Nt-3\}$. In addition, as $\XC_{\Nt-2} < 1$, we have $p_{\Nt-2}=1$ by definition. Now, as $p_0=p_j=1$ for all $j \in \{1, \hdots, \Nt-2\}$, it is easy to verify that \eqref{eq:eq-XC-a-1-new} is satisfied for all $j \in \{1, \hdots, \Nt-2\}$. Lastly, for $j = \Nt-1$, \eqref{eq:eq-XC-a-1-new} is satisfied by the choice of $(\XC_{\Nt-2}, \XC_{\Nt-1}, \XC_{\Nt})$. This proves the existence of a fixed point. Now, we prove the required bounds on any fixed point of \eqref{eq:odes-power-d-eq}. By \eqref{eq:eq-XC-a-2-new}, we have $\XC_{\Nt-1} \leq rT_B\lambda$ and so by \eqref{eq:eq-XC-a-3-new}, we have $\XC_{\Nt} \geq n - \lambda/\mu$. Thus, by \eqref{eq:eq-XC-a-1-new} for $j = \Nt-1$ and substituting the bounds on $\XC_{\Nt}$ and $\XC_{\Nt-1}$, we obtain that $\XC_{\Nt-1} \leq r/2$, similar to \eqref{eq: upper_bound_nt_minus_two}. Now, as $\XB_{\Nt-1} = \XC_{\Nt-2}/(rT_B\mu)$, we have $\XB_{\Nt-2} \leq \XB_{\Nt-1} = \XC_{\Nt-2}/(rT_B\mu) \leq \XC_{\Nt-2}/r \leq 1/2$. Thus, we have $\XC_{\Nt-2} + \XB_{\Nt-2} \leq (r+1)/2 \leq 1$ as $r<1$. Lastly, we have $\XC_0 \leq \XC_{\Nt-2} \leq r/2$ and so $p_0 = 1$ by definition. Also, we already showed that $I = 1$ and so, $p_0 I = 1$. This completes the proof by setting $\lambda_v \overset{\Delta}{=} \max\{\lambda_{v0}, \lambda_{v1}, \lambda_{v2}\}$. \hfill $\square$
\endproof
\section{Auxiliary Technical Results}\label{app:aux-tech-results}
 \begin{lemma}[Eq. 12 in \cite{srinivasa2009distance}]\label{lem:d-close}
Let $Z_1,Z_2,\dots,Z_k$ be a sequence of independent uniformly distributed random points in a 2-dimensional Euclidean ball $\mathcal{C}$. Let $Z^{(d)}(z_0)$ be the $d^{\text{th}}$ closest point in the sequence to a point $z_0$ in the interior of $\mathcal{C}$ then 
$$
\mathbb{E}\left[\|Z^{(d)}(z_0)-z_0\|\right] = \Theta\left(\sqrt{\frac{d}{k}}\right), \quad \text{as } k \uparrow \infty.
$$
\end{lemma}
\section{Additional Set-up Description and Results for Section~\ref{sim:synthetic_data}}
\subsection{Simulation Detailed Description}
We consider a $[0, 10 \text{ mi}] \times [0, 10 \text{ mi}]$ spatial region and generate random arrivals. The starting and ending coordinates of a trip request on each axis are sampled from a Uniform distribution on $[0, 10]$. The inter-arrival time between two arrivals is exponentially distributed with rate $\lambda$. In our simulations, we consider $\lambda \in \left\{5, 10, 20, 40, 80, 160, 320\right\}$ arrivals per min. We run our simulator for 1000 minutes (about 16 hours) and calculate parameters like service level by considering only the second half of the simulation to ensure that the system is in the steady state. We run each simulation instance for five randomly generated data sets, and the mean value of the output is reported.

For each EV in the system, a uniformly distributed initial location in $[0, 10 \text{ mi}] \times [0, 10 \text{ mi}]$ is sampled independent of all the other EVs. The SoC is initialized to be a uniform random variable on $[0.4, 0.6]$. Also, the initial state of the EVs is set to be Idle. We set the location of chargers uniformly at random on $[0, 10 \text{ mi}] \times [0, 10 \text{ mi}]$. Each charger has $m_p$ number of posts, i.e., charging ports, and there are $m_l = \lfloor m / m_p \rfloor$ different locations of chargers. We consider $m_p= 8$  \citep{no_of_posts} for simulations in this section.

The system's evolution will be governed by the parameters in \Cref{tab:my_label} unless otherwise stated \citep{bauer2018cost}.
\begin{table}[hbt!]
    \TABLE{Parameters for the Simulations.
    \label{tab:my_label}}{
    \centering
    \begin{tabular}{|l|c|} \hline
         \multicolumn{1}{|c|}{Parameter} & Value  \\ \hline
       Charge Rate $(r_c)$ (kilo-Watt)  & 20 kW  \\ \hline
        Discharge Rate $(r_d)$ (kilo-Watt) & 5 kW \\ \hline
        Average Velocity $(v_{\operatorname{avg}})$ (miles per hour) & 20 mi/hr \\ \hline
        Battery Pack Size $(p^\star)$ (kilo-Watt-hour) & 40 kWh  \\ \hline
        Target Service Level $(\alpha)$ & 0.9 \\ \hline
        Minimum SoC after finishing a trip ($s_{\min}$) & 0.2 \\ \hline
        Maximum SoC to dispatch an EV to a charger ($s_{\max}$) & 0.9 \\ \hline
    \end{tabular}}{}
\end{table}
For each trip request, the matching algorithm picks an EV to serve the customer. If the selected EV has sufficient battery to serve the trip and an additional buffer SoC of $s_{\min}$, then the customer is served, otherwise the customer leaves the system immediately. The buffer SoC of $s_{\min}$ allows the EV to be able to drive to a charger after finishing the trip. If the customer is served, then the EV drives from its current location to the customer's origin and then to the customer's destination. Once the trip ends, the EV's location is set to the customer's destination, and the SoC is reduced by $\Delta\text{SoC} = r_d \times \left(T_P + T_R\right) / p^\star$ [kWh], where $T_P$ and  $T_R$ are computed using $v_{\operatorname{avg}}$ and the corresponding euclidean distances. 
Next, the EV becomes idle, and if its SoC is less than $s_{\max}$,  the EV is dispatched to the closest available charger. In this case, the EV expels $r_d \times T_C / p^\star$ amount of SoC to reach the charger, where $T_C$ is computed using $v_{\operatorname{avg}}$ and the corresponding Euclidean distance.
While an EV is driving to a charger, another EV can reach the same charger and start charging earlier. Thus, an EV may reach an unavailable charger, in which case, the EV waits at the charger until a post becomes available. If a post is available at the charger, the EV starts charging immediately. In either case,  we allow for interrupting charging and dispatch an EV to pick up a customer. If the charging is not interrupted and the SoC reaches 1, the EV becomes idle. These dynamics continue for all customers and vehicles until the end of the simulation.

\subsection{Details on Asymptotic Simulations} \label{sec:app-sim-details}
We devise the following procedure to verify the theoretical scalings as in \Cref{thm:informal_upper-bound}. First, we fix a sequence of arrival rates $\lambda \in \{5, 10, 20, 40, 80, 160, 320\}$. Then, we fix a sequence of $m$ corresponding to a given value of $\beta \in \{0.7, 0.8, 0.9, 1\}$. In particular, we set 
\begin{equation} \label{eq: m}
    m = r\tilde{T}_R \alpha \lambda + c \left(\tilde{T}_R \lambda\right)^\beta
\end{equation}
where $\alpha = 0.9$ is the target service level, $r = 0.25$ is the ratio of discharge rate and charge rate, $c = 4$ is a constant, and $\tilde{T}_R$ is the average fulfilled trip time. Note that $\tilde{T}_R$ is not known a priori as it depends on $\alphasup$, which in turn depends on system parameters like $\lambda, n$, and $m$. In particular, Po$d$ is biased slightly towards shorter trips as longer trips are more likely to be rejected. Due to this bias, we get $\tilde{T}_R < T_R$. We observed empirically that $\tilde{T}_R \approx [15, 15.3]$ mins and $T_R \approx [15.6, 15.7]$ mins. Based on this empirical observation, we set $\tilde{T}_R = 15.14$, which allows us to approximate $m$ for a given $\beta$.

We make three comments on \eqref{eq: m}. First, we use $\tilde{T}_R$ as opposed to $T_R$ in \Cref{thm:informal_upper-bound} because the customers that are rejected do not affect the dynamics of the system at all. Hence, the trip time experienced by the EVs corresponds to $\tilde{T}_R$. Second, we set the second order term as $c\left(\tilde{T}_R \lambda\right)^\beta$ as opposed to $c\lambda^\beta$ in \Cref{thm:informal_upper-bound} mainly because the former is dimensionless and has a better interpretation. In particular, $\tilde{T}_R \lambda$ corresponds to the workload, which determines the number of EVs driving with a customer, in turn determining the minimum number of EVs required to charge to sustain the fleet SoC. The number of EVs required to charge then determines the charger requirements. Lastly, we set $c=4$ for two reasons. First, to ensure that the resultant average drive to the charger time is comparable to the average pickup time. Second, we set $c$ to be large enough to encourage achieving the asymptotic scaling as in \Cref{thm:informal_upper-bound}. For small values of $c$, we observe that the EVs have to wait at the chargers for small values of $\lambda$. In particular, the decision to select the closest available charger is blind to the EVs that are already driving to a charger. Thus, an EV dispatched to an available charger may have to wait at the charger before it starts charging. This effect is exacerbated for small values of $\lambda$ as the drive to the charger time is large. Due to EVs waiting at the chargers, the simulation starts departing from the theory, which warrants larger values of $\lambda$ to verify the theoretical scalings. 

After fixing the $(\lambda, m)$ sequence, we simulate various values of the fleet size to infer the fleet size corresponding to the 90\% service level. In particular, for each tuple $(\lambda, m)$, we simulate five different values of the fleet size. Each simulation is repeated five times with different seeds, and the average values are considered. As shown in Fig.~\ref{fig: fleet_size_90_percent}, we carry out linear regression between service level and fleet size buffer $(n - (1+r)\tilde{T}_R \alpha \lambda)$ to infer the fleet size and the fleet size buffer that corresponds to the 90\% service level. 
\begin{figure}[t]
\centering
    \FIGURE{
    \resizebox{\textwidth}{!}{
    \begin{tikzpicture}
        \node at (0, 0) {\input{figures_final/calculating_90_percent_fleet_size.pgf}};
        \node at (15, 0) {\input{figures_final/calculating_90_percent_fleet_size_beta_0.8.pgf}};
    \end{tikzpicture}}}{
    Fleet size buffer corresponding to 90\% service level. The markers correspond to simulated values, and the line is a linear regression fit. The plot corresponds to $\beta = 1$ (left) and $\beta = 0.8$ (right). 
    \label{fig: fleet_size_90_percent}}{}
\end{figure}
\begin{table}[hbt!]
    \centering
\TABLE{\revcolor{Fleet Size and Number of Chargers corresponding to the 90\% Service Level.} \label{tab: n_m_90_percent_data}}{
        \begin{tabular}{|c||c|c||c|c||c|c||c|c|} \hline
       $\lambda$  & \multicolumn{2}{c||}{Series A} & \multicolumn{2}{c||}{Series B}  & \multicolumn{2}{c||}{Series C}  & \multicolumn{2}{c|}{Series D} \\ \hline
         & $n$ & $m$ & $n$ & $m$ & $n$ & $m$ & $n$ & $m$ \\ \hline
    5 & 126 & 320 & 130 & 208 & 133 & 144 & 143 & 96 \\ \hline
    10 & 229 & 640 & 236 & 400 & 247 & 256 & 258 & 168 \\ \hline
    20 & 427 & 1280 & 437 & 752 & 451 & 456 & 472 & 288 \\ \hline
    40 & 806 & 2560 & 825 & 1416 & 851 & 808 & 889 & 488  \\ \hline
    80 & 1532 & 5120 & 1559 &  2656 & 1602 & 1448 & 1654 & 856 \\ \hline
    160 & 2958 & 10248 & 3003 & 5008 & 3072 & 2600 & 3188 & 1496 \\ \hline
    320 & 5769 & 20504 & 5841 & 9416 & 5956 & 4672 &  6123 & 2640 \\ \hline
    \end{tabular}}{}
\end{table}
For completeness, we report the fleet size and number of chargers that correspond to 90\% service level for different values of $\beta$ in Table~\ref{tab: n_m_90_percent_data}. The empirical $1-\gamma$ correspond to the slope of the linear regression line between $\log(\lambda)$ and $\log(n - (1+r)\tilde{T}_R \alpha \lambda)$. Similarly, the empirical $\beta$ corresponds to the slope of the linear regression line between $\log(\lambda)$ and $\log(m - r\tilde{T}_R \alpha \lambda)$. Note that, we recalculate the value of $\beta$ as $\tilde{T}_R$ is now known exactly. This recalculation slightly changes the value of $\beta$. For example, Series B now corresponds to $\beta = 0.906$ rather than the originally used $\beta=0.9$. We present the results in \Cref{tab: asymptotic_sim} in the main body of the paper and provide a graphical illustration in \Cref{fig: asymptotic_sim}.
\begin{figure}
\centering
\scalebox{0.8}{
\FIGURE{\begin{tikzpicture}    
\node at (0, 0) {\scalebox{0.5}{\includegraphics{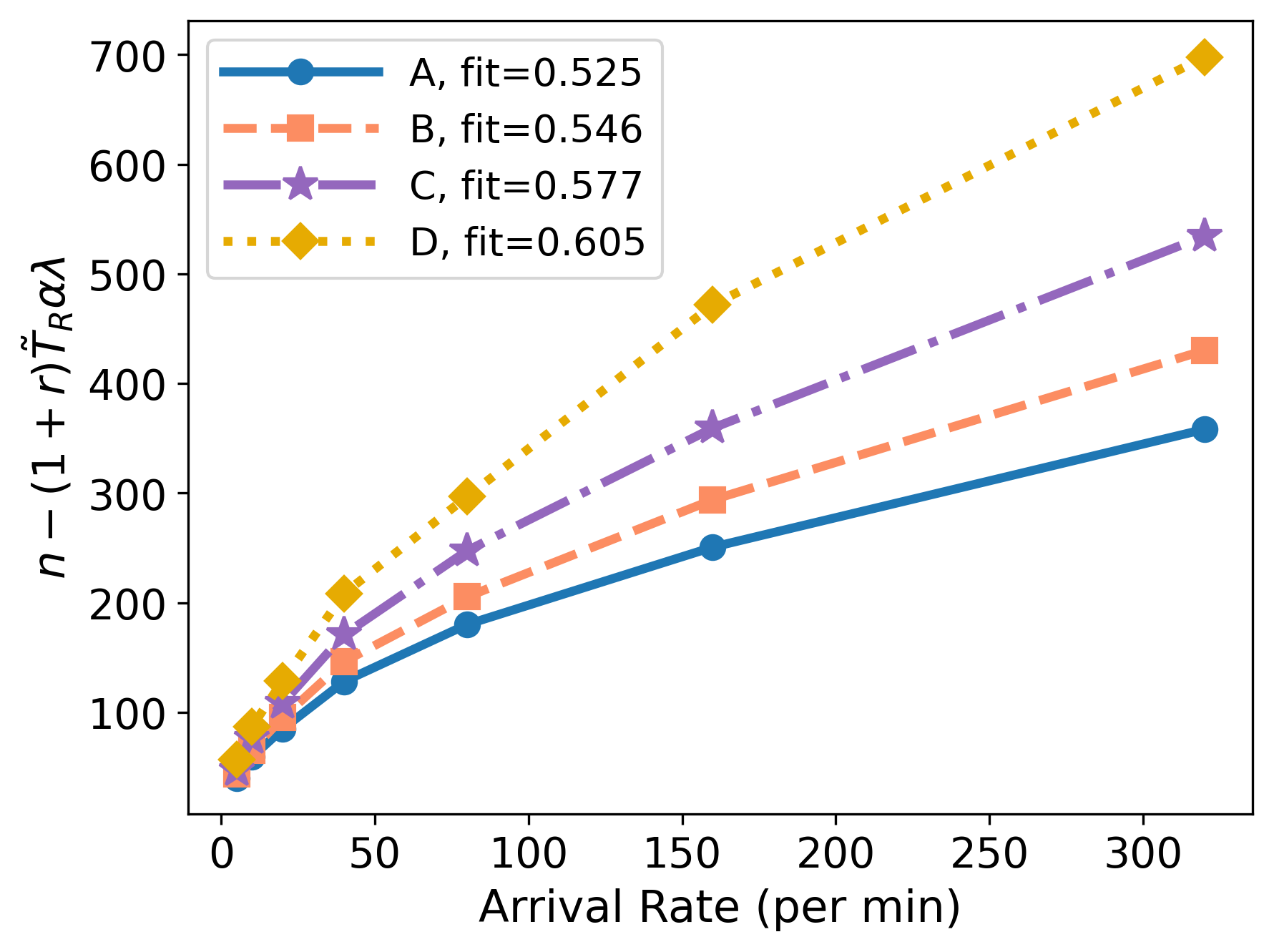}}};
\node at (8, 0) {\scalebox{0.5}{\includegraphics{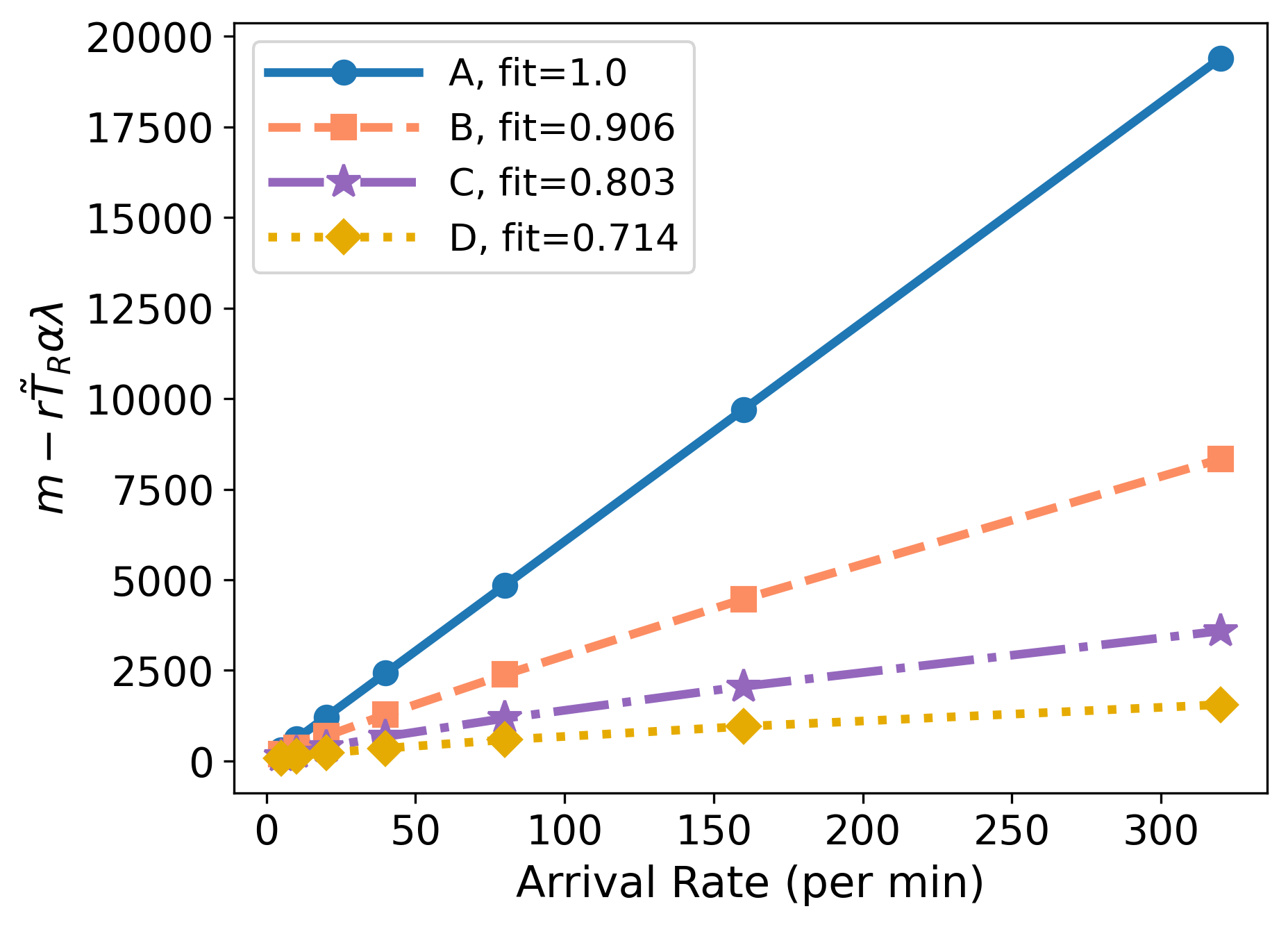}}};
\node at (0, -6.25) {\scalebox{0.5}{\includegraphics{figures_final/asymptotic_pickup_min.png}}};
\node at (8, -6.25) {\scalebox{0.5}{\includegraphics{figures_final/asymptotic_drive_to_charger_min.png}}};
\end{tikzpicture}}{
    \revcolor{Fleet size buffer (top left), number of chargers buffer (top right), average pickup time (bottom left), and average drive to charger time (bottom right) for Po2 as a function of arrival rate. Fit corresponds to the slope of the linear regression line of the log of $x$ and $y$ variables.}
    \label{fig: asymptotic_sim}}{}
    }
\end{figure}

Finally, in order to compute the theoretical values in Table~\ref{tab: asymptotic_sim}, we fix the value of $\beta$ and report the corresponding values of $1-\gamma$, pickup time scaling, and drive to charger time scaling as prescribed by \Cref{thm:informal_upper-bound}. In particular,
$$
\gamma = \min\left\{\frac{1}{2+1/\Nt}, \frac{\beta}{2}\right\} \implies 1-\gamma = \max\left\{0.508, 1-\frac{\beta}{2}\right\},\footnote{$\Nt = p^\star / (r_d \times T_B) \approx p^\star / (r_d \times T_R)$  = 32.}
$$
and that the pickup time and the drive to the charger time are approximately
$$
T_P \approx \frac{\tau_1\sqrt{d(1+r)}}{\sqrt{rn}} \sim n^{\gamma/(2\Nt) - 1/2} = n^{\gamma/64 - 1/2}, \quad \text{and}\quad
T_C \approx \frac{\tau_2}{\sqrt{m-A}} \sim n^{-\beta/2}.
$$

\subsection{Matching Algorithms Comparison} \label{app: matching_algo_comparison}
We compare the performance of three different matching algorithms: Closest Dispatch, Power-of-$d$ with $d=2$, and closest available dispatch (CAD). In CAD, we match a trip request to its closest available EV, i.e., the closest EV with sufficient battery to serve the trip. If no such EV exists, then the customer is dropped. That is, CAD is a natural extension to the Closest Dispatch, which keeps searching for a vehicle with enough SoC if the closest vehicle does not have enough battery to serve a request.

These matching algorithms are inspired by the long line of literature on routing in load balancing. In particular, one can interpret our model as a load balancing model by interpreting $\Nt-\text{SoC}_i(t)$ as queue lengths and dispatching an EV as routing an incoming customer, where $\text{SoC}_i(t)$ is the SoC of $i^{\text{th}}$ EV at time $t$. Then, CAD corresponds to joining an available queue at random, CD corresponds to random routing, and Po$d$ corresponds to Power-of-$d$ choices routing.

We start by qualitatively comparing the three algorithms in \Cref{fig: time_series_3_algos}.
\begin{figure}[thb!]
\def\sh{0}
\def\shb{0}
\def\shc{0}
    \centering
    \FIGURE{
    \begin{tikzpicture}[scale=1]
        \node[scale=0.35] at (0, -1.4) {\input{figures_final/CAD_stackplot.pgf}};
        \node[scale=0.35] at (5.6, -1.4) {\input{figures_final/CD_stackplot.pgf}};
        \node[scale=0.35] at (11.2, -1.4) {\input{figures_final/Po2_stackplot.pgf}};
          \draw[rectangle] (-1.2, -3.8) -- (12.5, -3.8) -- (12.5, -6) -- (-1.2, -6) -- (-1.2, -3.8);
        \filldraw[rectangle, blue] (-1, -4) -- (-0.5, -4) -- (-0.5, -4.5) -- (-1, -4.5) -- (-1, -4) node[black, right] at (-0.3, -4.25) {Driving with Customer};
        \filldraw[rectangle, orange] (4.2+\sh, -4) -- (4.7+\sh, -4) -- (4.7+\sh, -4.5) -- (4.2+\sh, -4.5) -- (4.2+\sh, -4)  node[black, right] at (4.9+\sh, -4.25) {Picking Up Customer};
        \filldraw[rectangle, darkgreen] (10+\shb, -4) -- (10.5+\shb, -4) -- (10.5+\shb, -4.5) -- (10+\shb, -4.5) -- (10+\shb, -4)  node[black, right] at (10.7+\shb, -4.25) {Idle};
        \filldraw[rectangle, purple] (-1, -5.4) -- (-0.5, -5.4) -- (-0.5, -5.9) -- (-1, -5.9) -- (-1, -5.4) node[black, right] at (-0.3, -5.65) {Driving to Charger};
        \filldraw[rectangle, yellow] (-1, -4.7) -- (-0.5, -4.7) -- (-0.5, -5.2) -- (-1, -5.2) -- (-1, -4.7) node[black, right] at (-0.3, -4.95) {Charging};
        \filldraw[rectangle, teal] (4.2+\sh, -4.7) -- (4.7+\sh, -4.7) -- (4.7+\sh, -5.2) -- (4.2+\sh, -5.2) -- (4.2+\sh, -4.7) node[black, right] at (4.9+\sh, -4.95) {Waiting for Charger};
        \draw[black, ultra thick] (10+\shb, -4.95) -- (10.5+\shb, -4.95) node[black, right] at (10.7+\shb, -4.95) {Avg SoC};
        \draw[brown, dashed, ultra thick] (4.2, -5.65) -- (4.7, -5.65) node[black, right] at (4.9, -5.65){Average Workload $(T_R\lambda)$};
    \end{tikzpicture}
}{Evolution of the state for CAD, CD, and Po2 with $n = 3072$, $m = 2600$ ($325$ locations with $8$ posts at each location), and $\lambda=160$ per min \label{fig: time_series_3_algos}}{}
\end{figure}
As the SoC of all the EVs is initialized away from 0, all three algorithms initially serve all incoming customers. The SoC of the fleet decreases with time, which then results in dropping incoming customers. Under CAD, all customers are served until the SoC of all the EVs reduces to the minimum allowable SoC.\footnote{In \Cref{fig: time_series_3_algos}, the SoC drops below $s_{\min}$ as EVs drive to the charger after finishing a trip.} As seen in \Cref{fig: time_series_3_algos} (left), initially, the blue region (almost) fills out the average workload, implying that all customers are initially served. As the fleet can only sustainably serve 90\% of the demand, the fleet average SoC (black line) consistently drops during this time. Now, with a lower fleet SoC, the density of available EVs reduces, which results in a larger pickup time (the orange region expands at $t \approx 100$). In turn, large pickup times negatively affect the fleet SoC, and this negative cycle continues until the fleet SoC drops to the minimum and stays there. In summary, in an attempt to serve all incoming customers, CAD experiences large pickup times in transience which, in turn, results in low fleet SoC in the steady state. On the other hand, CD and Po2 have a more stable fleet SoC in the steady state as these algorithms actively drop incoming customers if the nearby EVs do not have enough SoC.

\begin{figure}[thb!]
    \centering
    \FIGURE{
    \begin{tikzpicture}[scale=1, tight background]
        \node at (0, -0.7) {\includegraphics[scale=0.4]{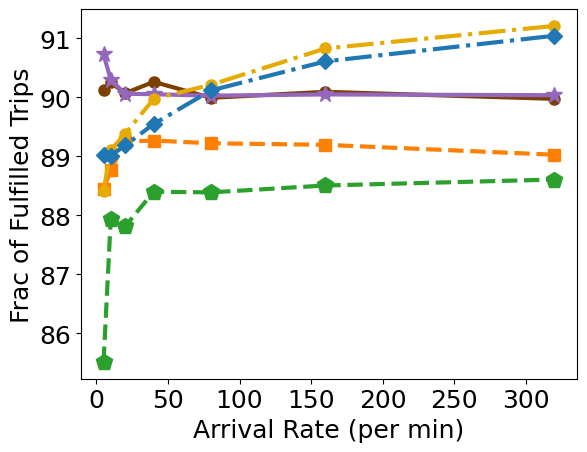}};
        \node at (6, -0.7) {\includegraphics[scale=0.4]{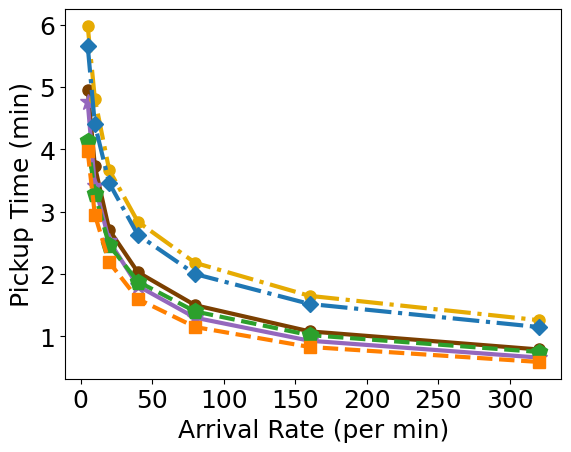}};
        \node at (12, -0.7) {\includegraphics[scale=0.4]{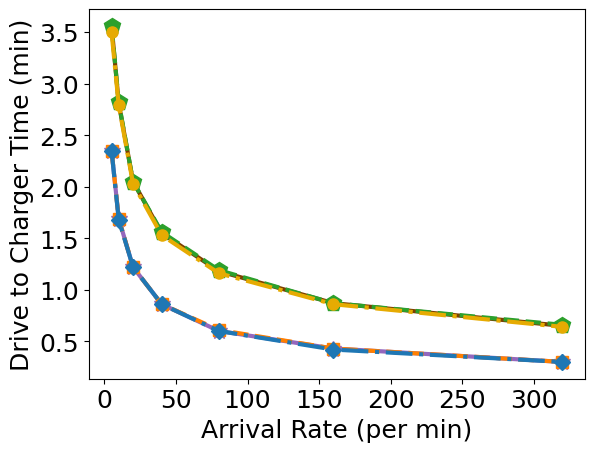}};
        \node at (3, -5.9) {\includegraphics[scale=0.4]{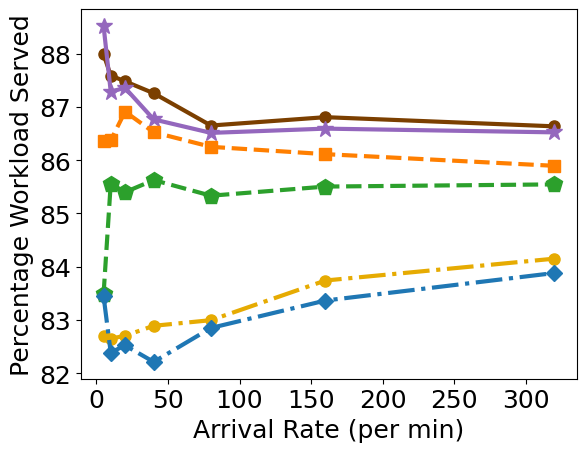}};
        \draw[rectangle] (6.2, -3.2) -- (11.4, -3.2) -- (11.4, -8.3) -- (6.2, -8.3) -- (6.2, -3.2);
        \node at (8.5, -3.5) {Legend};
        \draw[line_yellow, line width=0.8mm, dash pattern={on 7pt off 2pt on 1pt off 3pt}] (6.5, -4.25) -- (8, -4.25) node[black, right] at (8.2, -4.25) {CAD, $\beta = 0.803$};
        \draw[line_blue, line width=0.8mm, dash pattern={on 7pt off 2pt on 1pt off 3pt}] (6.5, -4.95) -- (8, -4.95) node[black, right] at (8.2, -4.95) {CAD, $\beta = 1$};
        \draw[line_brown, line width=0.8mm] (6.5, -5.65) -- (8, -5.65) node[black, right] at (8.2, -5.65) {Po2, $\beta = 0.803$};
        \filldraw[line_purple, line width=0.8mm] (6.5, -6.35) -- (8, -6.35) node[black, right] at (8.2, -6.35) {Po2, $\beta = 1$};
        \draw[line_green, line width=0.8mm, dash pattern={on 7pt off 2pt}] (6.5, -7.05) -- (8, -7.05)  node[black, right] at (8.2, -7.05) {CD, $\beta = 0.803$};
        \draw[line_orange, line width=0.8mm, dash pattern={on 7pt off 2pt}] (6.5, -7.75) -- (8, -7.75) node[black, right] at (8.2, -7.75) {CD, $\beta = 1$};
    \end{tikzpicture}
}{\revcolor{Fraction of fulfilled trips (top left), pickup time (top center), drive to charger time (top right), and percentage workload served (bottom left) as a function of the arrival rate for CAD, CD, and Po2. The fleet size and number of chargers correspond to Series A and C in Fig.~\ref{fig: asymptotic_sim}.}  \label{fig: comparing_algos}}{}
\end{figure}

We simulate the three matching algorithms for different values of $\lambda$. The fleet size and number of chargers correspond to Series A and Series C in \Cref{fig: asymptotic_sim}
, i.e., $\beta = 1$ and $\beta = 0.803$, respectively. %with the fleet size corresponding to 90\% service level for Po2. 
The results are reported in \Cref{fig: comparing_algos}. 
First, note that as the charging algorithm is identical across simulations, the drive to the charger time is identical for CAD, CD, and Po2. Next, we observe that the pickup time is the largest for CAD as it always tries to dispatch an EV irrespective of its distance to the customer. In addition, the pickup time for CD is the smallest as it either matches the closest EV or drops the customer. Even though the pickup times for Po2 are larger, it achieves a better service level compared to CD as it balances the SoC more efficiently across the fleet. Additionally, \Cref{prop:closest-fails} (Appendix~\ref{sec: closest_dispatch}) shows that CD has a fundamentally higher fleet requirement compared to Po$d$.  On the other hand, CAD can achieve a larger service level compared to Po2. However, observe that the percentage of workload served \footnote{Percentage of workload served is defined as the ratio of miles actually driven with a customer and the total miles requested by customers.
%, i.e. $\alphasup\tilde{T}_R/T_R$.
} is much smaller for CAD compared to Po2 as well as CD. A high service level and a low workload imply that CAD is biased toward shorter trips. In particular, as the fleet SoC is consistently low (as seen in Fig.~\ref{fig: time_series_3_algos}) for CAD, the longer the trip distance, the harder it is to find an available EV for it. Dropping long trip requests reduces the average fulfilled trip time, allowing the algorithm to serve more customers but with a lesser workload. In turn, Po2 (and more generally Po$d$) may be preferred over CD and CAD as it showcases good performance in terms of pickup times, service level, and workload, and it also nearly achieves the minimum infrastructure planning requirements.

\section{Additional Set-up Description and Results for Section~\ref{sec:sim_chicago}}\label{sec:sim_chicago_additional_setup}
\subsection{Simulation Detailed Description}
\subsubsection*{Trip data:}
For each trip, the dataset contains the start time, end time, origin location, destination location, and trip distance. 
In the dataset, the start time is rounded off to the nearest 15-minute mark. To avoid the arrival of customers in batches, we add a uniformly random time between $[0, 15]$ minutes to the start and end times while preserving the trip duration. We filter 95 percentile of the trips concerning their origin and destination locations, i.e., we only consider trips that are within 2.5 and 97.5th percentile of the population latitude and longitude. Lastly, to ensure fast simulation times, we uniformly filter 60\% of the trips (unless stated otherwise) and use the smaller dataset to run our simulations.

\subsubsection*{Vehicles and Chargers:} We consider \revcolor{four} representative EVs: Nissan Leaf ($\sim$\$30,000), \footnote{\url{https://ev-database.org/car/1656/Nissan-Leaf}} Tesla Model 3 Standard Range ($\sim$\$40,000), \footnote{\url{https://ev-database.org/car/1991/Tesla-Model-3}} Mustang Mach-E SR RWD ($\sim$\$50,000), \footnote{\url{https://ev-database.org/car/2034/Ford-Mustang-Mach-E-SR-RWD}}
and \revcolor{Hyundai IONIQ 5 ($\sim$\$50,000)}\footnote{\url{https://ev-database.org/car/2237/Hyundai-IONIQ-5-84-kWh-AWD}}. In the rest of the section, we denote them by Nissan, Tesla, Mustang, and Hyundai, respectively. After accounting for a nominal 10\% battery degradation typically experienced in 4-5 years of usage, we consider a battery capacity of 35.1 kWh for Nissan, 51.25 kWh for Tesla, 64.8 kWh for Mustang, and 75.6 kWh for Hyundai. While the consumption/mileage depends on the weather conditions and driving behavior, we consider a fixed average value. In particular, we consider a consumption of 270 Wh/mi for Nissan, 230 Wh/mi for Tesla, 250 Wh/mi for Mustang, and 260 Wh/mi for Hyundai. This consumption translates to a range of 130 miles for Nissan, 200 miles for Tesla, 260 miles for Mustang, and 290 miles for Hyundai. We consider a fixed charge rate of either 20 kW or 100 kW, corresponding to Level 2 chargers and DC fast chargers, respectively. For simplicity, we assume that the charge rate is agnostic to the SoC of the EV, which is reasonable for Level 2 chargers, but may not hold for DC fast chargers when the SoC increases beyond 80\%.

We set the charging station locations to be uniformly at random within the city of Chicago. We first independently sample the latitude and longitude within their bounds. We then set the sampled point as a charger location if it is at most a 20-minute driving distance away from one of the trips in the dataset. If not, we discard this location and resample. Such a procedure is computationally efficient and avoids hard-coding the boundary of the region of interest. Each charging station is endowed with 4 charging posts (chargers), and we denote this by $m_p = 4$. We initialize the location of each EV as the origin of a uniformly sampled trip from the dataset, and the SoCs to be uniformly sampled in $[0.7, 0.9]$. As each simulation instance is 3 days long, the initialization of the EVs does not affect our results.

\subsubsection*{System's dynamics:} For carrying out the simulation, we need access to pickup distances that are not available in the trip dataset, as it depends on the current locations of the EVs. We approximate the pickup times by the Manhattan distance between the trip origin and the EV's current location. Similarly, we approximate the drive to the charger time by the Manhattan distance between the selected charger's location and the EV's current location.

We implement the Power-of-$d$ vehicles dispatch policy as follows. For every customer, we consider the $d$ closest EVs that are either charging, waiting at the charger, or driving to the charger. We allow for fractional values of $d$ by sampling (independent of everything else) $\lceil d \rceil$ EVs with probability $\lceil d \rceil - d$ and $\lfloor d \rfloor$ otherwise. We dispatch the EV with the highest SoC among the ones in consideration if the following two conditions are met. First, the EV has enough SoC to drive to the customer's pickup location, then to the customer's destination, and then to the nearest charger while sparing at least 5\% SoC. Second, the pickup time is at most $T_{P, \max}$, which is equivalent to implementing the admission control policy as defined in Definition~\ref{def: admission_control}. We consider $T_{P, \max} \in \{30, 45, 60, \infty\}$ minutes for our simulations. If the selected EV does not meet the above two criteria, the customer is dropped, and it leaves the system immediately.    

After serving a customer, the EV is then dispatched to the nearest available charging station to replenish. We deem a charging station to be available if the number of available posts is more than 0.5 times the number of EVs already driving to that charging station. The contribution of the EVs driving to the charging station towards the availability ensures that we do not send too many EVs to the same charger. On the other hand, these EVs have not reached the charger yet and should not be assumed to be occupying a charging post, as this may result in poor utilization. Thus, we discount these EVs by a factor of 2.

\subsection{Detailed Description of Policies }\label{sec: detailed_policy_description}
We provide a more detailed description of the policies introduced in \Cref{sec:sim_chicago} including Closest Dispatch (CD), Closest Available Dispatch (CAD), Power-of-$d$ Vehicles Dispatch (Po$d$), Power-of-$T_{P, \max}$ (Po$T_{P, \max}$), and Charge at Night (CaN) + CAD/Po$d$ \revcolor{and its variants}:
\begin{itemize}
    \item \textbf{Closest Dispatch (CD):} The closest EV (charging, idle, waiting for charger, or driving to the charger) serves the customer if it has sufficient SoC and the pickup time is at most $T_{P, \max} \in \{30, 45, 60, \infty\}$ mins. On the other hand, if the closest EV does not have sufficient SoC or if the pickup time is more than $T_{P, \max}$, then the customer is dropped. After the EV serves the customer, it is immediately dispatched to the nearest charger.
    \item \textbf{Closest Available Dispatch (CAD):} The closest EV (charging, idle, waiting for charger, or driving to the charger) among the ones with sufficient SoC serves the customer if the pickup time is at most $T_{P, \max} \in \{30, 45, 60, \infty\}$ mins. If no such EV exists with pickup time at most $T_{P, \max}$, then the customer is dropped. After the EV serves the customer, it is immediately dispatched to the nearest charger.
    \item \textbf{Power-of-$d$ Vehicles Dispatch (Po$d$):} The EV (charging, idle, waiting for charger, or driving to the charger) with the highest SoC among the $d$ closest ones serves the customer if it has sufficient SoC and the pickup time is at most $T_{P, \max} \in \{30, 45, 60, \infty\}$. If all the $d$ closest EVs do not have sufficient SoC or the pickup time of the selected EV is more than $T_{P, \max}$, then the customer is dropped. After the EV serves the customer, it is immediately dispatched to the nearest charger.
    \item \textbf{Power-of-$T_{P, \max}$ (Po$T_{P, \max}$):} The EV (charging, idle, waiting for charger, or driving to the charger) with the highest SoC among the ones with pickup time at most $T_{P, \max} \in \{5, \hdots, 20\}$ serves the customer if it has sufficient SoC. If all EVs within $T_{P, \max}$ minutes away do not have sufficient SoC, then the customer is dropped. After the EV serves the customer, it is immediately dispatched to the nearest charger.
    \item \textbf{Charge at Night (CaN) + CAD/Po$d$:} CaN is a charging policy and can be combined with any of the above matching policies. We consider CaN+CAD and CaN+Po$d$ in the simulations. Under CaN, we do not allow interrupting a charging process, so only idle EVs are candidates for dispatching. The charging decisions are made as follows: During the day (5 am to 12 am), EVs are only sent to charge if their SoC falls below \revcolor{10\%/20\% (optimally chosen)}, and they only replenish the bare minimum needed to drive till 12 am. During the night (12 am to 5 am), all EVs are sent to charge and replenished to 100\%. \revcolor{We implement a logic to make sure that the charging is spread out uniformly through the night to ensure that not too many EVs are charging at any point in time.}
    \revcolor{\item \textbf{Charge at Night with Relocations (CaN-R) + Po$d$:} This algorithm is almost the same as CaN, with one key difference. After serving a trip, the EV is sent to the nearest available charging station, where it pretends to be charging, without actually increasing the SoC. Once the fake charging session is completed, the EV is deemed to be idle, at which point, it can again be sent to the nearest available charger for another fake charging session. In addition to the idle EVs, all EVs relocating to a charger or fake charging can be matched to an incoming customer. Such a policy mimics the relocation behavior of the vanilla Po$d$ policy.}
    \revcolor{\item \textbf{Charge at Night with Relocations (CaN-R) + Po$d$ at Night (Po$d$-N):} This policy is almost the same as CaN-R+Po$d$ with one key difference. During the night, the charging session of an EV can be interrupted to serve the incoming customer. Thus, idle, charging, driving to the charger, or waiting for the charger, vehicles are available to serve incoming trips during the night. The value of $d$ is allowed to be different for nighttime and daytime. This policy is inspired by the EV-aware policy as in Proposition~\ref{prop: pod_var}.}
\end{itemize}
\subsection{Infrastructure Planning Using Po\texorpdfstring{$d$}{d} ODEs}\label{app:infra_plan_odes}
To solve the system of ODEs \eqref{eq:odes-power-d}, we approximate the derivative by the finite difference and use a step size of 0.1 minutes. The arrival rate $\lambda(t)$ and the mean trip time $T_R(t)$ are allowed to be time-dependent and calculated by averaging the trips arriving within $t \pm 2$ minutes. The pickup time and the drive to the charger time are approximated by \revcolor{$T_P(t) = 99(n-\XB_{\Nt}(t))^{-0.57}$ and $T_{DC}(t) = 42(m - \XC_{\Nt-1}(t))^{-0.36}$} respectively. These functional forms are obtained by simulating Power-of-$d$ for the detailed model (the simulation in \Cref{sec:sim_chicago}) and fitting the data obtained from it. We set a nominal admission control for incoming customers and charging EVs to ensure that $T_P(t)$ and $T_{DC}(t)$ are bounded. In particular, we enforce $\XB_{\Nt} \leq n - 50$ and $\XC_{\Nt-1} \leq m - 20$. Lastly, the maximum service time $T_B$ is approximated by the mean cumulative trip, pickup, and drive to the charger time, which turns out to be $30.5$ minutes. 
Given the parameters for the Tesla Model 3, the value of $\Nt$ turns out to be 24 trips. 
As $T_P(t), T_{DC}(t)$, and
$T_B$ needs to be estimated from the real simulations, we also test the performance of the ODE by considering a parameter misspecification of $\pm 20\%$, i.e., $1.2(T_P(t), T_{DC}(t), T_B)$ and $0.8(T_P(t), T_{DC}(t), T_B)$. Fig~\ref{fig:ode_stackplot} compares the output of the ODE (with no parameter misspecification) to that of the detailed simulation by plotting the evolution of the system state with time. We highlight that the output of the ODE has a similar structure to that of the detailed simulation, providing a computationally efficient first-order approximation of the detailed simulation.
\begin{figure}
    \centering
    \FIGURE{\scalebox{0.7}{
    \begin{tikzpicture}
        \node at (0, 0) {\includegraphics[scale=0.65]{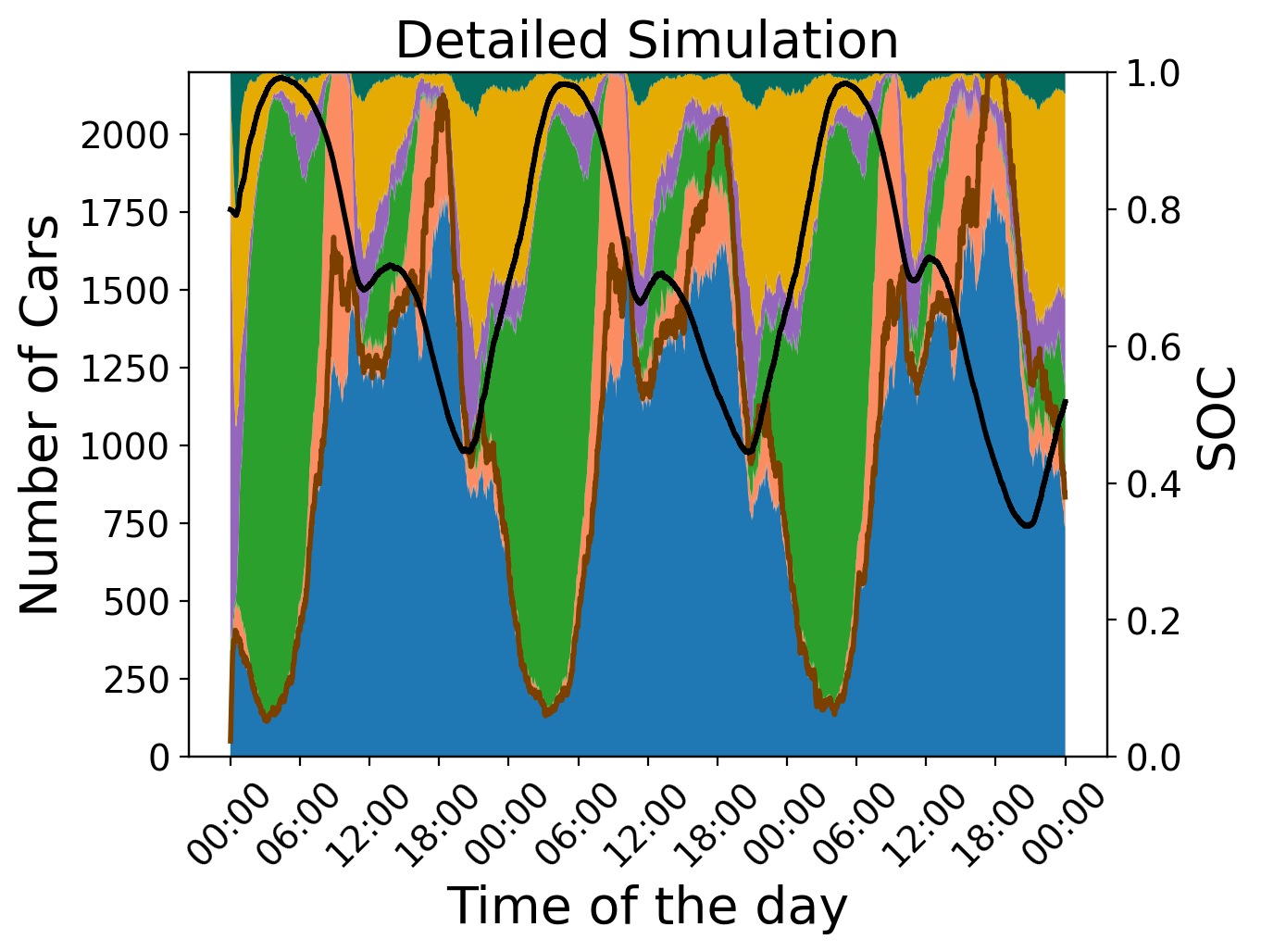}};
        \node at (11.5, 0) {\includegraphics[scale=0.65]{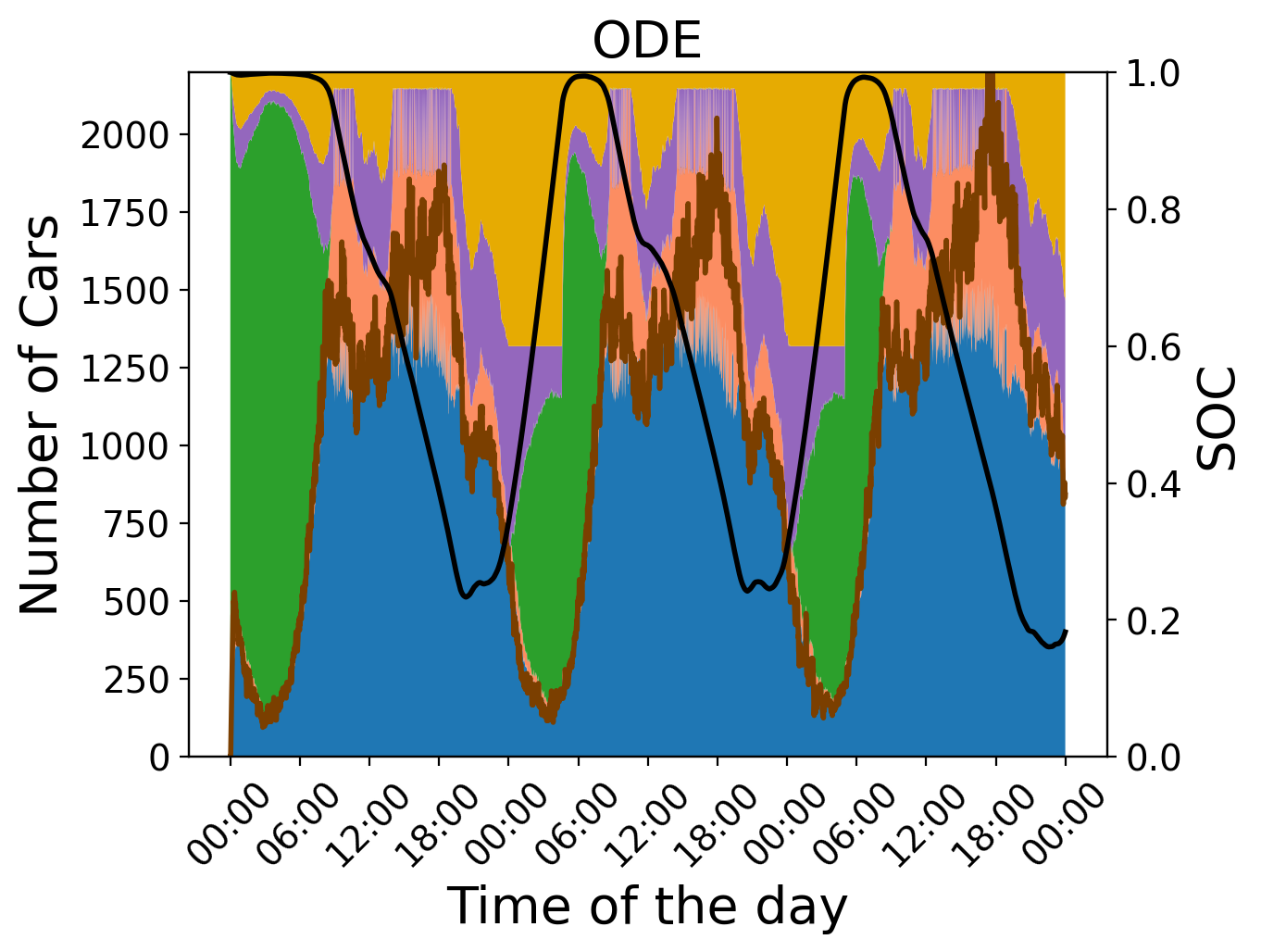}};
        \begin{scope}[yshift=-1cm,xshift=-2.3cm]
        \draw[rectangle] (-1.2, -3.8) -- (18, -3.8) -- (18, -5.4) -- (-1.2, -5.4) -- (-1.2, -3.8);
        \filldraw[rectangle, blue] (-1, -4) -- (-0.5, -4) -- (-0.5, -4.5) -- (-1, -4.5) -- (-1, -4) node[black, right] at (-0.3, -4.25) {Driving with Customer};
        \filldraw[rectangle, orange] (4.2, -4) -- (4.7, -4) -- (4.7, -4.5) -- (4.2, -4.5) -- (4.2, -4)  node[black, right] at (4.9, -4.25) {Picking up Customer};
        \filldraw[rectangle, darkgreen] (10, -4) -- (10.5, -4) -- (10.5, -4.5) -- (10, -4.5) -- (10, -4)  node[black, right] at (10.7, -4.25) {Idle};
        \filldraw[rectangle, purple] (12.8, -4) -- (13.3, -4) -- (13.3, -4.5) -- (12.8, -4.5) -- (12.8, -4) node[black, right] at (13.5, -4.25) {Driving to Charger};
        \filldraw[rectangle, yellow] (-1, -4.7) -- (-0.5, -4.7) -- (-0.5, -5.2) -- (-1, -5.2) -- (-1, -4.7) node[black, right] at (-0.3, -4.95) {Charging};
        \filldraw[rectangle, teal] (4.2, -4.7) -- (4.7, -4.7) -- (4.7, -5.2) -- (4.2, -5.2) -- (4.2, -4.7) node[black, right] at (4.9, -4.95) {Waiting for Charger};
        \draw[black, ultra thick] (10, -4.95) -- (10.5, -4.95) node[black, right] at (10.7, -4.95) {Avg SoC};
        \draw[brown, ultra thick] (12.8, -4.95) -- (13.3, -4.95) node[black, right] at (13.5, -4.95){(Potential) Active Trips};;
        \end{scope}
    \end{tikzpicture}}}{
    Stack plot that tracks the state of the system as a function of the simulation time for Tesla Model 3 with $n=2200$, $m=900$ under po$d$ with $d=1.4$ and $T_{P, \max} = 45$. The left and the right panels are the output of the detailed simulation and the ODE as described in Appendix~\ref{app:infra_plan_odes}, respectively.
    \label{fig:ode_stackplot}}{}
\end{figure}

\subsection{Additional Results and Discussions  }\label{sec:additional_result_sim}
\subsubsection*{Supplementary Results for Figure~\ref{fig: pod_vs_can_stackplot}}
In the left panel of \Cref{fig: pod_vs_can_stackplot}, we compare the performance of the policies by calculating the fleet size corresponding to a target 90\% workload served. To generate the plot, we filter $\{0.2, 0.3, 0.4, 0.5, 0.6\}$ fraction of the demand and set the number of charging stations to be $\{150, 225, 300, 375, 450\}$ respectively, corresponding to the setting of plenty of chargers. Then, we conduct several simulations to infer the fleet size corresponding to the 90\% workload. In \Cref{fig: pod_vs_can_stackplot_full}, we present a more complete version of the stack plots in \Cref{fig: pod_vs_can_stackplot} in the main body of the paper. 
\begin{figure}[tbh!]
    \centering
    \FIGURE{\scalebox{0.7}{
    \begin{tikzpicture}
        \node at (0, 0) {\includegraphics[scale=0.65]{figures_final/1_stackplot_tesla_can_po2_2700.png}};
        \node at (11.5, 0) {\includegraphics[scale=0.65]{figures_final/2_fleet_stackplot_pod_d_2_45_th_ev_2700_nc_450_tesla.png}};
        \node at (0, -9) {\includegraphics[scale=0.65]{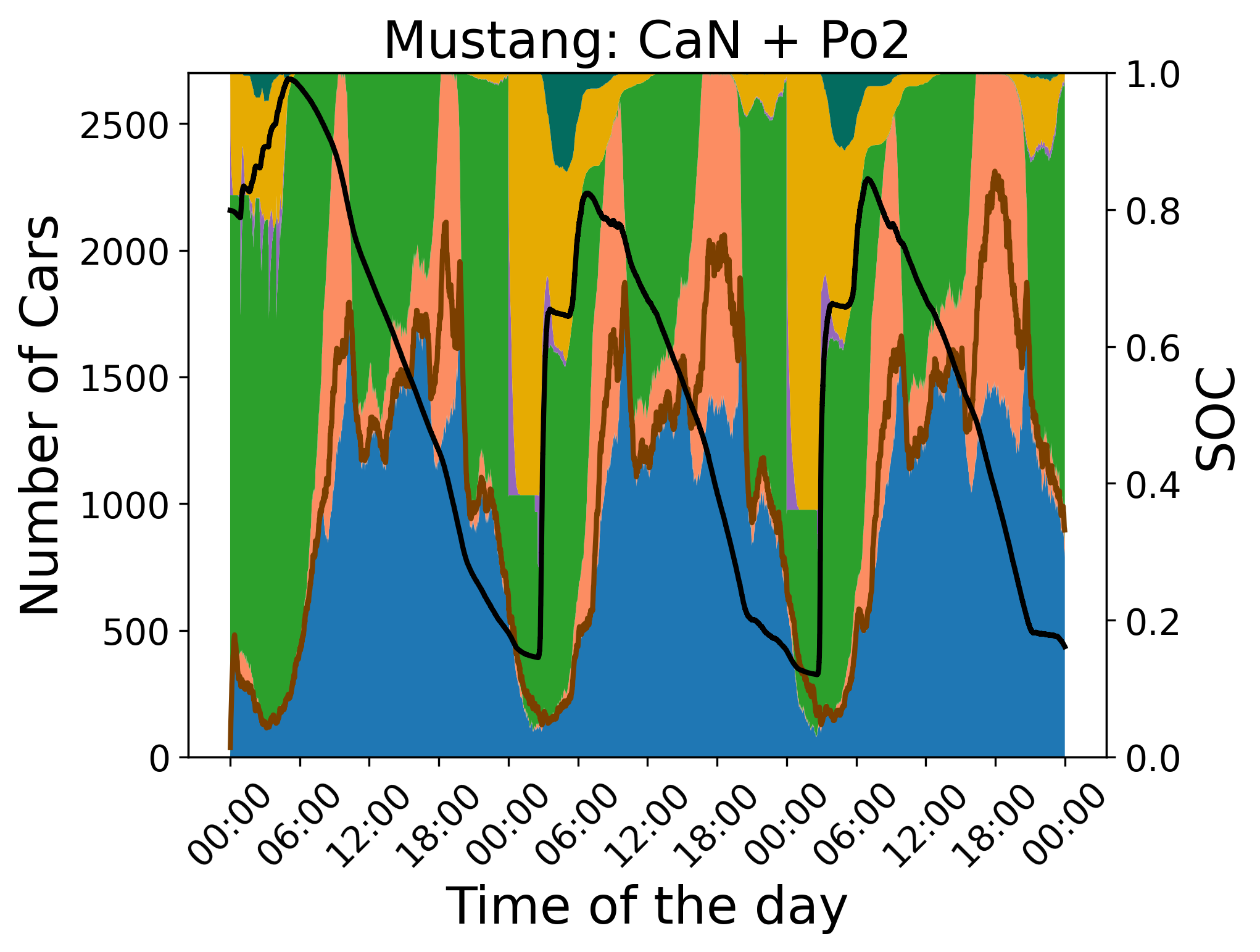}
        };
        \node at (11.5, -9) {\includegraphics[scale=0.65]{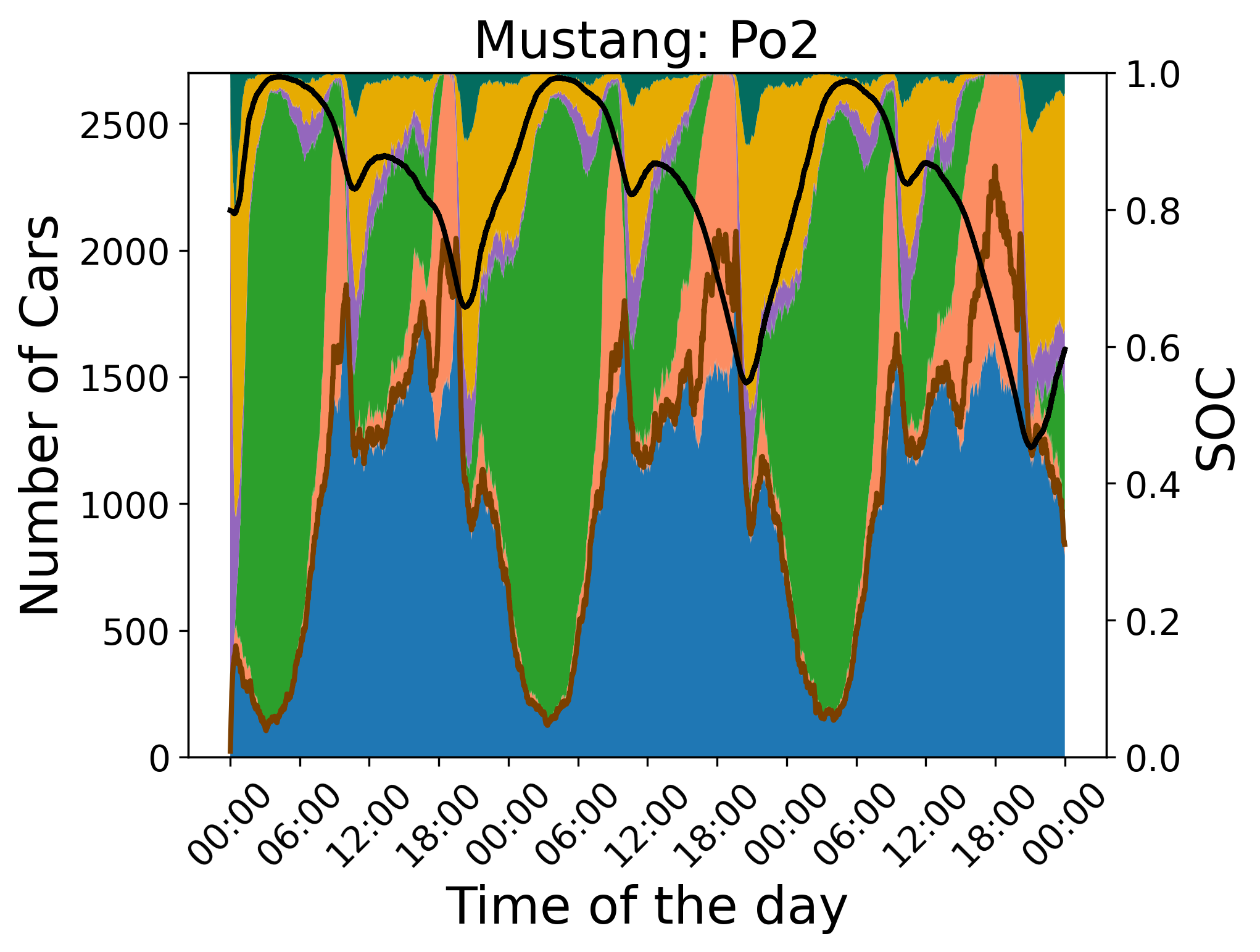}};
        \begin{scope}[yshift=-10cm,xshift=-2.3cm]
        \draw[rectangle] (-1.2, -3.8) -- (18, -3.8) -- (18, -5.4) -- (-1.2, -5.4) -- (-1.2, -3.8);
        \filldraw[rectangle, blue] (-1, -4) -- (-0.5, -4) -- (-0.5, -4.5) -- (-1, -4.5) -- (-1, -4) node[black, right] at (-0.3, -4.25) {Driving with Customer};
        \filldraw[rectangle, orange] (4.2, -4) -- (4.7, -4) -- (4.7, -4.5) -- (4.2, -4.5) -- (4.2, -4)  node[black, right] at (4.9, -4.25) {Picking up Customer};
        \filldraw[rectangle, darkgreen] (10, -4) -- (10.5, -4) -- (10.5, -4.5) -- (10, -4.5) -- (10, -4)  node[black, right] at (10.7, -4.25) {Idle};
        \filldraw[rectangle, purple] (12.8, -4) -- (13.3, -4) -- (13.3, -4.5) -- (12.8, -4.5) -- (12.8, -4) node[black, right] at (13.5, -4.25) {Driving to Charger};
        \filldraw[rectangle, yellow] (-1, -4.7) -- (-0.5, -4.7) -- (-0.5, -5.2) -- (-1, -5.2) -- (-1, -4.7) node[black, right] at (-0.3, -4.95) {Charging};
        \filldraw[rectangle, teal] (4.2, -4.7) -- (4.7, -4.7) -- (4.7, -5.2) -- (4.2, -5.2) -- (4.2, -4.7) node[black, right] at (4.9, -4.95) {Waiting for Charger};
        \draw[black, ultra thick] (10, -4.95) -- (10.5, -4.95) node[black, right] at (10.7, -4.95) {Avg SoC};
        \draw[brown, ultra thick] (12.8, -4.95) -- (13.3, -4.95) node[black, right] at (13.5, -4.95){(Potential) Active Trips};;
        \end{scope}
    \end{tikzpicture}}}{Stack plot that tracks the state of the system as a function of the simulation time. Charge at Night + Power-of-$d$ with \revcolor{$d=2, T_{P, \max} = 60$} (best choice among $\{30, 45, 60, \infty\}$) is simulated for Tesla (top left) and Mustang (bottom left). The performance is compared against Power-of-$d$ with $d=2, T_{P, \max}=45$ for Tesla (top right) and Mustang (bottom right). \revcolor{We set $n=2700$, $m=450 \times 4$, and $r_c = 20$ [kW].}
    \label{fig: pod_vs_can_stackplot_full}}{}
\end{figure}

\subsubsection*{Supplementary Results for \Cref{fig: pod_vs_can_20_100_opt_d}}
\revcolor{
The right panel in \Cref{fig: pod_vs_can_20_100_opt_d} shows that, in line with our discussion in \Cref{sec:var_arrival_rate}, the performance of long-range vehicles under CaN-R can be further improved by an EV-aware policy. The stack plots in \Cref{fig: pod_vs_can_stackplot__mustang_hyundai_canr} and \Cref{fig: pod_vs_can_stackplot_hyundai_canr} help to understand why. First, in \Cref{fig: pod_vs_can_stackplot__mustang_hyundai_canr} for Mustang and Hyundai under CaN-R, Closest Dispatch and 20[kW] (top and bottom left panels), the system is unable to serve a small fraction of the workload towards the end of the day, at around midnight, because the EVs that are not charging have low SoC which continue to decrease as they serve some of the night demand. That is, some EVs may not have enough SoC to serve trips at the start of the night. Moreover, some EVs may need to charge at the beginning of the next day. Although this effect is marginal in the presence of fast charging (top and bottom right panels in \Cref{fig: pod_vs_can_stackplot__mustang_hyundai_canr}), it can be completely overcome by an EV-aware policy. Indeed, a policy that utilizes Po$d$ (with an optimized $d$) at night (see \Cref{fig: pod_vs_can_stackplot_hyundai_canr}) is able to interrupt the charging of vehicles with higher SoC that would otherwise remain at a charging station, thereby increasing vehicle availability and allowing the fleet to uniformly raise its SoC, leading to improved performance.
}
\begin{figure}[tbh!]
    \FIGURE{
    \scalebox{0.7}{
    \begin{tikzpicture}
        \node at (0, 0) {\includegraphics[scale=0.65]{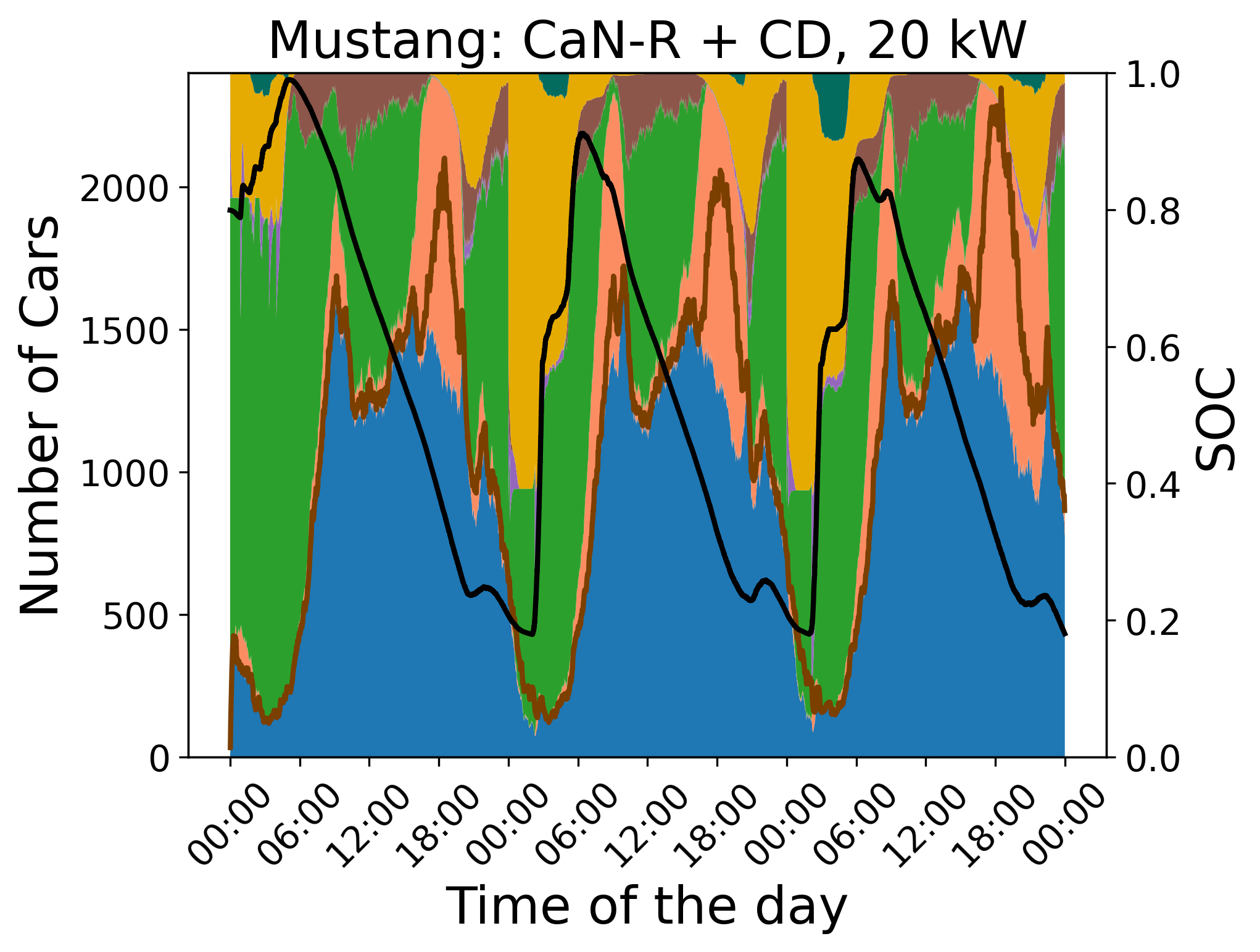}
        };
        \node at (11.5, 0) {\includegraphics[scale=0.65]{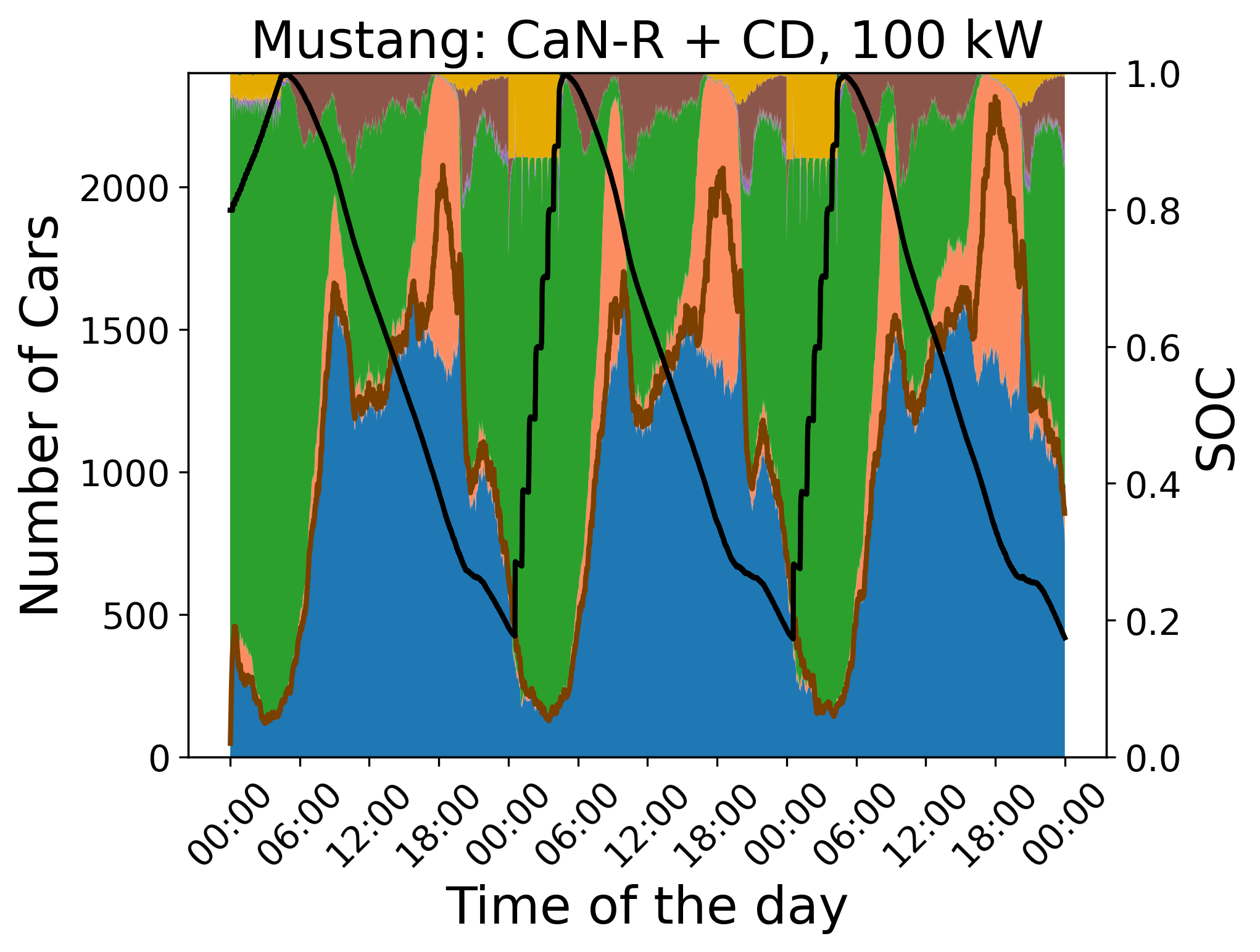}
        };
        \node at (0, -9.0) {\includegraphics[scale=0.65]{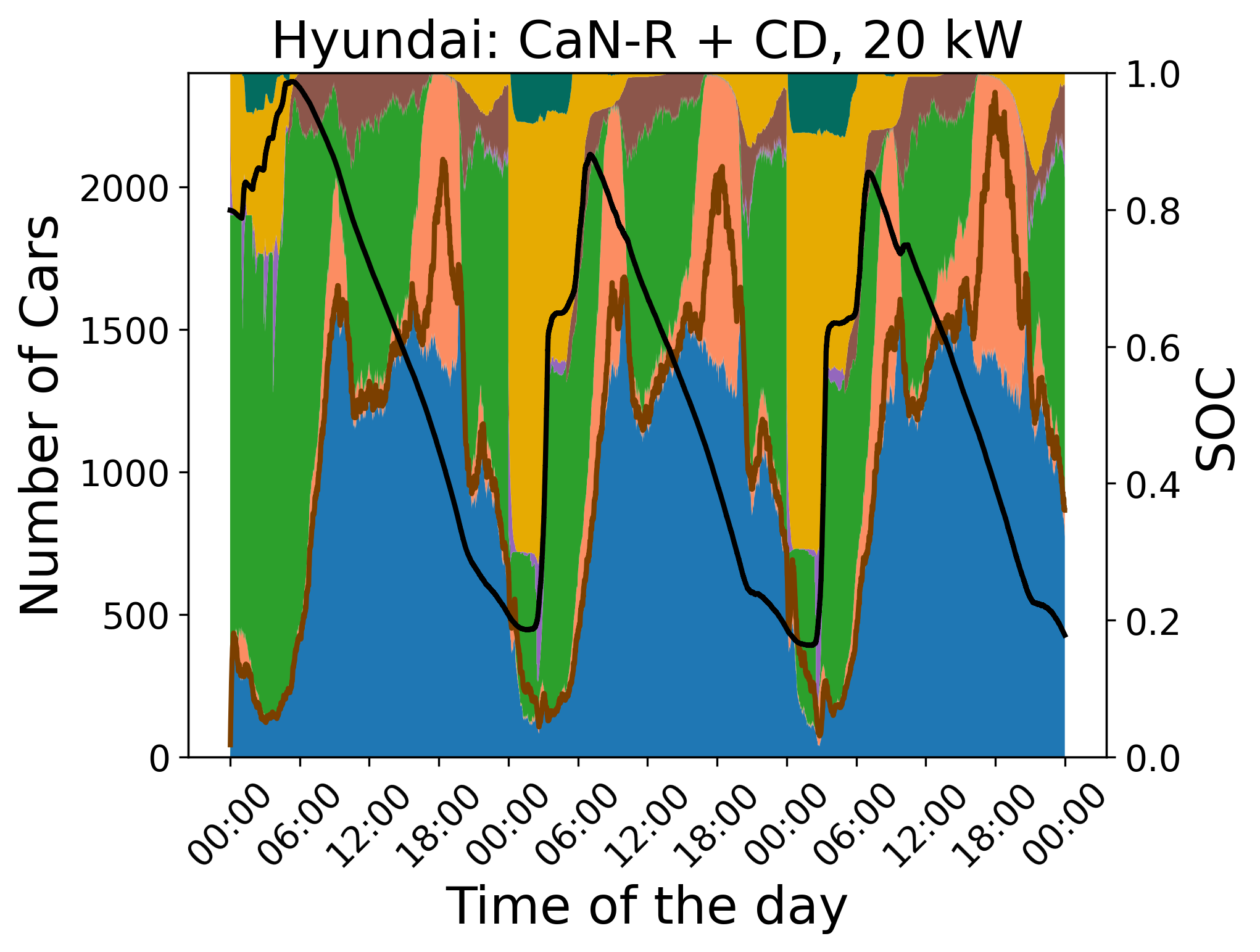}
        };
        \node at (11.5, -9.0) {\includegraphics[scale=0.65]{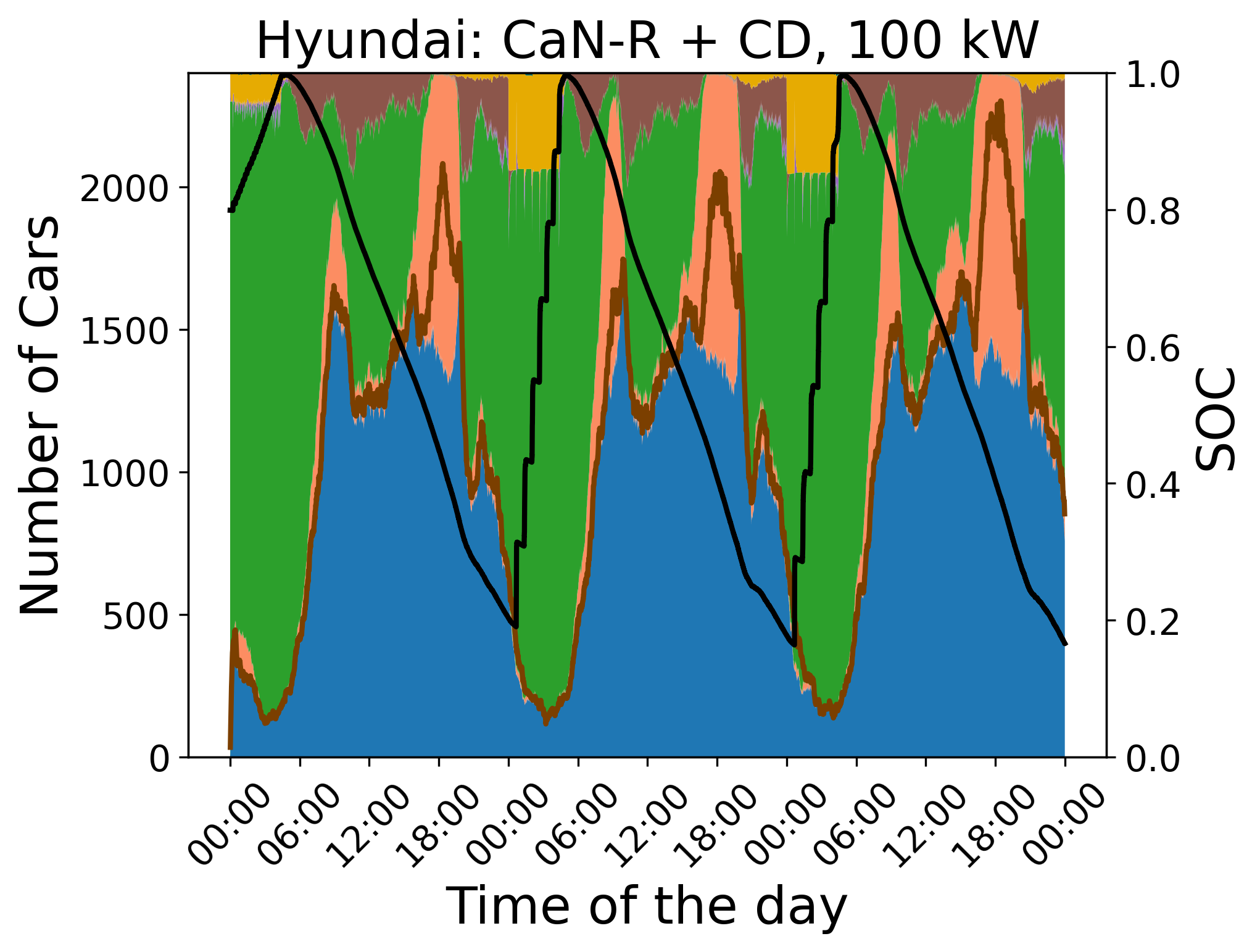}
        };
        \begin{scope}[yshift=-10cm,xshift=-1.2cm]
        \draw[rectangle] (-1.2, -3.8) -- (15.2, -3.8) -- (15.2, -6.2) -- (-1.2, -6.2) -- (-1.2, -3.8);
        \filldraw[rectangle, blue] (-1, -4) -- (-0.5, -4) -- (-0.5, -4.5) -- (-1, -4.5) -- (-1, -4) node[black, right] at (-0.3, -4.25) {Driving with Customer};
        \filldraw[rectangle, orange] (4.2, -4) -- (4.7, -4) -- (4.7, -4.5) -- (4.2, -4.5) -- (4.2, -4)  node[black, right] at (4.9, -4.25) {Picking up Customer};
        \filldraw[rectangle, darkgreen] (10, -4) -- (10.5, -4) -- (10.5, -4.5) -- (10, -4.5) -- (10, -4)  node[black, right] at (10.7, -4.25) {Idle};
        \filldraw[rectangle, purple] (-1, -4.7) -- (-0.5, -4.7) -- (-0.5, -5.2) -- (-1, -5.2) -- (-1, -4.7) node[black, right] at (-0.3, -4.95) {Driving to Charger};
        \filldraw[rectangle, brown] (4.2, -4.7) -- (4.7, -4.7) -- (4.7, -5.2) -- (4.2, -5.2) -- (4.2, -4.7) node[black, right] at (4.9, -4.95) {Relocating to Charger};

        \filldraw[rectangle, yellow] (10, -4.7) -- (10.5, -4.7) -- (10.5, -5.2) -- (10, -5.2) -- (10, -4.7) node[black, right] at (10.7, -4.95) {Charging};

        \filldraw[rectangle, teal] (-1, -5.4) -- (-0.5, -5.4) -- (-0.5, -5.9) -- (-1, -5.9) -- (-1, -5.4) node[black, right] at (-0.3, -5.65) {Waiting for Charger};
        
        \draw[black, ultra thick] (4.2, -5.65) -- (4.7, -5.65) node[black, right] at (4.9, -5.65) {Avg SoC};
        \draw[brown, ultra thick] (10, -5.65) -- (10.5, -5.65) node[black, right] at (10.7, -5.65){(Potential) Active Trips};;
        \end{scope}
    \end{tikzpicture}}}
    {\revcolor{Stack plots tracking the system state over time for Mustang under 
    CaN-R + CD + Level 2 charging (top left), and 
    CaN-R + CD + DC fast charging (top right); 
    and for Hyundai under 
    CaN-R + CD + Level 2 charging (bottom left), and 
    CaN-R + CD + DC fast charging (bottom right)
    with $n=2400$ and $m = 450\times 4$.}
    \label{fig: pod_vs_can_stackplot__mustang_hyundai_canr}}{}
    %\vspace{-0.25cm}
\end{figure}

\begin{figure}[tbh!]
    \FIGURE{
    \scalebox{0.7}{
    \begin{tikzpicture}
        \node at (0, -9.0) {\includegraphics[scale=0.65]{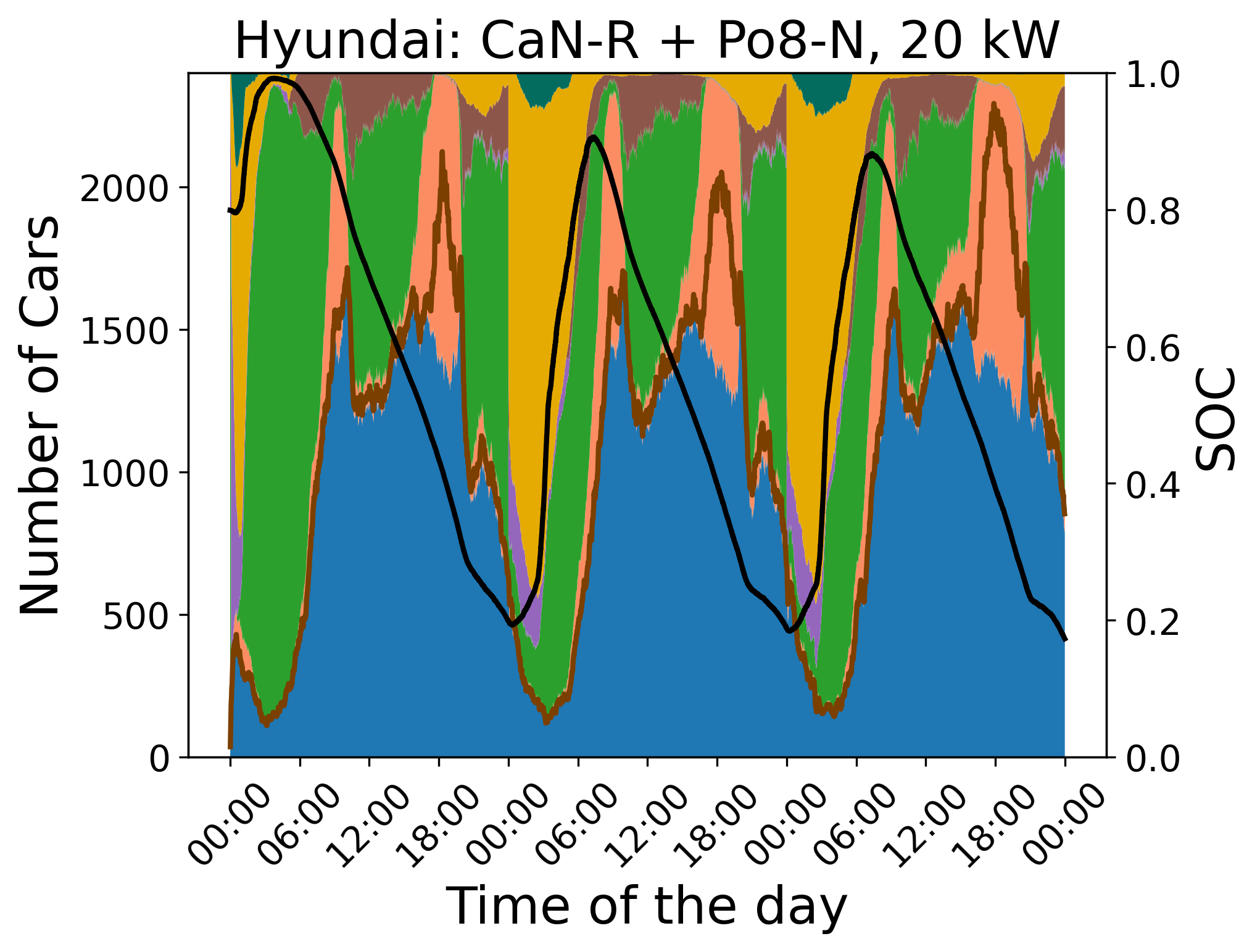}};
        \node at (11.5, -9.0) {\includegraphics[scale=0.65]{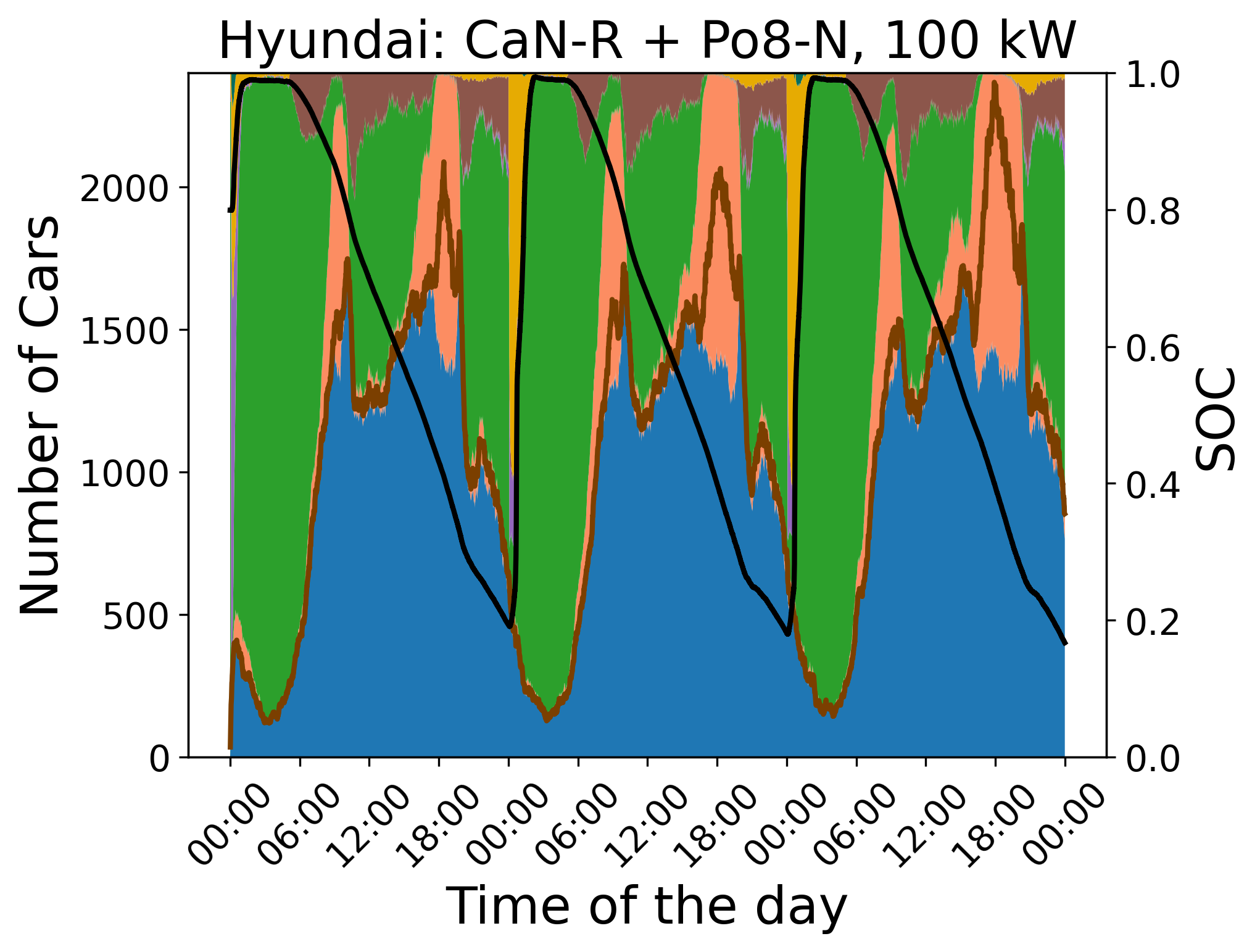}};
        \begin{scope}[yshift=-10cm,xshift=-1.2cm]
        \draw[rectangle] (-1.2, -3.8) -- (15.2, -3.8) -- (15.2, -6.2) -- (-1.2, -6.2) -- (-1.2, -3.8);
        \filldraw[rectangle, blue] (-1, -4) -- (-0.5, -4) -- (-0.5, -4.5) -- (-1, -4.5) -- (-1, -4) node[black, right] at (-0.3, -4.25) {Driving with Customer};
        \filldraw[rectangle, orange] (4.2, -4) -- (4.7, -4) -- (4.7, -4.5) -- (4.2, -4.5) -- (4.2, -4)  node[black, right] at (4.9, -4.25) {Picking up Customer};
        \filldraw[rectangle, darkgreen] (10, -4) -- (10.5, -4) -- (10.5, -4.5) -- (10, -4.5) -- (10, -4)  node[black, right] at (10.7, -4.25) {Idle};
        \filldraw[rectangle, purple] (-1, -4.7) -- (-0.5, -4.7) -- (-0.5, -5.2) -- (-1, -5.2) -- (-1, -4.7) node[black, right] at (-0.3, -4.95) {Driving to Charger};
        \filldraw[rectangle, brown] (4.2, -4.7) -- (4.7, -4.7) -- (4.7, -5.2) -- (4.2, -5.2) -- (4.2, -4.7) node[black, right] at (4.9, -4.95) {Relocating to Charger};

        \filldraw[rectangle, yellow] (10, -4.7) -- (10.5, -4.7) -- (10.5, -5.2) -- (10, -5.2) -- (10, -4.7) node[black, right] at (10.7, -4.95) {Charging};

        \filldraw[rectangle, teal] (-1, -5.4) -- (-0.5, -5.4) -- (-0.5, -5.9) -- (-1, -5.9) -- (-1, -5.4) node[black, right] at (-0.3, -5.65) {Waiting for Charger};
        
        \draw[black, ultra thick] (4.2, -5.65) -- (4.7, -5.65) node[black, right] at (4.9, -5.65) {Avg SoC};
        \draw[brown, ultra thick] (10, -5.65) -- (10.5, -5.65) node[black, right] at (10.7, -5.65){(Potential) Active Trips};;
        \end{scope}
    \end{tikzpicture}}}
    {\revcolor{Stack plots tracking the system state over time for Hyundai under 
    CaN-R + Po8 at night + Level 2 charging (left), and 
    CaN-R + Po8 at night + DC fast charging (right)
    with $n = 2400$ and $m = 450\times 4$.}
    \label{fig: pod_vs_can_stackplot_hyundai_canr}}{}
\end{figure}

\subsubsection*{Supplementary Results for \Cref{sec:sim_chicago_policy_comp}}
In \Cref{fig: real_sim_load_balancing_pickup_only}, we show the pickup times corresponding to \Cref{fig: real_sim_load_balancing}.  In line with the base model (cf. \eqref{eq: pickup_power_of_d}), we observe that they increase at a diminishing rate with $d$. Large pickup times dominate for large values of $d$, and the balancing of SoC across the fleet dominates for small values of $d$. 
\begin{figure}[hbt!]
    \centering    
    \FIGURE{\scalebox{0.75}{
    \begin{tikzpicture}
        \node at (0,0) {\includegraphics[scale=0.5]{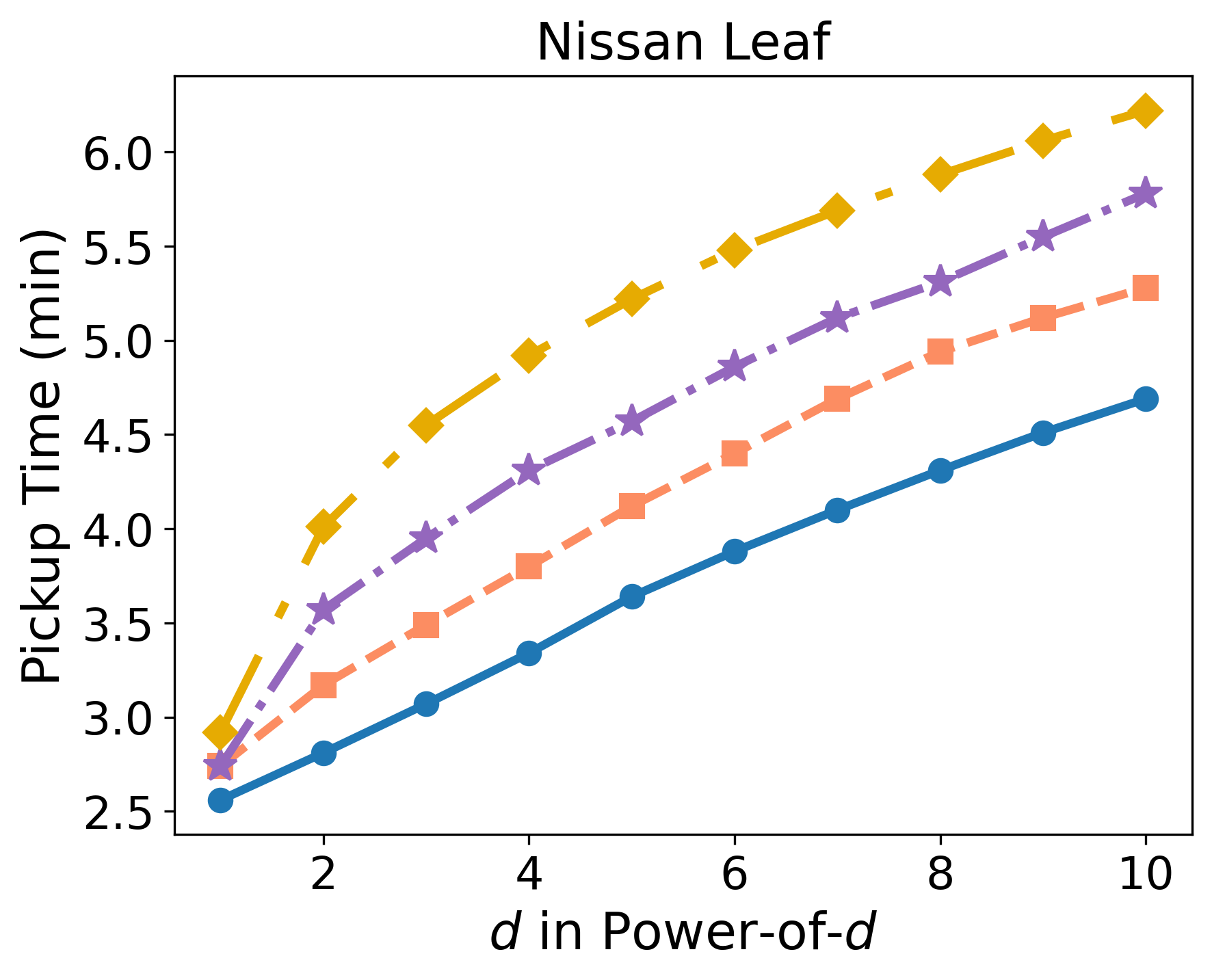}};
        \node at (8, 0) {\includegraphics[scale=0.5]{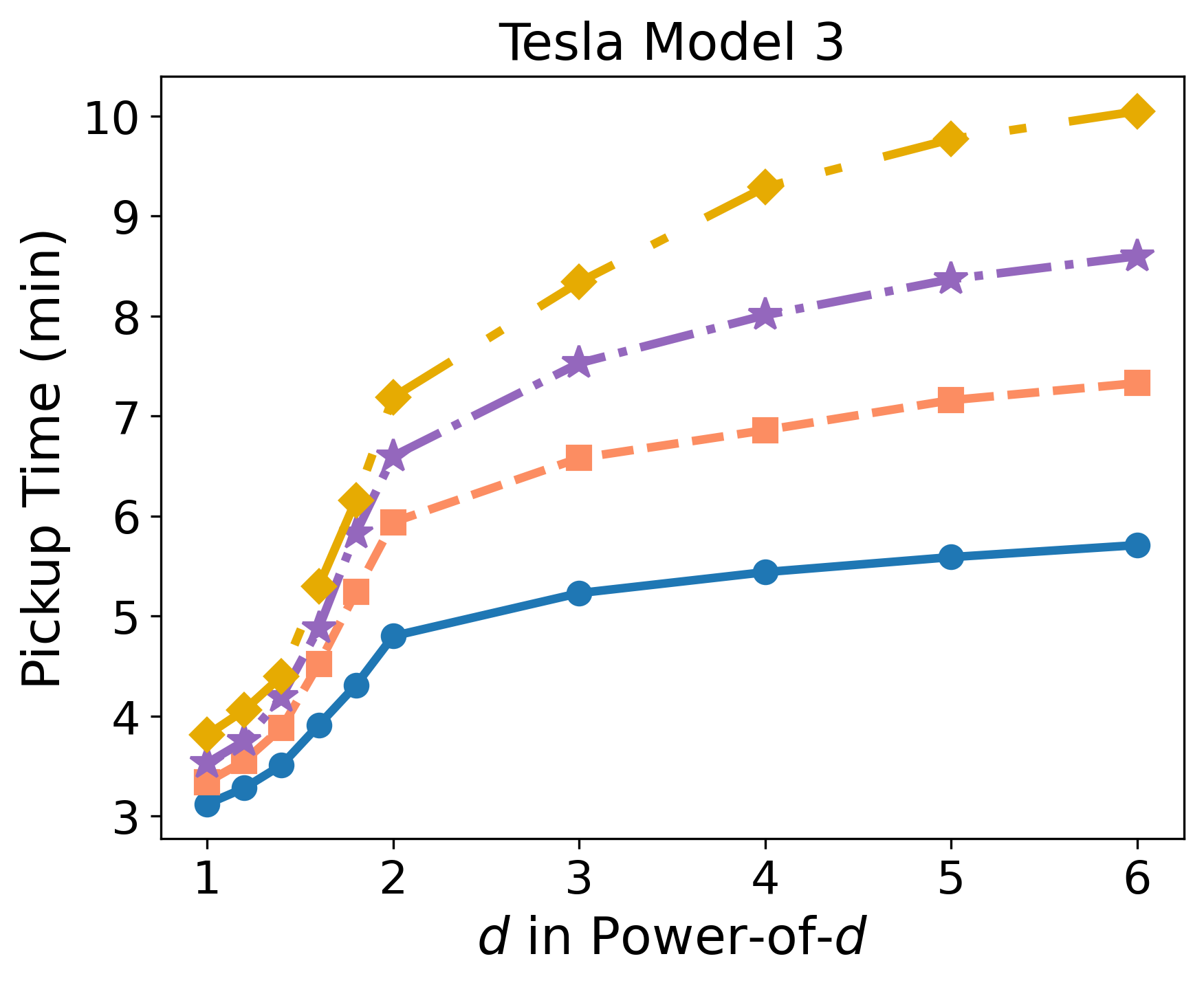  }};
        \begin{scope}[xshift = -0.5cm, yshift = 6.5cm]
        \draw (-1.3, -9.65) rectangle (10.7, -10.35);
        \draw[line_blue, ultra thick] (-1, -10) -- (0, -10) node[black, right] at (0.1, -10) {30 min};
        \filldraw[line_blue] (-0.5, -10) circle (3pt);
        \draw[line_orange, ultra thick, , dash pattern={on 5pt off 2pt on 5pt off 2pt}] (1.8, -10) -- (2.8, -10) node[black, right] at (2.9, -10) {45 min};
        \filldraw[line_orange] (2.2, -9.9) rectangle (2.4, -10.1); 
        \draw[line_purple, ultra thick, , dash pattern={on 7pt off 2pt on 1pt off 2pt}] (4.6, -10) -- (5.6, -10) node[black, right] at (5.7, -10) {60 min};
        \node[line_purple] at (5.1, -10) {$\bigstar$};
        \draw[line_yellow, ultra thick, dash pattern={on 9.5pt off 2pt on 1pt off 2pt}] (7.4, -10) -- (8.4, -10) node[black, right] at (8.5, -10) {No Bound};
        \filldraw[line_yellow] (7.9, -9.85) -- (8.05, -10) -- (7.9, -10.15) -- (7.75, -10) -- (7.9, -9.85);
        \end{scope}
    \end{tikzpicture}}
    }
    {
    \revcolor{Pickup Times (averaged over four runs) as a function of $d$ in Power-of-$d$ with $n = 2400$, $m = 225 \times 4$.}
    \label{fig: real_sim_load_balancing_pickup_only}}{}
\end{figure}
\Cref{tab:chicago_served_requests} complements the results in \Cref{tab: real_sim_cad_vs_pod}. As for the case of served workload, Po$d$ consistently delivers a competitive percentage of customers served even after optimizing the maximum allowable pickup time, $T_{P, \max}$. In all instances, it is \revcolor{better or comparable to CD, CAD, and Po$T_{P,\max}$.}
\begin{table}[hbt!]
    \centering
    \TABLE{Fraction of Customers Served \revcolor{(averaged over four runs)} for CD, CAD, Power-of-$T_{P, \max}$, and Power-of-$d$ with the optimal $T_{P, \max}$ and $d$. We set $n=2400$ and $m=225 \times 4$.
    \label{tab:chicago_served_requests}}{
    \begin{tabular}{|c|c|c|c|c|} \hline
       $T_{P, \max} = \infty$  & CD & CAD & Po$T_{P, \max}$ & Po$d$ $(d=\cdot)$ \\ \hline \hline
        Nissan & 86.2\% & 90.6\% & $\leq$ 65\% & 92.7\% $(5)$ \\ \hline
        Tesla & 93.5\% & 92.2\% & $\leq$ 75\% & 94.4\% $(1.2)$ \\ \hline
        Mustang & 94.5\% & 91.8\% & $\leq$ 75\% & 94.5\% $(1)$ \\ \hline \hline
        best $T_{P, \max}$ & CD $(T_{P, \max} = \cdot)$ & CAD $(T_{P, \max} = \cdot)$ & Po$T_{P, \max}$ $(T_{P, \max} = \cdot)$ & Po$d$ $(d=\cdot, T_{P, \max}=\cdot)$ \\ \hline \hline
        Nissan & 86.3\% (60) & 92.8\% (30) & 92.7\% (9) & 95.0\% (45, 6) \\ \hline
        Tesla & 93.9\% (45) & 95.6\% (30) & 96.7\% (11) & 95.5\% (45, 1.4) \\ \hline
        Mustang & 95.5\% (45) & 95.6\% (30) & 96.9\% (11) & 95.6\% (45, 1.05) \\ \hline
    \end{tabular}}{}
\end{table}
\subsubsection*{Supplementary Discussion for \Cref{sec:chicago-perf-pod}} \label{app: relocation_discussion}
As discussed in \Cref{sec:chicago-perf-pod}, Po$d$ dispatches an EV to the nearest available charger after finishing a trip, which naturally relocates the EVs to move away from destination-heavy locations, leading to reduced pickup times. While such a policy is simple and captures some of the benefits of relocation, one can design a more sophisticated \textit{relocation policy} in the context of EVs to further improve the system performance. More concretely, one shortcoming of the natural relocation under Po$d$ is that the relocation pattern depends on system parameters like the number of chargers, charge rate, and the density of chargers. For example, if we have too many chargers in any given location, then the relocation effect for that location is minimal. Similarly, if the charge rate is high, the chargers quickly become available, leading to several EVs accumulating at the same charger (EVs idle at the charger once fully charged until they are matched with an incoming customer). This endogenous relocation could lead to counterintuitive behavior, such as reduced workload, as the number of chargers increases much beyond the infrastructure planning prescription. These observations and our results suggest that a more comprehensive approach to optimizing relocations in a system with EVs---especially for longer-range, fast-charging systems---could be a fruitful avenue of future research.

\end{APPENDICES}

\end{document}